\documentclass[graybox,envcountchap]{svmono}

\usepackage{german}

\usepackage{amssymb}
\usepackage{mathtools}
\usepackage{amsfonts}
\usepackage{amsthm}

\usepackage{mathptmx}
\usepackage{helvet}
\usepackage{courier}
\usepackage{type1cm}         

\usepackage{makeidx}         % allows index generation
\usepackage{graphicx}        % standard LaTeX graphics tool
\usepackage{multicol}        % used for the two-column index
\usepackage[bottom]{footmisc}% places footnotes at page bottom

\usepackage[colorlinks=true]{hyperref}
\usepackage{pst-plot, pstricks, tabularx, array}
\usepackage{fancybox}
\usepackage{wasysym}
\usepackage{psfrag}
\usepackage{epsfig}
\usepackage{color}
\usepackage{subfig}
\usepackage{enumerate}
\usepackage{stmaryrd}
\usepackage{multirow}
\makeindex             % used for the subject index
\renewcommand{\D}{\displaystyle}
\newcolumntype{R}{>{\centering\arraybackslash}X} % neu
\newcommand{\N}{{\mathbb N}}

\newcommand{\Z}{{\mathbb Z}}
\newcommand{\R}{{\mathbb R}}
\newcommand{\C}{{\mathbb C}}
\newcommand{\Q}{{\mathbb Q}}

\newcommand{\F}{{\mathcal{F}}}
\newcommand{\M}{{\mathcal{M}}}

\newcommand{\sign}{{\text{sign}}}
\newcommand{\sgn}{{\mbox{sgn\,}}}
\numberwithin{equation}{chapter} 

\newcommand{\dokendDef}{\ensuremath{\hfill\square}}
\newcommand{\dokendSatz}{\ensuremath{\hfill\square}}
\newcommand{\dokendBem}{\ensuremath{\hfill\square}}
\newcommand{\dokendProof}{\ensuremath{\hfill\blacksquare}}
\newcommand{\ggT}{{\mbox{ggT\,}}}
\newcommand{\Det}{{\mbox{Det\,}}}
\newcommand{\kgV}{{\mbox{kgV\,}}}
\newcommand{\Id}{{\mbox{Id\,}}}

\newtheoremstyle{style}   
{0.5cm}                 %Space above    
{0.5cm}                 %Space below    
{\normalfont}           %Body font: original {\normalfont}    
{}  %Indent amount (empty = no indent,%\parindent = paraindent)    
{\normalfont\bfseries}  %Thm head font original 
{\normalfont}{\newline}%
{{\normalfont\bfseries \thmname{#1 }\thmnumber{#2:}\thmnote{~#3}}}
\theoremstyle{style}
\newtheorem{Satz}{Satz}[chapter]
\newtheorem{Def}[Satz]{Definition}
\newtheorem{DefSatz}[Satz]{Definition und Satz}

\newtheorem{Bem}[Satz]{Bemerkung}

\newtheorem{Beis}[Satz]{Beispiel}
\newtheorem{Folg}[Satz]{Folgerung}
\newtheorem{Zus}[Satz]{Zusammenfassung}
\newtheoremstyle{styleAuf}   
{0.8cm}                 %Space above    
{0.8cm}                 %Space below    
{\normalfont}           %Body font: original {\normalfont}    
{}  %Indent amount (empty = no indent,%\parindent = paraindent)    
{\normalfont\bfseries}  %Thm head font original 
{\normalfont}{\newline}%
{{\normalfont\bfseries \thmname{#1 }\thmnumber{#2:}\thmnote{~#3}}}
\theoremstyle{styleAuf}
\newtheorem{Auf}{Aufgabe}[chapter]
\theoremstyle{styleAuf}
\makeatletter
\let\@DefAlt\Def
\def\@DefAltI{\let\oldnewline\newline
	\let\newline\relax
	\@DefAlt\leavevmode\global\let\newline\oldnewline}
\def\@DefAltII{\@DefAlt}
\renewcommand{\Def}{%
	\@ifnextchar\begin\@DefAltI \@DefAltII}
\makeatother

\newcommand{\Mod}{{\mbox{mod\,}}}
\newcommand{\ord}{{\mbox{ord\,}}}

\begin{document}
	
	\author{Yauheniya Abramchuk, Alina Bondarava, Matthias Kunik}
	\title{Elementare Zahlentheorie}
	\subtitle{F\"ur Studierende der Mathematik, Informatik und Lehramt}

	\maketitle
	
	\frontmatter
	
	%%%%%%%%%%%%%%%%%%%%%%%%%%%%%%%%%%%%%%%%%%%%%%%%%%%%%%
	
	%\include{dedic}
	%\include{foreword}
	% !TEX root = buch_kunik.tex
\chapter*{Vorwort}

Gegenstand der elementaren Zahlentheorie sind vorrangig die nat\"urlichen Zahlen $\N=\{1,2,3,...\}$
sowie die ganzen Zahlen $\Z=\{0,\pm1,\pm2,\pm3,...\}$. 
Im Zusammenhang mit den Grundrechenarten in diesen Zahlenbereichen spielen dabei der Begriff der 
Teilbarkeit\index{Teilbarkeit} sowie der Euklidische Algorithmus \index{Euklidischer Algorithmus} eine zentrale Rolle.\\

Wir behandeln  im vorliegenden Lehrbuch klassische Themen der Zahlentheorie, 
die erstmals von Gau{\ss} in seinen "`{Disquisitiones Arithmeticae}"', 
Untersuchungen \"uber h\"ohere Arithmetik \cite{gauss}, 
zu einer systematischen Wissenschaft entwickelt worden sind. Der einfache Euklidische Algorithmus
in Lektion \ref{cha:2} ist die Grundlage f\"ur den Fundamentalsatz \ref{satz:2_12} der Arithmetik.
In nachfolgender Lektion \ref{cha:3} wird er f\"ur die Kettenbruchentwicklung 
reeller Zahlen zum erweiterten Euklidischen Algorithmus ausgebaut.
Die Anwendung der Kettenbruchentwicklung auf die reell quadratischen Irrationalzahlen
liefert wiederum mit den Grundlagen aus Abschnitt \ref{cha:8A} in Abschnitt \ref{cha:8B} und \ref{cha:8C} genau die perio\-dischen 
Ketten\-br\"uche sowie einen Algorithmus zur
Reduktion indefiniter quadratischer Formen. 
%Einen einfacheren Reduk\-tions\-algo\-rithmus tragen wir
%unabh\"angig von der Theorie der Kettenbr\"uche f\"ur positiv definite quadratische Formen in Abschnitt %\ref{cha:8D} nach.\\

Auch bei der Behandlung der Farey-Folgen in Abschnitt \ref{cha:4}
machen wir vom erweiterten Euklidischen Algorithmus Gebrauch, man vergleiche insbesondere
den Approximationssatz f\"ur Farey-Br\"uche \ref{satz:4_13} mit dem Satz \ref{satz:3_17} 
von den besten rationalen Approximationen in der Theorie der Kettenbr\"uche.
Damit zieht sich der Euklidische Algorithmus wie ein roter Faden durch das gesamte
Fundament der elementaren Zahlentheorie. \\

Die wichtigsten algebraischen Strukturen der elementaren Zahlentheorie sind Gruppen, 
Ringe und K\"orper,
mit deren Einf\"uhrung wir deshalb in Lektion \ref{cha:1} beginnen. Wir beschr\"anken uns dabei
auf die Untersuchung  derjenigen algebraischen Strukturen, 
die f\"ur sp\"atere Zwecke ben\"otigt werden. 
Wichtige Beispiele hierf\"ur sind die Permutationsgruppen, die
im Anhang \ref{cha:anhang2} und dem Aufgabenteil von Lektion \ref{cha:1} behandelt werden,
die Gruppe der multiplikativen zahlentheoretischen
Funktionen, die in Abschnitt \ref{cha:5} eingef\"uhrt und untersucht werden, aber auch
die Restklassenringe bzw. die prime Restklassengruppe
bez\"uglich eines Moduls in Lektion \ref{cha:6}. 
%
%Im Abschnitt \ref{cha:8E} betrachten wir noch
%ein tiefergehendes Resultat, n\"amlich die Gruppenstruktur einer bestimmten Art der Multiplikation 
%quadratischer primitiver Formenklassen mit fester, nicht quadratischer Diskriminante.
%Diese Theorie zur "`{Komposition der Formen}"' geht auf Gau{\ss} zur\"uck und wurde in \cite{gauss} 
%noch in einem weiter gespannten Rahmen behandelt.\\

Neben dem Euklidischen Algorithmus nimmt auch die Kongruenzrechnung  in Lektion \ref{cha:6}
einen bedeutenden Platz in der elementaren Zahlentheorie ein. Diese wird in  Lektion \ref{cha:7}
auf die Theorie der quadratischen Reste angewendet, deren wichtigstes Resultat 
das quadratische Reziprozit\"atsgesetz in Satz \ref{satz:7_9} darstellt. 
Gau{\ss} hat diesen Satz nicht nur als Erster
bewiesen, sondern in \cite{gauss} gleich mehrere Beweise geliefert.\\

Jede Lektion beginnt mit einer kurzen \"Ubersicht. Der nachfolgende theoretische Teil
wird durch zahlreiche Beispiele anschaulich gemacht, und  die Lektion
wird mit ausgew\"ahlten und vollst\"andig gel\"osten Aufgaben zur \"Ubung und Vertiefung des Stoffes abgeschlossen.
Im Anhang findet der Leser neben allgemeinen Grundlagen und Notationen zu  logischen Aussagen,
Mengen und Abbildungen die bereits erw\"ahnte kurze Einf\"uhrung der Permutationsgruppen
sowie eine Primzahltabelle.
\vspace{\baselineskip}
\begin{flushright}\noindent
Magdeburg, \today\hfill\\ {\it Yauheniya Abramchuk, Alina Bondarava und Matthias Kunik} \\
\end{flushright}

	\tableofcontents

	\mainmatter
	%%%%%%%%%%%%%%%%%%%%%%%%%%%%%%%%%%%%%%%%%%%%%%%%%%%%%%%
\chapter{Algebraische Grund\-lagen der elementaren Zahlen\-theorie}\label{cha:1}
Wir beginnen mit allgemeinen, aber unverzichtbaren algebraischen Grundlagen 
zu Gruppen und Ringen, zugeschnitten auf unsere sp\"ateren Anwendungen. 
Als Begleitlekt\"ure empfehlen wir van der Waerden's Lehrbuch \cite[Kapitel 2,\,3]{waerden1}
"`{Algebra~I}"' sowie das Lehrbuch \cite[Kapitel 2,\,3]{hornfeck} "`{Algebra}"' 
von Bernhard Hornfeck.\\
	\section{Gruppen}\label{cha:1A}
Beim Rechnen in Gruppen hat man nur eine (in der Regel multiplikativ geschriebene) 
assoziative Verkn\"upfung "`{$\circ$}"', die in einem gewissen Sinne "`{umkehrbar}"' ist:

\begin{Def}[Gruppen]\label{def:1_1}
Eine nichtleere Menge $G$ mit einer Verkn\"upfung $\circ:G\times G\to G$ hei{\ss}t Gruppe\index{Gruppe}\label{Gruppe}, wenn folgende Axiome erf\"ullt sind:
\begin{enumerate}[(G1)]
	\item  Assoziativgesetz\index{Assoziativgesetz}\label{Assoziativgesetz}: $(a\circ b)\circ c=a\circ(b\circ c)$ f\"ur alle $a,b,c\in G$.
	\item [(G2)] Existenz eines Einselementes\index{Einselement einer Gruppe}\label{Einselement_einer_Gruppe}: Es existiert eine Linkseins $e\in G$ mit $e\circ a=a$ f\"ur alle $a\in G$.
	\item [(G3)] Existenz eines inversen Elementes\index{inverses Element}\label{inverses_Element}: Zu jedem $a\in G$ existiert eine Linksinverse $a^{-1}\in G$ mit $a^{-1}\circ a=e$.
\end{enumerate}
Gilt zus\"atzlich
 \begin{enumerate}[(G4)]
 	\item  $a\circ b=b\circ a$ f\"ur alle $a,b\in G$,
 \end{enumerate}
 so erhalten wir einen Spezialfall und nennen die Gruppe $G$ abelsch \index{abelsche Gruppe}\label{abelsche_Gruppe}bzw. kommutativ.
 \index{kommutative Gruppe}\label{kommutative_Gruppe}~\dokendDef
\end{Def}
Die Elementeanzahl $|G|$ hei{\ss}t auch die Ordnung der Gruppe.

\begin{Satz}\label{satz:1_2} Es sei $(G,\circ, e)$ eine Gruppe mit der Linkseins $e$. Dann gilt: 
\begin{enumerate}[(a)]
	\item Ist $a^{-1}$ eine Linksinverse von $a\in G$, so auch eine Rechtsinverse:\\ $a^{-1}\circ a=e \Rightarrow a\circ a^{-1}=e$.
	\item Die Linkseins ist auch Rechtseins: $a\circ e=a$ f\"ur alle $a\in G$.
	\item F\"ur alle $a,b\in G$ sind die Gleichungen 
	$a \circ x=b$ bzw. $y\circ a=b$ in $G$ eindeutig l\"osbar.
	\item Das Einselement in $G$ ist eindeutig, ebenso die Inverse $a^{-1}$ zu jedem $a\in G$.
	\item $(a^{-1})^{-1}=a$ f\"ur alle $a\in G$.
	\end{enumerate}
\hfill\dokendSatz
\end{Satz}
{\bf Beweis:}~ 
\begin{enumerate}[(a)]
	\item Ist $a^{-1}$ Linksinverse zu $a\in G$, so auch Rechtsinverse, denn:
	Es sei $a^{-1}$ ein Linksinverses zu $a$ in $G$, d.h. $a^{-1}\circ a=e$, und $(a^{-1})^{-1}$ ein Linksinverses zu $a^{-1}$ in $G$, d.h. $(a^{-1})^{-1}\circ a^{-1}=e$. Dann gilt unter Verwendung der Gruppenaxiome
	$$	\begin{array}{lcl}
	a\circ a^{-1}& \underset{\text{(G2)}}{=} & e\circ (a\circ a^{-1}) \\
    & \underset{\text{(G3)}}{=} &((a^{-1})^{-1}\circ a^{-1})\circ (a\circ a^{-1}) \\
    & \underset{\text{(G1)}}{=} &(a^{-1})^{-1}\circ ((a^{-1}\circ a)\circ a^{-1}) \\
    & \underset{\text{(G3)}}{=} &(a^{-1})^{-1}\circ (e\circ a^{-1})\\
    & \underset{\text{(G2)}}{=} &(a^{-1})^{-1}\circ a^{-1}\underset{\text{(G3)}}{=}e.
	\end{array} $$
	
	\item Die Linkseins ist auch Rechtseins: Aus $e\circ a=a$ folgt mit (a) auch 
	$$a\circ e\underset{\text{(G3)}}{=}a\circ(a^{-1}\circ a)\underset{\text{(G1)}}{=}(a\circ a^{-1})\circ a\underset{\text{(a)}}{=}e\circ a\underset{\text{(G2)}}{=}a.$$
	\item Die Gleichung $a\circ x=b$ ist in $G$ l\"osbar:\\
	W\"ahle $a^{-1}\in G$ zu $a$ gem\"a{\ss} (G3) und setze  $x:=a^{-1}\circ b$. Dann wird nach (a) $$a\circ x=a\circ (a^{-1}\circ b)\underset{\text{(G1)}}{=}(a\circ a^{-1})\circ b\underset{\text{(a)}}{=}e\circ b\underset{\text{(G2)}}{=}b.$$
	
	Die Gleichung ist in $G$ eindeutig l\"osbar: \\
	Es sei $a\circ x=a\circ x'$ mit $a,x,x'\in G$. Dann folgen $a^{-1}\circ a\circ x=a^{-1}\circ a\circ x'$, also $e\circ x=e\circ x'$ und somit $x=x'$.\\ 
	
	Analog zeigt man die eindeutige L\"osbarkeit von $y\circ a=b$.\\
	
	\item Es folgt (d) sofort aus (c), da die Gleichungen $y\circ a=a$ bzw. $y\circ a=e$ in $G$ eindeutig l\"osbar sind.\\
	
	\item Nach (c) und (G3) hat $y\circ a^{-1}=e$ die eindeutige L\"osung $y=(a^{-1})^{-1}$, und nach~(a) darf $y=a$ gesetzt werden. Somit ist $(a^{-1})^{-1}=a$.
\end{enumerate}\dokendProof

{\it Merke:}~Es sei $(G,\circ,e)$ Gruppe.
\begin{enumerate}[1)]
	\item Bei der Auswertung eines Mehrfachproduktes $a_1\circ a_2\circ ... \circ a_n$ in $G$ k\"onnen wegen (G1) Klammern um je zwei Gruppenelemente beliebig gesetzt werden, so dass man Klammern nicht unbedingt mitschreiben muss. Auf die Reihenfolge der Faktoren ist jedoch zu achten, wenn (G4) nicht gilt.\\
	
	\item Es gibt genau ein $e\in G$ mit
	$$
	e\circ a=a\circ e=a \, \quad \mbox{ ~f\"ur~alle~} a\in G.
	$$
	
	\item Zu jedem $a\in G$ gibt es genau ein $a^{-1}\in G$ mit 
	$$
	a^{-1}\circ a=a\circ a^{-1}=e.
	$$
\end{enumerate}~~\dokendBem\\

\begin{Beis}\label{beis:1_3}
	\hspace*{0cm}\\\vspace{-1cm}
	\begin{enumerate}[(i)]
		\item $(\Z,+,0)$ bzw. $(\R,+,0)$ sind additiv geschriebene abelsche Gruppen, das "`{Neutralelement}"' 
		$0$ wird hierbei als Nullelement\index{Nullelement}\label{Nullelement} bezeichnet,
		und entsprechend das "`{inverse Element}"' $-a$ als die zu $a$ entgegengesetzte Zahl.
		Diese Begriffe verwendet man f\"ur alle additiv geschriebene abelsche Gruppen.\\
		
		\item F\"ur $\N_0=\{0,1,2,3,...\}$ ist $(\N_0,+,0)$ keine Gruppe, da (G3) verletzt ist.\\
		
		\item Die Menge aller $2\times 2$ Matrizen $M=\left( \begin{matrix}
		a & b  \\
		c & d\\
		\end{matrix} \right)$ mit $a,b,c,d\in\Z$ und Determinante $ad-bc=\pm1$ bildet die sogenannte Gruppe $GL(2,\Z)$ bzgl. der Matrizenmultiplikation \index{Matrizenmultiplikation}\label{Matrizenmultiplikation} $"`{\cdot}"'$ als Verkn\"upfung:
		$$
		\left( \begin{matrix}
		a & b  \\
		c & d\\
		\end{matrix} \right)\cdot\left( \begin{matrix}
		a' & b'  \\
		c' & d'\\
		\end{matrix} \right)=\left( \begin{matrix}
		aa'+bc' & ab'+bd'  \\
		ca'+dc' & cb'+dd'\\
		\end{matrix} \right).
		$$
		Dabei ist $\Det(M\cdot M')=\Det(M)\cdot \Det(M')\in \{+1,-1\}$ f\"ur $M, M'\in GL(2,\Z)$.\\
		Es gilt (G1), (G2) mit der Einheitsmatrix $\left( \begin{matrix}
		1 & 0  \\
		0 & 1\\
		\end{matrix} \right)$ als Neutralelement und (G3) mit der Inversen
		$$
		\left( \begin{matrix}
		a & b  \\
		c & d\\
		\end{matrix} \right)^{-1}=\dfrac{1}{ad-bc}
%		\begin{pmatrix*}[r]
%		d & -b  \\
%		-c & a\\
%		\end{pmatrix*}
		\left( \begin{matrix}
		d & -b  \\
		-c & a\\
		\end{matrix} \right)
		$$
		zu $\left( \begin{matrix}
		a & b  \\
		c & d\\
		\end{matrix} \right)\in GL(2,\Z)$. Es ist (G4) nicht erf\"ullt, z.B. 
		$$
		\begin{pmatrix*}[r]
		0 & 1  \\
		-1 & 0\\
		\end{pmatrix*} \cdot \begin{pmatrix*}[r]
		1 & 1  \\
		0 & 1\\
		\end{pmatrix*} =  \begin{pmatrix*}[r]
		0 & 1  \\
		-1 & -1\\
		\end{pmatrix*},\, \text{~aber~} \begin{pmatrix*}[r]
		1 & 1  \\
		0 & 1\\
		\end{pmatrix*} \cdot  \begin{pmatrix*}[r]
		0 & 1  \\
		-1 & 0\\
		\end{pmatrix*} = \begin{pmatrix*}[r]
		-1 & 1  \\
		-1 & 0\\
		\end{pmatrix*}.
		$$
		
		\item Die bijektiven Abbildungen $f:\N_n\to\N_n$ mit $\N_n=\{1,2,...,n\}$ und $n\in \N$ bilden bzgl. der Komposition $"`{\circ}"'$ dieser Abbildungen die Permuta\-tionsgruppe\index{Permutationsgruppe}\label{Permutationsgruppe}$\Sigma_n$ der Ordnung $|\Sigma_n|=n!$,
		siehe hierzu den entsprechenden Anhang \ref{cha:anhang2}.\\
		
		Die Identit\"at\index{Identit\"at}\label{Identitaet} $\Id : \N_n\to\N_n$ mit $\Id(k)=k$ f\"ur alle $k=1,...,n$ ist das Einselement von $\Sigma_n$, die inverse Abbildung $f^{-1}$ das zu $f\in\Sigma_n$ inverse Gruppenelement. Ab $n\geq3$ gilt nicht mehr (G4) f\"ur $\Sigma_n$.
	\end{enumerate}\dokendSatz
\end{Beis}
\begin{Def}[Untergruppe]\index{Untergruppe}\label{Untergruppe}\label{def:1_4}
	Es sei $(G,\circ,e)$ Gruppe und $U\subseteq G$ eine nicht leere Teilmenge von $G$. F\"ur \mbox{alle} \mbox{$a,b\in U$} gelte $a\circ b\in U$ sowie $a^{-1}\in U$. Dann hei{\ss}t $(U,\circ,e)$ Untergruppe von $(G,\circ,e)$. Wir sagen auch k\"urzer: $U$ ist Untergruppe von $G$. Hierbei ist $e\in U$ garantiert.~\dokendDef
\end{Def}
\begin{Beis}\label{beis:1_5}
	\hspace*{0cm}\\\vspace{-1cm}
	\begin{enumerate}[(i)]
		\item $(\Z,+,0)$ ist eine Untergruppe von $(\R,+,0)$.\\
		
		\item Die Menge aller Matrizen $M=\left( \begin{matrix}
		a & b  \\
		c & d\\
		\end{matrix} \right)\in GL(2,\Z)$ mit Determinante $$ad-bc=+1$$ liefert eine Untergruppe von $GL(2,\Z)$. Diese Untergruppe bezeichnet man mit $SL(2,\Z)$. Hierbei steht "`{G}"' f\"ur "`{general}"' und "`{S}"' f\"ur "`{special}"'.
	\end{enumerate}\dokendSatz
\end{Beis}

\begin{Satz}[Satz von Lagrange]\label{satz:1_6}\index{Satz von Lagrange}\label{Satz_von_Lagrange} Es sei $G$ eine Gruppe endlicher Ordnung und $U\subseteq G$ eine Untergruppe von $G$. Dann ist $\dfrac{|G|}{|U|}$ eine nat\"urliche Zahl, die auch Index \index{Index}\label{Index} von $U$ in $G$ genannt wird.\hfill\dokendSatz
\end{Satz}
{\bf Beweis:}~
Es seien $a,b,c\in G$. Wir definieren die Linksnebenklassen \index{Linksnebenklasse}\label{Linksnebenklasse}
$$a\circ U:=\{a\circ x:x\in U\},$$
 die wegen $e\in U$ ganz $G$ aussch\"opfen, und zeigen: Zwei Linksnebenklassen $a\circ U$, $b\circ U$ von $U$ sind entweder elementfremd oder identisch. Haben $a\circ U$ und $b\circ U$ ein Element $c=a\circ u_1=b\circ u_2$ mit $u_1, u_2\in U$ gemeinsam, so folgen $a=b\circ u_2\circ u^{-1}_1$ und $a\circ u=b\circ(u_2\circ u^{-1}_1\circ u)\in b\circ U$ f\"ur jedes $u\in U$, also $a\circ U\subseteq b\circ U$.
Durch Vertauschung  der Rollen von $a$ und $b$ folgt auch $b\circ U\subseteq a\circ U$, also $a\circ U= b\circ U$, wenn beide Linksnebenklassen nicht elementfremd sind. \\

Der Satz von Lagrange folgt nun aus Satz~\ref{satz:1_2} (c), der die eindeutige Aufl\"osbarkeit der Gleichungen $a\circ x=b$ garantiert, so dass jede Linksnebenklasse $a\circ U$ von $U$ genau $|U|$ Elemente besitzt.\dokendProof

\begin{Bem} Ist $G$ eine Gruppe endlicher Ordnung und $U$ eine Untergruppe von $G$, so zeigt obiger Beweis: Der Index $|G|/|U|$ von $U$ in $G$ stimmt mit der Anzahl aller Linksnebenklassen
	$$
	a\circ U=\{a\circ x: x\in U\},\quad (a\in G \text{ beliebig, aber fest})
$$
von $U$ \"uberein.
~~\dokendBem
\end{Bem}

Wir bilden nun die von einem Element $a$ der Gruppe $(G,\circ,e)$ 
erzeugten Potenzen $a^{k}$: Beginnend mit $a^0:=e$ definieren wir gem\"a{\ss} $a^{k+1}:=a\cdot a^k$ die Potenzen $a^k$ zun\"achst rekursiv f\"ur alle $k\in\N_0$, und dann mit 
$a^{-k}:=(a^{-1})^k$ auch f\"ur negative Exponenten $-k<0$. Da die Elemente $a,a^{-1},e$ in Mehrfachprodukten miteinander vertauschbar sind, gilt
\begin{equation}\label{eq:1_1}
 a^j\circ a^k=a^{j+k}\quad \mbox{f\"ur ~ alle ~}  j,k\in \Z.
\end{equation}	
	Wegen~\eqref{eq:1_1} ist
	\begin{equation}\label{eq:1_2} 
	U(a):=\{a^n:n\in\Z\}
	\end{equation}
	
	eine Untergruppe von $G$, die von $a$ erzeugte zyklische Untergruppe\index{zyklische Untergruppe}\label{zyklische_Untergruppe}. Wir nennen $|U(a)|$ die Ordnung von $a$ (in $G$).\\
	Wegen~\eqref{eq:1_1} ist $U(a)$ abelsche Untergruppe von $G$. Wir nehmen an, $G$ habe endliche Ordnung. Dann gibt es Exponenten $0\leq j<k$ mit $a^j=a^k$, und wegen~\eqref{eq:1_1} folgt hieraus $a^h=e$ f\"ur $h:=k-j\in\N$. Ist $h\geq1$ der kleinste nat\"urliche Exponent mit $a^h=e$ und $n\in\Z$, so gilt $n=k\cdot h+r$ mit $0\leq r\leq h-1$ f\"ur die gr\"o{\ss}te ganze Zahl $k\leq\frac{n}{h}$. Wir erhalten damit 
	$$
	a^n=a^{k\cdot h+r}=(a^h)^k\circ a^r=e\circ a^r=a^r,
	$$
	so dass unter Beachtung der Minimalit\"at von $h$ gilt:
	\begin{equation}
	U(a)=\{a^0,a^1,...,a^{h-1}\}, \quad |U(a)|=h.
\end{equation}

Nach Satz~\ref{satz:1_6}  ist $h$ ein Teiler von $|G|$. Somit gilt
\begin{Satz}\label{satz:1_8}$a^{|G|}=e$ f\"ur jedes $a$ aus einer endlichen Gruppe $(G,\circ,e)$.\hfill\dokendSatz
\end{Satz}
{\bf Beweis:}~
Nach Satz~\ref{satz:1_6} ist die Ordnung $h$ von $a$ ein Teiler von $|G|$. Es folgt
$$
a^{|G|}=(a^h)^{|G|/h}=e^{|G|/h}=e.
$$\dokendProof

	\section{Ringe}\label{cha:1B}

\begin{Def}[Ring]\label{def:1_9}
Eine algebraische Struktur $(R,+,\cdot)$ (oder kurz $R$) mit einer additiven Verkn\"upfung $+:R\times R\to R$ und einer multiplikativen Verkn\"upfung $\cdot :R\times R\to R$ hei{\ss}t ein Ring\index{Ring}\label{Ring}, wenn gilt:
\begin{enumerate}[(R1)]
\item $(R,+,\cdot)$ ist abelsche Gruppe\index{abelsche Gruppe}\label{abelsche_Gruppe2} mit dem Nullelement $0$ und dem zu $a\in R$ entgegengesetzten Element $-a$ mit $a+(-a)=0$.
\item $(a\cdot b)\cdot c=a\cdot(b\cdot c)$ \quad $\mbox{f\"ur ~ alle ~} a,b,c\in R$.
\item Es gelten die Distributivgesetze 
$$
a\cdot(b+c)=(a\cdot b)+(a\cdot c)=ab+ac \quad \text{sowie}
$$
$$
(b+c)\cdot a=(b\cdot a)+(c\cdot a)=ba+ca\quad \mbox{f\"ur ~ alle ~}  a,b,c\in R.
$$
\end{enumerate}
Das Zeichen "`{$\cdot$}"' bindet wie \"ublich st\"arker als "`{$+$}"' (Punkt- vor Strichrechnung) und wird nicht immer ausgeschrieben.\\

Gilt zus\"atzlich
\begin{enumerate}[(R4)]
	\item $ a\cdot b=b\cdot a\quad \mbox{f\"ur ~ alle ~}a,b\in R,$
\end{enumerate}
so wird der Ring kommutativ\index{kommutativer Ring}\label{kommutativer_Ring} genannt.  \dokendDef
\end{Def}

\begin{Bem}\label{bem:1_10}
\hspace*{0cm}\\\vspace{-1cm}
\begin{enumerate}[1)]
	\item Aus den Ringaxiomen (R1) bis (R3) folgert man m\"uhelos f\"ur alle $a,b,c \in R$ die Rechenregeln:
	\begin{enumerate}[$\bullet$]
		\item $a\cdot 0\underset{\text{(G2)}}{=}a\cdot (0+0)\underset{\text{(R3)}}{=}a\cdot 0+a\cdot 0$, und hieraus folgt 
		$a\cdot0=0$ nach Satz~\ref{satz:1_2} (c), angewendet auf die Gruppe $(R,+,0)$. Analog folgt \, $0\cdot a=0$ 
		\quad$\mbox{f\"ur~alle~} a\in R$.\\
		\item $0=a\cdot 0=a\cdot (b+(-b))=a\cdot b+a\cdot(-b)$, also $a\cdot (-b)=-a\cdot b$, und analog $(-a)\cdot b=-a\cdot b$.\\
		\item $a\cdot (b-c)=a\cdot (b+(-c))=a\cdot b+a\cdot(-c)=a\cdot b-a\cdot c$, und analog mit den Konventionen $b-c:=b+(-c)$ sowie "`{Punkt- vor Strichrechnung}"': 
		$(b-c)\cdot a=b\cdot a-c\cdot a$. 	\\
		\item $\left(\sum\limits_{j=1}^{n}a_j\right)\cdot\left(\sum\limits_{k=1}^{m}b_k\right)=\sum\limits_{j=1}^{n}\sum\limits_{k=1}^{m}a_j b_k$.
	\end{enumerate}

	\item Enth\"alt $R\neq\{0\}$ ein Element $1$ mit $1\cdot a=a\cdot 1=a$ 
	f\"ur alle $a\in R$, so nennen wir dieses Element Einselement \index{Einselement eines Ringes}\label{Einselement_eines_Ringes} von $R$. 
	Zwei Einselemente $1\neq 1'$ kann $R$ dann wegen $$1=1\cdot 1'=1'$$ nicht besitzen. 
	Da $a\cdot 0=0\cdot a=0$ f\"ur alle $a\in R$ gilt,
	ist \"uberdies $1\neq0$ garantiert. 
\end{enumerate}
~~\dokendBem
\end{Bem}

\begin{Def}[] %\label{def:}
	\hspace*{0cm}\\\vspace{-1cm}
	\begin{enumerate}[(a)]
	\item  Ein vom Nullring verschiedener kommutativer Ring $(R,+,\cdot,0)$ hei{\ss}t Integrit\"atsbereich\index{Integrit\"atsbereich}\label{Integritaetsbereich}, falls gilt:\\
	F\"ur alle $a,b\in R$ folgt aus $a\cdot b=0$ stets $a=0$ oder $b=0$.\\
	\item Ein kommutativer Ring $(R,+,\cdot,0,1)$ mit Einselement $1\neq0$ hei{\ss}t K\"orper\index{K\"orper}\label{Koerper}, wenn $(R\setminus\{0\},\cdot,1)$ (abelsche) Gruppe ist.
    \end{enumerate}\dokendDef
\end{Def}

\begin{Bem}	In einem Integrit\"atsbereich gilt die "`{K\"urzungsregel}"' $a\cdot b=a'\cdot b\Rightarrow$ $a=a'$ f\"ur alle $a,a',b\in R$ mit $b\neq0$, da man $a\cdot b=a'\cdot b$ nach Bemerkung~\ref{bem:1_10} $1)$ auch in der Form $(a-a')\cdot b=0$ schreiben kann.
~~\dokendBem
\end{Bem}

\begin{Beis}\label{beis:1_13}
\hspace*{0cm}\\\vspace{-1cm}	
	\begin{enumerate} [(a)] 
		\item Es ist $\Z$ mit der \"ublichen Addition $+$ und Multiplikation $\cdot$ ein Integrit\"atsbereich, ebenso 
		$$
		n\cdot\Z:=\{n\cdot k\,: \,k\in \Z\} \quad\text{f\"ur festes } n\in\N.
		$$
		Aber nur f\"ur $n=1$ ist $1\cdot\Z=\Z$ ein Integrit\"atsbereich mit (dem \"ublichen) Einselement $1$.\\
		\item Jeder K\"orper, wie z.B. $\Q, \R, \C$ mit den Grundrechenarten, ist auch ein Integrit\"atsbereich mit Einselement.	
	\end{enumerate}\dokendSatz
\end{Beis}

	\section{Aufgaben}\label{cha:1_A}
 
\noindent
{\bf Definition} zur Vorbereitung der Aufgabe~\ref{auf:1_1}:
\vspace{0.1cm}\\
Es seien $(G,\circ)$ und $(G',\circ')$ Gruppen
sowie $\varphi : G \to G'$ eine bijektive Abbildung.
Wir nennen die Abbildung $\varphi$ einen {\bf Isomorphismus}\index{Isomorphismus}\label{Isomorphismus}
zwischen den Gruppen $G$ und $G'$, wenn f\"ur alle
$a,b \in G$ folgendes gilt:
$$
\varphi(a \circ b) = \varphi(a) \circ' \varphi(b)\,.
$$
Die Gruppen $G$ und $G'$ hei{\ss}en in diesem Falle {\bf isomorph},
d.h. strukturgleich. Zur Bearbeitung der folgenden \"Ubungsaufgabe
beziehen wir uns auf die kurze Einf\"uhrung der Permutationsgruppen im Anhang \ref{cha:anhang2}.

\begin{Auf}[Permutationsgruppen]\label{auf:1_1}
%	\hspace*{0cm}\\\vspace{1cm}
Es sei $(G,\circ)$ eine beliebige Gruppe mit $|G|=n$ Elementen. 
Man zeige, dass $G$ dann einer Untergruppe der vollen
Permutationsgruppe\index{Permutationsgruppe}\label{Permutationsgruppe2} $\Sigma_n$ isomorph ist.\\

\underline{Hinweis:} Betrachte f\"ur beliebiges aber festes $b \in G$
die linksseitige Multi\-pli\-ka\-tion der Gruppenelemente 
$g_1,...,g_n$ von $G$ mit $b$.\\

\underline{Bemerkung:} \quad Isomorphe Gruppen unterscheiden sich nur hinsichtlich
der Bezeichnungsweise ihrer Elemente und ihrer Verkn\"upfung. Die Aufgabe 1
zeigt nun zus\"atzlich, dass die Untergruppen der Permutationsgruppen 
$\Sigma_n$ so allgemein sind, dass sie bereits alle endlichen Gruppen
beinhalten!\\	
\end{Auf}
{\bf L\"osung:}\\
Gegeben ist $G=\{g_1, g_2, ..., g_n\}$ mit $|G|=n$ Elementen. Wir zeigen: $G$ ist einer Untergruppe von $\Sigma_n$ isomorph.\\

Zun\"achst stellen nach Satz~\ref{satz:1_2} f\"ur festes $b\in G$ die $b\circ g_1$, $b\circ g_2$,..., $b\circ g_n$ eine Permutation der urspr\"unglichen $g_1$, $g_2$, ... , $g_n$ dar, d.h. es gibt zu jedem $b\in G$ eine Permutation $\pi_b\in\Sigma_n$ mit 
\begin{equation*}
b\circ g_j =g_{\pi_b(j)}\quad \mbox{f\"ur ~ alle ~ } j=1,...,n,
\end{equation*}
da in $G$ die Gleichung $b\circ g=a$ f\"ur alle $a,b\in G$ genau eine L\"osung $g$ besitzt, n\"amlich $g=b^{-1}\circ a$. Die Abbildung $\phi: G\rightarrow \Sigma_n$ mit $\phi(b):=\pi_b$ ist somit injektiv.\\

Betrachte $a,b\in G$. Dann gilt f\"ur alle $j=1,...,n$:
\begin{equation*}
\begin{array}{rcccl}
g_{\pi_{a\circ b}(j)} & = & (a\circ b)\circ g_j & = & a \circ(b\circ g_j) \\ 
                     & = & a\circ g_{\pi_b(j)}& =& g_{\pi_a(\pi_b(j))}\\
                     & = & g_{(\pi_a\circ\pi_b) (j)}\\
\end{array}
\end{equation*}
\begin{equation*}
\Rightarrow\quad\phi(a\circ b)=\pi_a\circ\pi_b=\phi(a)\circ\phi(b).
\end{equation*}
Die Untergruppe von $\Sigma_n$ ist das Bild
\begin{equation*}
\phi( G)=\{\phi(g):g\in G\}\subseteq\Sigma_n.
\end{equation*}
$\phi(G)$ ist Untergruppe von $\Sigma_n$ wegen 
\begin{equation*}
\phi(a\circ b^{-1})=\phi(a)\circ\phi(b)^{-1}\in\phi(G)\quad \mbox{f\"ur ~ alle ~ }  a,b\in G.
\end{equation*}

\begin{Auf}[Ein Ring mit Nullteilern]\label{auf:1_2}
Es werde $\mathcal{R}:=\left\{\,\begin{pmatrix}a & b\\ 0 & a\\  \end{pmatrix} \,:\,a,b \in \R\,\right\}$ mit der komponentenweisen Addition "`{$+$}"' zweier Matrizen und der \"ublichen Matrizenmultiplikation\index{Matrizenmultiplikation}\label{Matrizenmultiplikation2} "`{$\cdot$}"' versehen. Man zeige, dass dadurch ein kommutativer Ring mit Einselement entsteht, der kein Integrit\"atsbereich\index{Integrit\"atsbereich}\label{Integritaetsbereich2} ist. Hierzu bestimme man zwei Nullteiler\index{Nullteiler}\label{Nullteiler}, d.h. zwei von der Nullmatrix $\bf{0}$ verschiedene Matrizen $M, M' \in \mathcal{R}$ mit $M \cdot M' = \bf{0}$.\\

\end{Auf}
{\bf L\"osung:}\\
$\mathcal{R}:=\left\{\,\begin{pmatrix}a & b\\ 0 & a\\  \end{pmatrix} \,:\,a,b \in \R\,\right\}$ ist abgeschlossen unter den Rechenoperationen "`$+$"', "`$\cdot$"' im vollen Matrizenring $(\R^{2\times 2},+,\cdot)$, denn mit $\begin{pmatrix}a & b\\ 0 & a\\  \end{pmatrix}$, $\begin{pmatrix}a' & b'\\ 0 & a'\\  \end{pmatrix}\in\mathcal{R}$ folgt auch
\begin{equation*}
\begin{pmatrix}a & b\\ 0 & a\\  \end{pmatrix}+\begin{pmatrix}a' & b'\\ 0 & a'\\  \end{pmatrix}=\begin{pmatrix}a+a' & b+b'\\ 0 & a+a'\\  \end{pmatrix}\in\mathcal{R}
\end{equation*}
sowie 
\begin{equation*}
\begin{pmatrix}a & b\\ 0 & a\\  \end{pmatrix}\cdot \begin{pmatrix}a' & b'\\ 0 & a'\\  \end{pmatrix}=\begin{pmatrix}a' & b'\\ 0 & a'\\  \end{pmatrix}\cdot\begin{pmatrix}a & b\\ 0 & a\\  \end{pmatrix}=\begin{pmatrix}aa' & ab'+ba'\\ 0 & aa'\\  \end{pmatrix}\in\mathcal{R}
\end{equation*}
Die Matrizenmultiplikation ist bei Beschr\"ankung auf $\mathcal{R}$ kommutativ, auch ist sie assoziativ. $(\mathcal{R},+)$ ist abelsche Gruppe mit der entgegengesetzten Matrix
\begin{equation*}
-\begin{pmatrix}a & b\\ 0 & a\\  \end{pmatrix}=\begin{pmatrix}-a & -b\\ 0 & -a\\  \end{pmatrix}\in\mathcal{R}\quad \text{zu} \quad \begin{pmatrix}a & b\\ 0 & a\\  \end{pmatrix}\in\mathcal{R}
\end{equation*}
 und der Nullmatrix $\begin{pmatrix}0 & 0\\ 0 & 0\\  \end{pmatrix}\in\mathcal{R}$ als Nullelement. Die Distributivgesetze gelten schon allgemeiner in $(\R^{2\times 2},+,\cdot)$, und der Ring $\mathcal{R}$ hat $E=\begin{pmatrix}1 & 0\\ 0 & 1\\  \end{pmatrix}$ als Einselement. Da $\begin{pmatrix}0 & 1\\ 0 & 0\\  \end{pmatrix}\cdot \begin{pmatrix}0 & 1\\ 0 & 0\\  \end{pmatrix}=\begin{pmatrix}0 & 0\\ 0 & 0\\  \end{pmatrix}$ gilt, ist $\mathcal{R}$ kein Integrit\"atsbereich.\\
 
\hspace*{0cm}\\\vspace{0cm}	
{\bf Vorbereitung zur Bearbeitung der Aufgaben~\ref{auf:1_3} und \ref{auf:1_4}:}\\

Hier empfehlen wir f\"ur den Einstieg den ersten Teil des Anhangs \ref{cha:anhang1}
zu logischen Symbolen, Mengen und Abbildungen zu studieren.\\
Die Aussageform $\mathcal{A}(n)$ ordne jedem $n \in \N$ 
einen Wahrheitswert "`{wahr}"' oder "`{falsch}"' zu. Dann gilt das folgende
{\bf Induktionsprinzip:}\index{Induktionsprinzip}\label{Induktionsprinzip}\\

Wenn der Induktionsanfang $\mathcal{A}(1)$ wahr ist
und der Induktions\-schluss
\begin{equation*}
\mathcal{A}(n) \Rightarrow \mathcal{A}(n+1) 
\end{equation*}
f\"ur alle $n \in \N$ gilt, 
dann folgt bereits $\mathcal{A}(n)$ f\"ur alle $n \in \N$.\\

\begin{Auf}[Vollst\"andige Induktion\index{vollst\"andige Induktion}\label{vollstaendige_Induktion}]\label{auf:1_3}
	
Aus dem vorigen Induktionsprinzip sollen zwei Varianten 
hergeleitet werden.
\begin{enumerate}[(a)]
	\item Es sei  $\mathcal{B}(n)$ eine Aussageform f\"ur die nat\"urlichen Zahlen $n$
	und es be\-zeich\-ne $\N_n :=\{1,2,\ldots,n\}$ die Menge der ersten $n$ nat\"urlichen Zahlen.	
	Man zeige:\\
	
	Wenn der Induktionsanfang $\mathcal{B}(1)$ wahr ist
	und zudem
	\begin{equation*}
	\left(~\forall k \in \N_n \, : \, \mathcal{B}(k)  ~\right)
	\Rightarrow \mathcal{B}(n+1) 
	\end{equation*}
	f\"ur alle $n \in \N$ gilt, dann folgt $\mathcal{B}(n)$ f\"ur alle $n \in \N$.\\
	
	\item Es sei $k_0 \in \Z$ fest gew\"ahlt.
	Die Aussageform $\mathcal{B}(k)$ ordne jeder ganzen Zahl $k \geq k_0$ einen
	Wahrheitswert "`{wahr}"' oder "`{falsch}"' zu. Man zeige:\\
	
	Wenn der Induktionsanfang $\mathcal{B}(k_0)$ wahr ist
	und f\"ur alle ganzen Zahlen $k \geq k_0$ der Induktions\-schluss
	$
	\mathcal{B}(k) \Rightarrow \mathcal{B}(k+1)
	$
	gilt, dann folgt $\mathcal{B}(k)$ f\"ur alle ganzen Zahlen $k \geq k_0$.
\end{enumerate}	
\end{Auf}
{\bf L\"osung:}\\
Wir verwenden das eingangs formulierte Induktionsprinzip:
\begin{enumerate}[(a)]
	\item Es gelte $\mathcal{B}(1)$ und f\"ur alle $n\in\N$:
	\begin{equation}\label{eq:1_4}
	\mathcal{A}(n) \Rightarrow \mathcal{B}(n+1)
	\end{equation}
	mit $\mathcal{A}(n):\Leftrightarrow\bigwedge\limits_{k=1}^{n}\mathcal{B}(k) \Leftrightarrow \left(~\forall k \in \N_n \, : \, \mathcal{B}(k) ~\right)$ f\"ur $n\in\N$. Es gilt $\mathcal{A}(1)$ wegen $\mathcal{B}(1)$, und nach Definition von $\mathcal{A}$ f\"ur alle $n\in\N$ die \"Aquivalenz
	\begin{equation*}
	\mathcal{A}(n+1) \Leftrightarrow  \left(~\mathcal{A}(n) \wedge \mathcal{B}(n+1)  ~\right),	
	\end{equation*}
	so dass $\mathcal{A}(n+1)$ wegen \eqref{eq:1_4} f\"ur alle $n\in\N$ aus $\mathcal{A}(n)$ folgt. 
	Nach dem Induktionsprinzip gelten dann $\mathcal{A}(n)$ sowie 
	$\mathcal{B}(n)$ f\"ur alle $n\in\N$.\\
	
	\item folgt einfach, indem man die Aussageform $\mathcal{B}(k)$ durch die Aussageform \linebreak\mbox{$\mathcal{A}(n):\Leftrightarrow\mathcal{B}(k_0+n-1)$} mit $n\in\N$ ersetzt und dann auf $\mathcal{A}(n)$ Induktion anwendet.
\end{enumerate}

\begin{Auf}[Fibonacci-Folge, Teil 1]\label{auf:1_4}

Die Folge $(f_n)_{n \in \N_0}$ der Fibonacci-Zahlen\index{Fibonacci-Zahlen}\label{Fibonacci-Zahlen} ist rekursiv definiert
durch die beiden Anfangswerte $f_0=0$, $f_1 = 1$ sowie f\"ur alle $n \in \N_0$
durch die Rekursionsbeziehung $f_{n+2}=f_{n+1}+f_n$.
Zus\"atzlich definieren wir noch $f_{-1}:=1$.
\begin{enumerate}[(a)]
	\item Man zeige induktiv f\"ur alle $n \in \N_0$:
	$\begin{displaystyle}
	\begin{pmatrix}
	1 & 1\\
	1 & 0\\
	\end{pmatrix}^n =
	\begin{pmatrix}
	f_{n+1} & f_n\\
	f_n & f_{n-1}\\
	\end{pmatrix} \,, 
	\end{displaystyle}$\\
	und damit $f_{n+1}f_{n-1}-f_{n}^2=(-1)^n\,.$\\
	
	\item Mit den Eigenwerten $\begin{displaystyle}\lambda_{\pm} 
	:= \frac{1\pm \sqrt{5}}{2}\end{displaystyle}$ der Matrix
	$\begin{displaystyle}A := \begin{pmatrix}
	1 & 1\\
	1 & 0\\
	\end{pmatrix}\end{displaystyle}$
	und mit den Eigenvektoren 
	$\begin{displaystyle} \underline{x}_{\pm} = \begin{pmatrix}
	1 \\
	-\lambda_{\mp} \\
	\end{pmatrix}\end{displaystyle}$ zu den Eigenwerten $\lambda_{\pm}$ zeige man
	\begin{equation*}
	f_{n+1}=\lambda_{+}f_n+\lambda_{-}^n 
	=\lambda_{-}f_n+\lambda_{+}^n\quad  \mbox{f\"ur ~ alle ~ }  n \in \N_0\,.
	\end{equation*}
	
	\item Aus (b) leite man die Binetsche Formel\index{Binetsche Formel}\label{Binetsche_Formel} her:
	\begin{equation*}
	f_{n} =
	\frac{1}{\sqrt{5}}\,\left[\left(\frac{1+\sqrt{5}}{2}\right)^n
	-
	\left(\frac{1-\sqrt{5}}{2}\right)^n
	\, \right]\,~\mbox{f\"ur ~ alle ~ }  n \in \N_0\,.
	\end{equation*}
	
	\item Man zeige f\"ur alle $x \in \R$ mit $|x|<\frac{\sqrt{5}-1}{2}$:
	\begin{equation*}
	\sum \limits_{n=0}^{\infty}f_n x^n = \frac{x}{1-x-x^2}
	\end{equation*}
	mit absoluter Konvergenz der linksstehenden Reihe.
\end{enumerate}	
\end{Auf}
{\bf L\"osung:}\\
\begin{enumerate}[(a)]
	\item Wir zeigen induktiv:
	\begin{equation}\label{eq:1_5}
	\begin{pmatrix}
	1 & 1\\
	1 & 0\\
	\end{pmatrix}^n =
	\begin{pmatrix}
	f_{n+1} & f_n\\
	f_n & f_{n-1}\\
	\end{pmatrix}\quad \mbox{f\"ur~alle~} n\in\N_0.
	\end{equation}
	{\it Induktionsanfang}: 
	\begin{equation*}
	\begin{pmatrix}
	1 & 1\\
	1 & 0\\
	\end{pmatrix}^0=\begin{pmatrix}
	1 & 0\\
	0 & 1\\
	\end{pmatrix}=
	\begin{pmatrix}
	f_{1} & f_0\\
	f_0 & f_{-1}\\
	\end{pmatrix}.
	\end{equation*}
	{\it Induktionsannahme}: F\"ur ein $n\in\N_0$ sei $\begin{pmatrix}
	1 & 1\\
	1 & 0\\
	\end{pmatrix}^n=\begin{pmatrix}
	f_{n+1} & f_n\\
	f_n & f_{n-1}\\
	\end{pmatrix}$ bereits gezeigt. Dann folgt
		\begin{equation*}
		\begin{pmatrix}
		1 & 1\\
		1 & 0\\
		\end{pmatrix}^{n+1}=
		\begin{pmatrix}
		f_{n+1} & f_n\\
		f_n & f_{n-1}\\
		\end{pmatrix}\begin{pmatrix}
		1 & 1\\
		1 & 0\\
		\end{pmatrix}=\begin{pmatrix}
		f_{n+1}+f_n &\ f_{n+1}\\
		f_n+f_{n-1} &\  f_{n}\\
		\end{pmatrix}=\begin{pmatrix}
		f_{n+2} & f_{n+1}\\
		f_{n+1} & f_{n}\\
		\end{pmatrix},
		\end{equation*}
		wobei noch $f_1=f_0+f_{-1}$ zu beachten ist. Damit folgt \eqref{eq:1_5}. Aus dem Multiplikationssatz f\"ur Determinanten und \eqref{eq:1_5} folgt $(-1)^n=f_{n+1}f_{n-1}-f^2_n$ f\"ur alle $n\in\N_0$.\\
		
		\item Mit $\lambda_{\pm}=\frac{1\pm \sqrt{5}}{2}$ gilt 
		$\lambda_+ +\lambda_-=1$, $\lambda_+\cdot\lambda_-=-1$. Hieraus folgt f\"ur $\begin{displaystyle} \underline{x}_{\pm} = \begin{pmatrix}
		1 \\
		-\lambda_{\mp} \\
		\end{pmatrix}\end{displaystyle}$:
		\begin{equation*}
		A\,\underline{x}_{\pm}=
		\begin{pmatrix}
		1 & 1\\
		1 & 0\\
		\end{pmatrix}\underline{x}_{\pm}=\begin{pmatrix}
		1-\lambda_{\mp} \\
		1
		\end{pmatrix}=\begin{pmatrix}
		\lambda_{\pm} \\
		1\\
		\end{pmatrix}=\lambda_\pm\begin{pmatrix}
		1 \\
		-\lambda_{\mp} \\
		\end{pmatrix}=\lambda_\pm\underline{x}_{\pm},
		\end{equation*}
		und somit aus (a) f\"ur alle $n\in\N_0$:
		\begin{equation*}
		A^n\begin{pmatrix}
		1 \\
		-\lambda_{\mp} \\
		\end{pmatrix}=\begin{pmatrix}
		f_{n+1} & f_{n}\\
		f_n & f_{n-1}\\
		\end{pmatrix}\begin{pmatrix}
		1 \\
		-\lambda_{\mp} \\
		\end{pmatrix}=\begin{pmatrix}
		f_{n+1}-\lambda_\mp f_n\\
		f_n-\lambda_{\mp}f_{n-1}\\
		\end{pmatrix}=\lambda_{\pm}^n\begin{pmatrix}
		1\\
		-\lambda_\mp\\
		\end{pmatrix}.
		\end{equation*}
		Die Betrachtung der ersten Komponenten liefert $f_{n+1}=\lambda_\mp f_n+\lambda_\pm^n$.\\
		
		\item Aus (b) folgt $\lambda_+ f_n+\lambda_-^n=\lambda_- f_n+\lambda_+^n$, also wegen $\lambda_+-\lambda_-=\sqrt{5}$ und wegen $f_n=\dfrac{\lambda_+^n-\lambda_-^n}{\lambda_+ -\lambda_-}$ die Binetsche Formel.\\
		
		\item Aus der Binetschen Formel folgt
		\begin{equation*}
		\lim\limits_{n\rightarrow\infty}\sqrt[n]{f_n}=\lambda_+=\frac{1+\sqrt{5}}{2}\quad \text{mit} \quad \frac{1}{\lambda_+}=\frac{\sqrt{5}-1}{2}.
		\end{equation*}
		Somit ist $R=\frac{\sqrt{5}-1}{2}$ der Konvergenzradius der Potenzreihe $\sum\limits_{n=0}^\infty f_n x^n$, die f\"ur \mbox{$|x|<R$} absolut konvergiert. Es folgt f\"ur $|x|<R$:
		\begin{equation*}
		\begin{split}
		&(1-x-x^2)\sum\limits_{n=0}^\infty f_n x^n=
		\sum\limits_{n=0}^\infty (f_n x^n-f_nx^{n+1}-f_nx^{n+2})\\
		=&\sum\limits_{k=0}^\infty f_k x^k-\sum\limits_{k=1}^\infty f_{k-1} x^k-\sum\limits_{k=2}^\infty f_{k-2} x^k\\
		=&f_0 x^0+f_1 x^1-f_0 x^1 +\sum\limits_{k=2}^\infty (f_k-f_{k-1}-f_{k-2})x^k=x\,,\\
		\end{split}
		\end{equation*}
		wobei der letzte Schritt aus der Rekursionsformel $f_k=f_{k-1}+f_{k-2}$ folgt.
\end{enumerate}

\chapter{Euklidischer Algorithmus und Funda\-men\-tal\-satz der Arithmetik}\label{cha:2}
Die Berechnung des gr\"o{\ss}ten gemeinsamen Teilers zweier nat\"urlicher Zahlen 
mit Hilfe des Euklidischen Algorithmus\index{Euklidischer Algorithmus}\label{Euklidischer_Algorithmus} geht bis in die Antike zur\"uck.
Wie wir noch sehen werden, reicht die Bedeutung  des Euklidischen Algorithmus
weit \"uber diese einfache Aufgabenstellung hinaus. In dieser Lektion f\"uhren wir zun\"achst den
einfachen Euklidischen Algorithmus mit dem Ziel ein, den Fundamentalsatz der Arithmetik\index{Fundamentalsatz der Arithmetik}\label{Fundamentalsatz_der_Arithmetik}
zu beweisen. Dieser besagt, dass sich jede na\"urliche Zahl gr\"o{\ss}er als $1$ 
abgesehen von der Reihenfolge der Faktoren eindeutig in ein Produkt von Primzahlen\index{Primzahl}\label{Primzahl} zerlegen l\"asst.\\
	\section{Euklidischer Algorithmus}\label{cha:2A}

\begin{Def}[Gau{\ss}-Klammer]\label{def:2_1}
	Die Gau{\ss}-Klammer\index{Gau{\ss}-Klammer} \label{Gaussklammer}$\lfloor x\rfloor:=\max\{k\in\Z:k\leq x\}$ einer reellen Zahl $x$ bezeichnet die gr\"o{\ss}te ganze Zahl $k\leq x$. \\
        Die Gau{\ss}-Klammer einer reellen Zahl $x$ ist somit 
        diejenige ganze Zahl $k$, die durch die Ungleichungskette 
	\begin{equation}\label{eq:2_1}
	k\leq x<k+1
	\end{equation}
	eindeutig bestimmt ist.
  \dokendDef
\end{Def}

\begin{Bem}\label{bem:2_2}
	\hspace*{0cm}\\\vspace{-1cm}
	\begin{enumerate}[1)]
		\item Die Gau{\ss}-Klammer l\"asst die ganzen Zahlen unver\"andert, die nicht ganzen Zahlen werden dagegen stets abgerundet, z. B. ist 
	\begin{equation*}
	\lfloor 0.75\rfloor=0\quad \text{und} \quad\lfloor -0.5\rfloor=-1.
	\end{equation*}
		\item Entsprechend definiert man $\lceil x\rceil:=\min\{k\in\Z:k\geq x\}$ f\"ur $x\in\R$ durch Aufrunden, wobei $\lceil x\rceil=-\lfloor -x\rfloor$ gilt. 
	\end{enumerate}
	~~\dokendBem
\end{Bem}

\begin{center}
	\bf{Graphische Darstellung der Gau{\ss}-Klammer}
\end{center}
\begin{figure}
	\begin{center}
		\unitlength=0.9cm
		\begin{picture}(6,6)(-3,-3)
		\put(-3,0){\vector(1,0){6}}
		\put(3.1,-.2){\makebox(0,0)[t]{$x$}}
		\put(0,-3){\vector(0,1){6}}
		\put(-.2,3){\makebox(0,0)[r]{$y$}}
		\multiput(-3,-0.1)(1,0){6}{\line(0,1){0.2}}
		\multiput(-0.1,-3)(0,1){6}{\line(1,0){0.2}}
		\thicklines
		\multiput(-3,-3)(1,1){6}{\line(1,0){1}}
		\multiput(-3,-3)(1,1){6}{\circle*{.1}}
		\put(-3,-.2){\makebox(0,0)[t]{-3}}
		\put(-2,-.2){\makebox(0,0)[t]{-2}}
		\put(-1,-.2){\makebox(0,0)[t]{-1}}
		\put(.2,-.2){\makebox(0,0)[tl]{0}}
		\put(1,-.2){\makebox(0,0)[t]{1}}
		\put(2,-.2){\makebox(0,0)[t]{2}}
		\put(.2,-3){\makebox(0,0)[l]{-3}}
		\put(.2,-2){\makebox(0,0)[l]{-2}}
		\put(.2,-1){\makebox(0,0)[l]{-1}}
		\put(.2,1){\makebox(0,0)[l]{1}}
		\put(.2,2){\makebox(0,0)[l]{2}}
		\end{picture}
		\caption{Graphische Darstellung der Funktion $y = \lfloor x \rfloor$}
	\end{center}
\end{figure}
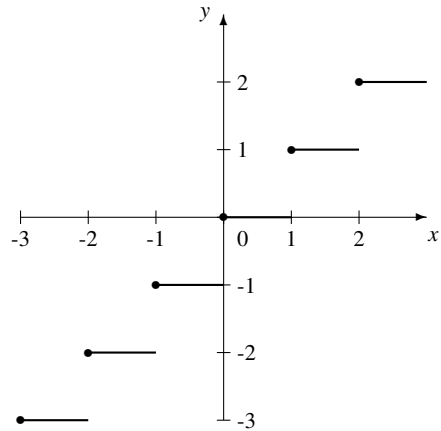
\begin{Def}[Teiler ]\label{def:2_3}
      \index{Teiler}\label{Teiler}
	\hspace*{0cm}\\\vspace{-1cm}	
	\begin{enumerate}[(a)]
		\item Es seien $d,k\in\Z$ mit $d\neq 0$. 
    Wir nennen $d$ einen Teiler von $k$ und schreiben daf\"ur $d|k$, 
wenn es ein $m\in\Z$ gibt mit $k=m\cdot d$. Es ist dann $m=\frac{k}{d}$ ganzzahlig.\\
		\item Es seien $a$ und $b$ ganze Zahlen, die nicht beide Null sind. 
Dann bezeichnen wir mit $\ggT(a,b)$ den gr\"o{\ss}ten gemeinsamen Teiler 
\index{gr\"o{\ss}ter gemeinsamer Teiler}\label{groesster_gemeinsamer _Teiler} von $a$ und $b$.
		Im Falle $\ggT(a,b)=1$ nennen wir $a$ und $b$ teilerfremd.
		\index{teilerfremde Zahlen}\label{teilerfremde_Zahlen}
	\end{enumerate}
	\dokendDef
\end{Def}
{\it Bemerkung}:~
Da $a$ und $b$ nicht beide verschwinden, gilt $|d|\leq\max(|a|,|b|)$ f\"ur jeden gemeinsamen Teiler $d\in\Z\setminus\{0\}$ von $a$ und $b$. Zudem ist $1$ ein gemeinsamer nat\"urlicher Teiler von $a$ und $b$. Somit ist die Menge aller gemeinsamer Teiler von $a$ und $b$ endlich und $\ggT(a,b)$ eine wohldefinierte nat\"urliche Zahl.\\

Zur Berechnung von $\ggT(a,b)$ beginnen wir mit dem
\begin{Satz}\label{satz:2_4}F\"ur je zwei Zahlen $a\in\Z$ und $b\in\N$ hat man eine eindeutige Darstellung der Form $a=q\cdot b+r$ mit $q\in\Z$ und $0\leq r<b$. Hierbei gilt $q=\left\lfloor\frac{a}{b}\right\rfloor$ .\hfill\dokendSatz
\end{Satz}
{\it Bemerkung}:~
Der Satz beschreibt die Division von $a$ durch $b$ 
mit Hilfe des Divi\-sions\-koef\-fizienten \index{Divisionskoeffizient}\label{Divisionskoeffizient} $q=\left\lfloor\frac{a}{b}\right\rfloor$ 
und des Divisions\-restes \index{Divisionsrest}\label{Divisionsrest} $r\in\N_0$.\\

{\bf Beweis des Satzes~\ref{satz:2_4}:}~
Wir zeigen zuerst, dass $q:=\left\lfloor\frac{a}{b}\right\rfloor$ und $r:=a-q\cdot b$ eine gew\"unschte Darstellung liefern:\\
Aus der Definition~\ref{def:2_1} der Gau{\ss}-Klammer folgt, siehe dort \eqref{eq:2_1}:

\begin{equation}\label{eq:2_2}
\left\lfloor\frac{a}{b}\right\rfloor\leq \frac{a}{b}<\left\lfloor\frac{a}{b}\right\rfloor+1\,.
\end{equation}
Die linke Ungleichung von \eqref{eq:2_2} ergibt
$$
r=a-\left\lfloor\frac{a}{b}\right\rfloor\cdot b\geq a-\frac{a}{b}\cdot b=0,
$$
und die rechte Ungleichung
$$
r<a-\left(\dfrac{a}{b}-1\right)\cdot b=b.
$$
Schlie{\ss}lich folgt eindeutig f\"ur jede Darstellung
$a=q'\cdot b+r'$ mit $q'\in\Z$ und \\$0\leq r'<b$:
$$
\left\lfloor\frac{a}{b}\right\rfloor=\left\lfloor q'+\frac{r'}{b}\right\rfloor=q'+\left\lfloor\frac{r'}{b}\right\rfloor=q'=q,
$$
$$
r'=a-\left\lfloor\frac{a}{b}\right\rfloor\cdot b=r.
$$\dokendProof

\begin{Satz}\label{satz:2_5}
	Die ganzen Zahlen $a$, $b$ m\"ogen nicht beide verschwinden. Dann gelten die folgenden Aussagen:
\begin{enumerate}[(a)]
	\item $\ggT(a,b)=\ggT(b,a)$.\\
	
	\item Die gemeinsamen Teiler des Zahlenpaares $a$, $b$ sind dieselben wie die des Zahlenpaares $a$, $|b|$. 
	\begin{center}
		Insbesondere gilt $\ggT(a,b)=\ggT(a,|b|)$.
	\end{center}	
	\item F\"ur $b\in\N$ setzen wir $r:=a-\left\lfloor\frac{a}{b}\right\rfloor\cdot b$. Dann sind die gemeinsamen Teiler des Zahlenpaares $a$, $b$ dieselben wie die des Zahlenpaares $b$, $r$.
	\begin{center}
		Insbesondere gilt $\ggT(a,b)=\ggT(b,r)$.
	\end{center}	
	\item F\"ur $b\in\N$ sind die gemeinsamen Teiler des Zahlenpaares $b$, $0$ genau die Teiler von $b$, und insbesondere gilt
	$$
	\ggT(b,0)=b.
	$$
\end{enumerate}\hfill\dokendSatz
\end{Satz}
{\bf Beweis:}~
Die Aussagen (a), (b) und (d) ergeben sich direkt aus Definition~\ref{def:2_3}, so dass wir uns auf den Nachweis von (c) beschr\"anken k\"onnen:\\

Es sei $d$ ein gemeinsamer Teiler von $a$ und $b$, $a=d \tilde a$,  $b=d \tilde b$ mit $\tilde a,\tilde b\in\Z$. Dann gilt auch 
$r=a-\left\lfloor a/b \right\rfloor b = 
d \left(\tilde a-\left\lfloor a/b \right\rfloor \tilde b\right)$ mit $d|r$.\\

Ist umgekehrt $d'$ gemeinsamer Teiler von $b$ und $r$, so gilt 
$d'| (r+\left\lfloor a/b \right\rfloor  b)$, d.h. $d'|a$. \\

Insbesondere stimmt der gr\"o{\ss}te gemeinsame Teiler von $a$ und $b$ mit dem von $b$ und $r$ \"uberein.
\hfill\dokendProof\\

\hspace*{0cm}\\\vspace{1cm}
\begin{center}
	\textbf{Formulierung des Euklidischen Algorithmus}\index{Euklidischer Algorithmus}\label{Euklidischer_Algorithmus2}
\end{center}
Die S"atze~\ref{satz:2_4} und \ref{satz:2_5} 
bilden das Fundament f\"ur den Euklidischen Algorithmus 
zur Berechnung von $\ggT(a,b)$:\\

Gegeben ist ein Paar $a''$, $b''$ ganzer Zahlen, nicht beide Null, 
wobei auf deren Reihenfolge zu achten ist. 
Wir beginnen mit zwei Startschritten gem\"a{\ss} Satz~\ref{satz:2_5}~(a),~(b):

\begin{enumerate}[\text{Schritt} 1:]
	\item Wir ersetzen das Paar $a''$, $b''$ durch $a'$, $b'$ 
         mit $a'=a''$, $b'=b''$ f\"ur $b''\neq0$ bzw. $a'=b''$, $b'=a''$ 
         f\"ur $b''=0$.\\	
	\item Wir ersetzen das Paar $a'$, $b'$ durch $a$, $b$ mit $a=a'$, $b=|b'|$. 
         Nun ist $\ggT(a'',b'')=\ggT(a,b)$ mit $a\in\Z$, $b\in\N$. 
         F\"ur den Startindex $j=1$ definieren wir nun 
         das aktuelle Zahlenpaar $r_{j-1}$, $r_j$ gem\"a{\ss} $r_{j-1}=r_0=b$, 
         $r_j=r_1=a-b\left\lfloor\frac{a}{b}\right\rfloor$. 
         F\"ur sp\"atere Zwecke geeignet sei $q_0=\left\lfloor\frac{a}{b}\right\rfloor$.\\
	\item Solange $r_j>0$ bleibt, berechnen wir rekursiv, 
        beginnend mit $j=1$, die Gr\"o{\ss}en
	$$
	q_j=\left\lfloor\frac{r_{j-1}}{r_j}\right\rfloor,\quad r_{j+1}=r_{j-1}-q_jr_j.
	$$	
\end{enumerate}\dokendSatz

Dann gilt nach Satz~\ref{satz:2_4}: Das Schema aus Schritt 3 endet f\"ur einen Abbruchindex $n_*\geq 1$ 
mit dem Divisionsrest $r_{n_*}=0$, denn die Folge der Divisionsreste
$r_j$ nimmt in jedem Schritt echt ab:
$$
r_0>r_1> ... >r_{n_*-1}>r_{n_*}=0.
$$
Nach Satz~\ref{satz:2_5}(c) gilt
$$
\ggT(a,b)=\ggT(r_{j-1},r_j)=r_{n_*-1}\quad \mbox{f\"ur~} j=1,...,n_*,
$$
und zudem stimmen die gemeinsamen Teiler von $a$, $b$ mit den Teilern von $r_{n_*-1}=\ggT(a,b)$ \"uberein. 
Somit gilt der
\begin{Satz}\label{satz:2_6}
	F\"ur $a,b\in\Z$ mit $|a|+|b|>0$ haben wir $d|\ggT(a,b)$ f\"ur jeden gemeinsamen Teiler $d$ von $a$, $b$.
%	\hfill
	\dokendSatz
\end{Satz}

%\hspace*{0cm}\\\vspace{1cm}
\begin{Beis}\label{beis:2_7}
	Wir bestimmen $\ggT(138,462)$ f\"ur $a=138$ und $b=462$, und 
        starten den Algorithmus mit 
        $r_0=b=462$, $\begin{displaystyle}q_0 
        = \left\lfloor \frac{a}{b} \right\rfloor\end{displaystyle}=0$
        sowie $r_1=a-q_0b=138$ unter Beachtung von
        $\ggT(138,462)=\ggT(r_0,r_1)=\ggT(462,138)$:
	\begin{center}
		\begin{minipage}[l]{5cm}
			\begin{align*}
			\\
			462 & = 3 \cdot 138 + 48\\
			138 & = 2 \cdot 48 + 42\\
			48  & = 1 \cdot 42 + 6\\
			42  & = 7 \cdot 6 + 0\, .
			\end{align*}
		\end{minipage}
		\begin{minipage}[r]{5cm}
			\begin{align*}
			\text{Schlusskette:}&\ \ggT(462, 138)\\
			=&\ \ggT(138, 48)\\
			=&\ \ggT(48, 42)\\
			=&\ \ggT(42, 6)\\
			=&\ \ggT(6, 0) = 6 \, .
			\end{align*}
		\end{minipage}
	\end{center}
	
	Da $6$ der letzte von $0$ verschiedene Divisionsrest \index{Divisionsrest}\label{Divisionsrest2}ist, 
        folgt $\ggT(138,462) = 6$.\\

	Wir stellen den Algorithmus als einfaches Rechenschema dar:
 
        Der Abbruchindex ist $n_*=5$, und f\"ur $1 \leq j < n_*-1=4$
        sind die Divisions\-koeffi\-zien\-ten \index{Divisionskoeffizient}\label{Divisionskoeffizient2}
        $q_{j} =\left \lfloor \frac{r_{j-1}}{r_j} \right \rfloor$
        mit der Rekursion $	r_{j+1}  =  r_{j-1}-q_{j}r_j$
        der Divi\-sions\-reste erkl\"art:

	\begin{equation*}
		\begin{tabular}{|l||c|c|} \hline
			$j$  & ~~ $q_j$ ~~  & ~~~$r_j$~~~\\
			\hline
			0  &    0    &  462 \\ \hline 
			1  &    3    &  138 \\ \hline
			2  &    2    &  48  \\ \hline
			3  &    1    &  42  \\ \hline
			4  &    7    &  6   \\ \hline
			5  &    ---  &  0   \\ \hline
		\end{tabular}\quad\quad\quad
		\begin{array}{rcccl}
		q_{0}& = & \left \lfloor 138/462 \right \rfloor\,& = & 0\,, \\
		q_1& = & \left \lfloor 462/138 \right \rfloor\,& = & 3\,, \\
		q_2& = & \left \lfloor 138/48 \right \rfloor\,& = & 2\,, \\
		q_3& = & \left \lfloor 48/42 \right \rfloor\,& = & 1\,, \\
		q_4& = & \left \lfloor 42/6 \right \rfloor\,& = & 7.\, \\
		\end{array}
	\end{equation*}

	\dokendSatz
\end{Beis}
In diesem Beispiel durchlaufen wir nun, 
beginnend mit der Darstellung des gr\"o{\ss}ten gemeinsamen Teilers 
im vorletzten Schritt, den Euklidischen Algorithmus in umgekehrter Reihenfolge, 
indem wir schrittweise den kleinsten auftretenden Rest $r_{j+1}$ 
mit dem gr\"o{\ss}ten Index $j+1$ durch $r_{j-1}-q_jr_j$ ersetzen. 
Auf diese Weise erhalten wir
$$
\begin{array}{rclcl}
6& = & 48-1\cdot 42=48-1\cdot(138-2\cdot 48)\\
& = &-1\cdot138+3\cdot48=-1\cdot138+3\cdot(462-3\cdot138) \\
& = & -10\cdot138+3\cdot462\,, %\text{ also:}
\end{array}
$$
also:
$$\ggT(138,462)=6=-10\cdot138+3\cdot462.
$$
F\"uhrt man dieses Verfahren allgemein durch, so erh\"alt man den
\begin{Satz}\label{satz:2_8}
	F\"ur $a,b\in\Z$ mit $|a|+|b|>0$ gibt es ganze Zahlen $\lambda$, 
$\mu$ mit $\ggT(a,b)=\lambda a+\mu b$.\hfill
	
	\dokendSatz
\end{Satz}
{\bf Beweis:}~
Die Menge $(a,b):=\{xa+yb:x,y\in\Z\}$ bildet einen eigenst\"andigen Unterring von $(\Z,+,\cdot,0)$, 
der $a$ und $b$ enth\"alt. Es ist nicht $(a,b)=\{0\}$, 
und folglich existiert die kleinste positive Zahl $g=\lambda a+\mu b$ in $(a,b)$ 
mit Koeffizienten $\lambda,\mu \in\Z$. Aus der Darstellung von $g$ folgt sofort
\begin{equation}\label{eq:2_3}
\ggT(a,b)\,|\, g.
\end{equation}
Nach Satz~\ref{satz:2_4} gibt es ganze Zahlen $q,r$ mit $a=q\cdot g+r$ und $0\leq r<g$. 
Mit \mbox{$a,g\in(a,b)$} ist aber auch $r=a-q\cdot g\in (a,b)$, und da $g$ minimal ist, 
folgt $r=0$, d.h. $a=q\cdot g$.
Entsprechend erhalten wir $b=q'\cdot g$ mit passendem $q'\in\Z$:
\begin{equation}\label{eq:2_4}
a=q\cdot g,\quad b=q'\cdot g.
\end{equation}
Aus Satz~\ref{satz:2_6} und \eqref{eq:2_4} folgt nun
$\begin{displaystyle}
g\, | \,\ggT(a,b)\,,
\end{displaystyle}$
und zusammen mit \eqref{eq:2_3}:
\begin{center}
	$\ggT(a,b)=g=\lambda a+\mu b$ mit passenden $\lambda,\mu\in\Z$.
\end{center} \dokendProof\\

{\it Bemerkung:}~
Dieser nichtkonstruktive Beweis ist dem eines allgemeineren
Resultates f\"ur sogenannte Euklidische Ringe nachempfunden, siehe hierzu
das Lehrbuch \cite[Kapitel 3, \S 17]{waerden1}. Der Beweis ist damit
eine interessante Alternative zu dem Ersetzungsverfahren 
im Anschluss an Beispiel~\ref{beis:2_7}. 
Mit Hilfe des sogenannten erweiterten Euklidischen Algorithmus 
\index{erweiterter Euklidischer Algorithmus}\label{erweiterter_Euklidischer_Algorithmus}
erhalten wir im folgenden Abschnitt 
noch eine konstruktive Beweisvariante. \dokendBem\\

Die folgende Zusammenfassung pr\"age man sich gut ein:
\begin{Satz}\label{satz:2_9}
	Es sei $a,b\in\Z$ mit $|a|+|b|>0$.\\
	Dann gelten die folgenden Aussagen:
	\begin{enumerate}[(a)]
		\item $d|\ggT(a,b)$ f\"ur jeden gemeinsamen Teiler $d$ von $a$, $b$.\\
		\item $\ggT(a,b)=\lambda a+\mu b$ mit passenden $\lambda,\mu\in\Z$.\\
		\item $\ggT(a,b)$ ist die kleinste positive Zahl, 
                die von der Linearform \index{Linearform}\label{Linearform}
                $x\cdot a+y\cdot b$ mit $x,y\in\Z$ dargestellt wird, und die Menge
		$$(a,b)=\{xa+yb:x,y\in\Z\}
		$$
		besteht genau aus den ganzzahligen Vielfachen von $\ggT(a,b)$.\\
		\item Sind speziell $a,b$ teilerfremd, also $\ggT(a,b)=1$, und gilt $a|b\cdot c$ f\"ur $a\in\Z\backslash\{0\}$ und $c\in\Z$, so folgt bereits $a|c$.
	\end{enumerate}
	\hfill	\dokendSatz
\end{Satz}

{\bf Beweis:}~
(a) ist der Satz~\ref{satz:2_6} und (b) der weiterreichende Satz~\ref{satz:2_8}.

Die erste Teilaussage von (c) haben wir im Beweis von Satz~\ref{satz:2_8} gezeigt. Es sei $k=xa+yb\in(a,b)$ mit $x,y\in\Z$. Nach Satz~\ref{satz:2_4} ist $k=q\cdot \ggT(a,b)+r$ mit 
$q\in\Z$ und $0\leq r<\ggT(a,b)$ sowie $r=k-q\cdot \ggT(a,b)\in(a,b)$, also muss $r=0$ und $k=q\cdot \ggT(a,b)$ sein. 

Wir zeigen (d): Bei $\ggT(a,b)=1$ gibt es nach (b) Zahlen $\lambda,\mu\in\Z$ mit $\lambda a+\mu b=1$. 
Es folgt mit $a|bc$, dass $c=\lambda ac+\mu bc$ durch $a$ teilbar ist.\hfill\dokendProof\\

	\section{Fundamentalsatz der Arithmetik}\label{cha:2B}

\begin{Def}[Primzahl, Einheit, Primelement]\label{def:2_10}\index{Primzahl}\label{Primzahl2}\index{Einheit}\label{Einheit}\index{Primelement}\label{Primelement}
	\hspace*{0cm}\\\vspace{-1cm}	
	\begin{enumerate}[(a)]
		\item Jede nat\"urliche Zahl $p>1$, die nur $1$ und $p$ als nat\"urliche Teiler besitzt, nennt man eine Primzahl.\\
		\item Eine Zahl $\varepsilon\in\{+1,-1\}$ hei{\ss}t Einheit in $\Z$.\\
		Die Zahlen $\varepsilon p=\pm p$ mit einer Primzahl $p$ nennt man die Primelemente von $\Z$.
	\end{enumerate}	
	\dokendDef
\end{Def}
\begin{Satz}\label{satz:2_11}
	Ist $p|ab$ mit $p$ als Primzahl und $a,b\in\N$, so gilt $p|a$ oder $p|b$. Allgemeiner: Gilt $p|a_1a_2...a_n$, dann teilt $p$ zumindest einen Faktor $a_j\in\N$ des Produktes.\hfill\dokendSatz
\end{Satz}
{\bf Beweis:}~
F\"ur $p\nmid a$ ist $\ggT(p,a)=1$ nach Definition der Primzahl $p$, und es gilt $p|a\cdot b$. Nach Satz~\ref{satz:2_9} (d) ist dann $p|b$. Die allgemeine Aussage folgt hieraus durch vollst\"andige Induktion nach der Anzahl $n$ der Faktoren.
\hfill\dokendProof\\
\begin{Satz}[Fundamentalsatz der Arithmetik]\label{satz:2_12}\index{Fundamentalsatz der Arithmetik}\label{Fundamentalsatz_der_Arithmetik2}
	Jede nat\"urliche Zahl $n>1$ kann als Produkt von Primzahlen dargestellt werden. Die Zerlegung in Primzahlen ist bis auf die Anordnung der Primfaktoren eindeutig.\hfill\dokendSatz
\end{Satz}
{\bf Beweis:}~
Unter allen Produktzerlegungen von $n>1$ mit nat\"urlichen Faktoren $\geq2$ existiert eine mit maximaler Anzahl $r$ von (m\"oglicherweise mehrfachen) Faktoren, etwa
\begin{equation}\label{eq:2_5}
	n=p_1p_2...p_r,\quad r\in\N,
\end{equation}
denn es gilt $n\geq 2^r$, und die Folge $\left(2^k\right)_{k\in\N}$ ist streng monoton wachsend und unbeschr\"ankt. Jedes $p_j\geq2$ in (\ref{eq:2_5}) muss Primzahl sein, da wir andernfalls $p_j=p_j'\cdot p_j''$ mit $p_j'\geq2$, $p_j''\geq2$ und einer Zerlegung von $n$ in $r+1$ Faktoren $\geq2$ h\"atten. Also ist~(\ref{eq:2_5}) eine Zerlegung von $n$ in Primfaktoren. Nun zeigen wir die Eindeutigkeit der Primfaktorzerlegung,
indem wir mittels Induktion f\"ur alle $n \in \N$ die folgende von $n$ abh\"angige Aussage $\mathcal{A}(n)$ beweisen: Wenn 
\begin{equation}\label{eq:2_6}
	n=q_1q_2...q_s,\quad s\in\N,
\end{equation}
und 
\begin{equation}\label{eq:2_6b}
	n=q'_1q'_2...q'_t,\quad t\in\N,
\end{equation}
zwei Primfaktorzerlegungen von $n$ sind, so stimmen diese bis auf die Reihenfolge der Faktoren \"uberein. Die Aussage stimmt f\"ur $n=1$. Wir nehmen an, dass bei einem gegebenem $n \geq 2$ die Aussage $\mathcal{A}(n')$ f\"ur alle $n'<n$ bereits stimmt, 
und m\"ussen $\mathcal{A}(n)$ zeigen. Hierf\"ur nehmen wir zwei Primfaktorzerlegungen \eqref{eq:2_6}, \eqref{eq:2_6b} von $n$ an. 
Nach Satz~\ref{satz:2_11} teilt die Primzahl $q'_1$ einen Primfaktor $q_j$ in \eqref{eq:2_6}, so dass $q'_1=q_j$ mit einem $j\in\{1,...,s\}$ gilt. 
Aus beiden Darstellungen l\"a{\ss}t sich die Primzahl $q'_1$ herausk\"urzen und hiernach die Induktionsannahme auf $n':=\frac{n}{q'_1}<n$ anwenden, d.h. es gilt $\mathcal{A}(n')$. Hieraus folgen weiter $s=t$ sowie die \"Uberstimmung von \eqref{eq:2_6} und \eqref{eq:2_6b} bis auf die Reihenfolge der Faktoren, und wir haben $\mathcal{A}(n)$ gezeigt. \hfill\dokendProof\\
\begin{Satz}\label{satz:2_13}
	Es gibt unendlich viele Primzahlen.\hfill\dokendSatz
\end{Satz}
{\bf Beweis nach Euklid:}~
H\"atte man nur endlich viele Primzahlen $p_1,...,p_s$, so w\"are $n=1+\prod\limits_{j=1}^{s}p_j>1$ 
durch kein $p_j$ teilbar, sondern durch eine "`{neue}"' Primzahl $p$, Widerspruch.\hfill\dokendProof\\

\begin{Bem}
	\hspace*{0cm}\\\vspace{-1cm}	
	\begin{enumerate}[(a)]
		\item Dem Fundamentalsatz\index{Fundamentalsatz der Arithmetik}\label{Fundamentalsatz_der_Arithmetik3} gem\"a{\ss} k\"onnen wir jede nat\"urliche Zahl $n>1$ in der kanonischen Form $n=p_{1}^{\alpha_1}\cdot p_{2}^{\alpha_2}\cdot ... \cdot p_{r}^{\alpha_r}$ mit paarweise verschiedenen Primzahlen $p_j$ und Exponenten $\alpha_j\in\N$ schreiben. Man darf dabei sogar $p_1<p_2<...<p_r$ voraussetzen, um die Reihenfolge eindeutig festzulegen.\\
		\item L\"asst man alternativ noch $\alpha_j\in\N_0$ sowie $r=0$ zu, 
		dann k\"onnen wir je zwei nat\"urliche Zahlen $a$, $b$ f\"ur passend gew\"ahltes $r \in \N_0$ in der Form
		$$a=p_{1}^{\alpha_1}\cdot p_{2}^{\alpha_2}\cdot ... \cdot p_{r}^{\alpha_r},\quad b=p_{1}^{\beta_1}\cdot p_{2}^{\beta_2}\cdot ... \cdot p_{r}^{\beta_r}
		$$
		mit $\alpha_j\geq 0$, $\beta_j\geq0$ und paarweise verschiedenen Primzahlen $p_j$ schreiben, $j=1,...,r$. Hiermit wird
		$$\ggT(a,b)=\prod\limits_{j=1}^{r}p_j^{\min(\alpha_j,\beta_j)}
		$$
		und
		$$\kgV(a,b):=\dfrac{a\cdot b}{\ggT(a,b)}=\prod\limits_{j=1}^{r}p_j^{\max(\alpha_j,\beta_j)}.
		$$ Wir nennen $\kgV(a,b)$ das kleinste gemeinsame Vielfache von $a$ und $b$. 
		Eine naheliegende Verallgemeinerung von $\ggT$ und $\kgV$ auf mehrere Argumente 
		findet der Leser in Aufgabe \ref{auf:5_2}.
		
		\begin{Beis}
			Man bestimme mittels Primfaktorzerlegung\index{Primfaktorzerlegung}\label{Primfaktorzerlegung}:
			$$\ggT(2520,1188)\quad \text{sowie}\quad \kgV(2520,1180).$$
			\textit{L\"osung}: Durch einfaches Probieren findet man
			$$\begin{array}{ll}
			2520=2^3\cdot3^2\cdot5^1\cdot7^1\cdot11^0&\text{sowie}  \\
			1188=2^2\cdot3^3\cdot5^0\cdot7^0\cdot11^1&\text{mit}
			\end{array}$$
			
			$$
			\begin{array}{l}
			\ggT(2520,1188)=2^2\cdot3^2\cdot5^0\cdot7^0\cdot11^0=36, \\
			\kgV(2520,1188)=2^3\cdot3^3\cdot5^1\cdot7^1\cdot11^1=83160.
			\end{array} $$

			F\"ur sehr gro{\ss}e Zahlen ist die Primfaktorzerlegung viel zu aufwendig oder undurchf\"uhrbar (Zahlen mit einigen Hundert Dezimalstellen), was man sich in der Kryptographie zu Nutze macht. Der Euklidische Algorithmus zur Berechnung des ggT ist dagegen sehr effizient!
						
			\dokendSatz
		\end{Beis}
		\item Ist $k\in\Z\backslash\{0\}$, so kann man auch 
		$$k=\varepsilon\cdot p_1\cdot p_2\cdot ... \cdot p_s\quad(s\geq0)
		$$
		mit einer Einheit $\varepsilon=\pm1$ und (nicht notwendigerweise verschiedenen) Primelementen $p_1 ,...,p_s$ von $\Z$ schreiben, die bis auf die Reihenfolge und das Vorzeichen eindeutig sind. Diese Form des Fundamentalsatzes\index{Fundamentalsatz_der_Arithmetik3}\label{Fundamentalsatz_der_Arithmetik4} findet eine nat\"urliche Verallgemeinerung in Euklidischen Ringen\index{Euklidischer Ring}\label{Euklidischer_Ring}. Man beachte, dass f\"ur $s=0$ das leere Produkt $\prod\limits_{j=1}^{s}p_j$ den Wert $1$ enth\"alt und $\varepsilon$ kein Primelement\index{Primelement}\label{Primelement2} in $\Z$ ist.
	\end{enumerate}		
	~~\dokendBem
\end{Bem}

	\section{Aufgaben}\label{cha:2_A}

\begin{Auf}[Euklidischer Algorithmus]\label{auf:2_1}

Mit Hilfe des Euklidischen Algorithmus bestimme man $\mbox{ggT}(462,390)$ und k\"urze anschliessend den
Bruch $390/462$.
\end{Auf}
{\bf L\"osung:}\\
Berechnung von $\ggT(462,390)$ und K\"urzung des Bruches $\D \frac{390}{462}$:

Der Algorithmus startet mit
\begin{equation*}
\begin{split}
	&r_0=b=390, \quad a=462, \quad q_0=\left\lfloor \frac{a}{b}\right\rfloor=1,\\
	&r_1=a-b q_0=462-390\cdot 1=72.
	\end{split}
\end{equation*}

Tabelle:
\vspace*{0cm}
\begin{equation*}
\begin{tabular}{|l||c|c|} \hline
$j$  & ~~ $q_j$ ~~  & ~~ $r_j~~$  \\
\hline
0  &    1    &  390 \\ \hline 
1  &    5    &  72 \\ \hline
2  &    2    &  30  \\ \hline
3  &    2    &  12  \\ \hline
4  &    2    &  6   \\ \hline
5  &    ---  &  0   \\ \hline
\end{tabular}\quad\quad\quad
\begin{array}{rl}
q_j &=  \left \lfloor \frac{r_{j-1}}{r_j} \right \rfloor \quad \text{f\"ur } j\geq 1, \\
r_{j+1}&= r_{j-1}-q_jr_j,\\
n_{\ast}&= 5 \quad \text{ist der Abbruchindex}.\\[0.3cm]
\text{K\"urzen}& \text{des Bruches mit }\ggT(462,390)=6:\\[0.1cm]
&\D\frac{390}{462}=\D\frac{390/6}{462/6}=\frac{65}{77}.
\end{array}
\end{equation*}

\begin{Auf}[Pythagoreische Zahlentripel]\label{auf:2_2}

Es sei $(a,b,c) \in \N^3$ ein {\it Pythagoreisches Zahlentripel}\index{Pythagoreisches Zahlentripel}\label{Pythagoreisches_Zahlentripel}, d.h. es gelte
\begin{equation*}
a^2+b^2=c^2\,.
\end{equation*}
Man zeige:
\begin{enumerate}[(a)]
	\item Genau dann sind $a$ und $c$ teilerfremd, wenn $b$ und $c$ teilerfremd sind.
	Wenn dies der Fall ist und zudem noch $a$ ungerade ist, 
	dann nennen wir $(a,b,c)$ ein {\it primitives Pythagoreisches Zahlentripel}.\\
			
	\item Man zeige, dass f\"ur die rationale Zahl $s := \frac{b}{a+c}$ die folgenden beiden
	Darstellungen gelten:
	\begin{equation*}
	\frac{a}{c}=\frac{1-s^2}{1+s^2}\,, \quad \frac{b}{c}=\frac{2s}{1+s^2}\,.
	\end{equation*}
	\item Man zeige mit Hilfe von (a) und (b):
	F\"ur je zwei teilerfremde nat\"urliche Zahlen $u$, $v$ mit $u>v$, 
	von denen nicht beide ungerade sind, erh\"alt man 
	ein primitives Pythagoreisches Zahlentripel $(u^2-v^2,2uv,u^2+v^2)$,
	und umgekehrt besitzt jedes primitive Pythagoreische Zahlentripel $(a,b,c) \in \N^3$
	eine solche Darstellung.
\end{enumerate}

\end{Auf}
{\bf L\"osung:}\\
F\"ur $a,b,c\in \N$ sei $a^2+b^2=c^2$. Dann ist $(a,b,c)$ Pythagoreisches Tripel.
\begin{enumerate}[(a)]
	\item  Es seien $a$ und $c$ teilerfremd. 
	Wir nehmen an, es sei $p\geq 2$ ein Primteiler von $b$ und $c$.
	Wegen $a^2=(c-b)(c+b)$ gilt dann auch $p|a^2$, denn $p$ ist ein Teiler von $c-b$.
	Wir erhalten $p|a$, da $p$ Primzahl ist, im Widerspruch zur Voraussetzung $\ggT(a,c)=1$. 
	Somit gilt $\ggT(b,c)=1$. Aus Symmetriegr\"unden folgt dann auch 
	$\ggT(a,c)=1$ aus $\ggT(b,c)=1$.
	\item F\"ur $s=b/(a+c)$ gilt
	\begin{equation*}
	\begin{split}
	\frac{1-s^2}{1+s^2}
	&=\frac{(a+c)^2-b^2}{(a+c)^2+b^2}=\frac{a^2+2ac+c^2-b^2}{a^2+2ac+c^2+b^2}\\
	&=\frac{a^2+2ac+a^2}{c^2+2ac+c^2}=\frac{2a(a+c)}{2c(a+c)}=\frac{a}{c}
	\end{split}
	\end{equation*}
	sowie mit einer Rechnung im Nenner wie oben:
	\begin{equation*}
	\begin{split}
	\frac{2s}{1+s^2}
	&=\frac{2b}{a+c}\cdot \frac{1}{1+\frac{b^2}{(a+c)^2}}=\frac{2b(a+c)}{(a+c)^2+b^2}\\
	&=\frac{2b(a+c)}{2c(a+c)}=\frac{b}{c}\,.
	\end{split}
	\end{equation*}
	
	\item F\"ur $u,v\in\N$ mit $\ggT(u,v)=1$ und $2|u\cdot v$ sowie $u>v$ sei 
	$$a':=u^2-v^2, \quad b':=2 uv, \quad c':=u^2+v^2.$$
	Dann gilt 
	$$a'\,^2+b'\,^2=(u^2-v^2)^2+4u^2v^2=u^4+2u^2v^2+v^4=c'\,^2,$$
	und $(a',b',c')$ ist Pythagoreisches Tripel.\\
	
	Aus $u>v$ und $u,v\in\N$ folgt dabei $a',b',c'\in\N$. Wegen $\ggT(u,v)=1$ und $2|u\cdot v$ ist $a'$ ungerade. Wir nehmen an, es sei $p\geq 3$ ein Primteiler von $a'$ und $c'$. Dann folgen $p|c'+a'$ und $p|c'-a'$, also $p|2u^2$, $p|2v^2$, und damit auch $p|\ggT(u,v)$, ein Widerspruch. Somit ist $(u^2-v^2,2uv,u^2+v^2)$ ein primitives Pythagoreisches Tripel.\\
	
	Nun sei umgekehrt das gegebene Pythagoreisches Tripel $(a,b,c)$ als primitiv vorausgesetzt. Dann gilt
	$$s=\frac{b}{a+c}=\frac{v}{u}\quad \text{mit}\quad  u,v\in\N,\ \ggT(u,v)=1.$$
	Aus (b) folgt 
	\begin{equation}\label{eq:2_7}
	\frac{a}{c}=\frac{1-s^2}{1+s^2}=\frac{u^2-v^2}{u^2+v^2}, \quad \frac{b}{c}=\frac{2s}{1+s^2}=\frac{2uv}{u^2+v^2}.
	\end{equation}
	W\"aren $a=2k+1$, $b=2m+1$ mit $k,m\in\N_0$ beide ungerade, so h\"atten wir
	$$a^2=1+8\frac{k(k+1)}{2}, \quad b^2=1+8\frac{m(m+1)}{2},$$
	und $c^2=a^2+b^2\equiv 2~(4)$ w\"are zwar gerade, aber nicht durch $4$ teilbar, ein Widerspruch. Also ist $b$ gerade, und $a,c$ m\"ussen ungerade sein. Aus der Darstellung \eqref{eq:2_7} folgt $cuv=b\frac{u^2+v^2}{2}$, und hieraus, 
	dass $u$ und $v$ nicht beide ungerade sind.\\
	
	Es folgt $2|uv$, und zusammen mit $u>v$, $\ggT(u,v)=1$
	wie zuvor gezeigt, dass $(u^2-v^2,2uv,u^2+v^2)$ ein {\it primitives} Pythagoreisches Tripel ist.
	nach Voraussetzung ist aber auch $(a,b,c)$ ein primitives  Pythagoreisches Tripel.
	Wir erhalten somit aus \eqref{eq:2_7}:
	$$a=u^2-v^2, \quad b=2uv, \quad c=u^2+v^2.$$	
	
\end{enumerate}
\begin{Auf}[Fibonacci-Folge, Teil 2]\label{auf:2_3}

Wir erinnern an die Definition der Fibonacci-Folge\index{Fibonacci-Folge}\label{Fibonacci-Folge}
$(f_n)_{n \in \N_0}$ mit den Fibonacci-Zahlen $f_0=0$, $f_1 = 1$ sowie 
$f_{n+2}=f_{n+1}+f_n$ f\"ur alle $n \in \N_0$, siehe Lektion~\ref{cha:1}, Aufgabe~\ref{auf:1_4}. %\ref{auf:1_4}
Zus\"atzlich definieren wir noch $f_{-1}:=1$.
\begin{enumerate}[(a)]
	\item Man zeige f\"ur alle $b \in \N_0$:
	$\begin{displaystyle}
	\mbox{ggT}(f_b,f_{b+1})=1\,.
	\end{displaystyle}$\\
	{\it Hinweis:} Lektion~\ref{cha:1}, Aufgabe~\ref{auf:1_4} (a).\\ %\ref{auf:1_4}
	
	\item Man zeige f\"ur alle $b,r \in \N_0$:
	$\begin{displaystyle}
	f_{b+r}=f_{b+1}f_r + f_b f_{r-1}\,.
	\end{displaystyle}$
	
	{\it Hinweis:}
	$\begin{displaystyle}
	\begin{pmatrix}
	1 & 1\\
	1 & 0\\
	\end{pmatrix}^{b+r} =
	\begin{pmatrix}
	1 & 1\\
	1 & 0\\
	\end{pmatrix}^b 
	\cdot
	\begin{pmatrix}
	1 & 1\\
	1 & 0\\
	\end{pmatrix}^r
	\,.
	\end{displaystyle}$\\
	
	\item Mit Hilfe von (a) und (b) zeige man f\"ur alle $b \in \N$ und $q,r \in \N_0$:
	\begin{equation*}
	\mbox{ggT}(f_b,f_{r})=\mbox{ggT}(f_b,f_{b+r})\,, \quad
	\mbox{ggT}(f_{qb+r},f_b)=\mbox{ggT}(f_b,f_{r})\,, 
	\end{equation*}
	und schlie{\ss}lich mit Hilfe des Euklidischen Algorithmus\index{Euklidischer Algorithmus}\label{Euklidischer_Algorithmus3}:
	$$
	\mbox{ggT}(f_a,f_{b})=f_{\mbox{ggT}(a,b)} \quad \mbox{f\"ur~alle~} a \in \N_0, b \in \N\,.
	$$
\end{enumerate}

\end{Auf}
{\bf L\"osung:}
$$f_{-1}=1, \quad f_0=0, \quad f_1=1, \quad f_{n+2}=f_{n+1}+f_n \quad \mbox{f\"ur~alle~} n\in\N_0.$$
Lektion~\ref{cha:1}, Aufgabe~\ref{auf:1_4}~(a) liefert
	$\begin{displaystyle}
	\begin{pmatrix}
	1 & 1\\
	1 & 0\\
	\end{pmatrix}^n =
	\begin{pmatrix}
	f_{n+1} & f_n\\
	f_n & f_{n-1}\\
	\end{pmatrix}\ \mbox{~f\"ur~alle~} n\in\N_0\,.
	\end{displaystyle}$

\begin{enumerate}[(a)]
	\item $f_{n+1} f_{n-1}-f_n^2=\Det \left(\begin{pmatrix}
	1 & 1\\
	1 & 0\\
	\end{pmatrix}^n\right)=(-1)^n\ \mbox{~f\"ur~alle~} n\in\N_0$. Ersetzen wir $n$ durch \mbox{$b\in\N_0$}, so folgt f\"ur $\lambda:=(-1)^bf_{b-1},\ \mu:=(-1)^{b+1}f_b\in\Z$:
	$$\lambda f_{b+1}+\mu f_b=1, \quad \text{d.h.} \quad \ggT(f_{b+1},f_b)=1.$$
	\item Wir haben f\"ur alle $b,r\in\N_0$ mit $A:=\begin{pmatrix}
		1 & 1\\
		1 & 0\\
	\end{pmatrix}$:
	\begin{eqnarray*}
	\begin{split}
	A^{b+r}&=
	\begin{pmatrix}
	f_{b+r+1} & f_{b+r}\\
	f_{b+r} & f_{b+r-1}\\
	\end{pmatrix}
	=A^bA^r\\
	&=\begin{pmatrix}
	f_{b+1} & f_b\\
	f_b & f_{b-1}\\
	\end{pmatrix}
	\begin{pmatrix}
	f_{r+1} & f_r\\
	f_r & f_{r-1}\\
	\end{pmatrix}\\
	&= \begin{pmatrix}
	f_{b+1}f_{r+1}+f_b f_r &\ \ f_{b+1}f_r+f_b f_{r-1}\\
	f_b f_{r+1}+f_{b-1}f_r &\ \ f_b f_r+ f_{b-1}f_{r-1}\\
	\end{pmatrix},
	\end{split}
	\end{eqnarray*}
	und hieraus $f_{b+r}=f_{b+1}f_r+f_b f_{r-1}$.\\
	
	\item Es seien $b\in\N$, $q,r\in\N_0$. Nach (b) ist jeder gemeinsame Teiler von $f_b$ und $f_r$ auch ein gemeinsamer Teiler von $f_b$ und $f_{b+r}$. Umgekehrt gelte $d|f_b$ und $d|f_{b+r}$. Aus $f_{b+1} f_r=f_{b+r}-f_b f_{r-1}$ folgt dann auch $d|f_{b+1}\cdot f_r$, und hieraus $d|f_r$, denn wegen $d|f_b$ und $\ggT(f_b,f_{b+1})=1$ ist $d$ auch zu $f_{b+1}$ teilerfremd. Damit haben $f_b$ und $f_r$ dieselben Teiler wie $f_b$ und $f_{b+r}$, insbesondere gilt
	$$ \ggT(f_b,f_r)=\ggT(f_b,f_{b+r}).$$
	Induktion bzgl. $q\in\N_0$ liefert nun
	\begin{equation}\label{eq:2_8}
	\ggT(f_{qb+r},f_b)=\ggT(f_b,f_r) \quad \mbox{f\"ur~alle~} b\in\N, q,r\in\N_0.
	\end{equation}
	Wir wenden auf $a\in\N_0$, $b\in\N$ den Euklidischen Algorithmus mit Abbruchindex $n_*$ an: 
	\begin{equation*}
	\begin{split}
	&r_0:=b,\ q_0:=\left\lfloor \frac{a}{b}\right\rfloor,\ r_1:=a-q_0b,\\
	&r_{j+1}=r_{j-1}-q_j r_j\ \ \text{f\"ur}\ q_j=\left\lfloor \frac{r_{j-1}}{r_j}\right\rfloor
	\ \text{und}\ j=1,...,n_{\ast}-1,\\
	&r_{n_\ast-1}=\ggT(a,b),\ r_{n_\ast}=0.
	\end{split}
	\end{equation*}
	Wir erhalten der Reihe nach aus \eqref{eq:2_8}:
	$$ \ggT(f_a,f_b)=\ggT(f_{q_0 b+r_1},f_b)=\ggT(f_b,f_{r_1})=\ggT(f_{r_0},f_{r_1}),$$
	sowie f\"ur $j=1,...,n_{\ast}-1$:
	$$\ggT(f_{r_{j-1}},f_{r_j})=\ggT(f_{q_jr_j+r_{j+1}},f_{r_j})=\ggT(f_{r_j},f_{r_{j+1}}).$$
	Hieraus folgt endlich
	$$\ggT(f_a,f_b)=\ggT(f_{r_{n_\ast-1}},\underbrace{f_{r_{n_\ast}}}_{=0})=f_{r_{n_\ast-1}}
	=f_{\ggT(a,b)}.$$

\end{enumerate}
\begin{Auf}[Eigenschaften der oberen und unteren Gau{\ss}-Klammer \index{Gaussklammer}\label{Gaussklammer2}]\label{auf:2_4}
	%	\hspace*{0cm}\\\vspace{1cm}

Man zeige, dass f\"ur jede reelle Zahl $x$ gilt:
\begin{enumerate}[(a)]
	\item $\begin{displaystyle}
	\lfloor x \rfloor \leq x < \lfloor x \rfloor + 1\,,\quad
	x-1  < \lfloor x \rfloor \leq x\,,\quad
	0 \leq x-\lfloor x \rfloor < 1\,,
	\end{displaystyle}$\\
	insbesondere ist $\begin{displaystyle}\lfloor x \rfloor \end{displaystyle}$
	die gr\"o{\ss}te ganze Zahl kleiner oder gleich $x$.\\
	
	\item $\begin{displaystyle}
	\lfloor x+k \rfloor= \lfloor x \rfloor + k \quad \mbox{f\"ur~alle~} k \in \Z\,,
	\end{displaystyle}$\\
	
	\item
	$\begin{displaystyle}
	\left\lfloor \frac{x}{n} \right\rfloor= 
	\left\lfloor{\frac{\lfloor x \rfloor }{n}}\right\rfloor \quad \mbox{f\"ur~alle~} n \in \N\,,
	\end{displaystyle}$\\
	
	\item F\"ur $\begin{displaystyle}\lceil x \rceil := -\lfloor -x \rfloor
	\end{displaystyle}$ ist\\
	$$
	\lceil x \rceil -1 < x \leq \lceil x \rceil\,,\quad
	x \leq \lceil x \rceil < x+1\,,\quad
	0 \leq \lceil x \rceil -x < 1\,,
	$$
	insbesondere ist $\begin{displaystyle}\lceil x \rceil \end{displaystyle}$
	die kleinste ganze Zahl gr\"o{\ss}er oder gleich $x$.
\end{enumerate}

\end{Auf}
{\bf L\"osung:}\\
Es sei $x\in\R$. Dann ist $\lfloor x \rfloor$ (ganzzahliger Anteil von $x$) diejenige ganze Zahl $j$, f\"ur die gilt: 	
\begin{equation}\label{eq:2_9}
j\leq x<j+1.
\end{equation}
\begin{enumerate}[(a)]
	\item Die erste Ungleichungskette entspricht \eqref{eq:2_9}, d.h. 
	$\lfloor x \rfloor \leq x<\lfloor x \rfloor+1$ f\"ur $j=\lfloor x \rfloor$, und die beiden anderen sind Umformulierungen dieser Ungleichungen.\\
	
	\item Aus (a) bzw. \eqref{eq:2_9} folgt f\"ur $k\in\Z$:
	$$\lfloor x \rfloor+k\leq x+k<(\lfloor x \rfloor+k)+1$$
	mit der ganzen Zahl $\lfloor x \rfloor+k$, d.h. $\lfloor x+k \rfloor=\lfloor x \rfloor+k\ \ \mbox{f\"ur~alle~} k\in\Z$.\\
	
	\item F\"ur $n\in\N$ ist zun\"achst nach (a):
	\begin{enumerate}[1)]
		\item 
		$\begin{displaystyle}
		\left\lfloor{\frac{\lfloor x \rfloor }{n}}\right\rfloor
		\leq \frac{\left\lfloor x \right\rfloor}{n} \leq \frac{x}{n},
		\end{displaystyle}$\\
		
	und aus $x<\lfloor x \rfloor+1$ in (a) erhalten wir\\	
	\item
	$\begin{displaystyle}
	\frac{x}{n} < \frac{\left\lfloor x \right\rfloor+1}{n}.
	\end{displaystyle}$\\
	
	Aus der zweiten Ungleichungskette in (a) folgt
	$$n	\left\lfloor{\frac{\lfloor x \rfloor }{n}}\right\rfloor >n\left(\frac{\lfloor x \rfloor}
	{n}-1\right)=\lfloor x \rfloor-n,$$
	und aus der Ganzzahligkeit von $\begin{displaystyle}\,n	\left\lfloor{\dfrac{\lfloor x \rfloor }{n}}\right\rfloor\end{displaystyle}$
	sowie $\lfloor x \rfloor-n$ die Ungleichung
	$$
	n\left\lfloor{\dfrac{\lfloor x \rfloor }{n}}\right\rfloor\geq \lfloor x \rfloor-n + 1\,.$$
	Die letzte Ungleichung schreiben wir in der \"aquivalenten  Form\\
	
	\item 
	$\begin{displaystyle}
	\frac{\left\lfloor x \right\rfloor+1}{n}
	\leq \left\lfloor{\frac{\lfloor x \rfloor }{n}}\right\rfloor+1
	\end{displaystyle}$.\\
	Wir erhalten endlich
	$\begin{displaystyle}
	\left\lfloor{\frac{\lfloor x \rfloor }{n}}\right\rfloor
	\leq \frac{x}{n}
	< \left\lfloor{\frac{\lfloor x \rfloor }{n}}\right\rfloor+1
	\end{displaystyle}$ aus 1) bis 3), d.h. \\
	
	$\begin{displaystyle}
	\left\lfloor{\frac{x}{n}}\right\rfloor
	=\left\lfloor{\frac{\lfloor x \rfloor }{n}}\right\rfloor
	\end{displaystyle}$.
\end{enumerate}
\item folgt aus (a), indem man dort $x$ durch $-x$ ersetzt.
\end{enumerate}

\begin{Auf}[Die h\"ochsten Primzahlpotenzteiler von $n!$
\index{Fakult\"at}\label{Fakultaet} \index{Primzahlpotenz}\label{Primzahlpotenz}]\label{auf:2_5}

Es sei $p$ eine Primzahl, $n$ eine nichtnegative ganze Zahl
und $\alpha_p(n)$ die gr\"o{\ss}te ganze Zahl 
$\alpha \geq 0$, f\"ur die $p^{\alpha}$
ein Teiler von $n!$ ist. Man zeige 
\begin{equation*}
\alpha_p(n) = \sum \limits_{k=1}^{\infty}\left\lfloor \frac{n}{p^k} \right\rfloor\,.
\end{equation*}
\end{Auf}
{\bf L\"osung:}\\
Die Summe in der zu beweisenden Formel mu{\ss} nur \"uber
die endlich vielen $k$ mit $p^k \leq n$ erstreckt werden. 
Wir beweisen die Formel durch Induktion nach $n$.

F\"ur $n=0$ stimmt sie mit dem Wert $\alpha_p(1)=0$,
wobei $0!=1$ zu beachten ist.
Nun nehmen wir an, es sei $n$ eine nat\"urliche Zahl mit 
\begin{equation*}
\alpha_p(m) = \sum \limits_{k=1}^{\infty}\left\lfloor \frac{m}{p^k} \right\rfloor
\end{equation*}
f\"ur alle ganzen Zahlen $m$ mit $0 \leq m <n$.
Da $p$ eine Primzahl ist, k\"onnen wir zur Berechnung von $\alpha_p(n)$
all diejenigen Faktoren $d$ aus dem Produkt
\begin{equation*}
n! = \prod \limits_{d=1}^{n}d
\end{equation*}
streichen, die nicht durch $p$ teilbar sind, so dass $p^{\alpha_p(n)}$
auch die h\"ochste Potenz von $p$ wird, welche das Produkt
\begin{equation*}
\prod \limits_{j \leq n/p} (pj) = p^{\lfloor n/p \rfloor} \cdot \lfloor n/p \rfloor!
\end{equation*}
teilt. Hieraus folgt 
\begin{equation*}
\alpha_p(n) = \left \lfloor \frac{n}{p} \right \rfloor +
\alpha_p\left( \left \lfloor \frac{n}{p} \right \rfloor\right)\,.
\end{equation*}
Nach der Induktionsannahme mit der Wahl von $m=\lfloor n/p \rfloor < n$
und der zuvor gel\"osten Aufgabe \ref{auf:2_4}(c) erhalten wir
\begin{equation*}
\alpha_p(n) = \left \lfloor \frac{n}{p} \right \rfloor +
\sum \limits_{k=1}^{\infty}\left\lfloor \frac{\lfloor n/p \rfloor}{p^k} \right\rfloor
=\left \lfloor \frac{n}{p} \right \rfloor +
\sum \limits_{k=1}^{\infty}\left\lfloor \frac{n}{p^{k+1}} \right\rfloor
= \sum \limits_{k=1}^{\infty}\left\lfloor \frac{n}{p^k} \right\rfloor\,,
\end{equation*}
so dass auch der Induktionsschritt gezeigt ist.

	\chapter[Kettenbr\"uche]{Erweiterter Euklidischer Algorithmus und Kettenbruchentwicklung reeller Zahlen}\label{cha:3}
	
Jede rationale Zahl l\"asst sich als endlicher Kettenbruch
\begin{equation*}
q_0+ \cfrac{1}{q_1 + \cfrac{1}{q_2 + \hspace{0.1cm}
\raisebox{-0.3cm}{$\ddots \hspace{0.1cm} \raisebox{-0.2cm}{$+ \,\cfrac{1}{q_{j-1} + \cfrac{1}{q_{j}}}$}$}}}
\end{equation*}
mit $q_0 \in \Z$ und $q_1,\ldots,q_j \in \N$ darstellen,
wie wir in diesem Abschnitt mit Hilfe des Euklidischen Algorithmus zeigen werden.
Mit der Zusatzforderung $q_j>1$ erh\"alt man \"uberdies die Eindeutigkeit dieser Darstellung.
Eine nat\"urliche Erweiterung des Euklidischen Algorithmus wird uns dabei sogar
unendliche Kettenbruchentwicklungen\index{Kettenbruchentwicklung}\label{Kettenbruchentwicklung} f\"ur alle reellen Irrationalzahlen liefern,
die sich ohne weitere Einschr\"ankung als eindeutig erweisen. Wir zeigen
in diesem Abschnitt auf, wie sich mit ihrer Hilfe die bestm\"oglichen rationalen Approximationen \index{beste rationale Approximation}\label{beste rationale Approximation}
der reellen Zahlen gewinnen lassen. Insbesondere werden wir im Abschnitt~\ref{cha:8} 
bei den rationalen Approximationen der reell quadratischen Irrationalzahlen \index{reell quadratische Irrationalzahl}\label{reell quadratische Irrationalzahl}
auf dieses Thema zur\"uckkommen. Letztere Thematik h\"angt eng mit der Theorie der
sogenannten indefiniten quadratischen Formen zusammen
und erweist sich f\"ur eine Einf\"uhrung in die elementare 
Zahlentheorie als interessant. \\

Als Lekt\"ure zur Vertiefung empfehlen wir die Lehrb\"ucher von Hardy/Wright 
\cite[Chapter X]{hw}, Niven/Zuckerman \cite[Band 47, Abschnitt 7]{nz},
Oswald und Steuding \cite{os}, Steuding \cite{st} 
sowie Perron's zeitlose Monographie \cite{perron}.
Diese Literaturquellen haben uns als Inspiration gedient.\\
        \section{Erweiterter Euklidischer Algorithmus und Kettenbr\"uche}\label{cha:3A}
\begin{Def}[endliche Kettenbr\"uche\index{endlicher Kettenbruch}\label{endlicher Kettenbruch}]\label{def:3_1}
	F\"ur $\lambda_0\in\R$ und positive reelle Zahlen $\lambda_1, ... ,\lambda_j$ definieren wir den Kettenbruch\index{Kettenbruch}\label{Kettenbruch}
	
	\begin{equation*}
	\langle\lambda_0,\lambda_1, ... ,\lambda_j\rangle= 
	\lambda_0+ \cfrac{1}{\lambda_1 + \cfrac{1}{\lambda_2 + \hspace{0.1cm}
			\raisebox{-0.3cm}{$\ddots \hspace{0.1cm} \raisebox{-0.2cm}{$+ \,\cfrac{1}{ \lambda_{j}}$}$}}}
	\end{equation*}

	induktiv gem\"a{\ss}
	$$\langle\lambda_0\rangle:=\lambda_0,\quad\langle\lambda_0,\lambda_1, ... ,\lambda_j\rangle:=\lambda_0+\dfrac{1}{\langle\lambda_1, ... ,\lambda_j\rangle},~\ j\in\N.
	$$
	\dokendDef
\end{Def}
\begin{Satz}\label{satz:3_2}
F\"ur alle $\lambda_0\in\R$ und alle $\lambda_1, ... ,\lambda_j>0$,~\ $j\in\N$, gilt 
$$\langle\lambda_0,\lambda_1, ... ,\lambda_j\rangle=\langle\lambda_0,\langle \lambda_1, ... ,\lambda_j\rangle\rangle,
$$
und f\"ur $j\geq2$ \"uberdies 
$$\langle\lambda_0,\lambda_1, ... ,\lambda_j\rangle=\langle\lambda_0,...,\lambda_{j-2},\lambda_{j-1}+\frac{1}{\lambda_j}\rangle.
$$
\hfill\dokendSatz
\end{Satz}
{\bf Beweis:}~ 
Wegen 
$$\langle\lambda_0,\lambda_1, ... ,\lambda_j\rangle=\lambda_0+\dfrac{1}{\langle\lambda_1, ... ,\lambda_j\rangle}=\langle\lambda_0,\langle\lambda_1, ... ,\lambda_j\rangle\rangle$$ 
folgt die erste Teilaussage des Satzes sofort aus der Definition~\ref{def:3_1}. 
Setzen wir speziell $j=2$, so erhalten wir wegen 
$$\langle\lambda_0,\lambda_1,\lambda_2\rangle=\langle\lambda_0,\langle\lambda_1,\lambda_2\rangle\rangle
=\langle\lambda_0,\lambda_1+\dfrac{1}{\lambda_2}\rangle$$
bereits den Induktionsanfang f\"ur die zweite Teilaussage. 
Wir nehmen an, die zweite Teilaussage sei f\"ur ein $j\geq 2$ bereits g\"ultig. Dann gilt sie auch f\"ur $j+1$ wegen 
$$\begin{array}{rclcl}
\langle\lambda_0,\lambda_1, ... ,\lambda_j,\lambda_{j+1}\rangle & = & \langle\lambda_0,\langle\lambda_1, ... ,\lambda_j,\lambda_{j+1}\rangle\rangle &  & \\
& = & \langle\lambda_0,\langle\lambda_1, ... ,\lambda_{j-1},\lambda_{j}+\dfrac{1}{\lambda_{j+1}}\rangle\rangle&  & (\text{Induktionsannahme})\\
& = &\langle\lambda_0,\lambda_1, ... ,\lambda_{j-1},\lambda_{j}+\dfrac{1}{\lambda_{j+1}}\rangle &  & 
\end{array} $$
f\"ur alle $\lambda_0\in \R$ und alle $\lambda_1, ... ,\lambda_{j+1}> 0$.
\dokendProof

Bei Verwendung der Klammer-Notation f\"ur Kettenbr\"uche lassen sich mit Hilfe von Satz \ref{satz:3_2}
endliche Kettenbr\"uche besonders einfach berechnen, z.B. erhalten wir
\begin{equation*}
\langle 1,2,3 \rangle =\langle 1,\, 2+\frac{1}{3} \rangle
=\langle 1,\, \frac{7}{3} \rangle= 1+\frac{3}{7} =\frac{10}{7}\,.
\end{equation*}

\begin{Satz}\label{satz:3_3}
	F\"ur $j\in\N$ seien $\lambda_0\in\R$ sowie $\lambda_1, ... ,\lambda_{j-1}>0$ reell. Setze $\underline{\lambda}:=(\lambda_0, ... ,\lambda_{j-1})\in\R^j$ f\"ur $j\geq 2$ bzw. $\underline{\lambda}:=(\lambda_0)$ f\"ur $j=1$ sowie $T_{\underline{\lambda}}:=T_{(\lambda_0)}\cdot\left( \begin{matrix}
	0 & 1  \\
	1 & \lambda_1\\
	\end{matrix} \right)\cdot...\cdot\left( \begin{matrix}
	0 & 1  \\
	1 & \lambda_{j-1}\\
	\end{matrix} \right)$ mit $T_{(\lambda_0)}:=\left( \begin{matrix}
	1 & \lambda_0  \\
	0 & 1\\
	\end{matrix} \right)$ f\"ur $j\geq 2$, $\sigma_0=1$, $\sigma_1=\lambda_0$, $\tau_0=0$, $\tau_1=1$ und iterativ $$\sigma_{k+1}=\sigma_{k-1}+\lambda_k \sigma_k,\quad \tau_{k+1}=\tau_{k-1}+\lambda_k \tau_k$$
	 f\"ur $1\leq k <j$. Dann gilt:
	\begin{enumerate}[(a)]
		\item $T_{\underline{\lambda}}=\left( \begin{matrix}
		\sigma_{j-1} & \sigma_j  \\
		\tau_{j-1} & \tau_j\\
		\end{matrix} \right)$.\\
		
		\item $\langle\lambda_0,\lambda_1, ... ,\lambda_{j-1},x\rangle=\dfrac{\sigma_j x+\sigma_{j-1}}{\tau_j x+\tau_{j-1}}$ f\"ur $x>0$.\\
		
		\item $\langle\lambda_0,\lambda_1, ... ,\lambda_{j-1},x\rangle-\langle\lambda_0,\lambda_1, ... ,\lambda_{j-1},x'\rangle=\dfrac{(-1)^j(x-x')}{(\tau_j x+\tau_{j-1})(\tau_j x'+\tau_{j-1})}$\\ f\"ur $x,x'>0$.
	\end{enumerate}
	\hfill\dokendSatz
\end{Satz}
{\bf Beweis:}~ 
Wir f\"uhren den Beweis von (a) und (b) durch vollst\"andige Induktion:
\begin{enumerate}[(a)]
	\item F\"ur $j=1$ haben wir $T_{(\lambda_0)}=\left( \begin{matrix}
	1 & \lambda_0  \\
	0 & 1\\
	\end{matrix} \right)=\left( \begin{matrix}
	\sigma_{0} & \sigma_1  \\
	\tau_{0} & \tau_1\\
	\end{matrix} \right)$ aufgrund der Startvorgaben. Ist die Aussage f\"ur einen Index $j\geq 1$ g\"ultig und setzen wir $\underline{\lambda}'=(\lambda_0,\lambda_1, ... ,\lambda_{j})$ mit $\lambda_j>0$, so wird 
	$$T_{\underline{\lambda}'}=\left( \begin{matrix}
	\sigma_{j-1} & \sigma_j  \\
	\tau_{j-1} & \tau_j\\
	\end{matrix} \right)\cdot\left( \begin{matrix}
	0 & 1  \\
	1 & \lambda_j\\
	\end{matrix} \right)=\left( \begin{matrix}
	\sigma_{j} & \sigma_{j-1}+\sigma_{j}\lambda_j  \\
	\tau_{j} & \tau_{j-1}+\tau_j \lambda_j\\
	\end{matrix} \right)=\left( \begin{matrix}
	\sigma_{j} & \sigma_{j+1}  \\
	\tau_{j} & \tau_{j+1}\\
	\end{matrix} \right).
	$$\\
	
	\item F\"ur $j=1$ ist
	$$\langle\lambda_0,x\rangle=\lambda_0+\dfrac{1}{x}=\dfrac{\sigma_1 x+\sigma_0}{\tau_1 x+\tau_0}.
	$$  Wir nehmen an, die zu beweisende Aussage sei f\"ur einen Index $j\geq1$ g\"ultig, w\"ahlen wieder $\lambda_j>0$ und setzen $\underline{\lambda}'=(\lambda_0,\lambda_1, ... ,\lambda_{j})$. Dann folgt aus Satz~\ref{satz:3_2} und unserer Induktionsannahme:
		
	\begin{equation*}
	\begin{split}	
	\langle\lambda_0,\lambda_1, ... ,\lambda_{j},x\rangle 
	=&  \langle\lambda_0,\lambda_1, ... ,\lambda_{j-1},\lambda_j+\frac{1}{x} \rangle\\
	=&  \dfrac{\sigma_{j}(\lambda_j +\frac{1}{x})+\sigma_{j-1}}{\tau_j(\lambda_j +\frac{1}{x})+\tau_{j-1}}	 =  \dfrac{(\sigma_{j-1}+\lambda_j \sigma_j)x+\sigma_j}{(\tau_{j-1}+\lambda_j \tau_j)x+\tau_j} =  \dfrac{\sigma_{j+1}x+\sigma_j}{\tau_{j+1}x+\tau_j},
	\end{split}
	\end{equation*}
	
	so dass die Aussage auch f\"ur $j+1$ stimmt.	 \\
	
	\item folgt aus (b) durch direktes Nachrechnen unter Beachtung von 
	$$\Det\, T_{\underline{\lambda}}=(-1)^{j-1}=\sigma_{j-1}\tau_j-\sigma_j \tau_{j-1}.$$
\end{enumerate}
\dokendProof\\

{\bf Einf\"uhrung des erweiterten Euklidischen Algorithmus}\\

Gegeben sind $a$, $b\in\R$ mit $b>0$.
 \begin{enumerate}[1)]
 	\item Wir definieren die Startwerte ~\ $x_0=\frac{a}{b}$,~\ $q_0=\left\lfloor\frac{a}{b}\right\rfloor$,~\ $r_0=b$,~\ $r_1=a-b\cdot\left\lfloor\frac{a}{b}\right\rfloor$,~\ $s_0=1$, $s_1=q_0$,~\ $t_0=0$,~\ $t_1=1$.\\
 	\item Solange $r_j\neq0$ ist, berechnen wir, beginnend mit $j=1$, schrittweise die Gr\"o{\ss}en
 	\begin{equation*}
        \begin{split}
 	x_j=\dfrac{r_{j-1}}{r_j},\quad q_j=\lfloor x_j\rfloor,\quad r_{j+1}=r_{j-1}-q_j \,r_j\,,\\
 	s_{j+1}=s_{j-1}+q_j \, s_j,\quad t_{j+1}=t_{j-1}+q_j \, t_j.\\
 	\end{split}
        \end{equation*} 
 	\item Falls $r_j=0$ f\"ur $j\in\N$ gilt, brechen wir den Algorithmus ab und nennen ihn terminierend mit Abbruchindex 
 	$n_*=n_*(a,b)=j$. Falls $r_j\neq0$ f\"ur alle $j\in\N_0$ definiert ist, 
        nennen wir den Algorithmus infinit und setzen $n_*=n_*(a,b)=\infty$\,.
 \end{enumerate}

Durch 1) bis 3) ist der erweiterte Euklidische Algorithmus mit 
{\it Eingabewerten} $a$ und $b$ erkl\"art.
Die {\it Ausgabewerte} sind $x_j$, $q_j$ mit $0\leq j<n_*$ 
sowie $r_j$, $s_j$, $t_j$ mit $j\in\N_0$ und $j\leq n_*$.
Ist der Algorithmus terminierend, so nennen wir auch $n_*\in\N$ einen Ausgabewert.

\begin{Bem}\label{bem:3_4}
	\hspace*{0cm}\\\vspace{-1cm}
	\begin{enumerate}[(a)]
		\item Ist $\lambda>0$ und ersetzen wir die Eingabewerte\index{Eingabewerte des erweiterten Euklidischen Algorithmus}\label{Eingabewerte des erweiterten Euklidischen Algorithmus} $a$, $b$ durch $\lambda a$, $\lambda b$, so bleiben der Abbruchindex $n_*\in\N\cup\{\infty\}$ und alle Ausgabewerte\index{Ausgabewerte des erweiterten Euklidischen Algorithmus}\label{Ausgabewerte des erweiterten Euklidischen Algorithmus} mit Ausnahme der $r_j$ erhalten; allein die $r_j$ m\"ussen durch die neuen Ausgabewerte $\lambda\cdot r_j$ ersetzt werden.\\
		
		\item Setzen wir in Satz~\ref{satz:3_3} \,$\lambda_k=q_k$ f\"ur $0\leq k<n_*$, so folgt dort $\sigma_k=s_k$, $\tau_k=t_k$ f\"ur $k\in\N_0$ mit $k\leq n_*$, wovon wir nun Gebrauch machen:
	\end{enumerate}
	~~\dokendBem
\end{Bem}

\begin{Satz}\label{satz:3_5}
	\hspace*{0cm}\\\vspace{-1cm}
	\begin{enumerate}[(a)]
		\item Es gilt $r_j>r_{j+1}\geq 0$ f\"ur $0\leq j<n_*$. Hierbei ist $r_{j+1}=0$ nur f\"ur $j+1=n_*<\infty$ m\"oglich.\\
		
		\item Es gilt $x_j>1$ und $q_j\in\N$ f\"ur $1\leq j<n_*$ sowie $x_j-q_j=\frac{1}{x_{j+1}}$ f\"ur $1\leq j+1<n_*$.\\
		
		\item $\langle q_0,...,q_{j-1},q_j\rangle=\frac{s_{j+1}}{t_{j+1}}$ mit $t_{j+1}\geq 1$ und $t_{j+1}\geq t_j$ f\"ur $0\leq j<n_*.$\\
		
		\item $\dfrac{a}{b}=\langle q_0,..,q_{j-1},x_j\rangle$ f\"ur $1\leq j<n_*.$\\
		
		\item Es ist $n_*$ genau dann endlich, wenn $\dfrac{a}{b}$ rational ist. In diesem Falle gilt die Beziehung $x_{n_*-1}=q_{n_*-1}$ mit 
		$$\dfrac{a}{b}=\langle q_0,...,q_{n_*-1}\rangle.
		$$
	\end{enumerate}
	\hfill\dokendSatz
\end{Satz}
{\bf Beweis:}~
	\begin{enumerate}[(a)]
		\item F\"ur $j=0$ haben wir $$r_0=b>a-b\left\lfloor\dfrac{a}{b}\right\rfloor=r_1\geq 0,$$
		 denn einerseits ist $$a-b\left\lfloor\dfrac{a}{b}\right\rfloor\geq a-b\frac{a}{b}=0,$$ 
		 und andererseits $$a-b\left\lfloor\dfrac{a}{b}\right\rfloor<a-b\left(\dfrac{a}{b}-1\right)=b,$$
		 da wir $b>0$ vorausgesetzt haben.\\
		 Solange noch $r_j>0$ f\"ur einen Index $j\geq 1$ ist, also $j<n_*$ gilt, erhalten wir zum einen 
		 $$r_{j+1}=r_{j-1}-\left\lfloor\dfrac{r_{j-1}}{r_j}\right\rfloor r_j\geq r_{j-1}-\dfrac{r_{j-1}}{r_j}\cdot r_j=0,$$
		 und zum anderen
		 $$r_{j+1}=r_{j-1}-\left\lfloor\dfrac{r_{j-1}}{r_j}\right\rfloor r_j<r_{j-1}-\left(\dfrac{r_{j-1}}{r_j}-1\right) r_j=r_j.
		 $$\\
		 Die Bedingung $r_{j+1}=0$ f\"ur $j+1=n_*<\infty$ ist genau die Abbruchbedingung f\"ur den erweiterten Euklidischen Algorithmus.\\
		
		\item Aus (a) folgt $x_j=\frac{r_{j-1}}{r_j}>1$ und $q_j=\lfloor x_j\rfloor\geq1$ f\"ur $1\leq j<n_*$. Nun sei $n_*>1$, $n_*$ endlich oder unendlich. Dann existiert $x_1>1$ mit 
		$$\dfrac{1}{x_1}=\dfrac{r_1}{r_0}=\dfrac{a-b\lfloor a/b\rfloor}{b}=\dfrac{a}{b}- \left \lfloor \frac{a}{b} \right \rfloor= x_0-\lfloor x_0\rfloor.
		$$
		F\"ur einen Index $j\in\N$ mit $2\leq j+1<n_*$ ist auch $x_{j+1}>1$ und
		$$\dfrac{1}{x_{j+1}}=\dfrac{r_{j+1}}{r_j}=\dfrac{r_{j-1}}{r_j}- q_j= x_j-\lfloor x_j\rfloor.
		$$ 
		
		\item Es sei $0\leq j<n_*$. Zun\"achst ist $t_1=1>t_0=0$. 
                Es gilt $t_{j+1}=t_{j-1}+q_j\cdot t_j\geq t_j$ 
                f\"ur $j\geq 1$ wegen $q_j \geq 1$. 
                Nun beachten wir Bemerkung~\ref{bem:3_4} (b) 
                und setzen $x=q_j$ in Satz~\ref{satz:3_3} (b). Es folgt
		$$\langle q_0,...,q_{j-1},q_j\rangle=\frac{s_j q_j+s_{j-1}}{t_j q_j+t_{j-1}}=\frac{s_{j+1}}{t_{j+1}}.$$\\
		
		\item Aus (b) folgt wegen $q_j=\lfloor x_j\rfloor$:
		\begin{equation}\label{eq:3_1}
		x_j=q_j+\frac{1}{x_{j+1}}\quad \text{f\"ur} \;1\leq j+1<n_*.
		\end{equation} 
			
			Wir zeigen mit vollst\"andiger Induktion:
		\begin{equation}\label{eq:3_2}	
			\frac{a}{b}=\langle q_0,..,q_{j-1},x_j\rangle\quad \text{f\"ur} \;1\leq j<n_*.
		\end{equation}
			 F\"ur $j=1$ und $n_*>1$ gilt $r_0>r_1>0$ sowie $x_1=\frac{r_0}{r_1}>1$ mit $$\langle q_0,x_1\rangle=q_0+\frac{1}{x_1}=x_0=\frac{a}{b}.$$
			  Wenn f\"ur einen Index $j$ mit $1\leq j<n_*$ die Gleichung (\ref{eq:3_2}) gilt und auch noch $j+1<n_*$ bleibt, so folgt mit (\ref{eq:3_1}) und Satz~\ref{satz:3_2}:
		$$\dfrac{a}{b}=\langle q_0,..,q_{j-1},q_j+\frac{1}{x_{j+1}}\rangle=\langle q_0,..,q_{j-1},q_j,x_{j+1}\rangle,
		$$	 	
		so dass (\ref{eq:3_2}) auch f\"ur $j+1$ stimmt.\\	 	
		
		\item Wenn $\frac{a}{b}$ rational ist, k\"onnen wir nach Bemerkung~\ref{bem:3_4} (a) voraussetzen, dass $a\in\Z$ und $b\in\N$ gilt. Damit ist der erweiterte Euklidische Algorithmus mit dem einfachen Euklidischen Algorithmus aus Lektion~\ref{cha:2} vertr\"aglich und terminiert mit $n_*< \infty$.\\
		Nun sei umgekehrt $n_*<\infty$ vorausgesetzt. 

                F\"ur $n_*=1$  ist $r_1=a-b\left\lfloor\frac{a}{b}\right\rfloor=0$ 
                und somit $x_0=q_0=\frac{a}{b}=\left\lfloor\frac{a}{b}\right\rfloor$ 
                mit $\frac{a}{b}=\langle q_0 \rangle$. 

                F\"ur $n_*\geq2$ ist dagegen $r_{n_*}=r_{n_*-2}-q_{n_*-1}r_{n_*-1}=0$ 
                mit $x_{n_*-1}=q_{n_*-1}$, und wegen (d) ist 
		$$\dfrac{a}{b}=\langle q_0,...,q_{n_*-1}\rangle\in\Q.
		$$
	\end{enumerate}
\dokendProof\\

\begin{Satz}\label{satz:3_6}
	\hspace*{0cm}\\\vspace{-1cm}
	\begin{enumerate}[(a)]
		\item Es gilt $s_j t_{j+1}-t_j s_{j+1}=(-1)^j$ f\"ur $0\leq j<n_*$ mit $\ggT(s_j,t_j)=1$ f\"ur $j\in\N_0$ mit $j\leq n_*$. Insbesondere sind f\"ur $0\leq j<n_*$ die Br\"uche $\frac{s_{j+1}}{t_{j+1}}$ aus Satz~\ref{satz:3_5} (c) bereits gek\"urzt.\\
		
		\item Es gilt $j\in\N_0$ mit $j\leq n_*$:
		$$b s_j-a  t_j=(-1)^j  r_j.$$
				
		\item Ist \"uberdies $n_*<\infty$ und gilt auch noch $a\in\Z$, $b\in\N$, so ist $s_{n_*}/t_{n_*}=a/b$ mit dem gek\"urzten Bruch $s_{n_*}/t_{n_*}$ sowie mit $\ggT(a,b)=r_{n_*-1}$:
		$$(-1)^{n_*-1}\cdot \ggT(a,b)=b s_{n_*-1}-a t_{n_*-1}.
		$$
			
	\end{enumerate}
	\hfill\dokendSatz
\end{Satz}
{\bf Beweis:}~
\begin{enumerate}[(a)]
	\item Nach Satz~\ref{satz:3_3} (a) und Bemerkung~\ref{bem:3_4} (b) gilt 
	$$T_{(q_0,...,q_j)}=\begin{pmatrix}
	s_j & s_{j+1}\\
	t_j & t_{j+1}
	\end{pmatrix}\quad \text{mit}\quad \Det T_{(q_0,...,q_j)}=(-1)^j,
	$$	
	$$s_j t_{j+1}-t_j s_{j+1}=(-1)^j\quad \text{f\"ur}\quad 0\leq j<n_*.
	$$
	Hieraus folgen auch die \"ubrigen Behauptungen von (a).\\
	
	\item Behauptung (b) stimmt f\"ur $j=0$, $1$ nach Wahl der Startwerte. Stimmt (b) bis zu einem $j\geq 1$ mit $r_j>0$ und $j+1\leq n_*$, so ist sie auch f\"ur $j+1$ erf\"ullt, denn
	\begin{multline*}
	b s_{j+1}-a t_{j+1}=b(s_{j-1}+s_j q_j)-a(t_{j-1}+t_j q_j)\\
	=(-1)^{j-1} r_{j-1}+q_j (-1)^j r_j=(-1)^{j+1} r_{j+1}.
	\end{multline*}
	
	\item folgt aus (b), indem man dort $j=n_*-1$, $\ggT(a,b)=r_{n_*-1}$ bzw. $j=n_*$, $r_{n_*}=0$ einsetzt und den Euklidischen Algorithmus aus Lektion~\ref{cha:2} verwendet.
\end{enumerate}
\dokendProof\\
\begin{Beis}\label{beis:3_7}
	Hier greifen wir das Beispiel~\ref{beis:2_7} zur Berechnung von $\ggT(138,462)$ noch einmal auf:
	F\"ur die Eingabewerte $a=138$, $b=462$ erhalten wir folgende Tabelle mit $n_*=5$:
	\newline
	\begin{center}
		\begin{tabular}{|l||c|c|c|c|c|} \hline
			$j$  & ~~ $q_j$ ~~  & ~~ $r_j$ ~~ & ~~ $s_j$ ~~ & ~~ $t_j$ ~~ & $bs_j-at_j=(-1)^j r_j$ \\
			\hline
			0  &    0    &  462   &  1   &  0  &   462  \\ \hline
			1  &    3    &  138   &  0   &  1  &  $-138$  \\ \hline
			2  &    2    &  48    &  1   &  3  &   48   \\ \hline
			3  &    1    &  42    &  2   &  7  &  $-42$   \\ \hline
			4  &    7    &  6     &  3   &  10 &   6    \\ \hline
			5  &    ---  &  0     &  23  &  77 &   0    \\ \hline
		\end{tabular}
	\end{center}
	\vspace{0.25cm}
	$$q_0=\left \lfloor \frac{a}{b} \right \rfloor=0,$$
	$$r_0=b=462,\; r_1=a-b\cdot q_0=138, $$
	$$s_0=1,\;t_0=0,\;s_1=0,\;t_1=1, $$
	und f\"ur $j=1,...,4$: $$q_j=\left\lfloor\dfrac{r_{j-1}}{r_j}\right\rfloor,$$
	$$ r_{j+1}=r_{j-1}-q_j r_j,$$
	$$ s_{j+1}=s_{j-1}+q_j s_j,$$
	$$ t_{j+1}=t_{j-1}+q_j t_j.$$\\
		
	F\"ur $n_*< \infty $ ist eine Spalte mit den Werten $x_j$ entbehrlich, stattdessen f\"ugen wir die Kontrollspalte $bs_j-at_j=(-1)^j r_j$ ein.\\	
	
	Hier haben wir $$\frac{a}{b}=\frac{138}{462}=\frac{s_5}{t_5}=\frac{23}{77}\; \text{~mit~}\;\ggT(23,77)=1,$$ $\ggT(138,462)=3\cdot 462 -10\cdot 138=6$, da $n_*-1=4$ gerade ist.
	\dokendSatz
\end{Beis}

F\"ur die Theorie der Kettenbr\"uche sind wir in Satz~\ref{satz:3_3} vor allem an nat\"urlichen Zahlen $\lambda_0,...,\lambda_{j-1}$ interessiert. Die einfachste Wahl ist $\lambda_0,...,\lambda_{j-1}=1$, was 
in Verbindung mit $\tau_0=0$, $\tau_1=1$ und der Rekursionsvorsicht $\tau_{k+1}=\tau_k+\tau_{k-1}$ zur Bildung der Fibonacci-Zahlen f\"uhrt, die in der Theorie der Kettenbr\"uche eine besondere Rolle spielen:
\begin{Def}
	Die Fibonacci-Zahlen\index{Fibonacci-Zahlen}\label{Fibonacci-Zahlen2} $f_k$ sind f\"ur $k\in\N_0$ rekursiv erkl\"art gem\"a{\ss} $f_0=0$, $f_1=1$ sowie $f_{k+1}=f_k+f_{k-1}$. Eine Liste der Anfangswerte lautet:
	\begin{center}
		\begin{tabular}{|r||c|c|c|c|c|c|c|c|c|c|c|c|c|c|c|}\hline
		
			  ~~ $k$~ ~& ~ 0 ~ & ~ 1 ~ & ~ 2 ~ & ~ 3 ~ & ~ 4 ~ & ~ 5 ~ & ~ 6 ~ & ~ 7 ~  
				& ~ 8 ~  & ~ 9 ~ & ~ 10 ~ & ~11~ & ~12~ & ~13~& ~14~\rule{0pt}{3.ex} \\\hline
			 ~~ $f_k$ ~~&  0  &  1 & 1 & 2 & 3 & 5 & 8 & 13 
			 & 21 & 34 & 55 & 89 & 144 & 233 & 377 \rule{0pt}{3.ex}\\ \hline
		\end{tabular}
	\end{center}
		\vspace{0.25cm}
		Da $f_k\geq 1$ f\"ur $k\geq1$ gilt, ist die Folge $(f_k)_{k\geq2}$ ab $k=2$ streng monoton wachsend und insbesondere unbeschr\"ankt.
		\dokendDef
\end{Def}
\begin{Beis}\label{beis:3_9}
	Wir wenden auf die positive der beiden L\"osungen $\lambda_{\pm}:=\frac{1\pm\sqrt{5}}{2}$ der quadratischen Gleichung $\lambda^2=\lambda +1$ den erweiterten Euklidischen Algorithmus an:\\	
	F\"ur $a:=\lambda_+=\frac{1+\sqrt{5}}{2}$, $b:=1$ erhalten wir die Tabelle
	\begin{center}
			\begin{tabular}{|c|c|c|c|c|c|} \hline
		 ~ $j$ ~ & ~~~ $x_j$ ~~~  & ~~~ $q_j$ ~~~ & ~~~ $r_j$ ~~~ & ~~~ $s_j$ ~~~ & ~~~ $t_j$ ~~~ \\
			\hline
			0  & \rule{0pt}{5ex}   $\dfrac{1+\sqrt{5}}{2}$    &  1   &  1   &  1  &   0  \\[0.4cm] \hline
			1  &  \rule{0pt}{5ex}  $\dfrac{1+\sqrt{5}}{2}$    &  1   &  $\dfrac{\sqrt{5}-1}{2}$  &  1  &  1  \\[0.4cm] \hline
			2  &  \rule{0pt}{5ex}  $\dfrac{1+\sqrt{5}}{2}$    &  1    &  $\left(\dfrac{\sqrt{5}-1}{2}\right)^2$  &  2  &   1  	  \\ [0.4cm]\hline
			3  &  \rule{0pt}{5ex}  $\dfrac{1+\sqrt{5}}{2}$    &  1    &  $\left(\dfrac{\sqrt{5}-1}{2}\right)^3$  &  3  &  2   \\[0.4cm] \hline
			...  &    ...    &  ...     &  ...   &  ... &   ...    \\ \hline
				\end{tabular}
	\end{center}

	\vspace{0.25cm}
	Hier ist $n_*=\infty$, denn nach Satz~\ref{satz:3_5} (b) ist 
	$$ x_{j+1}=\frac{1}{x_j-\lfloor x_j\rfloor}=x_j=\frac{\sqrt{5}+1}{2}$$
	f\"ur $j\in\N_0$ mit einer konstanten Folge $(x_j)_{j\in\N_0}$. Hierbei ist 
	$q_j=\lfloor x_j\rfloor=1$ f\"ur alle $j\in\N_0$. Aus $r_0=1$ und $\frac{r_{j-1}}{r_j}=x_j=\frac{\sqrt{5}+1}{2}$ f\"ur $j\in\N$ folgt:
	$$ r_j=\left(\frac{\sqrt{5}-1}{2}\right)^j	\quad \mbox{~f\"ur~alle~} j\in\N_0.$$
	Hier ist $s_j=t_{j+1}=f_{j+1}$ f\"ur alle $ j\in\N_0$, und Satz~\ref{satz:3_6} (b) liefert:
	\begin{equation}\label{eq:3_3}
	1\cdot s_j-\lambda_+\cdot t_j=f_{j+1}-\lambda_+ \cdot f_j=\lambda_-^j\quad \mbox{~f\"ur~alle~} j\in\N_0.
	\end{equation}
	Nun erhalten wir aus (\ref{eq:3_3}) f\"ur $j\in\N_0$:
	$$\begin{array}{rcl}
		(\lambda_+-\lambda_-)\cdot f_{j+1}+\lambda_-^{j+1}& = & (\lambda_+-\lambda_-)\cdot (\lambda_+\cdot f_j+\lambda_-^j)+\lambda_-^{j+1}  \\
				& = & (\lambda_+-\lambda_-)\cdot \lambda_+\cdot f_j+\lambda_+\cdot \lambda_-^j  \\
				& = & \lambda_+\cdot\left[(\lambda_+-\lambda_-)\cdot f_j+\lambda_-^j\right]. \\
	\end{array}$$
	Mit $(\lambda_+-\lambda_-)\cdot f_0+\lambda_-^0=1$ folgt hieraus sofort
	$$(\lambda_+-\lambda_-)\cdot f_{j}+\lambda_-^{j}=\lambda_+^j\quad \mbox{~f\"ur~alle~} j\in\N_0,$$
	und somit die Binetsche Formel\index{Binetsche Formel}\label{Binetsche_Formel2} f\"ur die Fibonacci-Zahlen:
		\begin{equation}\label{eq:3_4}
	f_j=\frac{1}{\sqrt{5}}\cdot\left[\left(\frac{1+\sqrt{5}}{2}\right)^j-\left(\frac{1-\sqrt{5}}{2}\right)^j\right]\quad \mbox{~f\"ur~alle~} j\in\N_0.
		\end{equation}
	\dokendSatz
\end{Beis}
\begin{DefSatz}
	Gegeben seien $\lambda_0\in\Z$ sowie eine unendliche Folge $(\lambda_j)_{j\in\N}$ nat\"urlicher Zahlen $\lambda_j$. Dann existiert der sogenannte unendliche Kettenbruch\index{Kettenbruch}\label{Kettenbruch2}\index{unendlicher Kettenbruch}\label{unendlicher Kettenbruch}
	$$
	\langle\lambda_0,\lambda_1, \lambda_2, ...\rangle:=\lim\limits_{j\rightarrow \infty}\langle\lambda_0,\lambda_1, ...,\lambda_j\rangle,$$
	den man auch in der Form
	$$
	\langle\lambda_0,\lambda_1, \lambda_2, ...\rangle=\lambda_0+ \cfrac{1}{\lambda_1 + \cfrac{1}{\lambda_2 + ...}}$$
	schreibt.	
\hfill\dokendSatz
\end{DefSatz}
{\bf Beweis:}~	
Wir setzen $x:=\lambda_j$ und w\"ahlen $x'\geq 1$ beliebig im Satz~\ref{satz:3_3}, und erhalten unter Beachtung von Satz~\ref{satz:3_2}:
\begin{equation}\label{eq:3_5}
	\langle\lambda_0,\lambda_1, ...,\lambda_j\rangle-\langle\lambda_0,\lambda_1, ...,\lambda_{j-1}+\frac{1}{x'}\rangle= \frac{(-1)^j(\lambda_j-x')}{\tau_{j+1}(\tau_j x'+\tau_{j-1})}	
\end{equation}
Im Limes $x'\rightarrow\infty$ erhalten wir aus (\ref{eq:3_5}):
\begin{equation}\label{eq:3_6}
d_j:=\langle\lambda_0,\lambda_1, ...,\lambda_j\rangle-\langle\lambda_0,\lambda_1, ...,\lambda_{j-1}\rangle= \frac{(-1)^{j+1}}{\tau_{j}\tau_{j+1}}	\quad \mbox{~f\"ur~alle~} j\in\N.
\end{equation}
Setzen wir noch $d_0:=\langle\lambda_0\rangle=\lambda_0$, so folgt aus (\ref{eq:3_6}):
\begin{equation}\label{eq:3_7}
\langle\lambda_0,\lambda_1, ...,\lambda_j\rangle=\sum\limits_{k=0}^{j}d_k=\lambda_0+\sum\limits_{k=1}^{j}\frac {(-1)^{k+1}}{\tau_{k}\cdot\tau_{k+1}}.
\end{equation}
Nun gilt f\"ur alle $k\in\N$ die Monotoniebeziehung
$$\frac {\tau_{k+1}\cdot \tau_{k+2}}{\tau_{k}\cdot\tau_{k+1}}=\frac{\tau_{k+2}}{\tau_{k}}=\frac {\tau_{k}+\lambda_{k+1} \tau_{k+1}}{\tau_{k}}>1
$$
sowie $\lim\limits_{k\rightarrow\infty}\frac{1}{\tau_k\cdot\tau_{k+1}}=0$. Das Leibniz-Kriterium liefert die Konvergenz der alternierenden Reihe in (\ref{eq:3_7}).\dokendProof\\
\begin{Bem}
	Es gilt $\tau_j \geq f_j$ f\"ur $j\in\N_0$, so dass die alternierende Reihe in (\ref{eq:3_7}) sogar absolut konvergiert. Nach Beispiel~\ref{beis:3_9} gilt insbesondere
	$$\langle1,1,1,...\rangle=1+ \cfrac{1}{1 + \cfrac{1}{1 + ...}}=\frac{\sqrt{5}+1}{2}.
	$$
~~\dokendBem
\end{Bem}
\begin{Satz}[Eindeutigkeit unendlicher Kettenbr\"uche]\label{satz:3_12}
	Es seien $\lambda_0\in\Z$ und $\lambda_j\in\N$ f\"ur $j\in\N$. Dann gilt:
	\begin{enumerate}[(a)]
	\item $\displaystyle\langle\lambda_0,\lambda_1, \lambda_2,...\rangle=\lambda_0+\frac{1}{\langle\lambda_1,\lambda_2,\lambda_3, ...\rangle}=\langle\lambda_0,\langle\lambda_1,\lambda_2,\lambda_3,...\rangle\rangle$.\\
	\item Wir setzen $y_j:=\langle\lambda_j,\lambda_{j+1},\lambda_{j+2},...\rangle$ f\"ur $j\in\N_0$. Dann gilt $\lambda_j<y_j<\lambda_j+1$ und $\lambda_j=\lfloor y_j\rfloor$ sowie $\displaystyle y_{j+1}=\frac{1}{y_j-\lfloor y_j \rfloor}$ f\"ur alle $j\in\N_0$.\\
	\item Wendet man den erweiterten Euklidischen Algorithmus auf die Eingabewerte $a=x_0:=\langle\lambda_0,\lambda_1,\lambda_2, ...\rangle$, $b:=1$ an, so folgen $n_*=\infty$, $x_j=y_j$ und $q_j=\lambda_j$ f\"ur alle $j\in\N_0$:\\
	Der unendliche Kettenbruch\index{unendlicher Kettenbruch}\label{unendlicher Kettenbruch2} hat eine eindeutige Darstellung
	$$\langle\lambda_0,\lambda_1, \lambda_2, ...\rangle=\lambda_0+ \cfrac{1}{\lambda_1 + \cfrac{1}{\lambda_2 + ...}}\;,$$
	und er liefert eine Irrationalzahl\index{Irrationalzahl}\label{Irrationalzahl} $x_0$.	
	\end{enumerate}
\hfill\dokendSatz	
\end{Satz}
{\bf Beweis:}~
\begin{enumerate}[(a)]
	\item $\langle\lambda_0,\lambda_1, ...,\lambda_j\rangle=\lambda_0+\frac{1}{\langle\lambda_1,...,\lambda_{j}\rangle}$ liefert im Limes $j\rightarrow\infty$ die Behauptung.\\
	\item Nach (a) gilt $y_j=\lambda_j+\frac{1}{y_{j+1}}$ mit 
	$$y_{j+1}=\lambda_{j+1}+\frac{1}{\langle\lambda_{j+2},\lambda_{j+3},...\rangle}>\lambda_{j+1}\geq 1,\quad j\in\N_0. $$
	Hieraus folgen $\displaystyle 0<\frac{1}{y_{j+1}}<1$, $\lambda_j=\lfloor y_j\rfloor$ sowie $\displaystyle y_{j+1}=\frac{1}{y_j-\lfloor y_j \rfloor}$ $\mbox{~f\"ur~alle~} j\in\N_0$.\\
	\item Die Startwerte und die Rekursionsvorschriften stimmen f\"ur beide Zahlenfolgen $(y_j)_{j\in\N_0}$ und $(x_j)_{j\in\N_0}$\"uberein. Somit gelten der Reihe nach 
        $$x_j=y_j\,, \quad q_j=\lfloor x_j\rfloor=\lfloor y_j\rfloor=\lambda_j \quad \mbox{~f\"ur~alle~} j\in\N_0\,.$$
Nach Satz~\ref{satz:3_5} (e) ist $x_0$ irrational mit $n_*=\infty$.
\end{enumerate}
\dokendProof\\

\begin{Bem}
	Es seien $\lambda\in\Z$, $\lambda'\in\N$ sowie $x\geq 1$, $y>1$ reell. Dann ist $\langle\lambda',x\rangle >1$, und es gilt $\lambda=\lfloor\langle\lambda,y\rangle\rfloor=\lfloor\langle\lambda,\lambda',x\rangle\rfloor$, so dass $\lambda$ sowohl durch $\langle\lambda,y\rangle$ als auch durch $\langle\lambda,\lambda',x\rangle$ eindeutig bestimmt ist. Wendet man diese Beziehungen schrittweise auf einen endlichen Kettenbruch
	$$ \rho=\langle\lambda_0, ...,\lambda_{j-1},1\rangle\quad \text{mit} \quad \lambda_0\in\Z,\;\lambda_1, ...,\lambda_{j-1}\in\N $$
	f\"ur $j\geq2$ an, so folgt, dass die rationale Zahl $\rho$ genau zwei Kettenbruchdarstellungen besitzt, n\"amlich
	\begin{equation}\label{eq:3_8}
	\rho=\langle\lambda_0, ...,\lambda_{j-1},1\rangle=\langle\lambda_0,...,\lambda_{j-2},\lambda_{j-1}+1\rangle, \quad j\geq2.
	\end{equation}
	Ebenso gilt
	\begin{equation}\label{eq:3_9}
	\langle\lambda_0,1\rangle=\langle\lambda_0+1\rangle,
	\end{equation}
	wobei sich jedes $\rho\in\Q$ entweder gem\"a{\ss} (\ref{eq:3_8}) oder (\ref{eq:3_9}) schreiben l\"asst. So liefern Satz~\ref{satz:3_5} (c) und Beispiel~\ref{beis:3_7}:
	$$\dfrac{3}{10}=\langle0,3,2,1\rangle=\langle0,3,3\rangle,$$
	$$\dfrac{23}{77}=\langle0,3,2,1,7\rangle=\langle0,3,2,1,6,1\rangle.$$
	Nur die Irrationalzahlen besitzen eine eindeutige Darstellung als (unendlicher) Kettenbruch.
	~~\dokendBem
\end{Bem}
\begin{Satz}\label{satz:3_14}
	F\"ur eine Irrationalzahl $x_0$ w\"ahlen wir $a:=x_0$, $b:=1$ 
als Eingabewerte des erweiterten Euklidischen Algorithmus. 
Seine Ausgabewerte bezeichnen wir wieder mit $x_j$, $q_j$ bzw. $r_j$, $s_j$, $t_j$ 
f\"ur $j\in\N_0$. Dann folgt $x_0=\langle q_0,q_1,q_2,...\rangle$, genauer f\"ur alle $j\in\N$:
\begin{enumerate}[(a)]
	\item Es gilt $$x_0-\langle q_0,..., q_{j-1}\rangle=x_0-\frac{s_j}{t_j}=(-1)^{j+1}\frac{r_j}{t_j}\quad 
	\mbox{~f\"ur~alle~} j\in\N$$
	mit der streng monoton fallenden positiven Nullfolge $\left(\frac{r_j}{t_j}\right)_{j\in\N}$.\\
	\item F\"ur die verallgemeinerten Divisionsreste gilt $$r_j=\frac{1}{t_j x_j+t_{j-1}}\quad \text{mit}\quad 0<r_j<\frac{1}{q_j\cdot t_j}$$
	und der streng monoton fallenden positiven Nullfolge $(r_j)_j\in\N$.\\
	\item Es besteht die Absch\"atzung $$\left|x_0-\frac{s_j}{t_j}\right|<\frac{1}{q_j t_{j}^2}$$
	mit den gek\"urzten N\"aherungsbr\"uchen $\frac{s_j}{t_j}=\langle q_0,...,q_{j-1}\rangle$ zu $x_0$.	
\end{enumerate}		
\hfill\dokendSatz	
\end{Satz}
{\bf Beweis:}~
\begin{enumerate}[(a)]
	\item folgt sofort aus Satz~\ref{satz:3_5} (a), (c) und Satz~\ref{satz:3_6} (b).\\
	\item In Satz~\ref{satz:3_3} (c) setzen wir $\lambda_k=q_k$ f\"ur $0\leq k<j$ bzw. $\tau_k=t_k$ f\"ur $0\leq k\leq j$ sowie $x=x_j$, und erhalten mit Satz~\ref{satz:3_5} (d):
	$$x_0-\langle q_0,...,q_{j-1},x'\rangle=\dfrac{(-1)^{j}(x_j-x')}{(t_j x_j+t_{j-1})(t_j x'+t_{j-1})}, $$
		und f\"ur $x'\rightarrow\infty$ im Limes:
	$$x_0-\langle q_0,...,q_{j-1}\rangle=\dfrac{(-1)^{j+1}}{t_j(t_j x_j+t_{j-1})}. $$
		Der Vergleich mit der Teilaussage (a) dieses Satzes liefert 
		$$0<r_j=\dfrac{1}{t_j x_j+t_{j-1}}\leq\dfrac{1}{t_j x_j}<\frac{1}{q_j t_j} \quad \mbox{~f\"ur~alle~} j\in\N,$$
		womit $(r_j)_{j\in\N}$ auch eine streng monoton fallende Nullfolge ist.\\
		\item folgt direkt aus (a) und (b).
\end{enumerate}
\dokendProof

\begin{Bem}
\hspace*{0cm}\vspace{0cm}	
\begin{enumerate}[(a)]
	\item Nach Satz~\ref{satz:3_14} (a) stellen die gek\"urzten Br\"uche $\frac{s_j}{t_j}$ N\"aherungsbr\"uche f\"ur $x_0$ dar, die abwechselnd kleiner bzw. gr\"o{\ss}er als $x_0$ sind. Diese Br\"uche liegen bei ungeradem Index $j$ links von $x_0$, und bei geradem Index $j$ rechts von $x_0$ gem\"a{\ss}
	$$\dfrac{s_1}{t_1}<\dfrac{s_3}{t_3}<\dfrac{s_5}{t_5}<...<x_0<\dfrac{s_6}{t_6}<\dfrac{s_4}{t_4}<\dfrac{s_2}{t_2}. $$
	Da allgemein $t_j\geq f_j$ mit der $j-$ten Fibonacci-Zahl\index{Fibonacci-Zahlen}\label{Fibonacci-Zahlen3} $f_j$ und $f_j\sim \frac{1}{\sqrt{5}}\left(\frac{1+\sqrt{5}}{2}\right)^j$ f\"ur $j \rightarrow\infty$ gilt, also 
	$$\lim\limits_{j\rightarrow\infty}\left(\dfrac{1}{f_j}\cdot\dfrac{1}{\sqrt{5}}\left(\frac{1+\sqrt{5}}{2}\right)^j\right)=1, $$
	erhalten wir:
	$$\left|x_0-\dfrac{s_j}{t_j}\right|<\dfrac{1}{f_j\cdot f_{j+1}}, $$
	und die $\dfrac{s_j}{t_j}=\langle q_0,...,q_{j-1}\rangle$ konvergieren zumindest exponentiell schnell gegen $x_0$.\\
	\item Satz~\ref{satz:3_14} (c) ist eine Versch\"arfung des klassischen Dirichletschen Approxima\-tionssatzes.\index{Dirichletscher Approximationssatz}\label{Dirichletscher Approximationssatz} Dieser besagt, dass f\"ur jede Irrationalzahl $x_0$ unendlich viele rationale Zahlen $\dfrac{s}{t}$ mit $\left|x_0-\dfrac{s}{t}\right|<\dfrac{1}{t^2}$ und $s\in\Z$, $t\in\N$ existieren.	
\end{enumerate}		
	~~\dokendBem
\end{Bem}
\begin{Satz}[Die Medianteneigenschaft\index{Mediantensatz}\label{Mediantensatz}]\label{satz:3_16}
	Gegeben sind $u,u'\in\Z$ und $v,v'\in\N$ mit $u' v-uv'=1$. Dann gilt $\frac{u}{v}<\frac{u'}{v'}$, und der sogenannte Mediant\index{Mediant}\label{Mediant} $\frac{u+u'}{v+v'}$ ist unter allen Br\"uchen $\frac{s}{t}$ mit $\frac{u}{v}<\frac{s}{t}<\frac{u'}{v'}$ und $s\in\Z$, $t\in\N$ der einzige mit dem kleinsten Nenner $t$.
\hfill\dokendSatz	
\end{Satz}
{\bf Beweis:}~
Wir zeigen, dass die Ungleichungen 
\begin{equation}\label{eq:3_10}
\dfrac{u}{v}<\dfrac{u+u'-\beta}{v+v'-\alpha}<\frac{u'}{v'}
\end{equation}
mit den Nebenbedingungen
\begin{equation}\label{eq:3_11}
\alpha,\;\beta\in\Z\quad\text{und}\quad 0\leq\alpha<v+v'
\end{equation}
die einzige L\"osung $\alpha=\beta=0$ besitzen.\\

Die linke Ungleichung von (\ref{eq:3_10}) ist \"aquivalent zu $u(v+v'-\alpha)<v(u+u'-\beta)$, also zu
\begin{equation}\label{eq:3_12}
v\beta-u\alpha<u'v-uv'=1.
\end{equation}
Entsprechend ist die rechte Ungleichung in (\ref{eq:3_10}) \"aquivalent zu 
\begin{equation}\label{eq:3_13}
u'\alpha-v'\beta<1.
\end{equation}
F\"ur $\alpha=\beta=0$ sind (\ref{eq:3_12}), (\ref{eq:3_13}) und somit auch (\ref{eq:3_10}) erf\"ullt.\\

Mit (\ref{eq:3_12}) und (\ref{eq:3_13}) folgt wegen $v,v'\geq1$ aus $\alpha=0$ schon $\beta=0$. Wir nehmen daher $\alpha\in\N$ an und m\"ussen diese Annahme zum Widerspruch f\"uhren:\\
Da $\ggT(u,v)=1$ ist, w\"urde aus $v\beta=u\alpha$ die Beziehung $v\arrowvert \alpha$ folgen, und hieraus $\alpha=\lambda v$, $\beta=\lambda u$ mit einem $\lambda\in\N$, was der Beziehung (\ref{eq:3_13}) widerspricht.\\

Mit der Ganzzahligkeit aller Gr\"o{\ss}en folgt somit $v \beta-u \alpha<0$ aus (\ref{eq:3_12}), d.h.
\begin{equation*}
\dfrac{\beta}{\alpha}<\dfrac{u}{v}.
\end{equation*}
Dies widerspricht wegen $\frac{u}{v}<\frac{u'}{v'}$ der aus (\ref{eq:3_13}) resultierenden Ungleichung
\begin{equation*}
\dfrac{u'}{v'}\leq\dfrac{\beta}{\alpha}.
\end{equation*}
\dokendProof\\
\begin{Satz}[Satz von den rationalen Bestapproximationen]\label{satz:3_17}
Mit den Voraussetzungen und Bezeichnungsweisen von Satz~\ref{satz:3_14}
gilt f\"ur alle $s\in\Z$ und $t\in\N$:\\
Aus $t\leq t_j$ sowie aus $\frac{s}{t}\neq \frac{s_j}{t_j}$ folgt
$$|t_j x_0-s_j|<|t x_0-s|\quad \text{f\"ur} \quad j\geq 2. $$
 Dies besagt, dass die N\"aherungsbr\"uche $\frac{s_j}{t_j}$ der Kettenbruchentwicklung von $x_0$ stets die besten Approximationen\index{rationale Bestapproximation}\label{rationale Bestapproximation} an $x_0$ mit rationalen Zahlen liefern. 
\hfill\dokendSatz	
\end{Satz}
{\bf Beweis:}~
Wir verwenden hier entscheidend die verallgemeinerten Divisionsreste $r_k$. Neben der Darstellung in Satz~\ref{satz:3_14} (b) brauchen wir 
\begin{equation}\label{eq:3_14}
r_k=|t_k x_0-s_k|\quad \mbox{~f\"ur~alle~} k\in\N_0.
\end{equation}
F\"ur $k=0$ ist dies klar und folgt f\"ur $k\in\N$ aus Satz~\ref{satz:3_14} (a). Den Fall $\frac{s}{t}=\frac{s_j}{t_j}$ haben wir ausgeschlossen.\\

Wegen $j\geq2$ k\"onnen wir $\frac{s_{j-1}}{t_{j-1}}$ bilden, denn es ist $t_{j-1}\geq1$.\\

Den Fall, dass $\frac{s}{t}$ echt zwischen $\frac{s_{j-1}}{t_{j-1}}$ und $\frac{s_{j}}{t_{j}}$ liegt, k\"onnen wir nach dem Mediantensatz~\ref{satz:3_16}\index{Mediantensatz}\label{Mediantensatz2} wegen $t\leq t_{j}$ ebenfalls ausschlie{\ss}en.\\

Nun betrachten wir den m\"oglichen Fall $\frac{s}{t}=\frac{s_{j-1}}{t_{j-1}}$ und beachten, dass $r_{j-1}>r_j$ nach Satz~\ref{satz:3_5} (a) gilt. In diesem Falle folgt bereits die Behauptung mit Verwendung von (\ref{eq:3_14}) wegen $s=\lambda\cdot s_{j-1}$, $t=\lambda\cdot t_{j-1}$ mit einem $\lambda\in\N$:
$$|t x_0-s|=\lambda\cdot r_{j-1}>r_j=|t_j x_0-s_j|. $$
Aus Satz~\ref{satz:3_3} (b), Bemerkung~\ref{bem:3_4} (b) und Satz~\ref{satz:3_5} (d) gewinnen wir folgende Darstellungsformel:
\begin{equation}\label{eq:3_15}
x_0=\dfrac{s_j x_j+s_{j-1}}{t_j x_j+t_{j-1}}.
\end{equation}
F\"ur die Position von $\frac{s}{t}$ m\"ussen wir nun nur noch zwei F\"alle unterscheiden:
\begin{enumerate}[\text{Fall} A:]
\item $\dfrac{s}{t}<\dfrac{s_{j-1}}{t_{j-1}}$ und $\dfrac{s}{t}<\dfrac{s_{j}}{t_{j}}$	bzw.\\
\item $\dfrac{s}{t}>\dfrac{s_{j-1}}{t_{j-1}}$ und $\dfrac{s}{t}>\dfrac{s_{j}}{t_{j}}$.
\end{enumerate}
In beiden F\"allen haben die beiden Terme $t s_j-s t_j$ und $t s_{j-1}-s t_{j-1}$ dasselbe Vorzeichen, und wir erhalten die entscheidende Absch\"atzung
$$\begin{array}{rclcl}
|t x_0-s|& = & \left|\dfrac{(t s_j-s t_j)x_j+(t s_{j-1}-s t_{j-1})}{t_j x_{j}+ t_{j-1}}\right| & (\text{wegen (\ref{eq:3_15})}) & \\
\\
& = & \dfrac{|t s_j-s t_j|x_j+|t s_{j-1}-s t_{j-1}|}{t_j x_{j}+ t_{j-1}} & (\text{Fall A bzw. B})\\
\\
& = & r_j x_j |t s_{j}-s t_{j}|+r_j|t s_{j-1}-s t_{j-1}|& (\text{Satz}~\ref{satz:3_14} \text{ (b)}) & \\
& > & r_j \quad (\text{wegen } x_j>1,\; |t s_{j}-s t_{j}|\geq 1 ), &
\end{array} $$
also auch im Falle A bzw. B:
$$|t x_0-s|>|t_j x_0-s_j|. $$
\dokendProof\\
\begin{Bem}
	Aus Satz~\ref{satz:3_17} folgt insbesondere f\"ur $j\geq2$:
	$$\left|x_0-\dfrac{s}{t}\right|\geq\dfrac{t}{t_j}\left|x_0-\dfrac{s}{t}\right|=\dfrac{1}{t_j}|t x_0-s|>\left|x_0-\dfrac{s_j}{t_j}\right|. $$
~~\dokendBem
\end{Bem}
Die folgende Version des Approximationssatzes von Hurwitz orientiert sich an Perrons Lehrbuch
\cite[\S 14]{perron}.
\begin{Satz}[Der Approximationssatz von Hurwitz\index{Approximationssatz von Hurwitz}\label{Approximationssatz von Hurwitz}]\label{satz:3_18}
	Es sei $x_0$ eine Irrationalzahl\index{Irrationalzahl}\label{Irrationalzahl2}. Hiermit w\"ahlen wir $a:=x_0$, $b:=1$ als Eingabewerte f\"ur den erweiterten Euklidischen Algorithmus. % Wir verwenden die Notationen und Resultate aus $\S4$.
	 Dann gilt:
	
	Von drei aufeinanderfolgenden N\"aherungsbr\"uchen zu $x_0$ hat mindestens einer, sagen wir $\frac{s_j}{t_j}$ mit $j\in\N$, die Eigenschaft 
	\begin{equation}\label{eq:3_16}
	\left|x_0-\frac{s_j}{t_j}\right|<\frac{1}{\sqrt{5}\,t_j^2}.
	\end{equation}
	Insbesondere gibt es zu jeder Irrationalzahl $x_0$ unendlich viele N\"aherungsbr\"uche $s_j/t_j$, die der Absch\"atzung \eqref{eq:3_16} gen\"ugen.
	\hfill\dokendSatz	
\end{Satz}
{\bf Beweis:}~
F\"ur jedes $j\geq 1$ finden wir nach Satz~\ref{satz:3_14}(a), (b) genau eine Zahl $\delta_j$ mit
\begin{equation}\label{eq:3_17}
x_0-\frac{s_j}{t_j}=(-1)^{j+1}\frac{\delta_j}{t_j^2},\quad 0<\delta_j <1.
\end{equation}
Aus \eqref{eq:3_17} folgt unter Beachtung von Satz~\ref{satz:3_6} (a) f\"ur $j\geq 2:$
$$
\frac{\delta_j}{t_j^2}+\frac{\delta_{j-1}}{t_{j-1}^2}=(-1)^{j+1}\left[x_0-\frac{s_j}{t_j}\right]-(-1)^{j+1}\left[x_0-\frac{s_{j-1}}{t_{j-1}}\right]$$
$$
=(-1)^{j}\left(\frac{s_j}{t_j}-\frac{s_{j-1}}{t_{j-1}}\right)=\frac{1}{t_j t_{j-1}},
$$
was wir auch in folgender Form schreiben k\"onnen:
\begin{equation}\label{eq:3_18}
\delta_{j-1}\left(\frac{t_j}{t_{j-1}}\right)^2-\frac{t_j}{t_{j-1}}+\delta_j=0,\quad j\geq2.
\end{equation}
Dies ist eine quadratische Gleichung in $t_j/t_{j-1}$, und deren Aufl\"osung 
ergibt mit $\sigma_j\in\{1,-1\}$:
\begin{equation}\label{eq:3_19}
\frac{t_j}{t_{j-1}}=\frac{1+\sigma_j\sqrt{1-4\delta_j\delta_{j-1}}}{2\delta_{j-1}}, \quad \frac{t_{j-1}}{t_j}=\frac{1-\sigma_j\sqrt{1-4\delta_j\delta_{j-1}}}{2\delta_{j}}.
\end{equation}
F\"ur nat\"urliche Zahlen $k\geq 2$ verwenden wir die erste Gleichung in \eqref{eq:3_19} f\"ur \mbox{$j=k+1$}, die zweite f\"ur $j=k$, und erhalten
\begin{equation*}
\frac{t_{k+1}}{t_k}-\frac{t_{k-1}}{t_k}=\frac{\sigma_{k+1}\sqrt{1-4\delta_{k+1}\delta_{k}}+\sigma_k\sqrt{1-4\delta_k\delta_{k-1}}}{2\delta_{k}}.
\end{equation*}
Unter Beachtung von $t_{k+1}=t_{k-1}+q_k t_k$ folgt
\begin{equation}\label{eq:3_20}
2\delta_k q_k=\sigma_{k+1}\sqrt{1-4\delta_{k+1}\delta_{k}}+\sigma_k\sqrt{1-4\delta_k\delta_{k-1}}, \quad k\geq2.
\end{equation}
Nun k\"onnen in \eqref{eq:3_18} keine zwei aufeinanderfolgenden Koeffizienten $\delta_{j-1}$, $\delta_j$ \"ubereinstimmen, da andernfalls
\begin{equation*}
\delta_{j-1}=\delta_j=\frac{t_j}{t_{j-1}}\cdot \frac{1}{1+\left(\frac{t_j}{t_{j-1}}\right)^2}
\end{equation*}
rational w\"are. Nehmen wir schlie{\ss}lich $\delta_{k-1}, \delta_{k}, \delta_{k+1}\geq\dfrac{1}{\sqrt{5}}$ in \eqref{eq:3_20} an, so erhalten wir aus dieser Gleichung den Widerspruch
	\begin{equation*}
	0<2\delta_k q_k<\sqrt{1-\frac{4}{5}}+\sqrt{1-\frac{4}{5}}=\frac{2}{\sqrt{5}},
	\end{equation*}
	denn es ist $q_k\geq 1$, und wenigstens eine der drei Zahlen $\delta_{k-1}$, $\delta_{k}$, $\delta_{k+1}$ 
	muss gr\"o{\ss}er als $\dfrac{1}{\sqrt{5}}$ sein.
\dokendProof\\

\begin{Bem}[~Zus\"atze zum Approximationssatz von Hurwitz]\label{bem:zusatz}\hspace*{0cm}\\\vspace{-1cm}
	\begin{enumerate}[(a)]
	\item Obiger Beweis zeigt auch, dass von je zwei aufeinanderfolgenden N\"aherungsbr\"uchen der Kettenbruchentwicklung von $x_0$ wenigstens einer, etwa $s_j /t_j$ mit $j\in\N$, die Eigenschaft
	\begin{equation*}
	\left|x_0-\frac{s_j}{t_j}\right|<\frac{1}{2 \, t_j^2}
	\end{equation*}
	besitzt, da $\delta_{j-1}=\delta_j$ ausgeschlossen ist und nach \eqref{eq:3_18} die quadratische Gleichung
	\begin{equation*}
	\delta_{j-1}x^2-x+\delta_j=0
	\end{equation*}
	die reelle L\"osung $x=\dfrac{t_j}{t_{j-1}}$ besitzt, so dass in \eqref{eq:3_19} immer Quadratwurzeln von nichtnegativen reellen Zahlen gebildet werden.
	\item In der Absch\"atzung \eqref{eq:3_16} 
kann die Konstante $\sqrt{5}$ durch keine gr\"o{\ss}ere Zahl ersetzt werden, 
wie das Beispiel $x_0:=\dfrac{\sqrt{5}+1}{2}$ lehrt: 
Hier ist $x_j=x_0$ und $t_j=f_j$, $s_j=f_{j+1}$ ($j$-te Fibonacci-Zahl $f_j$) 
f\"ur alle $j\in\N_0$, siehe Beispiel~\ref{beis:3_9}. 	
	Hier liefert der Satz~\ref{satz:3_14} (a), (b) f\"ur alle $j\in\N$:
	\begin{equation*}
	x_0-\frac{s_j}{t_j}=x_0-\frac{f_{j+1}}{f_j}
=\frac{(-1)^{j+1}}{f_j^2(x_0+\frac{f_{j-1}}{f_j})}
	\end{equation*}
	mit $\displaystyle\lim\limits_{j\rightarrow\infty}
(x_0+\frac{f_{j-1}}{f_j})=\frac{\sqrt{5}+1}{2}+\frac{2}{\sqrt{5}+1}=\sqrt{5}$.

\end{enumerate}
	~~\dokendBem
\end{Bem}

        \section{Historische Anmerkungen}\label{cha:3B}

\begin{enumerate}[$\bullet$]
	\item Der niederl\"andische Astronom, Mathematiker und Physiker \textit{Christiaan Huygens}\index{Huygens, Christiaan}\label{Christiaan Huygens} (1629-1695) 
verwendete die besten damals verf\"ugbaren Daten zur Konstruk\-tion eines mechanischen Modells unseres Sonnensystems. 
Dabei kamen Kettenbr\"uche zur Berechnung von Kalendern und Schaltjahren zum Einsatz.\\
	\item \textit{Leonard Euler}\index{Euler, Leonard}\label{Leonard Euler} (1707-1783) entwickelte 1737 
in seinem Werk "`{De Fractionibus Continuis Dissertatio}"' eine Theorie, die auch allgemeinere Kettenbr\"uche der Form
$$ a_0+ \cfrac{b_1}{a_1 + \cfrac{b_2}{a_2 + \cfrac{b_3}{a_3+...}}}
$$
		beinhaltet. F\"ur die nach ihm benannte Eulersche Zahl $e$ 
hat er die regelm\"a{\ss}ige Kettenbruchentwicklung $$e-1=<1,1,2,1,1,4,1,1,6,...>$$ angegeben. 
Auch wusste er, dass Kettenbruchentwicklungen, die in eine Periode einm\"unden, quadratische Irrationalzahlen darstellen.\\
	\item Neben Euler hat vor allem \textit{Joseph-Louis Lagrange}
	 \index{Lagrange, Joseph-Louis} \label{Joseph-Louis Lagrange}
(1736-1813) die  Theorie der Kettenbr\"uche vorangetrieben. 
Der Satz~\ref{satz:3_17} von den rationalen Bestapproximationen einer Irrationalzahl mit den endlichen Kettenbruchentwicklungen 
geht auf ihn aus dem Jahre 1770 zur\"uck. Er bewies, dass reell quadratische Irrationalzahlen eine Kettenbruchentwicklung besitzen, 
die in eine Periode m\"undet.\\
	\item \textit{Carl Friedrich Gau{\ss}}
	\index{Gau{\ss}, Carl Friedrich}\label{Carl Friedrich Gauss} (1777-1855) entwickelte in seinen "`{Disquisitiones Arithmeticae}"' von 1801,
  siehe \cite{gauss}, eine einheitliche Grundlage der Zahlentheorie. In seinem Werk nimmt die Theorie der quadratischen Formen\index{quadratische Form}\label{quadratische Form}
	$$F(x,y)=ax^2+bxy+cy^2\quad\text{mit}\;\; a,b,c\in\Z $$
	einen breiten Raum ein; in verkleideter Form hat diese Theorie, zumindest im indefiniten Fall $D:=b^2-4ac>0$, 
weitreichende Bez\"uge zur Kettenbruchentwicklung der reell quadratischen Irrationalzahlen\index{reell quadratische Irrationalzahl}\label{reell quadratische Irrationalzahl2} $\dfrac{\sqrt{D}-b}{2a}$,
siehe hierzu insbesondere noch den Abschnitt~\ref{cha:8} des vorliegenden Lehrbuches.\\
	 \item Das Lehrbuch von \textit{Oskar Perron}  \index{Perron, Oskar} \label{Oskar Perron}
\cite{perron} (1880-1975) "`{Die Lehre von den Kettenbr\"uchen}"' 
erschien 1913 im Teubner Verlag. Es ist bis heute eine wertvolle Einf\"uhrung in die Theorie geblieben. 
Im Vorwort dieses Buches werden die Arbeiten Eulers als Inspirationsquelle hervorgehoben:\\
	 
	 "`{Besonders die Arbeiten Eulers \"uber Kettenbr\"uche erweisen sich als wahre Fundgrube 
f\"ur h\"ochst interessante Beziehungen sowohl zwischen verschiedenen Kettenbr\"uchen 
als auch zwischen Kettenbr\"uchen und Reihen oder bestimmten Integralen; ...}"'.\\
	 
	 Die dritte verbesserte und erweiterte Auflage von Perrons Lehrbuch ist ab 1957 als Werk in zwei B\"anden erh\"altlich.
	
\end{enumerate}

        \section{Aufgaben}\label{cha:3_A}

\begin{Auf}[Erweiterter Euklidischer Algorithmus\index{erweiterter Euklidischer Algorithmus}\label{erweiterter_Euklidischer_Algorithmus2}]\label{auf:3_1}
Mit dem erweiterten Euklidischen Algorithmus ermittle man
zu den beiden teilerfremden Zahlen $a=7$, $b=81$ zun\"achst ein und dann
alle Zahlenpaare $(\lambda, \mu) \in \Z^2$ mit $81\lambda-7\mu=1$.\\
\end{Auf}

{\bf L\"osung:}\\
$$a=7,\  b=81:\quad r_0=b=81,\  r_1=a=7.$$
\begin{equation*}
\begin{tabular}{|l||c|c|c|c|} \hline
$j$  & ~~ $q_j$ ~~  & ~~ $r_j$~~  & ~~ $s_j$~~  & ~~ $t_j$~~\\
\hline
0  &    0    &  81   &    1    &  0 \\ \hline 
1  &    11   &  7    &    0    &  1 \\ \hline
2  &    1    &  4    &    1    &  11  \\ \hline
3  &    1    &  3    &    1    &  12  \\ \hline
4  &    3    &  1    &    2    &  23  \\ \hline
5  &    ---  &  0    &    7    &  81  \\ \hline
\end{tabular}\quad\quad\quad
\begin{array}{rl}
&\text{Es ist } n_\ast=n_\ast(a,b)=5.\\[0.3em]
&\text{Da } n_\ast-1 \text{ gerade ist, folgt}\\[0.3em]
&b s_4-a t_4=81\cdot 2-7\cdot 23=+1,\text{ also}\\[0.3em]
&81\lambda_0-7\mu_0=1 \text{ f\"ur } \lambda_0=2, \mu_0=23.\\
\end{array}
\end{equation*}
Nun m\"ogen $\lambda,\mu\in\Z$ eine weitere L\"osung von $81\lambda-7\mu=1$ ergeben.
Dann folgt durch Subtraktion:
\begin{equation*}
\begin{split}
&81 (\lambda-\lambda_0)-7(\mu-\mu_0)=0, \text{ d.h.}\\
&81 (\lambda-\lambda_0)=7(\mu-\mu_0).
\end{split}
\end{equation*} 
Wegen $\ggT(7,81)=1$ ist $81|\mu-\mu_0$ und $7|\lambda-\lambda_0$.
Setzen wir $\lambda-\lambda_0=7k$ mit $k\in\Z$, so folgt $\mu-\mu_0=81k$, und alle L\"osungen $(\lambda,\mu)\in\Z^2$ mit $81\lambda-7\mu=1$ sind gegeben durch 
\begin{equation*}
\lambda=2+7k, \quad \mu=23+81k \quad \text{mit}\ \ k\in\Z,
\end{equation*}
wie man durch eine Rechenprobe best\"atigt.

\begin{Auf}[Fibonacci-Folge\index{Fibonacci-Folge}\label{Fibonacci-Folge2}, Teil 3]\label{auf:3_2}

Wir wenden den erweiterten Euklidischer Algorithmus
auf die beiden Eingabewerte $a \in \Z$ und $b \in \N$ an
und verwenden dabei die Notationen im Hauptteil dieses Abschnitts, 
insbesondere sei $n_*(a,b) \in \N$ der zugeh\"orige
Abbruchindex.\\

Wir betrachten $f_0=0$, $f_1 = 1$,
$f_{n+2}=f_{n+1}+f_n$ f\"ur $n \in \N_0$.\\

Man zeige, dass f\"ur alle $n \in \N$ mit $n \geq 2$ folgende Aussage gilt:

F\"ur alle $a \in \Z$ und $b \in \N$ mit $b \leq f_n$ ist $n_*(a,b) \leq n-1$,
wobei $n_*(a,b)= n-1$ genau dann eintritt, wenn es ein $\tilde{q}_0 \in \Z$ gibt mit 
$a=f_{n-1}+\tilde{q}_0 f_n$ und $b=f_{n}$.\\

{\it Bemerkung:} Der Euklidische Algorithmus mit ganzen Zahlen als Einga\-be\-wer\-ten
hat insbesondere f\"ur $a=f_{n-1}$, $b=f_n$ und $n \geq 2$ die ung\"unstigste Laufzeit
mit $n_*(f_{n-1},f_n)= n-1$.
\end{Auf}

{\bf L\"osung:}\\
Wir zeigen vorab
\begin{equation}\label{eq:3_21}
1<f_{n+1}/f_n<2 \quad \text{f\"ur}\ n\geq 3. 
\end{equation}
Tabelle der ersten $5$ Fibonacci-Zahlen:
\begin{equation*}
\begin{tabular}{|l||c|c|c|c|c|} \hline
$n$  & ~~$0$~~  &  ~~$1$~~  &  ~~$2$~~  & ~~$3$~~ & ~~$4$~~ \\
\hline
$f_n$  &    0    &    1    &    1   &   2   &  3 \\ \hline 
\end{tabular}\quad\quad\quad
\begin{array}{rl}
&\text{F\"ur }n=3 \text{ stimmt die } \\[0.3em]
&\text{Behauptung wegen }\frac{f_4}{f_3}=\frac{3}{2}.
\end{array}
\end{equation*}

Wird sie f\"ur ein $n\geq 3$ als richtig angenommen, so folgt
\begin{equation*}
\frac{f_{n+2}}{f_{n+1}}=\frac{f_{n+1}+f_n}{f_{n+1}}=1+\frac{f_n}{f_{n+1}}
%\quad \text{mit}\quad 0<\frac{f_n}{f_{n+1}}<1
\end{equation*}
mit $0<\frac{f_n}{f_{n+1}}<1$. Die Behauptung gilt dann auch f\"ur $n+1$, womit \eqref{eq:3_21} bewiesen ist.\\

F\"ur $a\in\Z$ und $b\in\N$ gilt
\begin{equation}\label{eq:3_22}
n_\ast(a,b)=n_\ast(r_1,r_0)
\end{equation}
mit $r_0=b$, $r_1=a-b\left\lfloor\frac{a}{b}\right\rfloor$, $0\leq r_1<r_0$.\\

F\"ur $n=2$ gilt die Behauptung der Aufgabe mit $f_{n-1}=f_n=1$, $b=r_0=1$ mit $a\in\Z$ und $n_\ast(a,b)=n_\ast(a,1)=1$.\\

Wir m\"ussen die Behauptung nur noch f\"ur $n\geq 3$ zeigen. F\"ur jedes $n\geq 3$ ist die Behauptung nach~\eqref{eq:3_22} zur folgenden Aussage $\mathcal{A}(n)$ \"aquivalent:\\

F\"ur alle $r_1\in\N_0$ und alle $r_0\in\N$ mit $r_1<r_0\leq f_n$ ist $n_\ast(r_1,r_0)\leq n-1$, wobei $n_\ast(r_1,r_0)=n-1$ genau f\"ur $r_1=f_{n-1}$ und $r_0=f_n$ eintritt.\\

Diese Aussage $\mathcal{A}(n)$ beweisen wir induktiv f\"ur alle $n\geq 3$.\\

{\bf Induktionsanfang:}~F\"ur $n=3$ wird $0\leq r_1<r_0\leq 2$ wegen $f_3=2$ mit den ganzen Zahlen $r_1,r_0$ vorausgesetzt. $r_1=0$, $r_0=1,2$ liefert $n_\ast(r_1,r_0)=1<3-1$ mit \mbox{$r_1\neq f_{3-1}$}. 
F\"ur $r_1=1=f_{3-1}$ und $r_0=2=f_3$ ist $n_\ast(r_1,r_0)=2=3-1$. 
Insgesamt gilt damit $\mathcal{A}(3)$.\\

{\bf Induktionsschritt:}~Wir nehmen $\mathcal{A}(k)$ f\"ur $3\leq k\leq n$ und ein $n\geq 3$ an. Um damit $\mathcal{A}(n+1)$ zu zeigen, setzen wir 
\begin{equation}\label{eq:3_23}
0\leq r_1<r_0\leq f_{n+1}
\end{equation}
voraus. Gem\"a{\ss} \eqref{eq:3_23} unterscheiden wir drei F\"alle:\\

{\bf Fall A:}~ $r_1=0$.\\
Hier ist $n_\ast(r_1,r_0)=n_\ast (0,r_0)=1<n$, und auch f\"ur $r_0=f_{n+1}$ haben wir $r_1=0<f_n$.\\

Da  $\mathcal{A}(n+1)$ im Falle A gilt, werden wir im Folgenden $r_1>0$ voraussetzen, so dass gilt:
\begin{equation}\label{eq:3_24}
n_\ast(r_1,r_0)=n_\ast(r_2,r_1)+1.
\end{equation}
Hierbei ist $r_2$ der auf $r_0,r_1$ folgende Divisionsrest im Euklidischen Algorithmus. Es gilt
\begin{equation}\label{eq:3_25}
r_0=r_1\cdot \left\lfloor\frac{r_0}{r_1}\right\rfloor+r_2, \quad 0\leq r_2<r_1<r_0.
\end{equation}\\

{\bf Fall B:}~ $0<r_1<f_n$.\\
Wir wenden die Induktionsannahme $\mathcal{A}(n)$ auf das Zahlenpaar $r_2,r_1$ an, und erhalten aus \eqref{eq:3_24}:
$$n_\ast(r_1,r_0)<(n-1)+1=n.$$
Unter Beachtung von $r_1\neq f_n$ gilt hier $\mathcal{A}(n+1)$.\\

{\bf Fall C:}~ $f_n\leq r_1<r_0\leq f_{n+1}$.\\
Hier folgt $\left\lfloor\frac{r_0}{r_1}\right\rfloor=1$ aus \eqref{eq:3_21}, und weiter aus \eqref{eq:3_25}:
\begin{equation} \label{eq:3_26}
\left. \begin{tabular}{l}
$0<r_2=r_0-r_1\leq f_{n+1}-f_n=f_{n-1}$,\\
wobei $\ r_2=f_{n-1}\ $ genau f\"ur\\
$r_1=f_n\ $ und $\ r_0=f_{n+1}$ gilt.	 
\end{tabular}\right\}
\end{equation}
Hier k\"onnen wir den Divisionsrest $r_3\geq 0$ bilden mit 
\begin{equation}\label{eq:3_27}
n_\ast (r_2,r_1)=n_\ast(r_3,r_2)+1.
\end{equation}
F\"ur $n=3$ folgt $\mathcal{A}(4)$ aus $0<r_2=f_{3-1}=1$, $r_1=f_3=2$, $r_0=f_4=3$ gem\"a{\ss}~\eqref{eq:3_26}, so dass wir nun $n\geq 4$ voraussetzen. Aus~\eqref{eq:3_27}, \eqref{eq:3_24} erhalten wir
\begin{equation}\label{eq:3_28}
n_\ast(r_1,r_0)=n_\ast(r_3,r_2)+2.
\end{equation}
mit $r_2\leq f_{n-1}$. Wegen $n\geq 4$ gilt $\mathcal{A}(n-1)$, so dass
$$n_\ast(r_1,r_0)\leq(n-2)+2=n$$
mit \eqref{eq:3_28} folgt. Nehmen wir $n_\ast(r_1,r_0)=n$ an, d.h. 
$n_\ast(r_3,r_2)=n-2$, so erhalten wir $r_2=f_{n-1}$ aus $\mathcal{A}(n-1)$, und somit auch $r_1=f_n$, $r_0=f_{n+1}$ aus \eqref{eq:3_26}.
Aus $r_1=f_n$, $r_0=f_{n+1}$ folgen umgekehrt $r_2=f_{n-1}$ wegen \eqref{eq:3_26} und schlie{\ss}lich \mbox{$n_\ast(r_2,r_1)=n-1$} bzw. $n_\ast(r_1,r_0)=n$ wegen $\mathcal{A}(n)$.

\begin{Auf}[Quadratische Irrationalzahlen\index{quadratische Irrationalzahl}\label{quadratische Irrationalzahl}]\label{auf:3_3}

Das quadratische Polynom $P(x) := ax^2+bx+c$ mit $a,b,c \in \Z$
besitze die Diskriminante $D := b^2-4ac >0$, die keine Quadratzahl sei.
Es sei $f :=  \lfloor \sqrt{D} \rfloor$. 
Dem Polynom $P$ ordnen wir die Nullstelle 
$\begin{displaystyle}
x_P := \frac{\sqrt{D}-b}{2a}
\end{displaystyle}$
zu. Man zeige:
\begin{enumerate}[(a)]
	\item Es ist $x_P$ eine Irrationalzahl.\\
	
	\item F\"ur jedes $q \in \Z$ besitzt auch $Q(x):=ax^2+(b+2aq)x+(c+q(b+aq))$
	die positive Diskriminante $D$, und es gilt $x_Q=x_P-q=\frac{\sqrt{D}-(b+2aq)}{2a}$. \\
	
	\item $R(x):=-cx^2-bx-a$
	hat ebenfalls dieselbe Diskriminante $D$ wie $P$ und $Q$, und es gilt $x_R=1/x_P
	=\frac{\sqrt{D}+b}{-2c}$. \\
	
	\item
	\begin{equation*}
	\left \lfloor x_P \right \rfloor =
	\begin{cases}
	\displaystyle ~ \left \lfloor\frac{f-b}{2a}\right \rfloor\,,&
	\text{$a>0$}\,,\\[2ex]
	\displaystyle ~ \left \lfloor\frac{b-(f+1)}{2|a|}\right \rfloor\,,&
	\text{$a<0$}\,.\\
	\end{cases}
	%~\text{und}\quad
	%\left \lfloor\frac{1}{x_P} \right \rfloor =
	%\begin{cases}
	%\displaystyle ~ \left \lfloor\frac{f+b}{-2c}\right \rfloor\,,&
	%\text{$c<0$}\,,\\[2ex]
	%\displaystyle ~ \left \lfloor-\frac{f+b+1}{2c}\right \rfloor\,,&
	%\text{$c>0$}\,.\\
	%\end{cases}
	\end{equation*}
\end{enumerate}
\end{Auf}

{\bf L\"osung:}
\begin{enumerate}[(a)]
	\item Angenommen $\frac{\sqrt{D}-b}{2a}=\frac{s}{t}$ mit $s\in\Z$, $t\in\N$.
	Dann ist $t\sqrt{D}-tb=2as$ bzw. \mbox{$t^2 \cdot D=t'^2$} mit $t'=|tb+2as|>0$.
	Nach dem Satz von der eindeutigen Primfaktorzerlegung w\"are dann $D=p_1^{\alpha_1}\cdot ... \cdot p_j^{\alpha_j}$ mit paarweise verschiedenen Primzahlen und geraden Exponenten $\alpha_1,...,\alpha_j\in\N$, also $D$ im Widerspruch zur Annahme eine Quadratzahl. Somit ist $\frac{\sqrt{D}-b}{2a}$ eine Irrationalzahl.\\
	
	\item Die Diskriminante von $Q(x)=ax^2+(b+2aq)x+(c+q(b+aq))$ ist 
	\begin{equation*}
	\begin{split}
    (b+2aq)^2-4a(c+q(b+aq))
    &=b^2+4abq+4a^2q^2-4ac-4aqb-4a^2q^2\\
	&=b^2-4ac=D,
	\end{split}
	\end{equation*}
	damit gilt auch
	$$x_Q=\frac{\sqrt{D}-(b+2aq)}{2a}=\frac{\sqrt{D}-b}{2a}-q=x_P-q.$$
	\item Es ist klar, dass auch $R(x)=-cx^2-bx-a$ die Diskriminante $D$ hat mit
	$$x_P\cdot x_R=\frac{\sqrt{D}-b}{2a}\cdot \frac{\sqrt{D}+b}{-2c}=\frac{D-b^2}{(-4ac)}=1.$$
	
	\item $\D \left\lfloor x_P\right\rfloor
	=\left\lfloor \frac{\sqrt{D}-b}{2a}\right\rfloor
	=\left\lfloor \frac{\left\lfloor \sqrt{D}-b\right\rfloor}{2a}\right\rfloor
	=\left\lfloor \frac{\left\lfloor \sqrt{D}\right\rfloor-b}{2a}\right\rfloor
	=\left\lfloor \frac{f-b}{2a}\right\rfloor$\\
	folgt f\"ur $a>0$, d.h. $2a\in\N$, aus Aufgabe~\ref{auf:2_4} (c) und (b).\\
	
	F\"ur $a<0$ erhalten wir entsprechend
	$$ \left\lfloor x_P\right\rfloor
	=\left\lfloor \frac{b-\sqrt{D}}{2|a|}\right\rfloor
	=\left\lfloor \frac{ b+\left\lfloor-\sqrt{D}\right\rfloor}{2|a|}\right\rfloor
	=\left\lfloor \frac{b-\left\lceil \sqrt{D}\:\right\rceil}{2|a|}\right\rfloor
	=\left\lfloor \frac{b-(f+1)}{2|a|}\right\rfloor,$$
	denn $D$ ist keine Quadratzahl und somit $\left\lceil \sqrt{D}\:\right\rceil=f+1$, siehe Aufgabe~\ref{auf:2_4} (d).\\	 
\end{enumerate}

{\bf Vorbereitung zur Bearbeitung der Aufgabe~\ref{auf:3_4}}\\

Wir f\"uhren in tabellarischer Form die Kettenbruchentwicklung\index{Kettenbruchentwicklung}\label{Kettenbruchentwicklung2} \index{Kettenbruch}\label{Kettenbruch3}
einer qua\-dra\-tischen Irrational\-zahl an einem Beispiel vor, 
und verwenden die Resultate der Aufgabe~\ref{auf:3_3}. 
Analog soll dann bei der L\"osung der 
folgenden Aufgabe verfahren werden.\\

Wir entwickeln
$x_0 := \frac{\sqrt{2}+10}{14}$, 
indem wir die Folge $x_{j+1}=\frac{1}{x_j - q_j}$
mit den Divi\-sions\-koeffi\-zienten $q_j = \lfloor x_j \rfloor$ f\"ur $j \in \N_0$ bilden. 
In der folgenden Tabelle ist $x_j=x_{P_j}$ f\"ur
$P_j(x)=a_jx^2+b_jx+c_j$.
Aus den Koeffizienten $a_j, b_j, c_j $ des Poly\-noms $P_j(x)$
berechnen wir zun\"achst $q_j$ mit den Fallunterscheidungen $a_j>0$ bzw. $a_j<0$ gem\"a{\ss} 
der Formel mit den Gau{\ss}-Klammern aus Aufgabe \ref{auf:3_3}(d).
Mit den Notationen $x_j^*:=x_j - q_j=x_{P_j^*}$ f\"ur $P_j^*(x)=a_j^*x^2+b_j^*x+c_j^*$
berechnen wir hierauf die Koeffizienten $a_j^*, b_j^*, c_j^* $ des Poly\-noms $P_j^*(x)$ mit Hilfe von $q_j$
aus den Koeffizienten $a_j, b_j, c_j $ des Poly\-noms $P_j(x)$ gem\"a{\ss} Aufgabe \ref{auf:3_3}(b).
Schliesslich erhalten wir die neuen Koeffizienten $a_{j+1}=-c_j^*$, $b_{j+1}=-b_j^*$, $c_{j+1}=-a_j^*$
in der Folgezeile der Tabelle gem\"a{\ss} Aufgabe \ref{auf:3_3}(c), d.h. es gilt $x_{j+1}=1/x_j^*$
bzw. $x_{j+1}x_j^*=1$ f\"ur alle $j \in \N_0$, wobei die letzte Beziehung als einfache Rechenprobe 
bei der Erstellung dieser Tabellen dient.\\

Zu Beginn wird $P_0(x)=a_0x^2+b_0x+c_0$ mit $a_0,b_0,c_0 \in \Z$ und $x_{0}=x_{P_0}$ ermittelt:
Wir haben $(14x_0-10)^2=2$, d.h. $196x_0^2-280x_0+98=0$. Nach K\"urzung des
Faktors $14$ k\"onnen wir $a_0=14$, $b_0=-20$, $c_0=7$ mit $P_0(x)=14x^2-20x+7$ und $D=8$ w\"ahlen,
denn es gilt $x_0=x_{P_0}$. Wir erhalten $f = \lfloor \sqrt{D} \rfloor=2$.
\begin{center}
	\begin{tabular}{|c||c|c|c|c||c||c|c|c|c|c||c|} \hline
		~~$j$~~  & \rule{0pt}{2.7ex}   ~~$a_j$~~   &  ~~$b_j$~~  &  ~~$c_j$~~ &  ~~$x_j$~~ &  ~~$q_j$~~    &  ~~$a_j^*$~~   &  ~~$b_j^*$~~   &  ~~$c_j^*$~~  & ~~$x_j^*=x_j-q_j$~~   \\[0.3em] \hline \hline
		
		0  &   14   &  -20  &   7 &\rule{0pt}{3.5ex} $\frac{\sqrt{8}+20}{28}$ &  0        &  14   &  -20   &  7    &  $\frac{\sqrt{8}+20}{28}$ \\[0.3em] \hline %\hline
		
		1  &  -7    &  20   &  -14   &  \rule{0pt}{3.5ex}   $\frac{\sqrt{8}-20}{-14}$   &   1    &  -7    &  6   &  -1     &    $\frac{\sqrt{8}-6}{-14}$   \\[0.3em] \hline %\hline
		
		2      &  1    &  -6   &  7   &\rule{0pt}{3.5ex}    $\frac{\sqrt{8}+6}{2}$ &   4   &  1    &  2   &  -1   &  $\frac{\sqrt{8}-2}{2}$   \\[0.3em] \hline %\hline
		
		3      &  1    &  -2   &  -1    & \rule{0pt}{3.5ex}   $\frac{\sqrt{8}+2}{2}$ &  2  &  1    &  2   &  -1  
		&  $\frac{\sqrt{8}-2}{2}$  \\[0.3em] \hline
	\end{tabular}
\end{center}

Wir erhalten die periodische Kettenbruchentwicklung\index{periodische Kettenbruchentwicklung}\label{periodische Kettenbruchentwicklung}
\begin{equation*}
\frac{\sqrt{2}+10}{14}=\langle 0,1,4,\overline{2} \rangle\,.
\end{equation*}

\begin{Auf}[Kettenbruchentwicklung quadratischer Irrationalzahlen]\label{auf:3_4}
	
Man wende den obigen f\"ur qua\-dra\-tische Irra\-tio\-nalzahlen $x_0$ 
formulierten Ketten\-bruch\-algorithmus auf $x_0 := \sqrt{7}$ an. 
Damit zeige man, dass $\sqrt{7}$ eine perio\-dische
Ketten\-bruch\-ent\-wick\-lung besitzt, und gebe diese explizit an.
\end{Auf}

{\bf L\"osung:} Wir haben $D=28$, $f=5$ f\"ur $P_0(x)=x^2-7$ und
$x_0=\sqrt{7}=x_{P_0}$.
\begin{center}
	\begin{tabular}{|c||c|c|c|c||c||c|c|c|c|} \hline
		~~ $j$ ~~  &  ~~$a_j$~~   &  ~~$b_j$~~  &  ~~ $c_j$~~  & \rule{0pt}{3.5ex} ~~ $x_j$~~   &  ~~$q_j$ ~~
		&  ~~$a_j^*$~~   &  ~~$b_j^*$~~   & ~~ $c_j^*$~~  & ~~$x_j^*=x_j-q_j$~~  \\[0.3em]	\hline\hline
		0  &  1   &  0  &   -7 &   \rule{0pt}{3.5ex} $\frac{\sqrt{28}-0}{2}$   &  2 &  1   &  4   &  -3  & \rule{0pt}{3.5ex} $\frac{\sqrt{28}-4}{2}$ \\[0.3ex]\hline
		1 &  3    &  -4   &  -1  & \rule{0pt}{3.5ex}   $\frac{\sqrt{28}+4}{6}$ &  1 & 3    &  2   &  -2  & \rule{0pt}{3.5ex}    $\frac{\sqrt{28}-2}{6}$  \\[0.3ex]\hline
		2 &   2    &  -2   &  -3  &\rule{0pt}{3.5ex}    $\frac{\sqrt{28}+2}{4}$ & 1 &  2    &  2   &  -3   & \rule{0pt}{3.5ex}   $\frac{\sqrt{28}-2}{4}$   \\[0.3ex]\hline 
		3 &   3    &  -2   &  -2  &\rule{0pt}{3.5ex}   $\frac{\sqrt{28}+2}{6}$  & 1  &    3    &  4   &  -1  &\rule{0pt}{3.5ex}   $\frac{\sqrt{28}-4}{6}$\\[0.3ex]\hline
		4 &   1    &  -4   &  -3  & \rule{0pt}{3.5ex}   $\frac{\sqrt{28}+4}{2}$  & 4 & 1    &  4   &  -3  & \rule{0pt}{3.5ex}   $\frac{\sqrt{28}-4}{2}$  \\[0.3ex]\hline
	\end{tabular}
\end{center}
Die Tabelle liefert $\sqrt{7}=\langle 2,\overline{1,1,1,4} \rangle$.
Hierbei ist folgendes zu beachten: Da $x_4^*$ mit $x_0^*$ \"ubereinstimmt,
stimmt $x_5=1/x_4^*$ wieder mit $x_1=1/x_0^*$ \"uberein, so dass ab dem Index $j=1$ eine Periode der L\"ange 4
im Kettenbruch auftritt. Erst in Lektion \ref{cha:8} werden wir mit einem etwas einfacheren modifizierten Verfahren
zeigen, dass genau die quadratischen Irrationalzahlen\index{quadratische Irrationalzahl}\label{quadratische Irrationalzahl2}
eine periodische Kettenbruchentwicklung besitzen.

\begin{Auf}[Ein periodischer Kettenbruch mit zwei Parametern] \label{auf:3_4}
Gegeben sind zwei nat\"urliche Zahlen $a$ und $b$. Man berechne den periodischen
Kettenbruch $\langle \overline{a,b} \rangle$.
\end{Auf}

{\bf L\"osung:} Die zweimalige Anwendung von Satz \ref{satz:3_12}(a) auf 
$\xi =\langle \overline{a,b} \rangle>0$ liefert
\begin{equation*}
\xi =\langle a,  \langle \overline{b,a} \rangle \rangle 
= a + \frac{1}{\langle \overline{b,a} \rangle}
=a + \frac{1}{b + \frac{1}{\xi}}\,.
\end{equation*}
Hieraus erhalten wir f\"ur $\xi$ die quadratische Gleichung
\begin{equation*}
b\xi^2-ab\xi-a=0
\end{equation*}
mit den beiden L\"osungen
\begin{equation*}
\xi_{1,2}= \frac{a}{2} \pm \sqrt{\frac{a^2}{4}+\frac{a}{b}}\,.
\end{equation*}
Da $\xi$ positiv ist, folgt eindeutig
\begin{equation*}
\xi = \frac{a}{2} + \sqrt{\frac{a^2}{4}+\frac{a}{b}}\,.
\end{equation*}

\chapter{Farey-Sequenzen}\label{cha:4}

Wenn wir die gek\"urzten Br\"uche zwischen Null und Eins mit einem
vorgegeben maxi\-malen nat\"urlichen Nenner $n$ der Gr\"o\ss e nach ordnen, so
erhalten wir etwa f\"ur den maxi\-malen Nenner $n=5$:
\begin{equation*}
\frac{0}{1} < \frac{1}{5} < \frac{1}{4} < \frac{1}{3} <
\frac{2}{5} < \frac{1}{2} < \frac{3}{5} < \frac{2}{3} <
\frac{3}{4} < \frac{4}{5} < \frac{1}{1} \,.
\end{equation*}
Dem Geologen John Farey\index{Farey, John}\label{John Farey} (geboren 1766 in Woburn, Bedfordshire, England und
gestorben am 6. Januar 1826 in London, England) fiel beim Betrachten solcher
Folgen von Br\"uchen, die bis heute seinen Namen tragen, folgendes auf:
Bei drei aufeinanderfolgenden Br\"uchen ergibt sich
der Wert des mittleren Bruches als Quotient aus Z\"ahler- und Nennersumme 
von den beiden links und rechts benachbarten Br\"uchen, z.B.
\begin{equation*}
\frac{1}{3}=\frac{1+2}{4+5}\,, \quad 
\frac{2}{5}=\frac{1+1}{3+2}\,, \quad 
\frac{1}{2}=\frac{2+3}{5+5}\,.
\end{equation*}
Diese anhand von Beispielen erkannte Eigenschaft hielt Farey in einem 1816 erschienenen Artikel 
mit dem Titel ``On a curious property of vulgar fractions'' im {\em Philosophical Magazine} fest.
Der franz\"osische Mathematiker Cauchy las Fareys Aufsatz und lieferte noch im selben Jahr 
1816 den bei Farey fehlenden Beweis.

Farey war nicht der erste, der diese Eigenschaft erkannt hat. So schrieb Haros 1802
einen Artikel \"uber Dezimalbr\"uche, aus dem hervorgeht, dass er Fareys ``curious property'' f\"ur
$n=99$ verwendet hat. 

Im Folgenden untersuchen wir die nach Farey benannten Sequenzen von Br\"uchen.
Auch wenn wir hierf\"ur nur einfache
Rechenregeln f\"ur den Umgang mit Br\"uchen und Ungleichungen verwenden, 
wie sie aus dem Schulunterricht bekannt sind,
so erhalten wir dennoch eine F\"ulle interessanter zahlentheoretischer Eigenschaften.
Diese werden schlie\ss lich dazu verwendet, einfache Rechenschemata zu entwickeln, mit denen sich auch
beliebige Ausschnitte aus der $n$-ten Farey-Sequenz bzw. die besten Approximationen einer Irrationalzahl
mit Farey-Br\"uchen sehr effizient berechnen lassen.\\
       \section{Farey-Sequenzen}\label{cha:4T}
Wir beginnen mit der folgenden\\

{\bf Aufgabenstellung:}
~ F\"ur gegebenes $n\in\N$ finde und ordne man der Gr\"o{\ss}e nach alle gek\"urzten Br\"uche 
$\dfrac{a}{b}$ mit $0\leq\dfrac{a}{b}\leq 1$, $a\in\N_0$, $b\in\N$ und $b\leq n$.\\

Zur L\"osung konstruieren wir Zeile f\"ur Zeile folgende Tabelle:
\begin{enumerate} [(a)] 
	\item In der ersten Zeile steht nur $\dfrac{0}{1}$ und $\dfrac{1}{1}$.
	\item Die $n$-te Zeile bildet man, indem man die ($n-1$)-te Zeile noch einmal darunter schreibt und den Medianten $\dfrac{a+a'}{b+b'}$ zwischen die aufeinanderfolgenden Br\"uche $\dfrac{a}{b}$ und $\dfrac{a'}{b'}$ der abgeschriebenen ($n-1$)-ten Zeile setzt, falls $b+b'=n$ ist.
\end{enumerate}

\begin{Def}\label{def:4_1} Die $n$-te Zeile dieser Tabelle 
nennt man die Farey-Sequenz\index{Farey-Sequenz}\label{Farey-Sequenz} (Farey-Folge) $\F_n$ der Ordnung $n$.	
	\dokendDef
\end{Def} 

\begin{Beis}\label{beis:4_2} Konstruktion der Farey-Sequenzen bis zur Ordnung $n=5$:		
		\begin{center}	
			\begin{tabularx}{0.65\textwidth}{ X|*{11}{X}| }
				\hline
				\multicolumn{1}{|c||}{$n$}
				&  \multicolumn{11}{c|}{Br\"uche der $n$-ten Farey-Sequenz $\F_n$} \\ \hline
				\multicolumn{1}{|c||}{$1$} &
				\rule{0pt}{3.5ex}
				$\dfrac{0}{1}$ & & & &  &  & & & & &  $\dfrac{1}{1}$ \\[0.2cm]\hline
				
				\multicolumn{1}{|c||}{$2$} &
				\rule{0pt}{3.5ex} 
				$\dfrac{0}{1}$ & & & &  & $\dfrac{1}{2}$ & & & & &  $\dfrac{1}{1}$ \\ [0.2cm]\hline	
				
				\multicolumn{1}{|c||}{$3$} &
				\rule{0pt}{3.5ex} 
				$\dfrac{0}{1}$ & & & $\dfrac{1}{3}$ &  & $\dfrac{1}{2}$ & & $\dfrac{2}{3}$ & & &  $\dfrac{1}{1}$\\ [0.2cm]\hline
				
				\multicolumn{1}{|c||}{$4$} &
				\rule{0pt}{3.5ex}
				$\dfrac{0}{1}$ & & $\dfrac{1}{4}$  & $\dfrac{1}{3}$ &  &$\dfrac{1}{2}$ & & $\dfrac{2}{3}$ & $\dfrac{3}{4}$ & &  $\dfrac{1}{1}$ \\ 
				[0.2cm]\hline	
				
				\multicolumn{1}{|c||}{$5$} &
				\rule{0pt}{3.5ex}
				$\dfrac{0}{1}$ & $\dfrac{1}{5}$ & $\dfrac{1}{4}$ & $\dfrac{1}{3}$ & $\dfrac{2}{5}$ & $\dfrac{1}{2}$ & $\dfrac{3}{5}$ & $\dfrac{2}{3}$ & $\dfrac{3}{4}$ & $\dfrac{4}{5}$ &  $\dfrac{1}{1}$ 
				\\ [0.2cm]\hline
			\end{tabularx}
		\end{center}
\end{Beis}

{\bf Ziel:}~ Die Konstruktion der $n$-ten Farey-Sequenz $\F_n$ liefert der Gr\"o{\ss}e nach alle gek\"urzten Br\"uche $\dfrac{a}{b}$ von $\dfrac{0}{1}$ bis $\dfrac{1}{1}$ mit den Nennern $b\leq n$.

\newlength{\formel}
\settowidth{\formel}{$a'(b+b') - (a+a')b'$}

\begin{Satz}\label{satz:4_3}
		Sind $\dfrac{a}{b}$ und $\dfrac{a'}{b'}$ aufeinanderfolgende Br\"uche der $n$-ten
		Zeile, so gilt $a'b - ab' = 1$.
	\dokendSatz
\end{Satz}
{\bf Beweis:}~
	Wir beweisen den Satz mit Hilfe der vollst\"andigen Induktion.
	Der Satz gilt f\"ur $n = 1$ (Induktionsanfang). Wir nehmen an, er gilt bis zur
	Zeile $n-1$. Die aufeinanderfolgenden Br\"uche der Zeile $n$ sind dann
	\begin{equation*}
	\frac{a}{b}, \frac{a'}{b'}\quad \text{bzw.} \quad \frac{a}{b}, \frac{a + a'}{b + b'}
	\quad \text{bzw.} \quad \frac{a + a'}{b + b'}, \frac{a'}{b'}\, ,
	\end{equation*}
	wobei $\dfrac{a}{b}$, $\dfrac{a'}{b'}$ alle aufeinanderfolgenden Br\"uche der
	$(n-1)$-ten Zeile durchl\"auft. Wir erhalten in allen drei F\"allen
%	\begin{equation*}
%	\begin{split}
%	a'b - ab' & = 1 \quad \text{(Induktionsannahme)}\\	
%	\text{bzw.}\quad  (a+a')b - a(b+b') 	 =	 a'b - ab' & = 1  \\
%	\text{bzw.}\quad  a'(b+b') - (a+a')b'	 =	 a'b - ab' & = 1\, .	
%	\end{split}
%	\end{equation*}	  
	\begin{align*}
	a'b - ab' & = 1 \quad \text{(Induktionsannahme)}\\
	\text{bzw.}\quad  \makebox[\the\formel][c]{$(a+a')b - a(b+b')$} 	 =	 a'b - ab' & = 1  \\
	\text{bzw.}\quad  a'(b+b') - (a+a')b'	 =	 a'b - ab' & = 1\, .
	\end{align*}   
\dokendProof

\begin{Satz}\label{satz:4_4}
	Jeder Bruch\index{gek\"urzter Bruch}\label{gekuerzter Bruch} $\dfrac{a}{b}$ in der Tabelle ist gek\"urzt, d.h. $\ggT(a,b) = 1$.
	Die Br\"uche sind in jeder Zeile nach aufsteigender Gr\"o{\ss}e geordnet. Die Farey-Sequenz $\F_n$ in der $n$-ten Zeile ist \"uberdies vollst\"andig, d.h. sie enth\"alt alle gek\"urzten Br\"uche $\dfrac{a}{b}\in[0,1]$ mit~$1\leq b\leq n$.
	\dokendSatz
\end{Satz}
{\bf Beweis:}~ 
Sind $\dfrac{a}{b}$,  $\dfrac{a'}{b'}$ zwei aufeinanderfolgende Br\"uche von $\F_n$, so gilt 
\mbox{$a'b - ab' = 1$} nach Satz~\ref{satz:4_3}. Hieraus folgt zum einen $\ggT(a,b) = 1$, und zum anderen
$$\displaystyle \frac{a'}{b'}=\frac{a}{b}+\frac{1}{b b'}>\frac{a}{b},$$ so dass die Br\"uche der Farey-Sequenz $\F_n$ nach aufsteigender Gr\"o{\ss}e sortiert sind.\\

Die Vollst\"andigkeit von $\F_n$ beweisen wir unter Verwendung des Mediantensatzes~\ref{satz:3_16}\index{Mediantensatz}\label{Mediantensatz3} mit vollst\"andiger Induktion:

$\F_1$, bestehend aus den beiden Br\"uchen $\dfrac{0}{1}$, $\dfrac{1}{1}$, ist vollst\"andig (Induktionsanfang). 
Wir nehmen an, die Vollst\"andigkeit von $\F_{n-1}$ sei f\"ur ein $n\geq 2$ bereits gegeben.
Dann liegen alle gek\"urzten Br\"uche $\dfrac{a}{b}\in [0,1]$ mit $b\leq n-1$ bereits in $\F_{n-1}$, und somit auch in $\F_n$. Nun sei $\dfrac{A}{n}\in[0,1]$ ein beliebiger gek\"urzter Bruch. Wir m\"ussen $\dfrac{A}{n}\in\F_n$ zeigen. Wegen $n\geq 2$ folgt sogar $0<\dfrac{A}{n}<1$, und $\dfrac{A}{n}$ kann nicht in $\F_{n-1}$ liegen.
Damit gibt es eindeutig bestimmte und in $\F_{n-1}$ aufeinanderfolgende Br\"uche $\dfrac{a}{b}$,  $\dfrac{a'}{b'}$, so dass gilt:
$$
0\leq \frac{a}{b}<\frac{A}{n}<\frac{a'}{b'}\leq 1, \quad a'b-ab'=1.
$$
Mit Satz~\ref{satz:3_16} folgt $b+b'\leq n$, und aus der Vollst\"andigkeit von $\F_{n-1}$ erhalten wir $b+b'\geq n$, da andernfalls $b+b'\leq n-1$ w\"are und $\dfrac{a+a'}{b+b'}$ schon in $\F_{n-1}$ liegen w\"urde, im Widerspruch zur Wahl von $\dfrac{a}{b}$,  $\dfrac{a'}{b'}$.
Wir haben $b+b'= n$ gezeigt, und da $\dfrac{a+a'}{b+b'}$ nach Satz~\ref{satz:3_16} der einzige Bruch mit kleinstem Nenner ist, der echt zwischen $\dfrac{a}{b}$ und  $\dfrac{a'}{b'}$ liegt, folgt auch noch $A=a+a'$. Damit liegt $\dfrac{A}{n}=\dfrac{a+a'}{b+b'}$ in $\F_n$, und $\F_n$ ist vollst\"andig.

\dokendProof

\begin{Def}
	\label{def:4_5}
	Mit $\F_n^{ext}$ bezeichnen wir die erweiterte Farey-Sequenz\index{erweiterte Farey-Sequenz}\label{erweiterte Farey-Sequenz} der Ordnung $n$, bestehend aus allen gek\"urzten Br\"uchen $\dfrac{a}{b}$ mit $a\in\Z$, $b\in\N$ und $b\leq n$. Die Br\"uche von $\F_n^{ext}$ denken wir uns nach aufsteigender Gr\"o{\ss}e sortiert.
	\dokendDef
\end{Def}

\begin{Beis}\label{beis:4_5}
	Die erweiterte Farey-Sequenz $\F_3^{ext}$ lautet:
	$$
	\cdots < \frac{-1}{1}<\frac{-2}{3}<\frac{-1}{2}<\frac{-1}{3}<\frac{0}{1}<\frac{1}{3}<\frac{1}{2}<\frac{2}{3}
	<\frac{1}{1}<\frac{4}{3}<\frac{3}{2}<\cdots
	$$
	\dokendSatz
\end{Beis}

\begin{Satz}\label{satz:4_7}
	Es seien $a,a'\in\Z$, $b,b'\in\N$ und $\ggT(a,b)=\ggT(a',b')=1$. 
	Genau dann folgen die gek\"urzten Br\"uche $\dfrac{a}{b}<\dfrac{a'}{b'}$ in $\F_n^{ext}$ aufeinander, wenn gilt:
	\begin{equation}\label{eq:4_star1}
	a'b-ab'=1, \quad b\leq n, \quad b'\leq n \quad \text{und } \quad b+b'>n.
	\end{equation}
	\dokendSatz
\end{Satz}

{\bf Beweis:}~ 
	Wir nehmen an, dass $\dfrac{a}{b}<\dfrac{a'}{b'}$ in $\F_n^{ext}$  aufeinanderfolgen.
	Mit $q:=\left\lfloor\dfrac{a}{b}\right\rfloor$ bilden wir $\tilde{a}:=a-qb$, $\tilde{a}':=a'-qb'$.
	Dann ist $0\leq \tilde{a}<b$, und $\dfrac{\tilde{a}}{b}<\dfrac{\tilde{a}'}{b'}$ folgen bereits in $\F_n$ aufeinander. Insbesondere ist $\dfrac{\tilde{a}'}{b'}\in (0,1]$, und aus Satz~\ref{satz:4_3} folgt:
	$$
	1=\tilde{a}'b-\tilde{a}b'=(a'-qb')b-(a-qb)b'=a'b-ab'.
	$$
	Die Bedingungen $b\leq n$, $b'\leq n$ folgen aus der Definition von $\F_n^{ext}$. 
	Zudem ist $b+b'\leq n$ ausgeschlossen, da man sonst zwischen $\dfrac{a}{b}$ und $\dfrac{a'}{b'}$ den neuen Bruch $\dfrac{a+a'}{b+b'}\in \F_n^{ext}$ h\"atte.\\
	
	Nun setzen wir (\ref{eq:4_star1}) voraus. Dann folgt zun\"achst, dass $\dfrac{a}{b}<\dfrac{a'}{b'}$ in  $\F_n^{ext}$ liegen. Nach dem Mediantensatz~\ref{satz:3_16}\index{Mediantensatz}\label{Mediantensatz4} sind wegen $b+b'>n$ die beiden Br\"uche in  $\F_n^{ext}$ aufeinanderfolgend.
\dokendProof\\

Im Folgenden werden S\"atze hergeleitet, mit denen jede (erweiterte) Farey-Sequenz in einem beliebigen Abschnitt sehr effizient berechnet werden kann, ohne die vorhergehenden (erweiterten) Farey-Sequenzen kennen zu m\"ussen:

\begin{Satz}\label{satz:4_8}
	Es seien $\dfrac{a}{b}<\dfrac{a^*}{b^*}$ zwei gek\"urzte Br\"uche, die in $\F_b^{ext}$ aufeinanderfolgen. Ist dann $b\leq n$, so folgen die beiden Br\"uche
	$$
	\frac{a}{b}<\frac{a^*+a\left\lfloor \frac{n-b^*}{b}\right\rfloor}
	{b^*+b\left\lfloor \frac{n-b^*}{b}\right\rfloor}
	$$
	in  $\F_n^{ext}$ aufeinander.
	\dokendSatz
\end{Satz}

{\bf Beweis:}~ 
	Wir setzen $q^*:=\left\lfloor \dfrac{n-b^*}{b}\right\rfloor$ und verwenden Satz~\ref{satz:4_7}:\\
	
	Unter Beachtung von $a^*b-ab^*=1$ folgt auch 
	$$ (a^*+aq^*)\cdot b-a\cdot (b^*+bq^*)=a^*b-ab^*=1. $$
	Nach Voraussetzung ist $b\leq n$. Wir haben 
	$$ b^*+b\left\lfloor \dfrac{n-b^*}{b}\right\rfloor \leq b^*+b  \dfrac{n-b^*}{b} = n$$
	sowie
	$$ b+b^*+b\left\lfloor \dfrac{n-b^*}{b}\right\rfloor>b+b^*+b  \left(\dfrac{n-b^*}{b}-1\right)=n, $$
	womit der Satz schon bewiesen ist.
	\dokendProof
	
\begin{Bem}\label{bem:4_9}
Die linken Nachbarbr\"uche zu $a/b$ in $\F_n^{ext}$ werden analog zum Satz \ref{satz:4_8}
f\"ur $b \leq n$ in Aufgabe \ref{auf:4_3} berechnet.
Dazu sowie f\"ur die folgenden Betrachtungen merken wir folgendes an:
	Wenn $\dfrac{a}{b}<\dfrac{a^*}{b^*}$ in $\F_b^{ext}$ aufeinanderfolgen und $b\geq 2$ gilt (sonst w\"are $b=b^*=1$, $a^*=a+1$), so ist $1\leq b^*< b$, und nach Satz~\ref{satz:3_16} folgen die Br\"uche 
	\begin{equation}\label{eq:4_1}
		\frac{a-a^*}{b-b^*}<\frac{a}{b}<\frac{a^*}{b^*}
	\end{equation}
	in  $\F_b^{ext}$ aufeinander.\\
	
	F\"ur $b\geq 2$ setzen wir
	\begin{equation}\label{eq:4_2}
	a_*:=a-a^*, \quad b_*:=b-b^*,
	\end{equation}
	so dass der erweiterte Euklidische Algorithmus\index{erweiterter Euklidischer Algorithmus}\label{erweiterter_Euklidischer_Algorithmus3} mit den Eingabewerten $a,b$ f\"ur geraden Abbruchindex $n_*$ die Werte 
	\begin{equation}\label{eq:4_3}
	a_*=s_{n_*-1}, \quad b_*=t_{n_*-1}
	\end{equation}
	liefert, dagegen f\"ur ungerades $n_*>1$ die Werte
	\begin{equation}\label{eq:4_4}
	a^*=s_{n_*-1}, \quad b^*=t_{n_*-1},
	\end{equation} 
	siehe Satz~\ref{satz:3_6} (c), hier mit $\ggT(a,b)=1$.\\
	
	Der erweiterte Euklidische Algorithmus dient somit der Berechnung von $\dfrac{a^*}{b^*}$ aus $\dfrac{a}{b}$, bevor mit Satz~\ref{satz:4_8} der rechte Nachbarbruch von $\dfrac{a}{b}$ in $\F_n^{ext}$ berechnet werden kann.~	
	\dokendBem
\end{Bem}

\begin{Satz}\label{satz:4_10}
	Es seien $\displaystyle \frac{a}{b} < \frac{a'}{b'} < \frac{a''}{b''}$ drei aufeinanderfolgende
	Br\"uche von $\F_n^{ext}$, $n\in\N$. 
	Dann gilt mit $a''b-ab''>0$:
	\begin{equation*}
	a' = \frac{a+a''}{a''b - ab''}\,,\quad
	b' = \frac{b+b''}{a''b-ab''}\,,\quad
	\frac{a'}{b'} = \frac{a+a''}{b+b''}\,.  
	\end{equation*}
	\dokendSatz
\end{Satz}

{\bf Beweis:}~ 
	Aus Satz~\ref{satz:4_7} folgt:
	\begin{equation}\label{eq:4_5}
	a'b-ab'=1, \quad a''b'-a'b''=1
	\end{equation} 
	mit $\displaystyle \frac{a}{b}<\frac{a''}{b''}$ bzw. $a''b-ab''>0$.\\
	
	Es ist \eqref{eq:4_5} ein lineares Gleichungssystem f\"ur $a',b'$ mit der eindeutigen L\"osung
	\begin{equation}\label{eq:4_6}
	a'=\frac{a+a''}{a''b-ab''}, \quad b'=\frac{b+b''}{a''b-ab''}.
	\end{equation} 
	Durch Division folgt hieraus noch
	\begin{equation}\label{eq:4_7}
	\frac{a'}{b'}=\frac{a+a''}{b+b''}.
	\end{equation} 

\dokendProof

\begin{Satz}\label{satz:4_11}	
	Es seien $\dfrac{a}{b} < \dfrac{a'}{b'} < \dfrac{a''}{b''}$ drei aufeinanderfolgende Br\"uche
	von $\F_n^{ext}$, $n\in\N$. Dann gelten die folgenden
	Aussagen:
%	%\begin{itemize}
%	%\item[(i)]
%	\begin{equation}
%	b'' = b' \left \lfloor \frac{n+b}{b'} \right \rfloor - b, \tag{i}
%	\end{equation}
%	%\item[(ii)]
%	\begin{equation}
%	a'' = a' \left \lfloor \frac{n+b}{b'} \right \rfloor -a, \tag{ii}
%	\end{equation}
%	%\item[(iii)]
%	\begin{equation}
%	\left \lfloor \frac{n+b}{b'} \right \rfloor 
%	= \frac{b'' + b}{b'} 
%	= a''b - ab'' = \ggT (a+a'', b+b'').\tag{iii}
%	\end{equation}
$$ 	b'' = b' \left \lfloor \frac{n+b}{b'} \right \rfloor - b,\leqno \text{i)}$$
\vspace{-1.5em}
$$  a'' = a' \left \lfloor \frac{n+b}{b'} \right \rfloor -a,\leqno \text{ii)}$$
\vspace{-1.5em}
$$ 	\left \lfloor \frac{n+b}{b'} \right \rfloor
	=\frac{b'' + b}{b'}=a''b-ab''=\ggT (a+a'', b+b'').\leqno \text{iii)}$$
	\dokendSatz
\end{Satz}

{\bf Beweis:}~ 
	(ii) folgt aus (i) und Satz~\ref{satz:4_10}, angewendet auf die aufeinanderfolgenden Farey-Br\"uche $\dfrac{a}{b}$, $\dfrac{a'}{b'}$, $\dfrac{a''}{b''}$ aus $\F_n^{ext}$. Die letzten beiden Gleichungen von (iii) folgen aus Satz~\ref{satz:4_10}, w\"arend die erste Gleichung in (iii) zu (i) \"aquivalent ist.\\
	
	Wir m\"ussen nur noch die erste Gleichung von (iii) zeigen:
	
	Aus $\dfrac{a''}{b''}\in\F_n^{ext}$ folgen $b''\leq n$ sowie 
	
	\begin{equation}\label{eq:4_8}
	\frac{b''+b}{b'}\leq \frac{n+b}{b'}.
	\end{equation} 
	
	Aus Satz~\ref{satz:4_7} folgt $b''+b'>n$, wonach gilt:
	\begin{equation}\label{eq:4_9}
	\frac{b''+b}{b'}+1=\frac{b''+b'+b}{b'}>\frac{n+b}{b'}.
	\end{equation} 
	
	Wir fassen \eqref{eq:4_8} und \eqref{eq:4_9} zusammen:
	\begin{equation}\label{eq:4_10}
	\frac{b''+b}{b'}\leq \frac{n+b}{b'}<\frac{b''+b}{b'}+1.
	\end{equation} 
	
	Schlie{\ss}lich beachten wir, dass $\dfrac{b''+b}{b'}$ nach Satz~\ref{satz:4_10} eine nat\"urliche Zahl ist, so dass aus \eqref{eq:4_10} folgt:
	$$ \left\lfloor \frac{n+b}{b'}\right\rfloor=\frac{b''+b}{b'}. $$	
	\dokendProof

	Eine r\"uckw\"artslaufende Rekursion zweiter Ordnung zur
	Berechnung der Farey-Br\"uche der Ordnung $n$ findet der Leser in Aufgabe \ref{auf:4_2}.
	Dort erweist sie sich sogar als \"aquivalent zu der in Satz \ref{satz:4_11} .
Die Kombination der S\"atze~\ref{satz:4_8} und \ref{satz:4_11} gestattet nun eine sehr effiziente Berechnung von $\F_n^{ext}$ in einem vorgegebenen Abschnitt. Wir illustrieren dies in dem abschlie{\ss}enden

\begin{Beis}\label{beis:4_12}
	Wir berechnen den Abschnitt der Farey-Sequenz $\F_{24}$ im abgeschlossenen Intervall $\left[\dfrac{3}{8},\dfrac{7}{18}\right]$. Die Intervallr\"ander geh\"oren zu $\F_{24}$, und wir beginnen mit dem linken Randbruch $\dfrac{a}{b}$ f\"ur $a=3$, $b=8$, den Eingabewerten f\"ur den erweiterten Euklidischen Algorithmus:
%	\newline
	\begin{center}
		\begin{tabular}{|l||c|c|c|c|c|c|c|c|c|c|} \hline
			$j$  & ~~ $q_j$ ~~  & ~~ $r_j$ ~~ & ~~ $s_j$ ~~ & ~~ $t_j$ ~~  \\ \hline
			0  &    0    &  8   &  1   &  0   \\ \hline
			1  &    2    &  3   &  0   &  1   \\ \hline
			2  &    1    &  2   &  1   &  2   \\ \hline
			3  &    2    &  1   &  1   &  3   \\ \hline
			4  &    ---  &  0   &  3   &  8   \\ \hline
		\end{tabular}
	\end{center}
	
Der Abbruchindex $n_*=4$ ist gerade, und aus der letzten Zeile liest man ab, dass der linke Randbruch $\dfrac{3}{8}$ bereits gek\"urzt ist. 

Mit \eqref{eq:4_3} in Bemerkung~\ref{bem:4_9} erhalten wir aus der vorletzten Zeile der Tabelle, dass die drei Br\"uche
$$ 
\frac{1}{3}<\frac{3}{8}<\frac{3-1}{8-3}=\frac{2}{5} 
$$
in $\F_8$ aufeinanderfolgen. Wir setzen $n=24$, $a^*=2$, $b^*=5$ neben $a=3$, $b=8$ in Satz~\ref{satz:4_8}, und erhalten, dass die beiden Br\"uche
$$ 
\frac{3}{8}<\frac{2+3\cdot 2}{5+8\cdot 2}=\frac{8}{21} 
$$
in $\F_{24}$ benachbart sind, in \"Ubereinstimmung mit Satz~\ref{satz:4_7}.
Mit den beiden Startbr\"uchen $\dfrac{3}{8},\dfrac{8}{21}$ wenden wir noch zweimal den Satz~\ref{satz:4_11} an, und erhalten so den folgenden Abschnitt von $\F_{24}$:
$$ 
\frac{3}{8}<\frac{8}{21}<\frac{5}{13}<\frac{7}{18}.
$$
\dokendSatz
\end{Beis}

\begin{Satz}[Approximationssatz f\"ur Farey-Br\"uche\index{Approximationssatz f\"ur Farey-Br\"uche}\label{Approximationssatz fuer Farey-Brueche}]\label{satz:4_13}	
	Ist $x_0$ eine Irrationalzahl\index{Irrationalzahl}\label{Irrationalzahl3}, so wenden wir den erweiterten Euklidischen Algorithmus auf die beiden Eingabewerte $a=x_0$, $b=1$ an. Wir verwenden die Notationen von Lektion~\ref{cha:3}. Jeder nat\"urlichen Zahl $n\geq 2$ (Farey-Index) ordnen wir mit der Forderung $t_j<n\leq t_{j+1}$ genau einen Index $j\in\N$ zu. Der Zahl $n=1$ ordnen wir den Index $j=1$ zu. Damit setzen wir 
	$$
	s_{n, j+1}:=s_{j-1}+s_j\left \lfloor \frac{n-t_{j-1}}{t_j}\right\rfloor, \quad 
	t_{n, j+1}:=t_{j-1}+t_j\left \lfloor \frac{n-t_{j-1}}{t_j}\right\rfloor.
	$$
	Dann gilt f\"ur ungerades $j$: 
	$\displaystyle \frac{s_j}{t_j}<x_0<\frac{s_{n, j+1}}{t_{n, j+1}}$, und die gek\"urzten Br\"uche  $\displaystyle \frac{s_j}{t_j}<\frac{s_{n, j+1}}{t_{n, j+1}}$ sind in $\F_n^{ext}$ benachbart.\\
	
	F\"ur gerades $j$ gilt entsprechend: 
	$\displaystyle \frac{s_{n, j+1}}{t_{n, j+1}}<x_0<\frac{s_j}{t_j}$, und die gek\"urzten Br\"uche  \mbox{$\displaystyle \frac{s_{n, j+1}}{t_{n, j+1}}<\frac{s_j}{t_j}$} sind in $\F_n^{ext}$ benachbart.
	\dokendSatz
\end{Satz}

{\bf Beweis:}~
Wir setzen $\D q_{n,j}=\left \lfloor \frac{n-t_{j-1}}{t_j}\right\rfloor$, so dass gilt:
\begin{equation}\label{eq:4_11}
s_{n,j+1}=s_{j-1}+q_{n,j}\,s_j, \quad t_{n,j+1}=t_{j-1}+q_{n,j}\,t_j.
\end{equation}
Es gelten die Ungleichungen
\begin{equation}\label{eq:4_12}
0\leq q_{n,j}\leq q_j
\end{equation}
wegen 
$$
n\geq t_j \geq t_{j-1}, \quad \frac{n-t_{j-1}}{t_j}\leq \frac{t_{j+1}-t_{j-1}}{t_j}=q_j,
$$
sowie
\begin{equation}\label{eq:4_13}
1\leq t_j \leq n, \quad 1\leq t_{n, j+1} \leq n,
\end{equation}
denn $q_{n,j}=0$ ist nur f\"ur $t_j>1$, $t_{j-1}\geq 1$ m\"oglich, und es gilt 
$$
t_{n,j+1} \leq t_{j-1}+t_j \frac{n-t_{j-1}}{t_j}=n.
$$
Auch haben wir
\begin{equation}\label{eq:4_14}
t_j+t_{n,j+1}>n
\end{equation}
wegen 
$$
t_j+t_{n,j+1}>t_j+t_{j-1}+t_j\left(\frac{n-t_{j-1}}{t_j}-1\right)=n.
$$
Nun gilt nach \eqref{eq:4_11} und Satz~\ref{satz:3_6} (a):
\begin{equation}\label{eq:4_15}
\begin{split}
s_{n,j+1} t_j-s_j\,t_{n,j+1}
&=(s_{j-1}+q_{n,j}\,s_j)t_j-s_j(t_{j-1}+q_{n,j}\,t_j)\\
&=s_{j-1}t_j-s_j\,t_{j-1}\\
&=(-1)^{j-1}.
\end{split}
\end{equation}
F\"ur $q\geq 0$ definieren wir die Abbildung 
$M_{x_0,j}(q)=\dfrac{s_{j-1}+q s_j}{t_{j-1}+q t_j}$ mit der Ableitung 
$$ M_{x_0,j}'(q)=\dfrac{s_j t_{j-1}-s_{j-1}t_j}{(t_{j-1}+q t_j)^2}=\dfrac{(-1)^j}{(t_{j-1}+q t_j)^2}.$$ 
{\bf Fall A:}~ F\"ur ungerades $j$ ist $M_{x_0,j}$ monoton fallend, und wir erhalten mit \eqref{eq:4_12} sowie mit Satz~\ref{satz:3_6} (b):
$$
\frac{s_j}{t_j}<x_0<\frac{s_{j+1}}{t_{j+1}}=\frac{s_{j-1}+q_j s_j}{t_{j-1}+q_j t_j}
\leq\frac{s_{j-1}+q_{n,j} s_j}{t_{j-1}+q_{n,j} t_j}=\frac{s_{n,j+1}}{t_{n,j+1}}.
$$
Nach \eqref{eq:4_13}, \eqref{eq:4_14} und \eqref{eq:4_15} sind zudem die gek\"urzten Br\"uche $\displaystyle \frac{s_j}{t_j}<\frac{s_{n,j+1}}{t_{n,j+1}}$ in $\F_n^{ext}$ benachbart, siehe Satz~\ref{satz:4_7}.\\

{\bf Fall B:}~ F\"ur gerades $j$ ist $M_{x_0,j}$ monoton wachsend, und wir erhalten mit \eqref{eq:4_12} sowie mit Satz~\ref{satz:3_6} (b):
$$
\frac{s_{n,j+1}}{t_{n,j+1}}=\frac{s_{j-1}+q_{n,j} s_j}{t_{j-1}+q_{n,j} t_j}
\leq\frac{s_{j-1}+q_j s_j}{t_{j-1}+q_j t_j}=\frac{s_{j+1}}{t_{j+1}}<x_0<
\frac{s_j}{t_j}.
$$
Wieder sind nach \eqref{eq:4_13}, \eqref{eq:4_14}, \eqref{eq:4_15} und Satz~\ref{satz:4_7} die beiden gek\"urzten Br\"uche $\displaystyle \frac{s_{n,j+1}}{t_{n,j+1}}<\frac{s_j}{t_j}$ benachbart, jedoch in umgekehrter Reihenfolge.
\dokendProof\\\\

F\"ur die praktische Anwendung des Approximationssatzes mit vorgegebener Irrationalzahl $x_0$ und vorgegebenem Farey-Index $n\in\N$ ist es oft vorteilhaft, mit dem erweiterten Euklidischen Algorithmus die f\"unf Spalten $k$, $x_k$, $q_k$, $s_k$, $t_k$ f\"ur $k=0,1,\dots, j+1$ zu entwickeln, wobei $j\in\N$ derjenige Index ist, welcher der Farey-Ordnung $n$ zugeordnet ist. Das Schema hat dann die Startwerte
%
%\begin{equation} \label{eq:4_16}
%\left\{\begin{tabular}{l}
%$x_0,\; q_0=\lfloor x_0 \rfloor,\; s_0=1,\; t_0=0$,\\
%$s_1=q_0,\; t_1=1$,	 
%\end{tabular}\right.
%\end{equation}

\begin{equation*}% \label{eq:4_16}
\left\{\begin{tabular}{r l r}
$x_0,\; q_0=\lfloor x_0 \rfloor,\;s_0=$ & $1,\;$    & $t_0=0$, \\
								$ s_1=$ & $q_0,\;$  & $t_1=1$, 
\end{tabular}\right.
\end{equation*}

sowie f\"ur $k\geq 1$ die Iterationsvorschriften

\begin{equation*} %\label{eq:4_17}
\left\{\begin{tabular}{r l r l}
$ x_k$ & $=\dfrac{1}{x_{k-1}-q_{k-1}},\quad $  & $q_k$ & $=\lfloor x_k \rfloor$, \\
$ s_{k+1}$ & $=s_{k-1}+s_k q_k,\quad $  & $t_{k+1}$ & $=t_{k-1} +t_k q_k$. 
\end{tabular}\right.
\end{equation*}	
	
\begin{Beis}\label{beis:4_14}
		$x_0=\sqrt{2},\, n=20$. 	
		\begin{center}
			\begin{tabular}{|l||c|c|c|c|c|c|c|c|c|c|} \hline
				$k$  & ~~ $x_k$ ~~  & ~~ $q_k$ ~~ & ~~ $s_k$ ~~ & ~~ $t_k$ ~~  \\ \hline
				0  & 
				\rule{0pt}{2.3ex}   $\sqrt{2}$    &  1   &  1   &  0   \\ \hline
				1  & 
				\rule{0pt}{2.3ex}     $\sqrt{2}+1$  &  2   &  1   &  1   \\ \hline
				2  & 
				\rule{0pt}{2.3ex}     $\sqrt{2}+1$  &  2   &  3   &  2   \\ \hline
				3  & 
				\rule{0pt}{2.3ex}     $\sqrt{2}+1$  &  2   &  7   &  5   \\ \hline
				4  & 
				\rule{0pt}{2.3ex}     $\sqrt{2}+1$  &  2   &  17  &  12  \\ \hline
				5  & 
				\rule{0pt}{2.3ex}     $\sqrt{2}+1$  &  2   &  41  &  29  \\ \hline
			\end{tabular}
		\end{center}
		
		Hier ist $t_4<20\leq t_5$, also $j=4$. 
		Da $j$ gerade ist, folgt
		$\displaystyle \frac{s_{20,5}}{t_{20,5}}<\sqrt{2}<\frac{s_4}{t_4}$ mit den beiden Nachbarbr\"uchen $\displaystyle \frac{s_{20,5}}{t_{20,5}}<\frac{s_4}{t_4}$ in $\F_{20}^{ext}$,
		konkret
		\begin{equation*}
		\begin{split}
		s_{20,5}=s_3+s_4\left\lfloor\frac{20-t_3}{t_4}\right\rfloor&=7+17 \cdot 1=24,\\
		t_{20,5}=t_3+t_4\left\lfloor\frac{20-t_3}{t_4}\right\rfloor&=5+12 \cdot 1=17,
		\end{split}
		\end{equation*}
		und schlie{\ss}lich
		$\displaystyle \frac{24}{17}<\sqrt{2}<\frac{17}{12}$ f\"ur $\F_{20}^{ext}$.\\
		
		Mit dem hier entwickelten Rechenschema lassen sich allgemeiner 
                die besten rationalen Approximationen von $\sqrt{2}$ 
                in $\F_n^{ext}$ f\"ur $n \leq 29$ bestimmen.	
		So erhalten wir z.B. f\"ur $n=10$ den Index $j=3$ wegen $t_3<10\leq t_4$ mit ungeradem $j$,
		\begin{equation*}
		\begin{split}
		s_{10,4}=s_2+s_3\left\lfloor\frac{10-t_2}{t_3}\right\rfloor
		&=3+7\left\lfloor\frac{10-2}{5}\right\rfloor=10,\\
		t_{10,4}=t_2+t_3\left\lfloor\frac{10-t_2}{t_3}\right\rfloor
		&=2+5\left\lfloor\frac{10-2}{5}\right\rfloor=7
		\end{split}
		\end{equation*}
		und der besten rationalen Approximation 
		$\displaystyle \frac{7}{5}<\sqrt{2}<\frac{10}{7}$ in $\F_{10}^{ext}$.
	\dokendSatz
\end{Beis}

       \section{Aufgaben}\label{cha:4_A}

\begin{Auf}[Approximation einer Irrationalzahl\index{Approximation einer Irrationalzahl}\label{Approximation einer Irrationalzahl} mit Farey-Br\"uchen]\label{auf:4_1}
	Man bestimme die besten Approximationen an $x_0=\sqrt{7}$ von links und rechts
	mit Br\"uchen aus $\mathcal{F}_{200}^{ext}$.
\end{Auf}

{\bf L\"osung:}\\
Zun\"achst wenden wir den erweiterten Euklidischen Algorithmus auf die Eingabewerte $a=x_0=\sqrt{7}$, $b=1$ an, und bestimmen f\"ur $n=200$ einen Index $j\in\N$ mit $t_j<n\leq t_{j+1}$.

\begin{center}
	\begin{tabular}{|l||c|c|c|c|c|c|c|c|c|c|} \hline
		$k$  & ~~ $x_k$ ~~  & ~~ $q_k$ ~~ & ~~ $s_k$ ~~ & ~~ $t_k$ ~~  \\ \hline
		0  & 
		\rule{0pt}{2.8ex}   $\sqrt{7}$    &  2   &  1   &  0   \\[0.05cm] \hline
		1  & 
		\rule{0pt}{2.8ex}     $\frac{\sqrt{28}+4}{6}$  &  1   &  2   &  1   \\[0.05cm] \hline
		2  & 
		\rule{0pt}{2.8ex}     $\frac{\sqrt{28}+2}{4}$  &  1   &  3   &  1   \\[0.05cm] \hline
		3  & 
		\rule{0pt}{2.8ex}     $\frac{\sqrt{28}+2}{6}$  &  1   &  5   &  2   \\[0.05cm] \hline
		4  & 
		\rule{0pt}{2.8ex}     $\frac{\sqrt{28}+4}{2}$  &  4   &  8  &  3  \\[0.05cm] \hline
	\end{tabular}
	\quad
	\begin{tabular}{|l||c|c|c|c|c|c|c|c|c|c|} \hline
		$k$  & ~~ $x_k$ ~~  & ~~ $q_k$ ~~ & ~~ $s_k$ ~~ & ~~ $t_k$ ~~  \\ \hline
		5  & 
		\rule{0pt}{2.8ex}     $\frac{\sqrt{28}+4}{6}$  &  1   &  37  &  14  \\[0.05cm] \hline
		6  & 
		\rule{0pt}{2.8ex}     $\frac{\sqrt{28}+2}{4}$  &  1   &  45   &  17   \\[0.05cm] \hline
		7  & 
		\rule{0pt}{2.8ex}     $\frac{\sqrt{28}+2}{6}$  &  1   &  82   &  31   \\[0.05cm] \hline
		8  & 
		\rule{0pt}{2.8ex}     $\frac{\sqrt{28}+4}{2}$  &  4   &  127  &  48  \\[0.05cm] \hline
		9  & 
		\rule{0pt}{2.8ex}     $\frac{\sqrt{28}+4}{6}$  &  1   &  590  &  223  \\[0.05cm] \hline
	\end{tabular}
\end{center}
	Wir haben $t_8=48<200\leq 223=t_9$ mit geradem Index $j=8$.\\
	
	Es ist \quad \quad \quad $\D\left\lfloor\frac{200-t_7}{t_8}	\right\rfloor
	=\left\lfloor\frac{200-31}{48}	\right\rfloor=3$,  \quad
	$\D\frac{s_8}{t_8}=\frac{127}{48},$
	
	\begin{equation*} s_{200,9}=s_7+3\cdot s_8=463,
	\quad t_{200,9}=t_7+3\cdot t_8=175\,,
	\end{equation*}
	und die besten Approximationen von $\sqrt{7}$ in $\mathcal{F}_{200}^{ext}$ von links und rechts sind gegeben durch
	$$\frac{463}{175}<\sqrt{7}<\frac{127}{48}.$$

\begin{Auf}[R\"uckl\"aufige Rekursion f\"ur Farey-Br\"uche\index{r\"uckl\"aufige Rekursion f\"ur Farey-Br\"uche}\label{ruecklaeufige Rekursion fuer Farey-Brueche}]\label{auf:4_2}
	
	Es seien $\frac{a}{b} < \frac{a'}{b'} < \frac{a''}{b''} $ drei aufeinanderfolgende Br\"uche
	aus $\mathcal{F}_{n}^{ext}$, $n \geq 1$ eine beliebige nat\"urliche Zahl.\\
	
	Man zeige die folgenden Darstellungsformeln:
	\begin{equation*}
	b = b' \left \lfloor \frac{n+b''}{b'} \right \rfloor - b'',\quad
	a = a' \left \lfloor \frac{n+b''}{b'} \right \rfloor  -a'',\quad
	\left \lfloor \frac{n+b''}{b'} \right \rfloor = \left \lfloor \frac{n+b}{b'} \right \rfloor\,.
	\end{equation*}

\end{Auf}

{\bf L\"osung:}\\
Voraussetzung: $\frac{a}{b} < \frac{a'}{b'} < \frac{a''}{b''}$ sind gek\"urzte Br\"uche, die f\"ur $n\in\N$ in $\mathcal{F}_n^{ext}$ aufeinanderfolgen. 
Nach Satz~\ref{satz:4_11}~(i) gilt
\begin{equation}\label{eq:4_20}
\frac{b''+b}{b'}=\left\lfloor\frac{n+b}{b'}\right\rfloor.
\end{equation}
Nach Satz~\ref{satz:4_11}~(ii) muss nur noch
$$ \left\lfloor\frac{n+b''}{b'}\right\rfloor=\left\lfloor\frac{n+b}{b'}\right\rfloor$$
gezeigt werden: Aus $b\leq n$ folgt zun\"achst unter Beachtung der Ganzzahligkeit von 
$\frac{b''+b}{b'}$ in \eqref{eq:4_20}:
\begin{equation}\label{eq:4_21}
\frac{b''+b}{b'}\leq\left\lfloor\frac{n+b''}{b'}\right\rfloor.
\end{equation}
Nach Satz~\ref{satz:4_7} ist $b+b'>n$, und hieraus folgt
$$\frac{b''+b}{b'}+1=\frac{b''+b+b'}{b'}>\frac{n+b''}{b'}
\geq\left\lfloor\frac{n+b''}{b'}\right\rfloor,$$
d.h.
\begin{equation}\label{eq:4_22}
\left\lfloor\frac{n+b''}{b'}\right\rfloor<\frac{b''+b}{b'}+1.
\end{equation}
Aus \eqref{eq:4_20}-\eqref{eq:4_22} erhalten wir schlie{\ss}lich:
$$\frac{b''+b}{b'}=\left\lfloor\frac{n+b}{b'}\right\rfloor
=\left\lfloor\frac{n+b''}{b'}\right\rfloor,$$
was noch zu zeigen war.

\begin{Auf}\label{auf:4_3}
	Es seien $\dfrac{a_*}{b_*}<\dfrac{a}{b}$ zwei gek\"urzte Br\"uche, die in $\F_b^{ext}$ aufeinanderfolgen. 
	Man zeige: f\"ur $b\leq n$ folgen die beiden Br\"uche
	$$
	\frac{a_*+a\left\lfloor \frac{n-b_*}{b}\right\rfloor}
	{b_*+b\left\lfloor \frac{n-b_*}{b}\right\rfloor}< \frac{a}{b}
	$$
	in  $\F_n^{ext}$ aufeinander.
\end{Auf}
{\bf L\"osung:}\\
	Wir setzen $q_*:=\left\lfloor \dfrac{n-b_*}{b}\right\rfloor$ und verwenden Satz~\ref{satz:4_7}:\\
	
	Unter Beachtung von $ab_*-a_*b=1$ folgt auch 
	$$a\cdot (b_*+bq_*) -(a_*+aq_*)\cdot b=1. $$
	Nach Voraussetzung ist $b\leq n$. Wir haben 
	$$ b_*+b\left\lfloor \dfrac{n-b_*}{b}\right\rfloor \leq b_*+b  \dfrac{n-b_*}{b} = n$$
	sowie
	$$ b+b_*+b\left\lfloor \dfrac{n-b_*}{b}\right\rfloor>b+b_*+b  \left(\dfrac{n-b_*}{b}-1\right)=n, $$
	womit alles gezeigt ist.
Diese Aufgabe ist eine Erg\"anzung zum Satz \ref{satz:4_8}, siehe auch Bemerkung \ref{bem:4_9},
um die hier gemachte Voraussetzung mit der vom Satz \ref{satz:4_8} in Einklang zu bringen.

\begin{Auf}\label{auf:4_4}

Es sei $n$ eine nat\"urliche Zahl.
	\begin{itemize}
	\item[(a)] Man zeige, dass die ersten $1+\lfloor n/2 \rfloor$ Nachbarbr\"uche in  $\F_n$  rechts von $0/1$ gegeben sind durch
	$$
	 \frac{1}{n} <  \frac{1}{n-1} < \ldots < \frac{1}{\lceil n/2 \rceil}\,.
	$$
	\item[(b)] Man berechne f\"ur $n \geq 2$ den linken und rechten Nachbarbruch von $1/2$ in $\F_n$.
	\end{itemize}
\end{Auf}
{\bf L\"osung:}\\
Zum Beweis von (a) verwenden wir den Satz \ref{satz:4_7} mit den dortigen Notationen,
setzen zun\"achst $a=0$, $b=1$, $a'=1$, $b'=n$ und erhalten $a'b-ab'=1$ sowie $b+b'=n+1>n$ neben
$a/b, a'/b' \in \F_n$. Damit ist gezeigt, dass $1/n$ der rechte Nachbarbruch von $0/1$ in $\F_n$ ist.
F\"ur $n=1$ ist (a) schon bewiesen, so dass wir $n \geq 2$ voraussetzen d\"urfen.
Damit betrachten wir f\"ur nat\"urliche Zahlen $k \leq \lceil n/2 \rceil$ die beiden Br\"uche
\begin{equation}\label{auf:4_4nachbarn}
\frac{1}{k+1}<\frac{1}{k}\,.
\end{equation}
F\"ur $a=a'=1$, $b=k+1$, $b'=k$ gilt wieder $a'b-ab'=1$, und wegen $n \geq 2$ haben wir $k+1 \leq  \lceil n/2 \rceil+1 \leq n$
sowie $k \leq n$. Schliesslich ist $b+b'=2k+1\geq 2\frac{n}{2}+1>n$, so dass die beiden Br\"uche 
in \eqref{auf:4_4nachbarn} in $\F_n$ benachbart sind.\\

F\"ur die Teilaufgabe (b) setzen wir $a=1$, $b=2$, und erhalten aus Aufgabe \ref{auf:4_3} mit $a_*=0$, $b_*=1$
sowie aus Satz \ref{satz:4_8} mit $a^*=1$, $b^*=1$, dass die folgenden drei B\"uche f\"ur $n \geq 2$ in $\F_n$ benachbart sind:
$$
\frac{\left \lfloor \frac{n-1}{2} \right\rfloor}{1+2\left \lfloor \frac{n-1}{2} \right\rfloor}
< \frac12
<\frac{1+\left \lfloor \frac{n-1}{2} \right\rfloor}{1+2\left \lfloor \frac{n-1}{2} \right\rfloor}
\,.
$$

\chapter{Zahlentheoretische Funktionen}\label{cha:5}
Zahlentheoretische Funktionen	sind zun\"achst nichts anderes als reell- oder komplexwertige
Zahlenfolgen. Motiviert durch die Einschr\"ankung auf sogenannte multiplikative
zahlentheoretische Funktionen	werden grundlegende spezielle
zahlentheoretische Funktionen wie die M\"obius-Funktion $\mu$ und die Eulersche Funktion $\varphi$
eingef\"uhrt und studiert. Im Rahmen dieser Untersuchungen
wird man neben der punktweisen Multiplikation
von Zahlenfolgen noch auf die allgemeine Dirichletsche Faltung von Zahlenfolgen gef\"uhrt,
eine interessante weitere Art der Multiplikation, welche die Einf\"uhrung
von zwei f\"ur die Zahlentheorie wichtigen abelschen, multiplikativen Gruppen 
in Satz \ref{satz:5_5} erm\"oglicht. \\
		\section{Zahlentheoretische Funktionen}\label{cha:5T}
Ganz allgemein nennt man eine Abbildung $f:\N\rightarrow\C$ zahlen\-theo\-re\-tische Funk\-tion\index{zahlentheoretische Funktion}\label{zahlentheoretische Funktion}. Man schreibt sie auch als Zahlenfolge $(a_n)_{n\in\N}$ mit $a_n=f(n)$.
Von besonderem Interesse sind dabei multiplikative bzw. vollst\"andig multiplikative Funktionen:

\begin{Def}\label{def:5_1}	
	Die zahlentheoretische Funktion $f:\N\rightarrow\C$ hei{\ss}t multiplikativ\index{multiplikativ}\label{multiplikativ}, wenn $f(1)=1$ ist sowie $f(n_1\cdot n_2)=f(n_1)\cdot f(n_2)$ f\"ur alle teilerfremden nat\"urlichen Zahlen $n_1, n_2$ gilt. Wenn $f(1)=1$ ist und \"uberdies $f(n_1\cdot n_2)=f(n_1)\cdot f(n_2)$ f\"ur alle $n_1, n_2 \in \N$ gilt, dann wird $f$ sogar vollst\"andig multiplikativ\index{vollst\"andig multiplikativ}\label{vollstaendig multiplikativ} genannt.	
	\dokendDef
\end{Def} 

\begin{Bem}\label{bem:5_2}
	Nach dem Satz von der eindeutigen Primfaktorzerlegung\index{eindeutige Primfaktorzerlegung}\label{eindeutige Primfaktorzerlegung} der nat\"urlichen Zahlen $n$ ist eine multiplikative Funktion $f$ durch ihre Werte an allen Primzahlpotenzen $p^\alpha$ mit $\alpha\in\N$ eindeutig festgelegt:\\
	
	Aus der Zerlegung $n=p_1^{\alpha_1}\cdot \ldots \cdot p_m^{\alpha_m}$ mit paarweise verschiedenen Primzahlen $p_1,\dots, p_m$ und Exponenten $\alpha_1,\dots, \alpha_m\in\N$ folgt ja
	
	\begin{equation}\label{eq:5_1}
	f(p_1^{\alpha_1}\cdot \ldots \cdot p_m^{\alpha_m})=
	f(p_1^{\alpha_1})\cdot \ldots \cdot f(p_m^{\alpha_m}).
	\end{equation}
	
	Soll $f$ sogar vollst\"andig multiplikativ sein, so gen\"ugt es, f\"ur jede Primzahl $p$ und jedes $\alpha\in\N$ neben $f(1)=1$ und neben \eqref{eq:5_1} noch Folgendes zu fordern:
	
	\begin{equation}\label{eq:5_2}
	f(p^{\alpha})=f(p)^\alpha.
	\end{equation}
	
	\dokendBem
\end{Bem}

\begin{Def}\label{def:5_3}
	\hspace*{0cm}\\\vspace{-1cm}
%  ~\vspace{-0.6cm}
	\begin{enumerate} [(a)]
		\item Wir definieren $\varepsilon,\, 1,\, \Id:\N\rightarrow\N_0$ mit	
			
		$$		\varepsilon(n):=   \left\{\begin{tabular}{l}
									$1,\quad n=1$,\\
									$0, \quad n\geq 2$,	 
									\end{tabular}\right.\quad   
		1(n):= 1, \quad \Id(n):=n. 
		$$
		Dies sind vollst\"andig multiplikative Funktionen.
		\item Ist $p$ irgendeine Primzahl und $\alpha\in\N$, so definiert man gem\"a{\ss} Bemerkung~\ref{bem:5_2} durch die Festlegungen		
		\begin{align*}
			\mu(p^\alpha):=  & \left\{\begin{tabular}{rl}
			$-1,\quad $ & $\alpha=1$,\\
			$0, \quad$ & $ \alpha\geq 2$,	 
			\end{tabular}\right.		&  \text{bzw.}\quad
			\varphi(p^\alpha):=&p^{\alpha-1}\cdot (p-1)
		\end{align*}
		die multiplikative M\"obius-Funktion\index{M\"obius-Funktion}\label{Moebius-Funktion} $\mu:\N\rightarrow \{0,\pm 1\}$ bzw. die multiplikative Euler-Funktion\index{Eulersche Funktion}\label{Eulersche Funktion} $\varphi:\N \rightarrow \N $.
		Weder $\mu$ noch $\varphi$ sind vollst\"andig multiplikativ, da \eqref{eq:5_2} in Bemerkung~\ref{bem:5_2} nicht allgemein gilt.
	\end{enumerate}	
	Tabelle:

\begin{center}%\resizebox{\textwidth}{!} {
	\begin{tabularx}{\textwidth}
		{|r||R|R|R|R|R|R|R|R|R|R| R|R|R|R|R|R|R|R|R|R|}
		\hline
		$n$  &  $1$  &  $2$  &  $3$   &   $4$   &   $5$   &   $6$    &   $7$   
		&   $8$   &   $9$   &   $10$  &  $11$   &   $12$   &   $13$  &   $14$   &   $15$  &  $16$  &   $17$   &  $18$   &   $19$   &   $20$\\
		\hline
		$\varepsilon(n)$  &    1    &  0   &  0   &  0  &  0   &  0   &  0 &  0   &  0   &  0  &  0   &  0   &  0 &  0   &  0   &  0 &  0   &  0 &  0   &  0 \\ \hline
		
		$1(n)$  &    1    &  1   &  1   &  1  &  1   &  1   &  1 &  1   &  1   &  1  &  1   &  1   &  1 &  1   &  1   &  1 &  1   &  1 &  1   &  1 \\ \hline
		
		$\Id(n)$  &    1    &  2   &  3   &  4  &  5   &  6   &  7 &  8  &  9   &  10  &  11   &  12   &  13 &  14   &  15   &  16 &  17   &  18 &  19   &  20 \\ \hline
		
		$\mu(n)$  &    1    &  -1   &  -1   &  0  &  -1   & 1   &  -1 &  0   &  0   &  1  &  -1   &  0   & -1 &  1   &  1   &  0 &  -1   &  0 &  -1   &  0 \\ \hline
		
		$\varphi(n)$  &    1    &  1   &  2   &  2  &  4   &  2   &  6 &  4   &  6   &  4  &  10   &  4   &  12 &  6   &  8   &  8 &  16   &  6 &  18   &  8 \\ \hline
	\end{tabularx}
\end{center}

Wir haben 
$$\varepsilon(1)=1(1)=\Id(1)=\mu(1)=\varphi(1)=1,$$ 
$$f(n_1\cdot n_2)=f(n_1)\cdot f(n_2) \quad \text{f\"ur alle}~ n_1, n_2 \in \N$$ 
%$\forall n_1, n_2 \in \N$ 
f\"ur $f=\varepsilon,1, \Id$\!;
schlie{\ss}lich 
$$\ggT(n_1,n_2)=1\Rightarrow g(n_1\cdot n_2)=g(n_1)\cdot g(n_2)$$
f\"ur $g=\mu$ bzw. $g=\varphi$ und jeweils f\"ur alle teilerfremden $n_1, n_2 \in \N$.
	\dokendDef
\end{Def}

\begin{Def}[Dirichlet-Faltung]\label{def:5_4}
	Je zwei zahlentheoretischen Funktionen $f,g:\N \rightarrow\C$ ordnen wir ihre Dirichlet-Faltung\index{Dirichlet-Faltung}\label{Dirichlet-Faltung} $f\ast g:\N \rightarrow\C$ zu mit
	$$
	(f\ast g)(n):=\sum_{d|n}f(d)g\left(\frac{n}{d}\right),
	$$
	wobei $d$ alle nat\"urlichen Teiler von $n$ durchl\"auft. Da mit $d$ auch $\dfrac{n}{d}$ alle nat\"urlichen Teiler von $n$ durchlaufen werden, ist die Dirichlet-Faltung kommutativ:
	$$
	(f\ast g)(n)=(g\ast f)(n),
	$$
	was man auch mit folgender symmetrischer Kurzschreibweise ausdr\"uckt: 
	$$\D (f\ast g)(n)=\sum_{d_1 d_2=n}f(d_1)f(d_2).$$
	\dokendDef
\end{Def} 

\begin{Satz}\label{satz:5_5}
	\hspace*{0cm}\\\vspace{-1cm}
%	~\vspace{-0.6cm}
	\begin{enumerate}[(a)]
		\item Die Dirichlet-Faltung zahlentheoretischer Funktionen ist kommutativ und assoziativ mit der Funktion $\varepsilon$ aus Definition~\ref{def:5_3} (a) als Einselement:
		$$ \varepsilon\ast f=f \quad \text{f\"ur jedes } f:\N \rightarrow \C.
		$$
		
		\item Jedes $f:\N \rightarrow \C$ mit $f(1)\neq 0$ besitzt bzgl. der Dirichlet-Faltung eine Inverse $f_\ast^{-1}:\N \rightarrow \C$ mit $f_\ast^{-1}(1)\neq0$ und $f_\ast^{-1} \ast f=\varepsilon$. Mit der Dirichlet-Faltung ist 
		$$ \F_\ast:=\left\{f:\N \rightarrow \C: f(1)\neq 0\right\}$$ eine abelsche Gruppe, die gro{\ss}e Faltungsgruppe\index{gro{\ss}e Faltungsgruppe}\label{grosse Faltungsgruppe}.
		
		\item Die Menge $\M$ aller multiplikativen zahlentheoretischen Funktionen ist bzgl. ,,$\ast$'' eine Untergruppe  der gro{\ss}en Faltungsgruppe $\F_\ast$.
		Wir nennen $\M$ die Faltungsgruppe der multiplikativen Funktionen\index{Faltungsgruppe der multiplikativen Funktionen}\label{Faltungsgruppe der multiplikativen Funktionen}.
	\end{enumerate}	
	\dokendSatz
\end{Satz}
{\bf Beweis:}~ 
\begin{enumerate}[(a)]
	\item
	Die Kommutativit\"at von ,,$\ast$'' wurde schon gezeigt, und die Assoziativit\"at folgt f\"ur alle $n\in\N $ und je drei Funktionen $f,g,h:\N \rightarrow\C$ aus
	\begin{equation*}
	\begin{split}
		\left(\left(f\ast g\right)\ast h\right)(n)
		&=\sum_{d_3|n} (f\ast g) \left(\frac{n}{d_3}\right)h(d_3)\\
		&=\sum_{d_3|n}\,	\sum_{d_2|\frac{n}{d_3}} f\left(\frac{n}{d_2 d_3}\right)g(d_2)h(d_3)\\
		&=\sum_{\underset{d_1\cdot d_2 \cdot d_3=n}{(d_1,d_2,d_3)\in\N^3:}} f(d_1)g(d_2)h(d_3)\\
%		&\overset{\underset{\text{form}}{\text{Kurz}}}{=}
%		&\overset{\text{Kurzform}}{=}
		&=
		\sum_{d_1 d_2 d_3=n} f(d_1)g(d_2)h(d_3)=(f\ast(g\ast h))(n).
		\end{split}	
	\end{equation*}
	
	Auch ist $\D (\varepsilon \ast f)(n)=\sum_{d|n}\varepsilon(d)f\left(\frac{n}{d}\right)
	=\varepsilon(1)f\left(\frac{n}{1}\right)=f(n)$ f\"ur alle $n\in\N$ klar. 
	
	\item F\"ur (b) beachten wir zun\"achst $(f\ast g)(1)=f(1)g(1)\neq 0$ f\"ur alle $f,g\in\F_{\ast}$, so dass auch $f\ast g\in\F_{\ast}$ ist.
	Zu jedem $f\in\F_{\ast}$ konstruieren wir nun $f_{\ast}^{-1}\in\F_{\ast}$ aus den Rekursionsformeln
	$$ 
	f_{\ast}^{-1}(1)=\frac{1}{f(1)}, \quad 
	f_{\ast}^{-1}(n)=-\frac{1}{f(1)}
	\sum_{\underset{d<n}{d|n:}}f\left(\frac{n}{d}\right)f_{\ast}^{-1}(d)\,,
	$$
	wobei $n>1$ ist. Dann folgt zun\"achst f\"ur $n=1$:
	$$
	(f_{\ast}^{-1}\ast f)(1)=f_{\ast}^{-1}(1)\cdot f(1)=1=\varepsilon(1).
	$$
	F\"ur $n>1$ haben wir dagegen
	$$
	(f_{\ast}^{-1}\ast f)(n)
	=\sum_{\underset{d<n}{d|n:}}f_{\ast}^{-1}(d)f\left(\frac{n}{d}\right)
	+f_{\ast}^{-1}(n)f\left(\frac{n}{n}\right)=0=\varepsilon(n),
	$$
	und insgesamt $f_{\ast}^{-1}\ast f=\varepsilon$.
	Zusammen mit (a) folgt, dass $\F_{\ast}$ bzgl. ,,$\ast$'' eine abelsche Gruppe ist.\\
	
	\item Mit $\M$ haben wir die Menge aller multiplikativen zahlentheoretischen Funktionen bezeichnet. Es seien $f$, $g\in\M$. Dann ist $(f\ast g)(1)=f(1)g(1)=1$. Die nat\"urlichen Zahlen $m$, $n$ seien teilerfremd. Es gilt 
	$$(f\ast g)(mn)=\sum_{d|mn}f(d)g\left(\frac{mn}{d}\right). $$
	Wegen $\ggT(m,n)=1$ entspricht jedem nat\"urlichen Teiler $d$ von $m\cdot n$ umkehrbar eindeutig ein Zahlenpaar $(d', d'') \in\N^2$ mit $d'|m$, $d''|n$,
	so dass $d=d'\cdot d''$ wird. Aus der Multiplikativit\"at von $f$ und $g$ folgt somit
		\begin{equation*}
			\begin{split}
				\left(f\ast g\right)(mn)
				&=\sum_{d'|m}\,	\sum_{d''|n}f\left(d'd''\right) g\left(\frac{m}{d'}\cdot\frac{n}{d''}\right)\\
				&=\sum_{d'|m}\,	\sum_{d''|n}f\left(d'\right)g\left(\frac{m}{d'}\right)f\left(d''\right) g\left(\frac{n}{d''}\right)\\
				&=\left(\sum_{d'|m}\,	f\left(d'\right)g\left(\frac{m}{d'}\right)\right)\cdot\left(\sum_{d''|n}f\left(d''\right) g\left(\frac{n}{d''}\right)\right)\\
				&=\left(f\ast g\right)(m)\cdot \left(f\ast g\right)(n).
			\end{split}	
		\end{equation*}
		Damit ist wieder $f\ast g\in\M$. \\
		
		Schlie{\ss}lich m\"ussen wir noch die Abgeschlossenheit von $\M$ unter der Dirichletschen Inversion zeigen.
	  Wir setzen hierf\"ur $f\in\M$ voraus und m\"ussen $f_{\ast}^{-1}\in\M$ zeigen:
	Nach (b) ist $f_{\ast}^{-1}\in\F_{\ast}$ eindeutig konstruierbar. Auch haben wir im ersten Beweisteil von Satz~\ref{satz:5_5} (c) bereits
	$$ f,g\in\M\Rightarrow f\ast g\in\M$$
	gezeigt. In $\F_{\ast}$ gilt $\D f_{\ast}^{-1}(1)=\frac{1}{f(1)}=1$ wegen $f(1)=1$.\\
	
	Wir definieren $g\in\M$ an Primzahlpotenzstellen $p^\alpha$ gem\"a{\ss} $g(p^\alpha)=f_{\ast}^{-1}(p^\alpha)$ und setzen dann $g$ unter Beachtung von $g(1)=1$ multiplikativ fort. Dann gilt in~$\F_{\ast}$ f\"ur alle $\alpha \in\N$ und alle Primzahlen $p$:
	$$
	(f\ast g)(p^\alpha)=\sum_{d|p^\alpha} f(d)f_{\ast}^{-1}\left(\frac{p^\alpha}{d}\right)=\varepsilon(p^\alpha)=0 \,.
	$$
	Wegen $f\ast g\in\M$ folgt hieraus $f\ast g=\varepsilon$ mit $g=f_{\ast}^{-1} \in \M$.
\end{enumerate}	
\dokendProof

\begin{Def}\label{def:5_6}
	\hspace*{0cm}\\\vspace{-1cm}
	\begin{enumerate}[(a)]
		\item Eine nat\"urliche Zahl $n$ hei{\ss}t quadratfrei, wenn $k^2|n$ f\"ur keine nat\"urliche Zahl $k>1$ gilt.
		\item Ist $n=p_1^{\alpha_1}\cdot ... \cdot p_j^{\alpha_j}$ mit Primzahlen $p_1<...<p_j$ und $\alpha_1,...,\alpha_j\in\N$ die Primfaktorzerlegung von $n\in\N$ f\"ur $n>1$ mit $j$ verschiedenen Primzahlen, so setzen wir $\omega (n):=j$. Zudem setzen wir $\omega(1):=0$.
	\end{enumerate}
	\dokendDef
\end{Def} 

\begin{Satz}\label{satz:5_7}
	\hspace*{0cm}\\\vspace{-1cm}
%	~\vspace{-0.6cm}
	\begin{enumerate}[(a)]
		\item Die multiplikative M\"obius-Funktion\index{M\"obius-Funktion}\label{Moebius-Funktion2} $\mu$ berechnet sich nach der Formel
		\begin{align*}
			\mu(n)=  & \left\{\begin{tabular}{cl}
			$(-1)^{\omega(n)},\quad $ & $\text{falls } n\in\N \text{ quadratfrei ist}$,\\
			$0, \quad$ & \text{sonst}.	 
			\end{tabular}\right.	
		\end{align*}
		Es gilt f\"ur alle $n \in \N$:
		\begin{align*}
			\sum_{d|n}\mu(d)=  
			& \left\{\begin{tabular}{rl}
			$1,\quad $ & $n=1$,\\
			$0, \quad$ & $ n> 1$,	 
			\end{tabular}\right.		&  \text{d.h.}\quad
			\mu\ast1=&\varepsilon \quad \text{bzw.}\quad \mu=1_{\ast}^{-1}.
		\end{align*}
		
		\item Die multiplikative Eulersche Funktion\index{Eulersche Funktion}\label{Eulersche Funktion2} $\varphi$ erf\"ullt die Beziehungen 
		$$
		\sum_{d|n} \varphi(d)=n\quad \text{f\"ur~} n\in\N,\quad \text{d.h.}\quad \varphi\ast1=\Id,\quad
		\text{sowie}\quad \varphi=\mu\ast \Id.
		$$
		Es ist $\varphi(n)$ die Anzahl der zu $n$ teilerfremden Zahlen $k\in\N$ mit $1\leq k\leq n$.
	\end{enumerate}	
	\dokendSatz
\end{Satz}

{\it Bemerkung zu Satz~\ref{satz:5_7}:}~ Mit der Berechnungsformel f\"ur $\mu (n)$ in (a) kann man die Beziehung $\varphi=\mu\ast \Id$ in (b) f\"ur alle $n\in\N$ auch in der folgenden Form schreiben:
$$\varphi (n)=n\underset{\underset{p\text{ prim}}{p|n:}}{\Pi}\left(1-\frac{1}{p}\right)\,.$$ 

{\bf Beweis} von Satz~\ref{satz:5_7}:~ 
\begin{enumerate}[(a)]
	\item Die Berechnungsformel f\"ur $\mu(n)$ mit Hilfe von $\omega(n)$ ergibt sich direkt aus Definition~\ref{def:5_3} (b).
	Nach Satz~\ref{satz:5_5}~(c)   ist mit $\mu,1\in\M$ auch $\mu\ast 1\in\M$,
	so dass wir die Beziehung $\mu\ast 1=\varepsilon$ nur an Primzahlpotenzstellen $p^{\alpha}$ mit $\alpha\in\N$ zeigen m\"ussen:
	$$
	(\mu\ast 1)(p^\alpha)=\sum_{d|p^\alpha}\mu(d)
	=\underset{\underset{p\text{ quadratfrei}}{d|p^\alpha:}}{\sum}\mu(d)=1+\mu(p)=0=\varepsilon(p^\alpha).
	$$
	$\mu= 1_\ast^{-1}$ folgt damit ebenfalls aus Satz~\ref{satz:5_5}~(c).\\
	
	\item Wir verwenden Satz~\ref{satz:5_5}~(c): Es ist $\varphi\ast 1\in\M$ sowie 
	$$
	(\varphi\ast 1)(p^\alpha)=\sum_{d|p^\alpha}\varphi(d)
	=\sum_{\beta=0}^\alpha\varphi(p^\beta)
	=1+\sum_{\beta=1}^\alpha\left(p^\beta-p^{\beta-1}\right)=p^\alpha=\Id(p^\alpha)
	$$
	an jeder Primzahlpotenzstelle $p^\alpha$ mit $\alpha\in\N$, siehe auch Definition~\ref{def:5_3}~(b). Daher gelten $	\varphi\ast 1=\Id$ bzw. $\varphi=(\varphi\ast 1)\ast 1_{\ast}^{-1}=\Id\ast\mu$ allgemein.
	Zur Interpretation von $\varphi(n)$ definieren wir f\"ur jeden nat\"urlichen Teiler $d$ von $n$ die Mengen 
	$$A_{d,n}=\left\{k\in\N: 1\leq k\leq n \quad\text{und}\quad \ggT(k,n)=d\right\}.$$
	Deren Elementeanzahl ist $|A_{d,n}|=\tilde{\varphi}\left(\frac{n}{d}\right)$, wenn $\tilde{\varphi}(j)$ f\"ur $j\in\N$ die Anzahl der nat\"urlichen Zahlen $k\leq j$ mit $\ggT(k,j)=1$ bezeichnet.\\
	
	Die Mengen $A_{d,n}$ sind f\"ur festes $n\in\N$ elementfremd mit 
	$$\bigcup_{d|n} A_{d,n}=\{k\in\N:1\leq k\leq n\} \quad \text{und}
	\quad \sum_{d|n} |A_{d,n}|=\sum_{d|n} \tilde{\varphi}\left(\frac{n}{d}\right)=n,$$
	also gilt $1\ast\tilde{\varphi}=\Id$.
	Aus Satz~\ref{satz:5_5}~(a) und Satz~\ref{satz:5_7}~(a) folgt endlich
	$$ \varphi=\mu\ast\Id=\mu\ast(1\ast\tilde{\varphi})
	=(\mu\ast1)\ast\tilde{\varphi}=\varepsilon\ast\tilde{\varphi}=\tilde{\varphi}.$$
\end{enumerate}	
\dokendProof	

\begin{Satz}[M\"obiussche Umkehrformel\index{M\"obiussche Umkehrformel}\label{Moebiussche Umkehrformel}]\label{satz:5_8}
	Zu jedem $g:\N\rightarrow \C$ gibt es genau ein $f:\N\rightarrow \C$ mit der Eigenschaft $$g(n)=\sum_{d|n}f(d)\quad \text{f\"ur~alle~} n\in\N,\quad \text{d.h.}\quad g=f\ast 1.$$
	F\"ur dieses gilt 
	$$f(n)=\sum_{d|n}\mu(d)g\left(\frac{n}{d}\right)\quad \text{f\"ur~alle~} n\in\N,\quad \text{d.h.}\quad f=\mu\ast g.$$
	\dokendSatz
\end{Satz}

{\bf Beweis:}~Zu $g:\N\rightarrow \C$ definieren wir $f:\N\rightarrow \C$ mit $f:=\mu\ast g$. Dann gilt nach Satz~\ref{satz:5_5}~(a):
$$ f=g\ast\mu, \quad f\ast1=(g\ast\mu)\ast1=g\ast(\mu\ast1),$$
und weiter mit Satz~\ref{satz:5_7}~(a): $$f\ast1=g\ast\varepsilon=g.$$

Zur Eindeutigkeit von $f$ nehmen wir noch $\tilde{f}\ast1=f\ast1$ mit einem $\tilde{f}:\N\rightarrow \C$ an. Wie zuvor folgt
$$ 
 (\tilde{f}\ast1)\ast\mu=(f\ast1)\ast\mu, \quad \text{also wegen}\quad 1\ast\mu=\varepsilon:
$$
$$
 \tilde{f}\ast(1\ast\mu)=f\ast(1\ast\mu), \quad \text{und}\quad \tilde{f}=\tilde{f}\ast\varepsilon=f\ast\varepsilon=f.
 $$	
% 
% \begin{align*}
% &(\tilde{f}\ast1)\ast\mu=(f\ast1)\ast\mu, &\text{also wegen}\quad &1\ast\mu=\varepsilon:\\
% &\tilde{f}\ast(1\ast\mu)=f\ast(1\ast\mu), &\text{und}\quad &\tilde{f}=\tilde{f}\ast\varepsilon=f\ast\varepsilon=f.
% \end{align*}
\dokendProof	

\begin{Bem}\label{bem:5_9}
	In Satz~\ref{satz:5_8} m\"ussen weder $f$ noch $g$ multiplikativ sein, nicht einmal $f\in\F_{\ast}$ oder $g\in\F_{\ast}$ mu{\ss} gelten. 
	Dagegen gilt dort nach Satz~\ref{satz:5_5}~(c) die \"Aquivalenz
	$$ f\in\M \Leftrightarrow g\in\M.$$
	\dokendBem
\end{Bem}

\begin{Satz}\label{satz:5_10}
	Es sei $f:\N\rightarrow\C$ vollst\"andig multiplikativ\index{vollst\"andig multiplikativ}\label{vollstaendig multiplikativ2}, siehe Definition~\ref{def:5_1}. 
	Dann gilt \mbox{$f_{\ast}^{-1}=f\cdot \mu$}, d.h f\"ur alle $n\in\N$ ist
	$$f_{\ast}^{-1}(n)=f(n)\cdot\mu(n) \quad \text{f\"ur~alle~} n\in\N.$$
	\dokendSatz
\end{Satz}

{\bf Beweis:}~Es ist $f$ auch multiplikativ, also $f\in\M$ und $f\cdot\mu\in\M$, so dass wir gem\"a{\ss} Satz~\ref{satz:5_5}~(c) die Beziehung $(f\cdot\mu)\ast f=\varepsilon$ nur an Primzahlpotenzstellen $p^\alpha$ mit $\alpha\in\N$ zeigen m\"ussen. Dort gilt in der Tat
\begin{equation*}
\begin{split}
((f\cdot \mu)\ast f)(p^\alpha)
=& \sum_{d|p^\alpha}f(d)\mu(d)f\left(\frac{p^\alpha}{d}\right)\\
=& f(1)\mu(1)f(p^\alpha)+f(p)\mu(p)f(p^{\alpha-1})\\
=& 1\cdot f(p^\alpha)-f(p)f(p^{\alpha-1})=f(p^\alpha)-f(p^\alpha)=0=\varepsilon(p^\alpha).
\end{split}
\end{equation*}
\dokendProof

\begin{Bem}\label{bem:5_11}
	Die Funktion $f=1$ ist vollst\"andig multiplikativ, dagegen $f_{\ast}^{-1}=\mu$ nur noch multiplikativ. Auch sieht man leicht, dass die Dirichlet-Faltung $f\ast g$ zweier vollst\"andig multiplikativer Funktionen $f,g$ im Allgemeinen nur noch multiplikativ ist:\\
	
	Als Beispiel betrachten wir $f=g=1$ mit der Divisorfunktion $\tau:=1\ast1$, wobei $\tau(n) $ die Anzahl der nat\"urlichen Teiler von $n$ liefert. Es ist zwar $\tau\in\M$, aber $\tau$ ist nicht vollst\"andig multiplikativ, denn an Primzahlpotenzstellen $p^\alpha$ mit $\alpha\in\N$ gilt: $\tau(p^\alpha)=\alpha+1$. Zum Beispiel ist $\tau(p^2)=3$, dagegen $\tau(p)\cdot \tau(p)=2\cdot 2=4$.
	\dokendBem
\end{Bem}

		\section{Aufgaben}\label{cha:5_A}

\begin{Auf}[Inversion zahlentheoretischer Matrizen]\label{auf:5_1}
Es sei $\lambda : \N \to \C$ vollst\"andig multiplikativ\index{vollst\"andig multiplikativ}\label{vollstaendig multiplikativ3},
also $\lambda(1)=1$ und $$\lambda(mk)=\lambda(m)\lambda(k)\quad \forall m, k \in \N\,. $$ 

Definiere $A_{\lambda,n} = (a_{jk})_{j,k=1,\ldots,n} \in \C^{n \times n}$,
$B_{\lambda,n} = (b_{km})_{k,m=1,\ldots,n} \in \C^{n \times n}$ mit den Matrix-Eintr\"agen
\begin{align*}
a_{jk} &= \begin{cases}
\lambda\left( \frac{k}{j} \right) \, ,&\ \text{f\"ur}\ j | k \, ,\\
0 \, ,&\ \text{f\"ur}\ j \not{|} ~ k\,,\\
\end{cases}\\
b_{km} &= \begin{cases}
\mu\left( \frac{m}{k} \right) \lambda\left( \frac{m}{k} \right)\, ,&\ \text{f\"ur}\ k | m \, ,\\
0 \, ,&\ \text{f\"ur}\ k \not{|} ~ m\,.\\
\end{cases}\\
\end{align*}
Man zeige: Die Matrix $B_{\lambda,n}$ ist invers zu $A_{\lambda,n}$, d.h.
\begin{equation*}
\sum \limits_{k=1}^{n} a_{jk}b_{km}=\delta_{jm}=
\begin{cases}
1 \, ,&\ \text{f\"ur}\ j = m \, ,\\
0 \, ,&\ \text{f\"ur}\ j \neq m\,.\\
\end{cases}
\end{equation*}
\end{Auf}
{\bf L\"osung:}\\
Die Eintr\"age der Produktmatrix $A_{\lambda,n}$, $B_{\lambda,n}$ lauten 
\begin{equation*}
c_{jm}=\sum \limits_{k=1}^{n} a_{jk}b_{km}
=\sum \limits_{\underset{j|k,\; k|m}{k=1:}}^{n}\lambda\left(\frac{k}{j}\right)
\mu\left(\frac{m}{k}\right)
\lambda\left(\frac{m}{k}\right).
\end{equation*}
F\"ur $j\nmid m$ ist $c_{jm}=0$ (leere Summe), w\"ahrend wir f\"ur $j|m$ erhalten:
\begin{equation*}
c_{jm}=\lambda\left(\frac{m}{j}\right)\sum \limits_{\underset{j|k,\; k|m}{k=1:}}^{n}\mu\left(\frac{m}{k}\right)=\lambda\left(\frac{m}{j}\right)\sum \limits_{l|\frac{m}{j}}\mu\left(\frac{m}{l\cdot j}\right)
\end{equation*}
mit $$c_{jm}=
\begin{cases}
\lambda\left(\frac{m}{j}\right)=1 \, ,&\ \text{f\"ur}\ m = j \; (\text{siehe Satz~\ref{satz:5_7}~(a)}),\\
0 \, ,&\ \text{sonst}.\\
\end{cases}$$

Insgesamt ist $c_{jm}=\delta_{jm}=
\begin{cases}
1 \, ,& j = m \, ,\\
0 \, ,& j \neq m\,.\\
\end{cases}$

\begin{Auf}[Eigenschaften des ggT und seine Verallgemeinerung auf mehrere] {\bf Argumente}\label{auf:5_2}
	
	\begin{enumerate}[(a)]
		\item F\"ur jedes feste $n \in \N$ zeige man die Multiplikativit\"at der zahlentheoretischen Funktion
		$\mbox{ggT}(n,\cdot) : \N \to \N$.\\
		
		\item  Mit $\mbox{ggT}(a_1,\ldots,a_n)$ bezeichnen wir den gr\"o{\ss}ten gemeinsamen Teiler\index{gr\"o{\ss}ter gemeinsamer Teiler} \label{groesster_gemeinsamer _Teiler2}
		von $n \in \N$ ganzen Zahlen $a_1,\ldots,a_n$, die nicht alle zugleich verschwinden.
		Man zeige, dass es ganze Zahlen $\lambda_1,\ldots,\lambda_n$ gibt mit
		$$\mbox{ggT}(a_1,\ldots,a_n)=\sum \limits_{k=1}^{n}\lambda_k a_k$$
		und dass $d \, | \,\mbox{ggT}(a_1,\ldots,a_n)$  f\"ur jeden gemeinsamen Teiler $d$ von $a_1,\ldots,a_n$ gilt.
	\end{enumerate}
\end{Auf}

{\bf L\"osung:}
\begin{enumerate}[(a)]
	\item
Betrachte $\ggT(n,\cdot):\N\rightarrow\N$ f\"ur festes $n\in\N$. Die Zahlen $j,k\in\N$ seien teilerfremd. 
Dann gelten die Darstellungen
\begin{equation*}
\begin{split}
&j=\prod_{\rho=1}^r p_\rho^{\alpha_\rho},\;\qquad  k=\prod_{\sigma=1}^{s}p_\sigma'^{\:\beta_\sigma}\quad \text{und}\\
&n=\prod_{\rho=1}^r p_\rho^{\gamma_\rho}
\cdot\prod_{\sigma=1}^{s}p_\sigma'^{\:\delta_\sigma}
\cdot\prod_{\tau=1}^{t}p_\tau''^{\:\varepsilon_\tau},\quad 
\end{split}
\end{equation*}
mit $r,s,t\in\N_0$, mit paarweise verschiedenen Primzahlen $p_1,...,p_r$; $p'_1,...,p'_s$; $p''_1,...,p''_t$ und Exponenten $\alpha_\rho,\beta_\sigma,\varepsilon_\tau\in\N$, $\gamma_\rho,\delta_\sigma\in\N_0$. 
F\"ur $r=0$ bzw. $s=0$ bzw. $t=0$ erhalten die leeren Produkte jeweils den Wert 1. Zun\"achst beachten wir $\ggT(n,1)=1$.
Nun erhalten wir allgemein f\"ur $\ggT(j,k)=1$:
\begin{equation*}
\begin{split}
\ggT(n,j)&=\prod_{\rho=1}^r p_\rho^{\text{min}(\alpha_\rho,\gamma_\rho)},\\
\ggT(n,k)&=\prod_{\sigma=1}^{s}p_\sigma'^{\:\text{min}(\beta_\sigma,\delta_\sigma)},\\
\ggT(n,jk)&=\prod_{\rho=1}^r p_\rho^{\text{min}(\alpha_\rho,\gamma_\rho)}
\cdot\prod_{\sigma=1}^{s}p_\sigma'^{\:\text{min}(\beta_\sigma,\delta_\sigma)}\\
&=\ggT(n,j)\cdot\ggT(n,k), 
\end{split}
\end{equation*}
d.h. $\ggT(n,\cdot)$ ist multiplikativ.\\

\item
Wir zeigen durch Induktion bzgl. $n\in\N$:\\
Wenn $a_1,...,a_n\in\Z$ nicht alle verschwinden, dann gibt es ganze Zahlen $\lambda_1,...,\lambda_n$ mit 
\begin{equation}\label{eq:a5_2}
\ggT(a_1,...,a_n)=\sum_{k=1}^n\lambda_k a_k.
\end{equation}
Ist dann $d$ ein gemeinsamer Teiler von $a_1,...,a_n$, also $a_k=d\cdot a_k'$ mit $a'_k\in\Z$ f\"ur $k=1,...,n$, so folgt aus
$$\ggT(a_1,...,a_n)=d\cdot\sum_{k=1}^n\lambda_k a'_k$$
auch $d|\ggT(a_1,...,a_n)$.\\

{\bf Induktionsanfang:}~Es gilt \eqref{eq:a5_2} f\"ur $n=1$ wegen $\ggT(a_1)=|a_1|=\lambda_1a_1$ mit dem Vorzeichen
$$\lambda_1=\sgn(a_1)=\begin{cases}
1 &\text{f\"ur}\quad a_1>0\\
-1 &\text{f\"ur}\quad a_1<0.
\end{cases}$$
{\bf Induktionsannahme:}~Wir nehmen an, f\"ur ein $n\geq 2$ m\"oge es zu jedem\linebreak \mbox{$j\in\{1,...,n-1\}$} ganze Zahlen $\lambda_1,...,\lambda_j$ geben mit
$$ \ggT(a_1,...,a_j)=\sum_{k=1}^j\lambda_k a_k.$$
Hierbei k\"onnen wir voraussetzen, dass alle Zahlen $a_1,...,a_n$ von Null verschieden sind, da man Argumente $a_k=0$ in $\ggT(a_1,..,a_j)$ einfach streichen kann, um danach auf die reduzierte Liste die Induktionsannahme anwenden zu k\"onnen.\\

Nach der Induktionsannahme gibt es zu $d_\ast:=\ggT(a_1,...,a_{n-1})$ ganze Zahlen \mbox{$\lambda_1,...,\lambda_{n-1}$} mit 
\begin{equation*}
d_\ast=\sum_{k=1}^{n-1}\lambda_k a_k.
\end{equation*}
Weiterhin gibt es nach Satz \ref{satz:2_8} Zahlen $\alpha, \beta\in\Z$ mit
\begin{equation}\label{eq:4_25}
\ggT(d_\ast,a_n)=\alpha d_\ast+\beta a_n=\sum_{k=1}^{n-1}\alpha\lambda_k a_k+\beta a_n.
\end{equation}
Nach \eqref{eq:4_25} ist $d^\ast:=\ggT(a_1,...,a_n)$ ein Teiler von $\ggT(d_\ast,a_n)$ und $\ggT(d_\ast,a_n)$ gem\"a{\ss} Konstruktion ein gemeinsamer Teiler von $a_1,...,a_n$.
Aus der Maximalit\"at von $d^\ast$ folgt $d^\ast=\ggT(d_\ast,a_n)$, und wiederum nach~\eqref{eq:4_25} ist $d^\ast$ eine Linearkombination von $a_1,...,a_n$ mit ganzzahligen Koeffizienten, womit der Induktionsschritt gezeigt ist.
\end{enumerate}

\begin{Auf}[Variante der Umkehrformel von M\"obius\index{M\"obiussche Umkehrformel}\label{Moebiussche Umkehrformel2}] \label{auf:5_3}

Wir betrachten Funktionen $F, G : [1,\infty) \to \C$. Man zeige die \"Aquivalenz der beiden Bedingungen
$$ G(x)=\sum \limits_{n \leq x} F\left( \frac{x}{n}\right) \quad \mbox{f\"ur~alle~} x \geq 1 $$
sowie
$$ F(x)=\sum \limits_{n \leq x} \mu(n) G\left( \frac{x}{n}\right) \quad \mbox{f\"ur~alle~} x \geq 1 \,.$$
\end{Auf}
{\bf L\"osung:}\\
Aus der ersten Bedingung folgt
$$ \sum \limits_{n \leq x} \mu(n) G\left( \frac{x}{n}\right)=
\sum \limits_{n \leq x} \mu(n)\sum \limits_{k \leq x/n}  F\left( \frac{x}{nk}\right)  =
\sum \limits_{m \leq x}\,\sum \limits_{n \,|\,m}  \mu(n) F\left( \frac{x}{m}\right) = F(x)
$$
f\"ur alle $x \geq 1$ unter Beachtung von $\begin{displaystyle} \sum \limits_{n \,|\,m}  \mu(n) 
= \left \lfloor \frac1m \right\rfloor\,.\end{displaystyle}$
Umgekehrt erhalten wir aus der zweiten Bedingung die erste, denn
$$ \sum \limits_{n \leq x} F\left( \frac{x}{n}\right)=
\sum \limits_{n \leq x} \,\sum \limits_{k \leq x/n} \mu(k) G\left( \frac{x}{nk}\right)  = 
\sum \limits_{m \leq x}\,\sum \limits_{k \,|\,m}  \mu(k) G\left( \frac{x}{m}\right) =G(x)\,.
$$

\begin{Auf}[Addition und Multiplikation zahlentheoretischer Funktionen] \label{auf:5_4}
Mit $R$ bezeichnen wir hier die Menge aller zahlentheoretischer Funktionen. Wir versehen $R$ mit der punktweisen
Addition  $+$zweier zahlentheoretischer Funktionen sowie mit der Dirichlet-Faltung $*$ als Multiplikation.
Man zeige, dass $(R,+,*)$ ein kommutativer Ring\index{kommutativer Ring} \label{kommutativer_Ring2} mit Einselement ist.
\end{Auf}
{\bf L\"osung:}\\
Wir \"uberpr\"ufen die Ringaxiome (R1)-(R4) aus Definition \ref{def:1_9}.
Bezeichnen wir die Nullfunktion mit ${\bf 0}$, also ${\bf 0}(n)=0$ f\"ur alle $n \in \N$,
so folgt aus den Gruppeneigenschaften von $(\R,+,0)$ sofort, dass auch $(R,+,{\bf 0})$ abelsche Gruppe ist.
Somit gilt (R1). Mit der Dirichletschen Faltung $*$ als Multiplikation gilt (R2) und (R4) nach Satz \ref{satz:5_5}(a).
Wegen (R4) m\"ussen wir f\"ur den Nachweis von (R3) nur das erste Distributivgesetz \"uberpr\"ufen.
Es seien hierzu drei zahlentheoretische Funktionen $f,g,h : \N \to \C$ gegeben. Dann gilt f\"ur jedes $n \in \N$:
\begin{equation*}
\begin{split}
(f*(g+h))(n)&=\sum \limits_{d\,|\,n} f(d)\left( g\left(\frac{n}{d}\right) + h\left(\frac{n}{d}\right)\right)\\
=& \sum \limits_{d\,|\,n} f(d) g\left(\frac{n}{d}\right)+\sum \limits_{d\,|\,n} f(d) h\left(\frac{n}{d}\right)\\
=&(f*g)(n)+(f*h)(n)\,.\\
\end{split}
\end{equation*}
Wir erhalten einen kommutativen Ring, der nach Satz \ref{satz:5_5}(a) die zahlentheoretische Funktion $\varepsilon$
aus Definition \ref{def:5_3}(a) als Einselement besitzt.

	\chapter{Kongruenzen, Restklassenringe und prime Restklassengruppen}\label{cha:6}
	Zwei ganze Zahlen nennt man kongruent bzgl. eines sogenannten Moduls $n \in \N$, wenn deren Differenz durch $n$ teilbar ist. Man stellt dann nicht nur fest, dass die Kongruenz eine 
	\"Aquivalenzrelation ist, sondern auch, dass man die zugeh\"origen 
	\"Aquivalenzklassen, von denen es nur endlich viele gibt, wie ihre Representanten
	addieren und multiplizieren kann. So wird man auf die Restklassenringe sowie die prime Restklassengruppe modulo $n$ gef\"uhrt. Dies sind endliche algebraische Strukturen,
die in den darauffolgenden Abschnitten wichtige neue Anwendungen erm\"oglichen.
		\section{Kongruenzrechnung}\label{cha:6T}
Grundlage dieses Abschnittes ist die
\begin{Def}\label{def:6_1}
	Es sei $n$ eine nat\"urliche Zahl, hier auch Modul\index{Modul}\label{Modul} genannt. Zwei ganze Zahlen $a$, $a'$ mit $n|a-a'$ werden kongruent bzgl. des Moduls $n$ genannt. Wir schreiben
\begin{equation*}
a\equiv a'\mod\, n\quad \text{bzw.}\quad a\equiv a'~\ (n)
\end{equation*}
oder einfach nur $a\equiv a'$, wenn der entsprechende Modul eindeutig aus dem Kontext hervorgeht.\hfill\dokendDef
\end{Def}
\begin{Satz}\label{satz:6_2}
	Die Kongruenz $\equiv$ $\Mod n\in\N$ ist eine \"Aquivalenzrelation. Es gibt genau $n$ verschiedene \"Aquivalenzklassen $[a]_n:=\{a'\in\Z:a\equiv a'~\ (n)\}=a+n\Z$, gegeben durch die Repr\"asentanten $a\in\Z$ mit $0\leq a\leq n-1$.\hfill\dokendSatz
\end{Satz}
{\bf Beweis:}~
Wegen $n|a-a$ ist die Kongruenz\index{Kongruenz}\label{Kongruenz} $\equiv$ $\Mod n$ reflexiv, wegen $n|a-a'$ $\Rightarrow$ $n|a'-a$ symmetrisch, und die Transitivit\"at folgt aus der Implikation
\begin{equation*}
n|a-a'\;\wedge\;n|a'-a''\;\Rightarrow\;n|a-a''\;\text{wegen}\;a-a''=(a-a')+(a'-a'').
\end{equation*} 
F\"uhrt man f\"ur jedes $a\in\Z$ die Division durch $n$ mit Rest $r$ aus gem\"a\ss
\begin{equation*}
a=r+n\left\lfloor\dfrac{a}{n}\right\rfloor,
\end{equation*}
so erh\"alt man nach Satz~\ref{satz:2_4} jeweils ein $r\equiv a~\ (n)$ mit genau $n$ paarweise $\Mod n$ inkongruenten Resten $0\leq r\leq n-1$.
\dokendProof
\begin{DefSatz}\label{defsatz:6_3}
	F\"ur je zwei Restklassen $[a]_n$, $[b]_n$ wird gem\"a\ss 
	\begin{equation*}
	[a]_n+[b]_n:=[a+b]_n,\quad [a]_n\cdot[b]_n:=[a\cdot b]_n
	\end{equation*}
eine von der Auswahl der Repr\"asentanten $a$, $b$ unabh\"angige Addition bzw. Multiplikation eingef\"uhrt. Damit ist
\begin{equation*}
\Z_n:=\{[a]_n:a\in\Z\}\;\text{f\"ur jedes feste}\;n\in\N
\end{equation*}
ein kommutativer Ring\index{kommutativer Ring} \label{kommutativer_Ring3} mit $n$ Restklassen, der sogenannte Restklassenring $\Mod n$. F\"ur $n>1$ besitzt dieser das Einselement $[1]_n$.\\

Schlie{\ss}lich ist $\Z_n$ genau dann ein Integrit\"atsbereich\index{Integrit\"atsbereich} \label{Integritaetsbereich3}, wenn $n=p$ eine Primzahl\index{Primzahl}\label{Primzahl3} ist. In diesem Fall ist $\Z_p$ sogar ein K\"orper\index{K\"orper}\label{Koerper2} (mit $p$ Elementen).\hfill\dokendSatz
\end{DefSatz}
{\bf Beweis:}~
Hat man $a'\equiv a~\ (n)$, $b'\equiv b~\ (n)$, so gibt es ganze Zahlen $k$, $l$ mit $a=a'+kn$, $b=b'+ln$. Hieraus folgt die Wohldefiniertheit der Addition und Multiplikation von Restklassen, denn
\begin{equation*}
(a+b)-(a'+b')=n\cdot(k+l)\quad\text{und}\quad a\cdot b-a'\cdot b'=n\cdot(a'l+b'k+nkl).
\end{equation*}
Nun erinnern wir uns an Definition~\ref{def:2_1}, die axiomatischen Beschreibung eines Ringes:\\

Es gilt das assoziative Gesetz
\begin{multline*}
([a]_n+[b]_n)+[c]_n=[a+b]_n+[c]_n=\\
=[a+b+c]_n=[a]_n+[b+c]_n=[a]_n+([b]_n+[c]_n),
\end{multline*}
und das kommutative Gesetz folgt noch einfacher:
\begin{equation*}
[a]_n+[b]_n=[a+b]_n=[b+a]_n=[b]_n+[a]_n.
\end{equation*}
Es ist $[0]_n$ das Nullelement und $[-a]_n$ das zu $[a]_n$ entgegengesetzte Element. Damit ist $(\Z_n,+,[0]_n)$ abelsche Gruppe, die additive Restklassengruppe $\Mod n$, und es gilt (R1). Das assoziative Gesetz (R2) der Restklassenmultiplikation beweist man v\"ollig analog wie bei der Addition, ebenso die Kommutativit\"at der Restklassenmultiplikation, so dass sich (R3) schon folgenderma{\ss}en ergibt:
\begin{equation*}
[a]_n\cdot([b]_n+[c]_n)=[a]_n\cdot[b+c]_n=[a(b+c)]_n=[a\cdot b+a\cdot c]_n=[a]_n\cdot[b]_n+[a]_n\cdot[c]_n.
\end{equation*} 
Damit ist $\Z_n$ ein kommutativer Ring, f\"ur $n>1$ mit dem Einselement $[1]_n$. F\"ur $n=1$ besteht $\Z_n$ dagegen nur aus dem Nullelement $[0]_n$, was wir in Definition~\ref{def:2_3} f\"ur einen Integrit\"atsbereich ausgeschlossen haben.\\

Ist $n=a\cdot b>1$ mit nat\"urlichen Zahlen $1<a<n$, $1<b<n$ keine Primzahl, so wird
\begin{equation*}
[0]_n=[n]_n=[a]_n\cdot[b]_n
\end{equation*}
mit den Nullteilern $[a]_n$, $[b]_n$, siehe Definition~\ref{def:2_3}, und $\Z_n$ ist wieder kein Integrit\"atsbereich. Schlie{\ss}lich sei $n=p\geq2$ eine Primzahl und $[a]_p\cdot[b]_p=[0]_p$. Dann gilt $p|a\cdot b$, also $p|a$ oder $p|b$ und mithin $[a]_p=[0]_p$ oder $[b]_p=[0]_p$, so dass $\Z_p$ Integrit\"atsbereich ist.\\

Wir zeigen, dass $\Z_p$ ein K\"orper ist. Hierzu sei $[a]_p\neq [0]_p$, d.h. $p\nmid a$ und somit auch $\ggT(a,p)=1$. Der Euklidische Algorithmus liefert ganze Zahlen $b$, $k$ mit $a\cdot b+k\cdot p=1$, und hieraus folgt
\begin{equation*}
[a]_p\cdot[b]_p=[1]_p,
\end{equation*}
so dass $[b]_p$ die multiplikative Inverse\index{multiplikative Inverse}\label{multiplikative Inverse} zu $[a]_p$ liefert. Damit ist Satz~\ref{defsatz:6_3} bewiesen.
\dokendProof
\begin{Beis}\label{defsatz:6_4}
	Betrachte den K\"orper $\Z_3$ zum Primzahlmodul $n=3$. In folgenden Tabellen rechnet man gem\"a{\ss} Satz~\ref{satz:6_2} nur noch mit den eindeutig bestimmtem Divisionsresten $\Mod n$, $0\leq a\leq n-1$, d.h. man ersetzt $[a]_n$ durch diese Werte von $a$.
\begin{equation*}
		\begin{tabular}{c|ccc}
		 ~~ $+$ ~~  & ~~ $0$ ~~  & ~~ $1$ ~~ & ~~ $2$ ~~  \\ \hline
			0  & \rule{0pt}{2.3ex}   0  &  1   &  2   \\ 
			1  & \rule{0pt}{2.3ex}   1  &  2   &  0   \\ 
			2  & \rule{0pt}{2.3ex}   2  &  0   &  1   \\ 
		\end{tabular}\quad\quad\quad
		\begin{tabular}{c|ccc}
		~~ $\cdot$ ~~  & ~~ $0$ ~~  & ~~ $1$ ~~ & ~~ $2$ ~~  \\ \hline
		0  & \rule{0pt}{2.3ex}   0  &  0   &  0   \\ 
		1  & \rule{0pt}{2.3ex}   0  &  1   &  2   \\ 
		2  & \rule{0pt}{2.3ex}   0  &  2   &  1   \\ 
		\end{tabular}
	\end{equation*}
	Tabellen: Addition und Multiplikation modulo $3$.\\
	
	F\"ur $\Z_4$ erhalten wir dagegen keinen K\"orper:
\begin{equation*}
\begin{tabular}{c|cccc}
~~ $+$ ~~  & ~~ $0$ ~~  & ~~ $1$ ~~ & ~~ $2$ ~~ & ~~ $3$ ~~ \\ \hline
0  & \rule{0pt}{2.3ex}   0  &  1   &  2  & 3  \\ 
1  & \rule{0pt}{2.3ex}   1  &  2   &  3  & 0 \\ 
2  & \rule{0pt}{2.3ex}   2  &  3   &  0  & 1 \\ 
3  & \rule{0pt}{2.3ex}   3  &  0   &  1  & 2 \\ 
\end{tabular}\quad\quad\quad
\begin{tabular}{c|cccc}
~~ $\cdot$ ~~  & ~~ $0$ ~~  & ~~ $1$ ~~ & ~~ $2$ ~~ & ~~ $3$ ~~ \\ \hline
0  & \rule{0pt}{2.3ex}   0  &  0   &  0  & 0  \\ 
1  & \rule{0pt}{2.3ex}   0  &  1   &  2  & 3 \\ 
2  & \rule{0pt}{2.3ex}   0  &  2   &  0  & 2 \\ 
3  & \rule{0pt}{2.3ex}   0  &  3   &  2  & 1 \\ 
\end{tabular}
\end{equation*}	
Tabellen: Addition und Multiplikation modulo $4$.\hfill\dokendSatz
\end{Beis}
\begin{DefSatz}\label{defsatz:6_5}
	Es sei $n>1$ ein Modul und $\Z^*_n:=\{[a]_n:\ggT(a,n)=1\}$. Dann ist $(\Z^*_n,\cdot,[1]_n)$ bzgl. der Restklassenmultiplikation eine abelsche Gruppe mit dem Einselement $[1]_n$ als Neutralelement, die sogenannte prime Restklassengruppe\index{prime Restklassengruppe}\label{prime Restklassengruppe} $\Mod n$. Diese besitzt genau $\varphi(n)$ Elemente, und f\"ur jedes $a\in\Z$ mit $\ggT(a,n)=1$ gilt $a^{\varphi(n)}\equiv1~\ (n)$.\hfill\dokendSatz
\end{DefSatz}
{\bf Beweis:}~
Da aus $\ggT(a,n)=\ggT(b,n)=1$ auch $\ggT(ab,n)=1$ folgt, ist $\Z^*_n$ unter der assoziativen Multiplikation abgeschlossen. Nun betrachten wir eine ganze Zahl $a$, die zu $n$ teilerfremd ist. Der Euklidische Algorithmus liefert dann ganze Zahlen $b$, $k$ mit $ab+kn=1$, so dass $[a]_n\cdot[b]_n=[1]_n$ mit dem Einselement $[1]_n\in\Z^*_n$ wird. Damit ist $(\Z^*_n,\cdot, [1]_n)$ eine abelsche Gruppe. Deren Elementeanzahl ist $|\Z^*_n|=\varphi(n)$, da $\varphi(n)$ die Anzahl der nat\"urlichen Zahlen $a$ mit $a\leq n$ liefert, die zu $n$ teilerfremd sind, siehe Satz~\ref{satz:5_7}~(b). Aus dem Satz~\ref{satz:1_8} folgt nun
\begin{equation*}
a^{\varphi(n)}\equiv 1~\ (n)\quad\text{f\"ur}\quad\ggT(a,n)=1.
\end{equation*}
\dokendProof
\begin{Bem}\label{zusatz:6_6}
	Auch f\"ur $n=1$ definieren wir noch $\Z^*_1:=\{[0]_1\}$ als einelementige abelsche Gruppe mit $a^{\varphi(n)}\equiv a\equiv1\equiv0~\ (1)$ $\mbox{~f\"ur~} a\in\Z$.\hfill\dokendSatz
\end{Bem}
\begin{Beis}\label{beis:6_7}
	\hspace*{0cm}\\\vspace{-1cm}
	\begin{enumerate}[(a)]
		\item $\Z_{12}^*=\{[1]_{12},[5]_{12},[7]_{12},[11]_{12}\}$. Die prime Restklassengruppe modulo $12$ hat $\varphi(12)=\varphi(3)\cdot\varphi(4)=2\cdot 2=4$ Elemente. Ihre Gruppentafel lautet:
		\begin{equation*}
		\begin{tabular}{c|cccc}
		~~ $\cdot$ ~~  & ~~ $1$ ~~  & ~~ $5$ ~~ & ~~ $7$ ~~ & ~~ $11$ ~~ \\ \hline
		1  & \rule{0pt}{2.3ex}   1   &  5   &  7   & 11  \\ 
		5  & \rule{0pt}{2.3ex}   5   &  1   &  11  & 7 \\ 
		7  & \rule{0pt}{2.3ex}   7   &  11  &  1   & 5 \\ 
		11  & \rule{0pt}{2.3ex}  11  &  7   &  5   & 1 \\ 
		\end{tabular}
		\end{equation*}		
		\item F\"ur die Primzahl $n=7$ haben wir $\varphi(7)=6$ und $\Z_{7}^*=\{[1]_7,[2]_7,[3]_7,[4]_7,[5]_7,[6]_7\}$. Tabelle der inversen Elemente in $\Z^*_7$:
			\begin{equation*}
			\begin{tabular}{c|cccccc}
			~~ $a$ ~~  & ~~ $1$ ~~  & ~~ $2$ ~~ & ~~ $3$ ~~ & ~~ $4$ ~~ & ~~ $5$ ~~ & ~~ $6$ ~~\\ \hline
			$a^{-1}$  & \rule{0pt}{2.3ex}   1   &  4   &  5   & 2 & 3 & 6  \\			
			\end{tabular} \quad\text{modulo}\; 7.
			\end{equation*}	\\
	\end{enumerate}
	\hfill\dokendSatz
\end{Beis}
\begin{Satz}[Chinesischer Restsatz\index{Chinesischer Restsatz}\label{Chinesischer Restsatz}]\label{satz:6_8}
	Es seien $n_1,n_2,...,n_r$ nat\"urliche Zahlen, die paarweise teilerfremd sind:
	\begin{equation*}
	\ggT(n_j,n_k)=1\quad\text{f\"ur alle}\; j,k\in\{1,...,r\}\;\text{mit}\;j\neq k.
	\end{equation*}
Sind dann $a_1,a_2,...,a_r$ weitere $r$ ganze Zahlen, dann haben die Kongruenzen
\begin{equation*}
x\equiv a_j~\ (n_j)\quad\text{f\"ur}\;j=1,...,r
\end{equation*}	
gemeinsame L\"osungen $x\in\Z$. Je zwei dieser L\"osungen sind einander modulo 
\mbox{$ n:=n_1n_2...n_r$} kongruent.
	\hfill\dokendSatz	
\end{Satz}
{\bf Beweis:}~
F\"ur $n=n_1 n_2...n_r$ ist $n'_j:=\dfrac{n}{n_j}$ eine nat\"urliche Zahl, und es gilt $$\ggT(n_j,n'_j)=1 \quad\text{f\"ur~alle~}\;j=1,...,r.$$
W\"are n\"amlich $p$ ein gemeinsamer Primteiler von $n_j$ und $n'_j$, so m\"usste $r>1$ sein und $p$ schon einen der Faktoren $n_k$ mit $k\neq j$ teilen, aus denen $n'_j$ zusammengesetzt ist. Man h\"atte dann den Widerspruch $p|n_j$ $\wedge$ $p|n_k$. Nach Satz~\ref{defsatz:6_5} gibt es ganze Zahlen $b_j$ mit
\begin{equation}\label{eq:6_1}
n'_j b_j\equiv1~\ (n_j),\quad j=1,2,...,r,
\end{equation}
was trivialerweise auch f\"ur Indizes $j$ mit $n_j=1$ gilt.\\
 Offenbar gilt f\"ur alle $j,k=1,...,r$ mit $j\neq k$:
\begin{equation}\label{eq:6_2}
n'_k b_k\equiv 0~\ (n_j).
\end{equation}
Nun ist $x:=\sum\limits_{k=1}^{r}n'_k b_k a_k$ eine gesuchte L\"osung, denn es gilt nach \eqref{eq:6_1} und \eqref{eq:6_2} f\"ur alle $j=1,...r:$
\begin{equation*}
x=n'_j b_j a_j+\sum\limits_{\underset{k\neq j}{k=1}}^{r}n'_k b_k a_k\equiv a_j~\  (n_j).
\end{equation*}
Es sei $\tilde x$ eine weitere L\"osung des Kongruenzsystems $\tilde x\equiv a_j~\ (n_j)$, $j=1,...,r$. Dann ist $n_j|\tilde x-x$, und f\"ur $j\neq k$ haben $n_j$, $n_k$ stets verschiedene Primteiler. Nach dem Satz von der eindeutigen Primfaktorzerlegung folgt dann auch $n| \tilde x-x$, d.h. $x\equiv \tilde x~\  (n)$.
\dokendProof\\

Nach Satz~\ref{defsatz:6_5} gilt $a^{\varphi(n)}\equiv1~\ (n)$ f\"ur $\ggT(a,n)=1$, siehe auch Bemerkung~\ref{zusatz:6_6} und Satz~\ref{satz:1_8}. Dies rechtfertigt die
\begin{Def}\label{def:6_9}
	F\"ur $n\in\N$ und $a\in\Z$ sei $\ggT(a,n)=1$. Dann bezeichnet $h=\ord_n(a)$ die kleinste nat\"urliche Zahl mit $a^h\equiv 1~\ (n)$. Wir nennen $h$ die Ordnung\index{Ordnung}\label{Ordnung} oder auch den Exponenten\index{Exponent}\label{Exponent} von $a$ modulo $n$.\dokendDef
\end{Def}
In dieser Definition ist $h\geq1$. Bei $\ggT(a,n)=1$ k\"onnen wir f\"ur alle $j\in\Z$ Potenzen $a^j$ modulo $n$ berechnen: Nach Satz~\ref{defsatz:6_5} gibt es ein $b\in\Z$ mit $\ggT(b,n)=1$ und $a\cdot b\equiv1~\ (n)$, wobei $b$ modulo $n$ eindeutig ist. Auch f\"ur $j<0$ erhalten wir dann aus $a^j\equiv b^{-j}~\ (n)$ einen modulo $n$ zu $a^j$ kongruenten Divisionsrest $r\in\N_0$. Dann gilt\\
\begin{Satz}\label{satz:6_10}
	F\"ur $h=\ord_n(a)$ gilt $h|\varphi(n)$. Des weiteren ist die Kongruenz $a^j\equiv a^k~\ (n)$ f\"ur alle $j,k\in\Z$ genau dann erf\"ullt, wenn $j\equiv k~\ (h)$ gilt.
\hfill\dokendSatz	
\end{Satz}
{\bf Beweis:}~
Wir dividieren $j$ und $k$ mit Rest jeweils durch $h$, also  $$j=r+\left\lfloor\dfrac{j}{h}\right\rfloor h, \quad k=r'+\left\lfloor\dfrac{k}{h}\right\rfloor h$$ mit ganzzahligen Resten $r,r'\in\{0,...,h-1\}$. Dann ist wegen \mbox{$a^{\lambda h}\equiv 1~\ (n)$} $\text{~f\"ur~}\lambda\in\Z$ die Kongruenz $a^j\equiv a^k~\  (n)$ zu $a^r\equiv a^{r'}~\  (n)$ \"aquivalent.\\
W\"are $r\neq r'$, etwa $r<r'$, so w\"urde aus $a^r\equiv a^{r'}~\ (n)$ die Kongruenz $a^{r'-r}\equiv1~\ (n)$ mit $1\leq r'-r<h$ folgen, im Widerspruch zur Minimalit\"at von $h\in\N$. Es muss also $r=r'$ gelten, und die Kongruenzen $a^j\equiv a^k~\  (n)$ sowie $j\equiv k~\ (h)$ sind \"aquivalent.
\dokendProof
\begin{Satz}\label{satz:6_11}
	Aus $\ggT(a,n)=1$ und $\ord_n(a)=h$ folgt
	\begin{equation*}
	\ord_n(a^k)=\dfrac{h}{\ggT(h,k)}\quad \text{f\"ur~alle~} k\in\Z.
	\end{equation*}
	\hfill\dokendSatz	
\end{Satz}
{\bf Beweis:}~
Es ist $s:=\ord_n(a^k)$ die kleinste Zahl $s\in\N$ mit $a^{k s}\equiv1~\ (n)$, d.h. mit $h|ks$ nach Satz~\ref{satz:6_10}. Die letzte Bedingung ist \"aquivalent dazu, dass
$\dfrac{h}{\ggT(h,k)}$ ein Teiler von $\dfrac{k}{\ggT(h,k)}\cdot s$ 
und somit auch von $s$ ist,  da $\dfrac{h}{\ggT(h,k)}$ und $\dfrac{k}{\ggT(h,k)}$ teilerfremd sind. 
Aus der Minimalit\"at von $s$ folgt nun $s=\dfrac{h}{\ggT(h,k)}$.
\dokendProof
\begin{Satz}\label{satz:6_12}
	Es sei $q(x)=\sum\limits_{k=0}^{n}a_k x^k$ ein Polynom mit ganzzahligen Koeffizienten\index{Polynom mit ganzzahligen Koeffizienten}\label{Polynom mit ganzzahligen Koeffizienten} $a_0,...,a_n$ f\"ur $n\in\N_0$, und es sei $p$ eine Primzahl, die kein Teiler von $a_n$ ist. Dann gibt es modulo~$p$ nicht mehr als $n$ zueinander inkongruente L\"osungen $y\in\Z$ von $q(y)\equiv0~\ (p)$.
	\hfill\dokendSatz	
\end{Satz}
{\bf Beweis:}~
Wenn der Satz nicht gilt, gibt es ein Polynom $\tilde q(x)=\sum\limits_{k=0}^{m}b_k x^k$ mit minimalem Grad $m\in\N$, zu dem wir eine Primzahl $p$ mit $p\nmid b_m$ sowie $m+1\mod p$ zueinander inkongruente L\"osungen 
\begin{equation*}
x_1,...,x_m,x_{m+1}\in\Z
\end{equation*}
mit $\tilde q(x_\alpha)\equiv0~\ (p)$ f\"ur $\alpha =1,...,m+1$ finden. Es muss dann $m\geq 1$ gelten, da f\"ur $m=0$ niemals $\tilde q(x)=b_0$ durch $p$ teilbar ist. Nun setzen wir 
\begin{equation*}
Q(x):=\tilde q(x)-b_m\cdot\prod\limits_{k=1}^{m}(x-x_k).
\end{equation*}
Hierbei k\"urzt sich der Anteil der gr\"o{\ss}ten Potenz $x^m$ heraus. Es ist $$Q(x_{m+1})\equiv -b_m\prod\limits_{k=1}^{m}(x_{m+1}-x_k)~\ (p),$$
 wobei kein Faktor $b_m$ bzw. $x_{m+1}-x_k$ durch $p$ teilbar ist. Daher ist $Q(x)=\sum\limits_{j=0}^{N}a_j x^j$ ein Polynom vom Grad $N <m$, dessen Koeffizienten $a_j$ nicht alle durch $p$ teilbar sind. Ist $n\leq N$ der gr\"o{\ss}te Index, f\"ur den $a_n$ nicht durch $p$ teilbar ist, so setzen wir
\begin{equation*}
q(x):=\sum\limits_{j=0}^{n}a_j x^j.
\end{equation*}
Nun gilt $p\nmid a_n$ mit $n<m$, $q(y)\equiv Q(y)~\ (p)$ f\"ur alle $y\in\Z$ sowie f\"ur $\alpha=1,..,m$ $q(x_\alpha)\equiv Q(x_\alpha)\equiv\tilde q(x_\alpha)\equiv0~\  (p)$, und wieder muss $n\geq 1$ sein. Damit hat $\tilde q(x)$ nicht den geforderten Minimalgrad $m$.\dokendProof
\begin{Def}\label{def:6_13}
Es sei $\ggT(a,n)=1$. Wir nennen $a$ eine Primitivwurzel\index{Primitivwurzel}\label{Primitivwurzel} $\Mod$ $n$, wenn gilt: $$\ord_n(a)=\varphi(n).$$
\hfill\dokendDef	
\end{Def}
Primitivwurzeln haben demnach die gr\"o{\ss}tm\"ogliche Ordnung. Jedoch gibt es nicht zu jedem Modul $n$ Primitivwurzeln:	

\begin{Satz}\label{satz:6_14}
	\hspace*{0cm}\\\vspace{-1cm}
	\begin{enumerate}[(a)]
		\item Ist $p$ eine Primzahl, so gibt es $\varphi(p-1)$ Primitivwurzeln $\Mod$ $p$.\\
		
		\item Die einzigen Moduln $n\in\N$ , zu denen es Primitivwurzeln gibt, sind gegeben durch
		$$n=1, 2, 4, p^\beta, 2p^\beta,$$
		wobei $\beta$ eine beliebige nat\"urliche Zahl ist und $p$ eine ungerade Primzahl.
		\end{enumerate}
	\hfill\dokendSatz	
\end{Satz}
{\bf Beweis:}~
(a) Jedes $a\in\N$ mit $1\leq a\leq p-1$ ist zu $p$ teilerfremd und besitzt eine Ordnung $h:=\ord_p (a)\mod p$ mit $h|p-1$. F\"ur $k=0,...,h-1$ gilt dann 
	\begin{equation*}
	(a^k)^h\equiv 1~\ (p),
	\end{equation*}
	und $1,a^1,...,a^{h-1}$ sind nach Satz~\ref{satz:6_10} zueinander $\Mod p$ inkongruent. Somit stellen nach Satz~\ref{satz:6_12} die $a^k$ mit $0\leq k\leq h-1$ alle ganzzahligen L\"osungen von $x^h-1\equiv0~\ (p)$ dar. Davon gibt es nach Satz~\ref{satz:6_11} und Satz~\ref{satz:5_7}~(b) genau $\varphi(h)$ L\"osungen, die exakt die Ordnung $h\mod p$ besitzen, n\"amlich die $a^k$ mit $0\leq k \leq h-1$ und $\ggT(h,k)=1$. Insbesondere ist jedes $a$ mit $p\nmid a$ und $\ord_p(a)=h$ L\"osung von $x^h-1\equiv0 ~\ (p)$. \\
	
	Jedes $a\in\N$ mit $1\leq a\leq p-1$ besitzt $\Mod p$ eine Ordnung $h$, die $p-1$ teilt. Daher gilt, wenn $\varphi_p(h)$ die Anzahl dieser Zahlen $a$ mit Ordnung $h\mod p$ bezeichnet:
	\begin{equation}\label{eq:6_3}
	\sum\limits_{h|p-1}\varphi_p(h)=p-1.
	\end{equation} 
	Nach den vorigen Betrachtungen ist entweder $\varphi_p(h)=\varphi(h)$ oder $\varphi_p(h)=0$, je nachdem, ob es \"uberhaupt ein $a\in\{1,...,p-1\}$ mit Ordnung $h$ gibt oder nicht. Es folgt
	\begin{equation}\label{eq:6_4}
	\varphi_p(h)\leq \varphi(h)\quad\text{f\"ur jedes}~\ h\in\N ~\ \text{mit} ~\ h|p-1.
	\end{equation}
	Nach Satz~\ref{satz:5_7}~(b) gilt zudem
	\begin{equation}\label{eq:6_5}
	\sum\limits_{h|p-1}\varphi(h)=p-1.
	\end{equation}
	Aus \eqref{eq:6_3} bis \eqref{eq:6_5} folgt sofort $\varphi_p(h)=\varphi(h)$ f\"ur jedes $h|p-1$. Insbesondere ist $\varphi_p(p-1)=\varphi(p-1)\geq 1$, womit (a) bewiesen ist.\\
	
 (b) Den Nachweis von (b) f\"uhren wir in vier Schritten durch:
\begin{enumerate}[\text{Schritt} 1:]
	\item Wir zeigen, dass die Moduln $n=2^k$ f\"ur $k=0,1,2$ Primitivwurzeln besitzen, f\"ur $k\geq 3$ dagegen nicht. Wir beginnen mit $n=1,2,4$:
	\begin{center}
		\begin{tabular}{|c|l|} \hline
		\rule{0pt}{3ex}	~~ Modul $n$ ~~ & ~~  Primitivwurzeln $\Mod n$: ~~\\ \hline
			\rule{0pt}{2.5ex}	1 & \quad ~ 1 (bzw. 0)\\ 
			\rule{0pt}{2.5ex}	2 & \quad ~ 1 \\ 
			\rule{0pt}{2.5ex}	4 & \quad ~ 3 \\ \hline
		\end{tabular}
		\end{center}			
	Nun sei $n=2^k$ mit $k\geq 3$ und $a\in\N$ ungerade. Wir zeigen induktiv:
	\begin{equation}\label{eq:6_6}
	a^{2^{k-2}}\equiv1~\ (2^k)\quad\text{f\"ur alle}~\ k\geq 3.
	\end{equation}
	Wegen $2^{k-2}=\dfrac{1}{2}\varphi(2^k)=\dfrac{1}{2}\varphi(n)$ ist dann Schritt 1 getan. Mit $a=2t+1$ wird $a^2=4t(t+1)+1\equiv 1~\ (8)$, da $t(t+1)$ gerade ist. Somit gilt \eqref{eq:6_6} f\"ur $k=3$. Nun nehmen wir an, \eqref{eq:6_6} sei f\"ur ein $k\geq3$ schon gezeigt, d.h. $a^{2^{k-2}}=1+2^k\cdot u$ mit einem $u\in\N_0$.\\
	Durch Quadrieren folgt hieraus der Induktionsschritt:
	\begin{equation*}
	a^{2^{k-1}}=1+2^{k+1}u+2^{2k} u^2\equiv 1~\ (2^{k+1}).
	\end{equation*} 
	
	\item Nun sei $n=2^k\cdot\prod\limits_{j=1}^{m}p_j^{\alpha_j}>2$ mit $k\in\N_0$, $m,\alpha_j\in\N$ und paarweise verschiedenen ungeraden Primzahlen $p_j$, $j=1,...,m$. Wir zeigen, dass $n$ f\"ur $m\geq2$ oder f\"ur $m=1$, $k\geq2$ keine Primitivwurzel besitzt:\\
	Die Zahlen $n_1:=p_1^{\alpha_1}$ und $n_2:=\frac{n}{n_1}$ sind teilerfremd mit $n=n_1\cdot n_2$ und $\varphi(n)=\varphi(n_1)\varphi(n_2)$. Es ist $\varphi(n_1)=p_1^{\alpha_1-1}(p_1-1)$ gerade. Nun sei $m\geq2$ oder $m=1$, $k\geq2$. Dann ist auch $\varphi(n_2)$ gerade. Es sei $a\in\N$ zu $n$ teilerfremd. Dann ist $a$ auch zu $n_1$ und $n_2$ teilerfremd, und nach Satz~\ref{defsatz:6_5} gilt
	\begin{equation*}
	a^{\varphi(n_1)}\equiv 1~\ (n_1), \quad a^{\varphi(n_2)}\equiv 1~\ (n_2),
	\end{equation*}
	also wegen $\ggT(n_1,n_2)=1$:
	\begin{equation*}
	a^{\frac{1}{2}\varphi(n)}=(a^{\varphi(n_1)})^{\frac{1}{2}\varphi(n_2)}=(a^{\varphi(n_2)})^{\frac{1}{2}\varphi(n_1)}\equiv 1~\ (n),
	\end{equation*}
	und $a$ ist keine Primitivwurzel $\Mod n$.\\
	
	\item Es sei $n=p^\beta$ mit ungerader Primzahl $p$. F\"ur die Suche nach Primitivwurzeln $\Mod n$ d\"urfen wir $\beta\geq2$ nach der bereits bewiesenen Aussage (a) voraussetzen und eine Primitivwurzel $b\mod p$ als gegeben betrachten. Dann ist f\"ur jedes $t\in\Z$ mit $b$ auch $\tilde a:=b(1+tp)$ Primitivwurzel $\Mod p$. Es gilt $b^{p-1}=1+sp$ f\"ur ein $s\in\Z$, und weiter mit Hilfe des binomischen Lehrsatzes:
	\begin{equation*}
    \begin{array}{lcl}
		\tilde a^{p-1} & = & (1+sp)(1+tp)^{p-1}  \\
				 & \equiv & (1+sp)(1+tp(p-1))~\ (p^2) \\
			  	 & \equiv & 1+(s-t)p ~\ (p^2) .
	\end{array}
	\end{equation*}
	F\"ur $t\not\equiv s~\ (p)$	ist $\tilde a$ eine Primitivwurzel $\Mod p$ mit $\tilde a^{p-1}\not\equiv1~\ (p^2)$.\\
	Es gibt also immer eine Primitivwurzel $a$  $\Mod p$ mit $a^{p-1}\not\equiv1~\ (p^2)$.
	Wir zeigen nun, dass dieses $a$ bereits eine Primitivwurzel $\Mod p^\beta$ ist: \\
	Aus der Darstellung
	\begin{equation}\label{eq:6_7}
	a^{p-1}=1+\lambda p\quad \text{mit}~\ \lambda\not\equiv0~\ (p)
	\end{equation}
	folgt mit vollst\"andiger Induktion f\"ur alle $j\in\N_0$:
	\begin{equation*}
a^{p^j(p-1)}\equiv1+\lambda p^{j+1}~\ (p^{j+2}).	
	\end{equation*}
	Setzen wir $j:=\beta-2$, so erhalten wir
	\begin{equation}\label{eq:6_8}
	a^{p^{\beta-2}(p-1)}\equiv1+\lambda p^{\beta-1}~\ (p^{\beta}).	
	\end{equation}
	Es sei $d$ die Ordnung von $a$ modulo $p^\beta$. Es ist $\varphi(p^\beta)=p^{\beta-1}(p-1)$, und Satz~\ref{satz:6_10} liefert
	\begin{equation}\label{eq:6_9}
	d|p^{\beta-1}(p-1).
	\end{equation} 
	Es gilt $a^d\equiv1~\ (p^\beta)$, und hieraus folgt $a^d\equiv a^0~\ (p)$. Nun wenden wir Satz~\ref{satz:6_10} auf die letzte Kongruenz an, und beachten, dass die Primitivwurzel $a$ modulo $p$ die Ordnung $p-1$ hat. Wir erhalten $d\equiv0~\ (p-1)$, d.h. 
	\begin{equation}\label{eq:6_10}
	p-1|d.
	\end{equation}
	Aus \eqref{eq:6_9} und \eqref{eq:6_10} folgt mit einem Exponenten $k\leq\beta-1$:
	\begin{equation}\label{eq:6_11}
	d=p^k(p-1).
	\end{equation}
	Wegen \eqref{eq:6_8} und $\lambda\not\equiv0~\ (p)$ in \eqref{eq:6_7} ist $k\leq\beta-2$ ausgeschlossen, da $a^d\equiv1~\ (p^\beta)$ gelten muss. Somit gilt $k=\beta-1$, und aus \eqref{eq:6_11} folgt $d=\varphi(p^\beta)$, so dass $a$ in der Tat Primitivwurzel $\Mod p^\beta$ ist.\\
	
	\item Es sei $n=2 p^\beta$ mit einer ungeraden Primzahl $p$ und $\beta\in\N$, sowie $b$ eine Primitivwurzel $\Mod p^\beta$. Nun setzen wir 
	\begin{equation}\label{eq:6_12}
	a:=\begin{cases}
	b, &\text{falls $b$ ungerade ist},\\
	b+p^\beta, &\text{falls $b$ gerade ist}.
	\end{cases}
	\end{equation}
	Dann liefert \eqref{eq:6_12} eine ungerade Primitivwurzel $a$ $\Mod p^\beta$, und es gilt $\ggT(a,n)=1$. Es sei $d$ die Ordnung von $a$ modulo $n$. Dann gilt $d|\varphi(n)$ mit $\varphi(n)=\varphi(p^\beta)$, also $d|\varphi(p^\beta)$. Aus $a^d\equiv1~\ (2p^\beta)$ folgen aber auch $a^d\equiv1~\ (p^\beta)$ sowie $\varphi(p^\beta)|d$, da $a$ Primitivwurzel $\Mod p^\beta$ ist. Es folgt endlich $d=\varphi(p^\beta)=\varphi(2p^\beta)$, so dass $a$ Primitivwurzel $\Mod 2p^\beta$ ist.
	%
	% Damit ist der Satz bewiesen.
\end{enumerate}
\dokendProof
\begin{Def}\label{def:6_15}
	Es sei $n\in\N$.
	\begin{enumerate}[(a)]
		\item Wir sagen, die Zahlen $a_1, a_2,...,a_n$ bilden ein vollst\"andiges Restsystem\index{vollst\"andiges Restsystem}\label{vollstaendiges Restsystem} $\Mod n$, wenn $\Z_n=\{[a_1]_n,[a_2]_n,...,[a_n]_n\}$ gilt. Man beachte, dass dann die $a_j$ f\"ur \mbox{$j=1,...,n$} 
		zueinander modulo $n$ inkongruent sind.\\
		
		\item Wir sagen, die Zahlen  $a_1, a_2,...,a_{\varphi(n)}$ bilden ein reduziertes Restsystem\index{reduziertes Restsystem}\label{reduziertes Restsystem} $\Mod n$, wenn $\Z^*_n=\{[a_1]_n,[a_2]_n,...,[a_n]_{\varphi(n)}\}$ gilt. Man beachte, dass dann die $a_j$ f\"ur $j=1,...,\varphi(n)$ zueinander modulo $n$ inkongruent und alle zum Modul $n$ teilerfremd sind.\\
	\end{enumerate}\hfill\dokendDef
\end{Def}
\begin{Bem}\label{bem:6_16}
	Ist $n=1,2,4,p^\beta, 2p^\beta$ ein Modul aus Satz~\ref{satz:6_14} mit einer ungeraden Primzahl $p$ und $\beta\in\N$ und $a$ eine Primitivwurzel $\Mod n$, so ist die Gruppe $\Z^*_n$ zyklisch, denn die Potenzen
	\begin{equation*}
	a^1,a^2,...,a^{\varphi(n)}
	\end{equation*}
	bilden ein reduziertes Restsystem $\Mod n$. Von diesen sind nach Satz~\ref{satz:6_11} genau die Potenzen $a^k$ mit $1\leq k\leq \varphi(n)$ und $\ggT(\varphi(n),k)=1$ Primitivwurzeln $\Mod n$, so dass es genau $\varphi(\varphi(n))$ Primitivwurzeln $\Mod n$ gibt.\\
	
	Aus dem Beweis von Satz~\ref{satz:6_14} geht klar hervor, dass aus der Kenntnis der Primitivwurzeln zu Primzahl-Moduln\index{Primzahl-Modul}\label{Primzahl-Modul} $p$ sofort die Primitivwurzeln zu obigen Moduln $n$ gewonnen werden k\"onnen. Aus diesem Grund tabelliert man meist nur die Primitivwurzeln\index{Primitivwurzel}\label{Primitivwurzel2} $\Mod p$.\hfill\dokendSatz
\end{Bem}
\begin{Beis}\label{satz:6_17}
	Ist $a$ eine Primitivwurzel $\Mod n$ und $\lambda_1,\lambda_2,...,\lambda_{\varphi(n)}$ ein vollst\"andiges Restsystem $\Mod \varphi(n)$, so bilden die Potenzen $a^{\lambda_1},a^{\lambda_2},...,a^{\lambda_{\varphi(n)}}$ ein reduziertes Restsystem $\Mod n$:\\
	Es sei $p=3$, $n=p^2=9$, $\varphi(n)=9-3=6$. Dann ist $b=2$ Primitivwurzel $\Mod 3$, und wegen $2^{3-1}=4\not\equiv 1~\ (3^2)$ auch Primitivwurzel $\Mod 9$. Die Zahlen $0$, $\pm1$, $\pm2$, $3$ bilden ein vollst\"andiges Restsystem $\Mod 6$, und somit 
	\begin{equation*}
	1=2^0,\quad 2=2^1,\quad 5\equiv 2^{-1}~\ (9), \quad 4=2^2,\quad 7\equiv2^{-2}~\ (9),\quad 8=2^3
	\end{equation*}
	ein reduziertes Restsystem $\Mod 9$. Davon gibt es nur $\varphi(\varphi(9))=\varphi(6)=2$ Primitivwurzeln $\Mod 9$, n\"amlich $2$ und $5$.\\
	Tabelle der Ordnungen: Die zyklische Gruppe $\Z^*_9$ wird von $[2]_9$ und $[5]_9$ erzeugt. 
	\begin{equation*}
	\begin{tabular}{c|cccccc}	
	~~ $k$ ~~  & ~~ $1$ ~~  & ~~ $2$ ~~ & ~~ $4$ ~~ & ~~ $5$ ~~ & ~~ $7$ ~~ & ~~ $8$ ~~\\ \hline
	~~ 	$\ord_9(k)$ ~~  & \rule{0pt}{2.3ex}   1   &  6   &  3   & 6 & 3 & 2  \\		
	\end{tabular}\hfill\dokendSatz	
	\end{equation*}	
\end{Beis}

		\section{Aufgaben}\label{cha:6_A}

\begin{Auf}[Der Wilsonsche Satz\index{Wilsonscher Satz}\label{Wilsonscher Satz}]\label{auf:6_1}
Man zeige: F\"ur jede Primzahl $p$ und nur f\"ur Primzahlen $p$ gilt bei $p>1$:
$$ (p-1)! \equiv -1 ~\,(p)\,.$$
\end{Auf}
{\bf L\"osung:}\\
F\"ur $p=2$ ist $(2-1)!=1\equiv -1~\ (2)$, und f\"ur $p=3$ haben wir $(3-1)!=2\equiv -1~\ (3)$. Nun sei $p\geq 5$ eine Primzahl. Dann besitzt die Kongruenz $x^2 \equiv 1~\ (p)$ f\"ur $1\leq x< p$ genau die beiden L\"osungen $x=1$ bzw. $x=p-1$. Dies folgt aus Satz \ref{satz:6_12} mit $q(x)=x^2-1$. Die \"ubrigen von Null verschiedenen Divisionsreste $2,3,...,p-2 \mod p$ lassen sich f\"ur $\alpha=1,2,...,\dfrac{p-3}{2}$ zu paarweise disjunkten Mengen $\{x_\alpha, y_\alpha \}$ mit $x_\alpha\cdot y_\alpha\equiv 1~\ (p)$ und $x_\alpha\neq y_\alpha$ zusammenfassen. Es folgt $(p-1)!=1\cdot(p-1)\prod\limits_{\alpha=1}^{\frac{p-3}{2}}(x_\alpha\cdot y_\alpha)\equiv -1 ~\ (p)$. Ist schlie{\ss}lich $n=ab$ aus den nat\"urlichen Zahlen $a>1$, $b>1$ zusammengesetzt, so enth\"alt das Produkt $(n-1)!$ die Faktoren $a$ und $b$, und folglich wird 
\begin{equation*}
(n-1)!\equiv 0\not\equiv-1~\ (n).
\end{equation*}

\begin{Auf}[Ordnungen in der primen Restklassengruppe $\Z^*_{17}$]\label{auf:6_2}
F\"ur die Zahlen $a \in \{\pm 1, \ldots,\pm 8\}$ fertige man eine Tabelle der Ordnungen 
von $a\mod 17$ an. Welche davon sind Primitivwurzeln $\Mod 17$?
\end{Auf}
{\bf L\"osung:}\\
Die Ordnungen von $k$ und $-k \mod 17$ stimmen f\"ur $k=2,3,...,8$ \"uberein, da sie Teiler von $\varphi(17)=16=2^4$ gr\"o{\ss}er als Eins und somit gerade Zahlen sind.\\

	Tabelle der Ordnungen:
	\begin{equation*}
	\begin{tabular}{c|c|c|c|c|c|c|c|c|c|}	
	~~ $a$ ~~  & ~~ $1$ ~~ & ~~ $-1$ ~~ & ~~ $\pm2$ ~~ & ~~ $\pm3$ ~~ & ~~ $\pm4$ ~~ & ~~ $\pm5$ ~~ & ~~ $\pm6$ ~~ & ~~ $\pm7$ ~~ & ~~ $\pm8$ ~~\\ \hline
	~~ 	$\rule{0pt}{2.3ex}\ord_{17}(a)$ ~~  &  1 &  2   & 8   & 16 & 4 & 16 & 16 & 16 & 8 \\		
	\end{tabular}
	\end{equation*}	
	Davon sind Primitivwurzeln $\Mod 17$:
	$
	\pm3,\; \pm5,\; \pm6,\; \pm7.
	$
	
\begin{Auf}[Die $b$-adische Darstellung\index{$b$-adische Darstellung}\label{b-adische Darstellung} nat\"urlicher Zahlen]\label{auf:6_3}
Wir verwenden die nat\"urliche Zahl $b>1$ als Basis f\"ur $b$-adische Zahldarstellungen.
\begin{itemize}
\item[(a)] ~Man zeige f\"ur alle $n \in \N_0$: Jede ganze Zahl $a$ mit $0 \leq a<b^{n+1}$ besitzt genau eine 
$b$-adische Darstellung $$a=a_n b^n+a_{n-1} b^{n-1}+\ldots+a_0$$
mit den Ziffern $a_0,\ldots,a_{n-1},a_{n} \in \{0,1,\ldots,b-1\}$. 
\item[(b)] ~Die nat\"urliche Zahl $a$ besitze die Dezimaldarstellung $$a=a_n 10^n+a_{n-1} 10^{n-1}+\ldots+a_0$$
mit den Ziffern $a_0,\ldots,a_{n-1},a_{n} \in \{0,1,\ldots,9\}$. Man zeige:
F\"ur die Quer\-summe $\begin{displaystyle}Q(a):=\sum \limits_{j=0}^n a_j\end{displaystyle}$ 
 bzw. f\"ur die alter\-nierende Quer\-summe $\begin{displaystyle}Q_{-}(a):=\sum \limits_{j=0}^n (-1)^j a_j\end{displaystyle}$ 
haben wir $Q(a) \equiv a~\,(9)\,$
bzw. $Q_{-}(a) \equiv a~\,(11)\,$.
Durch iterierte Bildung von Quersummen bzw. von alternierenden Quersummen
erh\"alt man so einfache Rechenproben modulo 9 bzw. modulo 11.
\end{itemize}
\end{Auf}
{\bf L\"osung:}\\
(a) folgt durch Induktion bzgl. $n \in \N_0$. Der Induktionsanfang f\"ur $n=0$ ist mit der eindeutigen Darstellung
$a=a_0 \in \{0,1,\ldots,b-1\}$ f\"ur jedes ganze $a$ mit $0 \leq a <b$ erf\"ullt.
Wir nehmen an, die Behauptung sei f\"ur ein $n \in \N_0$ richtig, und betrachten eine beliebige
ganze Zahl $a$ mit $0 \leq a<b^{n+2}$. Nun dividieren wir $a$ durch $b^{n+1}$ mit Rest,
und erhalten nach Satz \ref{satz:2_4} eindeutig bestimmte ganze Zahlen $q \geq 0$ und $r$ mit $0 \leq r < b^{n+1}$, so dass
$a = qb^{n+1}+r$ gilt. Es ist aber
$$ q = \left\lfloor \frac{a}{b^{n+1}} \right\rfloor \leq \frac{a}{b^{n+1}} < b\,,$$
und somit $q \in \{0,1,\ldots,b-1\}$ eindeutig bestimmt. Auf den Divisionsrest $r$ wenden wir die Induktions\-annahme an,
und erhalten eindeutig bestimmte Ziffern\\ \mbox{$a_0,\ldots,a_{n-1},a_{n} \in \{0,1,\ldots,b-1\}$} mit
$$r=a_n b^n+a_{n-1} b^{n-1}+\ldots+a_0\,.$$
Hieraus folgt mit $a_{n+1}=q$ f\"ur $a$ die eindeutige Darstellung
$$a=a_{n+1} b^{n+1} + a_nb^n+a_{n-1} b^{n-1}+\ldots+a_0\,,$$
was zu zeigen war.\\
(b) ist ein Anwendung der Kongruenzrechnung und des vorigen Resultates:
Aus $10 \equiv 1 ~(9)$ bzw. $10 \equiv -1 ~(11)$ folgen $10^j \equiv 1 ~(9)$ bzw. $10^j \equiv (-1)^j ~(11)$
f\"ur alle $j \in \N_0$, und somit nach der Teilaufgabe (a) auch $Q(a) \equiv a ~(9)$ bzw. $Q_-(a) \equiv a ~(11)$,
da unabh\"angig von der Wahl der Repr\"asentanten modulo 9 bzw. 11 addiert und multipliziert werden darf,
siehe Definition und Satz \ref{defsatz:6_3}.\\

\begin{Auf}[Kongruenzen\index{Kongruenz}\label{Kongruenz2} mit den Fibonacci-Zahlen\index{Fibonacci-Zahlen}\label{Fibonacci-Zahlen4}]\label{auf:6_4}
Wir betrachten f\"ur $k \in \N_0$ die Fibonacci-Zahlen $f_k$ aus Aufgabe \ref{auf:1_4}.
Man zeige:
\begin{itemize}
\item[(a)] ~Es gilt
$\begin{displaystyle}
f_{12n} \equiv 0 ~\,(144) \mbox{~bzw.~}f_{12n} \equiv 0~\,(9) \mbox{~f\"ur~ alle~} n \in \N_0\,.
\end{displaystyle}$ 
\item[(b)] ~F\"ur alle $n,j \in \N_0$ gilt
$
f_{12n+j} \equiv (-1)^n f_j ~\,(9)\,.
$
\end{itemize}
\end{Auf}
{\bf L\"osung:}\\
$\mbox{ggT}(f_{12n},f_{12})=f_{\mbox{ggT}(12n,12)}=f_{12}=144$ liefert  $f_{12n} \equiv 0 ~\,(144)$
bzw. $f_{12n} \equiv 0~\,(9)$  f\"ur  alle $n \in \N_0\,,$ siehe Aufgabe \ref{auf:2_3}(c).
Aus Aufgabe  \ref{auf:2_3}(b) erhalten wir weiter unter Beachtung von $f_{13}=233\equiv -1~\,(9)$:
$$
f_{12+j}=f_{13} f_j + f_{12}f_{j-1} \equiv -f_j+0 \equiv -f_j~\,(9)\,,
$$
und somit gilt $f_{12n+j} \equiv (-1)^n \, f_j~\,(9)$ f\"ur alle $n \in \N_0$\,.\\

\begin{Auf}[Eine vollst\"andig multiplikative\index{vollst\"andig multiplikative Funktion}\label{vollstaendig multiplikative Funktion},  periodische Funktion\index{periodische Funktion}\label{periodische Funktion}] \label{auf:6_5}
Betrachte die 3-periodische Zahlenfolge $\chi: \N \to \Z$ mit  $\chi=(\overline{ 1,-1,0})\,,$
$\chi(1)=1$, $\chi(2)=-1$, $\chi(3)=0$ usw. Man zeige, dass $\chi$ vollst\"andig multiplikativ ist,
und berechne die Dirichlet-Inverse $\chi_*^{-1}$ zu $\chi$. F\"ur die Werte $\chi_*^{-1}(n)$
mit $1 \leq n \leq 20$ fertige man zudem eine Tabelle an.
\end{Auf}
{\bf L\"osung:}\\
Die beiden Zahlen $\pm 1$ bilden ein reduziertes Restsystem modulo 3, und f\"ur die Funktion $f : \Z_3^* \to \{+1,-1\}$
mit $f([1]_3)=1$ und $f([-1]_3)=-1$ gilt offenbar
$$
f([n \cdot m]_3)=f([n ]_3)\cdot f([m]_3)
$$
f\"ur alle ganzen Zahlen $n,m$, die nicht durch 3 teilbar sind. Speziell f\"ur nicht durch 3 teilbare {\it nat\"urliche} Zahlen 
$n,m$ erhalten wir aufgrund der 3-Periodizit\"at von $\chi$:
$$
\chi(n\cdot m)=f([n \cdot m]_3)=f([n]_3)\cdot f([m]_3) =\chi(n) \cdot \chi(m)\,.
$$
Ist dagegen zumindest eine der beiden nat\"urlichen Zahlen $n,m$ durch 3 teilbar, so gilt 
$\chi(n\cdot m)= \chi(n) \cdot \chi(m)=0$. Damit ist $\chi$ vollst\"andig multiplikativ, und Satz \ref{satz:5_10}
liefert $\chi_*^{-1} = \chi \cdot \mu$ f\"ur die Dirichlet-Inverse von $\chi$. Wir erhalten f\"ur die ersten
20 Funktionswerte von $\chi_*^{-1}$ die folgende Tabelle:

\begin{center}%\resizebox{\textwidth}{!} {
	\begin{tabularx}{\textwidth}
		{|r||R|R|R|R|R|R|R|R|R|R| R|R|R|R|R|R|R|R|R|R|}
		\hline
		$n$  &  $1$  &  $2$  &  $3$   &   $4$   &   $5$   &   $6$    &   $7$   
		&   $8$   &   $9$   &   $10$  &  $11$   &   $12$   &   $13$  &   $14$   &   $15$  &  $16$  &   $17$   &  $18$   &   $19$   &   $20$\\
		\hline
		$\chi_*^{-1}(n)$  &    1    &  1   & 0   &  0  & 1   & 0   &  -1 &  0   &  0   &  1  &  1   &  0   & -1 &  -1   &  0   &  0 &  1   &  0 &  -1   &  0 \\ \hline
	\end{tabularx}
\end{center}

	\chapter{Quadratische Reste}\label{cha:7}
Quadratische Reste bzgl. eines Moduls, auch kurz Reste genannt, 
sind die einfachsten Potenzreste, die auf mathematisch 
anspruchsvolle Fragestellungen f\"uhren.
Deren Untersuchung geht schon auf Euler, Fermat und Lagrange zur\"uck, doch erst Gau{\ss}
gab in seinem Buch "'Disquisitiones Arithmeticae"' eine systematische Theorie
an. Er war der erste, der das sogenannte
quadratische Reziprozit\"atsgesetz nicht nur bewiesen hat,
siehe Satz \ref{satz:7_9},
sondern gleich mehrere unterschiedliche Beweiszug\"ange geliefert hat.
Heute kennt man etwa $200$ Beweise dieses grundlegenden Satzes,
allerdings sind die meisten davon nur leichte Varianten von vorausgegangenen Beweisen.
Wir zitieren Gau{\ss} aus der deutschen \"Ubersetzung \cite[Art. 131]{gauss} seines 
lateinischen Originalwerkes:

"`Ist $p$ eine Primzahl von der Form $4n+1$, so wird $+p$, ist
dagegen $p$ eine solche von der Form $4n+3$, so wird $-p$ Rest
oder Nichtrest jeder Primzahl sein, welche, positiv genommen,
Rest oder Nichtrest von $p$ ist.\\
Da fast alles, was sich \"uber die quadratischen Reste sagen l\"asst, auf
diesem Satze beruht, so wird die Bezeichnung "`Fundamentalsatz"', die wir
im Folgenden gebrauchen werden, f\"ur denselben nicht unpassend sein."'

	\section{Quadratische Reste\index{quadratische Reste}\label{quadratische Reste}}\label{cha:7T}

Hier untersuchen wir f\"ur einen Modul $m=p_1^{\alpha_1}\cdot p_2^{\alpha_2}\cdot... \cdot p_k^{\alpha_k}\geq 2$ mit paarweise verschiedenen Primzahlen $p_1$, $p_2$, ..., $p_k$ und Exponenten $\alpha_1$, $\alpha_2$, ..., $\alpha_k\in\N$ f\"ur gegebenes $a\in\Z$ die L\"osungen der quadratischen Kongruenz\index{quadratischer Kongruenz}\label{quadratischer Kongruenz}
\begin{equation}\label{eq:7_1}
x^2\equiv a~\ (m).
\end{equation}
Aus \eqref{eq:7_1} folgt $x^2\equiv a~\ (p_j^{\alpha_j})$ f\"ur alle $j=1,...,k$. Sind umgekehrt die $x_j\in\Z$ f\"ur $j=1,...,k$ L\"osungen der Kongruenzen
\begin{equation*}
x_j^2\equiv a~\ (p_j^{\alpha_j}),
\end{equation*}
so liefert der chinesische Restsatz eine modulo $m$ eindeutige L\"osung $x\in\Z$ des Kongruenzsystems $x\equiv x_j~\ (p_j^{\alpha_j})$, $j=1,...,k$, so dass $x$ auch L\"osung von \eqref{eq:7_1} ist. Somit gen\"ugt es, anstelle von \eqref{eq:7_1} die Kongruenz
\begin{equation}\label{eq:7_2}
x^2\equiv a~\ (p^\alpha)
\end{equation}
nur f\"ur Primzahlpotenz-Moduln $p^\alpha$ zu l\"osen, also mit Primzahlen $p$ und Exponenten $\alpha\in\N$.\\

F\"ur $a\equiv 0~\ (p^\alpha)$ erh\"alt man nur triviale L\"osungen, n\"amlich genau die ganzen Zahlen $x$ mit $x\equiv0~\ (p^{\lceil\frac{\alpha}{2}\rceil})$, wobei $\lceil y\rceil=-\lfloor -y\rfloor=\min\{k\in\Z: k\geq y\}$ f\"ur $y\in\R$ ist.\\

F\"ur $a\not\equiv 0~\ (p^\alpha)$ ist dagegen \eqref{eq:7_2} h\"ochstens dann l\"osbar, wenn $a=p^{2\beta}\cdot \tilde a$ und $x\equiv 0~\ (p^\beta)$ f\"ur ganzzahlige $\tilde a$, $\beta$ ist mit $0\leq \beta <\frac{\alpha}{2}$ und $\tilde a \not\equiv 0~\ (p)$.\\

Mit $\tilde x:=\frac{x}{p^\beta}$ muss dann nur noch $\tilde x^2\equiv\tilde a ~\ (p^{\alpha-2\beta})$ gel\"ost werden. Im Folgenden sei daher $a\not\equiv 0~\ (p)$.\\

Wir betrachten zun\"achst $p=2$. Dann muss $a$ ungerade sein, so dass nur ungerade L\"osungen $x=2k+1$ in Frage kommen mit 
\begin{equation}\label{eq:7_3}
x^2=(2k+1)^2=1+8\cdot\frac{k(k+1)}{2}\equiv 1~\ (8).
\end{equation}
F\"ur $\alpha=1$ muss $a$ nur ungerade sein und f\"ur $\alpha=2$ \"uberdies $a\equiv 1~\ (4)$ erf\"ullen, und jedes ungerade $x$ ist L\"osung von \eqref{eq:7_2}.\\

Wir zeigen, dass \eqref{eq:7_2} mit $p=2$, $\alpha\geq 3$ genau f\"ur $a\equiv1~\ (8)$ l\"osbar ist: Nach \eqref{eq:7_3} ist $a\equiv1~\ (8)$ f\"ur $\alpha\geq3$ notwendig. Die Umkehrung zeigen wir induktiv: Nach \eqref{eq:7_3} gilt der Induktionsanfang f\"ur $\alpha=3$. Wir nehmen $x_0^2\equiv a~\ (2^\alpha)$ f\"ur ein $\alpha \geq 3$ an.\\

Hiermit w\"ahlen wir ein $\lambda\in\Z$ so, dass gilt:
$$(x_0+\lambda 2^{\alpha-1})^2=x_0^2+x_0\lambda 2^\alpha+\lambda^2 2^{2\alpha-2}\equiv a~\ (2^{\alpha+1}). $$
Das ist m\"oglich, da $2\alpha-2\geq \alpha+1$ f\"ur $\alpha\geq3$ gilt und $\frac{x_0^2-a}{2^\alpha}+x_0\cdot\lambda\equiv0~\ (2)$ l\"osbar ist. Die Behauptung ist bewiesen. \\

Von nun ab betrachten wir nur noch Primzahlen $p\geq 3$ in \eqref{eq:7_2} und beachten dabei $a\not\equiv 0~\ (p)$. Damit \eqref{eq:7_2} l\"osbar ist, muss es ein $x_0\in\Z$ mit $x_0^2 \equiv a~\ (p)$ geben. Diese Bedingung ist auch hinreichend zur L\"osbarkeit von \eqref{eq:7_2}. Genauer zeigen wir induktiv:

Es gibt eine rekursiv konstruierte Folge $(x_n)_{n\in\N_0}$ ganzer Zahlen $x_n$, so dass f\"ur alle $n\in\N_0$ gilt:
\begin{equation}\label{eq:7_4}
x_n^2\equiv a~\ (p^{2^n}),\quad 2x_nx_{n+1}\equiv x_n^2+a~\ (p^{2^{n+1}}).
\end{equation}
F\"ur $n=0$ haben wir $x_0^2\equiv a~\ (p)$ vorausgesetzt, und k\"onnen wegen $2x_0\not\equiv0~\ (p)$ die Kongruenz $2x_0x_1\equiv x_0^2+a~\ (p^2)$ nach $x_1$ aufl\"osen. Ist \eqref{eq:7_4} f\"ur ein $n\geq0$ bereits gezeigt, so folgt $2x_n x_{n+1}\equiv 2x_n^2~\ (p^{2^n})$, also $x_{n+1}\equiv x_n~\ (p^{2^n})$ durch K\"urzen des Faktors $[2x_n]_{p^{2^n}}$ in der Gruppe $\Z^*_{p^{2^n}}$. Hieraus erhalten wir $$0\equiv(x_{n+1}-x_n)^2\equiv x_{n+1}^2-2x_{n+1}x_n+x_n^2\equiv  x_{n+1}^2-a~\ (p^{2^{n+1}}).$$
Indem wir noch eine L\"osung $x_{n+2}$ der Kongruenz $2x_{n+1}x_{n+2}\equiv x_{n+1}^2+a~\ (p^{2^{n+2}})$ ermitteln, was wegen $2x_{n+1}\not\equiv 0~\ (p)$ m\"oglich ist, folgt \eqref{eq:7_4} f\"ur alle $n\in\N_0$.\\

Nun definieren wir quadratische Reste:
\begin{Def}\label{def:7_1}
Es sei $m\in\N$ mit $m\geq 2$. Eine ganze Zahl $a$ mit $\ggT(a,m)=1$ hei{\ss}t quadratischer Rest $\Mod m$, wenn es ein $x\in\Z$ gibt mit 
$$x^2\equiv a~\ (m).$$ 
Damit ist notwendigerweise auch $\ggT(x,m)=1$.\hfill\dokendDef
\end{Def}
Wir erhalten nun den
\begin{Satz}\label{satz:7_2}
Bei $\ggT(a,m)=1$ und $m=2^\alpha\cdot m'\geq 2$ mit ungeradem $m'$ und $\alpha\in\N_0$ ist die Kongruenz $x^2\equiv a~\ (m)$ genau dann l\"osbar, wenn gilt:

Die Kongruenz $x^2\equiv a~\ (p)$ ist f\"ur jeden Primteiler $p$ von $m'$ l\"osbar, und \"uberdies gilt $a\equiv 1~\ (4)$ f\"ur $\alpha=2$ bzw. $a\equiv 1~\ (8)$ f\"ur $\alpha \geq3$.

Die Anzahl der L\"osungen $x$ $\Mod m$ von $x^2\equiv a~\ (m)$ ist in diesem Falle gegeben durch
$$\min\left(4,2^{\max(\alpha,1)-1}\right)\cdot 2^{\omega(m')}=\min(4,\varphi(2^\alpha))\cdot 2^{\omega(m')} $$
mit der Anzahl $\omega(m')$ der verschiedenen Primfaktoren von $m'$ (ohne Vielfachheiten).\\

\underline{Beachte}: Die Anzahl der L\"osungen von $x^2\equiv a~\ (m)$ ist insbesondere f\"ur alle quadratischen Reste $a$ $\Mod m$ dieselbe.
		\hfill\dokendSatz	
\end{Satz}
{\bf Beweis:}~
Wir m\"ussen gem\"a{\ss} den vorausgegangenen Betrachtungen nur noch die Formel f\"ur die L\"osungsanzahl zeigen: Gilt
\begin{equation*}
x^2_0\equiv a~\ (m), \quad x^2\equiv a~\ (m),
\end{equation*}
so finden wir ein $x_0^*\in\Z$ mit $x_0\cdot x_0^*\equiv 1~\ (m)$. Damit gilt 
$$(x\cdot x_0^*)^2 \equiv a\cdot x_0^{*^2}\equiv x_0 ^2\cdot x_0^{*^2}\equiv 1~\ (m)$$
sowie $x\equiv x_0\cdot y~\ (m)$ f\"ur die L\"osung $y:= x\cdot x_0^*$ der Kongruenz $y^2\equiv 1~\ (m)$. Umgekehrt liefert jede L\"osung $y$ von $y^2\equiv 1~\ (m)$ bei festem $x_0$ ein $x\equiv x_0\cdot y~\ (m)$, das L\"osung von $x^2\equiv a~\ (m)$ ist. Damit gen\"ugt es, die Anzahlformel f\"ur den einfachsten quadratischen Rest $a=1$ zu zeigen:\\

Wir beginnen mit dem Spezialfall $m=2^\alpha$:

F\"ur $\alpha=0$ bzw. $\alpha=1$ erhalten wir jeweils nur eine L\"osung von $x^2\equiv 1~\ (2^\alpha)$ modulo $2^\alpha$. F\"ur $\alpha=2$ haben wir die beiden L\"osungen $x\equiv\pm 1~\ (4)$, und schlie{\ss}lich f\"ur $\alpha\geq 3$ genau vier L\"osungen von $x^2\equiv1~\ (2^\alpha)$ modulo $2^\alpha$, n\"amlich
$$x_{1}\equiv -1,~~x_{2}\equiv 1,~~ x_{3}\equiv 2^{\alpha-1}-1,~~x_{4}\equiv2^{\alpha-1}+1~\ (2^\alpha). $$
Dies sind die einzigen, denn f\"ur sie gilt
$$x^2-1=(x-1)(x+1)\equiv0~\ (2^\alpha) $$
mit den beiden geraden Faktoren $x\pm1$, von denen jeweils genau einer nicht durch $4$ teilbar ist. Somit erh\"alt man f\"ur $m=2^\alpha$ in jedem Fall genau $\min\left(4,2^{\max(\alpha,1)-1}\right)$ L\"osungen.\\

Nun betrachten wir den Spezialfall $m=p^\alpha$ mit einer Primzahl $p\geq 3$. Dann hat die Kongruenz $x^2\equiv1~\ (p^\alpha)$ modulo $p^\alpha$ die beiden L\"osungen
$$x_{1,2}\equiv\pm1~\ (p^\alpha), $$
und wegen $x^2-1=(x-1)\cdot(x+1)\equiv 0~\ (p^\alpha)$ sind dies modulo $p^\alpha$ die einzigen, da in jedem Produkt $(x-1)(x+1)$ nur jeweils ein Faktor durch $p$ und damit schon durch $p^\alpha$ teilbar ist.\\

Da gem\"a{\ss} dem chinesischen Restsatz die L\"osungen von \eqref{eq:7_1} f\"ur einen aus paarweise teilerfremden Primzahlpotenzen $p_j^{\alpha_j}$ zusammengesetzten Modul $m$ aus den L\"osungen der Kongruenzen $x_j^2\equiv a~\ (p_j^{\alpha_j})$ hervorgehen, hier mit $a=1$, folgt die Anzahlformel durch Produktbildung.
\dokendProof\\
\begin{Beis}$$x^2\equiv13~\ (324).$$
	Hier ist $m=324=4\cdot 81=2^2\cdot 3^4$, $\alpha=2$, $m'=81$, $a=13$. Da $x^2\equiv13\equiv1~\ (3)$ modulo $3$ die L\"osungen $x=\pm 1$ hat und $x^2\equiv13\equiv1~\ (4)$ modulo $4$ die L\"osungen $x=\pm 1$, besitzt die Ausgangskongruenz genau vier L\"osungen:
	$$\min\left(4,2^{\max(\alpha,1)-1}\right)\cdot 2^{\omega(m')}=2\cdot 2=4. $$
	Mit \eqref{eq:7_4} bestimmen wir die L\"osungen von $z^2\equiv13~\ (81)$, beginnend mit $z_0=1$:
	$$	2z_0z_1 \equiv  z_0^2+13~\ (9)\ \text{ liefert }\ z_1\equiv7~\ (9). $$
	 $$ 2z_1z_2 \equiv  z_1^2+13~\ (81)\  \text{ f\"uhrt auf }\ 14z_2\equiv62~\ (81)\ \text{ bzw. }\ 7z_2\equiv31~\ (81). $$
	Wir bestimmen das multiplikative Inverse zu $7$ $\Mod 81$:
\begin{equation*}
	\begin{tabular}{|c||c|c|c|c|} \hline
	 $j$  & ~~ $q_j$ ~~  & ~~ $r_j$ ~~  & ~~ $s_j$ ~~  & ~~ $t_j$ ~~\\
		\hline
	\rule{0pt}{2.4ex} 0~  &    0    &  81   &    1    &  0 \\ \hline 
	\rule{0pt}{2.4ex} 1~  &    11   &  7    &    0    &  1 \\ \hline
	\rule{0pt}{2.4ex} 2~  &    1    &  4    &    1    &  11  \\ \hline
	\rule{0pt}{2.4ex} 3~  &    1    &  3    &    1    &  12  \\ \hline
	\rule{0pt}{2.4ex} 4~  &    3    &  1    &    2    &  23  \\ \hline
	\rule{0pt}{2.4ex} 5~  &    ---  &  0    &    7    &  81  \\ \hline
	\end{tabular}\quad\quad\quad
	\begin{array}{rl}
		& \text{Wir erhalten }\\[0.3em]
		& 2\cdot 81-23\cdot 7=1,\\[0.3em]
		&-23\cdot 7\equiv1~\ (81),\\[0.3em]
		& z_2\equiv-23\cdot 31\equiv-65~\ (81).\\[0.3em]
		& \text{Nun ist sogar }(\pm65)^2\equiv13~\ (324).\\
	\end{array}
\end{equation*}
Wir l\"osen jeweils vier simultane Kongruenzsysteme:
$$\begin{array}{llll}
	1) \quad x\equiv 1~\ (4),&\quad x\equiv -65~\ (81)&\quad \text{ liefert } &\quad x\equiv 97~\ (324),\\
	2)\quad x\equiv -1~\ (4),&\quad x\equiv -65~\ (81)&\quad \text{ liefert } &\quad x\equiv -65~\ (324),\\
	3)\quad x\equiv 1~\ (4),&\quad x\equiv 65~\ (81)&\quad \text{ liefert } &\quad x\equiv 65~\ (324),\\
	4)\quad x\equiv -1~\ (4),&\quad x\equiv 65~\ (81)&\quad \text{ liefert } &\quad x\equiv -97~\ (324).
\end{array}$$
Wir erhalten die L\"osungen $x\equiv\pm65$ bzw. $x\equiv\pm97$ $\Mod 324$ von $x^2\equiv 13~\ (324)$.
\hfill\dokendSatz
\end{Beis}
Nun charakterisieren wir quadratische Reste\index{quadratische Reste}\label{quadratische Reste2} (Q-Reste) bzw. quadratische Nichtreste\index{quadratische Nichtreste}\label{quadratische Nichtreste} (Q-Nichtreste) bzgl. eines Primzahlmoduls $p\geq3$:
\begin{Satz}[Eulersches Kriterium\index{Eulersches Kriterium}\label{Eulersches Kriterium}]\label{satz:7_4}
F\"ur Primzahlen $p\geq3$ und $a\in\Z$ definieren wir das Legendre-Symbol\index{Legendre-Symbol}\label{Legendre-Symbol}
$$(a|p):=\begin{cases}
+1,& \text{wenn } x^2\equiv a~\ (p) \text{ mit } x\not\equiv0~\ (p) \text{ l\"osbar ist},\\
-1,& \text{wenn } x^2\equiv a~\ (p) \text{ nicht l\"osbar ist},\\
\ \ \ 0,& \text{wenn  } a\equiv 0~\ (p) \text{ gilt}.\\
\end{cases} $$	
Dann ist
$$(a|p)\equiv a^{\frac{p-1}{2}}~\ (p).$$
\hfill\dokendSatz
\end{Satz}	
{\bf Beweis:}~
Da f\"ur $a\equiv 0~\ (p)$ die Behauptung stimmt, d\"urfen wir $a\not\equiv 0~\ (p)$ voraussetzen. Ist $x^2\equiv a~\ (p)$ mit $x\not\equiv 0~\ (p)$ l\"osbar, so ist $x^{p-1}\equiv 1~\ (p)$ nach Satz~\ref{defsatz:6_5}. In diesem Falle folgt
\begin{equation*}
(a|p)=1\equiv (x^2)^{\frac{p-1}{2}}\equiv a^\frac{p-1}{2}~\ (p).
\end{equation*}
Nun sei $x^2\equiv a~\ (p)$ nicht l\"osbar, d.h. $(a|p)=-1$. Wegen $$a^{p-1}-1=\left(a^\frac{p-1}{2}-1\right)\left(a^\frac{p-1}{2}+1\right)\equiv 0~\ (p)$$ kommt nur $a^\frac{p-1}{2}\equiv\pm1~\ (p)$ in Frage, so dass wir nur $a^\frac{p-1}{2}\equiv1~\ (p)$ ausschliessen m\"ussen: Die Quadrate 
\begin{equation}\label{eq:7_5}
1^2, 2^2, ... , \left(\frac{p-1}{2}\right)^2
\end{equation}
sind L\"osungen der Kongruenz $u^\frac{p-1}{2}-1\equiv 0~\ (p)$, und modulo $p$ voneinander verschieden, da $j^2-k^2=(j+k)(j-k)\not\equiv0~\ (p)$ gilt wegen $0<j+k<p$, $0<j-k<p$ f\"ur $j>k$ und $j,k\in\{1,2,...,\frac{p-1}{2}\}$. Nach Satz~\ref{satz:6_12} sind dies modulo $p$ alle L\"osungen von $u^\frac{p-1}{2}-1\equiv 0~\ (p)$.\\

Da $a$ als Q-Nichtrest vorausgesetzt wurde, ist er modulo $p$ keiner der Zahlen aus \eqref{eq:7_5} kongruent, und es folgt $a^\frac{p-1}{2}\equiv -1~\ (p)$.
\dokendProof\\
\begin{Folg}\label{folg:7_5}
Modulo einer Primzahl $p\geq 3$ gibt es genau $\frac{p-1}{2}$ Q-Reste, die einer der Zahlen
\begin{equation*}
1^2, 2^2, ... , \left(\frac{p-1}{2}\right)^2
\end{equation*}
kongruent sind, und damit auch ebensoviele Q-Nichtreste. Dabei gilt $$(a|p)\cdot(a'|p)=(aa'|p)\quad\forall a,a'\in\Z.$$
\hfill\dokendSatz
\end{Folg}
\underline{Merke}: F\"ur Primzahl $p\geq3$ gilt
$$\left.\begin{array}{ccccc}
\text{Q-Rest}& \cdot & \text{Q-Rest}& =&\text{Q-Rest},\\
\text{Q-Rest}& \cdot & \text{Q-Nichtrest}&=&\text{Q-Nichtrest},\\
\text{Q-Nichtrest}& \cdot & \text{Q-Nichtrest}&=&\text{Q-Rest}.\\
\end{array}\right\}\Mod p.$$
{\bf Beweis:}~	
Dies folgt sofort aus dem Eulerschen Kriterium und seinem Beweis.
\dokendProof\\

Setzen wir  $a:=-1$ in Satz~\ref{satz:7_4}, so erhalten wir die
\begin{Folg}\label{folg:7_6}
	F\"ur jede Primzahl $p\geq 3$ gilt $$(-1|p)=(-1)^{\frac{p-1}{2}}.$$ 
	Somit ist $-1$ genau f\"ur $p\equiv 1~\ (4)$ ein Q-Rest $\Mod p$.
	\dokendDef
\end{Folg}

Neben Satz~\ref{satz:7_4} dient auch das folgende Kriterium der Bestimmung des Restsymboles $(a|p)$:

\begin{Satz}[Gau{\ss}sches Lemma\index{Gau{\ss}sches Lemma}\label{Gausssches Lemma}, erweiterte Version]\label{satz:7_7}
	F\"ur jede Primzahl $p\geq 3$ und $\ggT(a,p)=1$ gilt:\\
	Wenn $t$ die Anzahl derjenigen kleinsten positiven Reste der Zahlen $a,2a,3a,...,\frac{p-1}{2} a$ modulo $p$ ist, die gr\"o{\ss}er als $\frac{p}{2}$ sind, dann gilt $(a|p)=(-1)^t$.
	Hierbei ist $$
	t\equiv \sum_{j=1}^{\frac{p-1}{2}}\left\lfloor\frac{ja}{p}\right\rfloor+(a-1)\frac{p^2-1}{8} \mod 2.$$
	\dokendSatz
\end{Satz}

{\bf Beweis:}~
	Wir k\"onnen die kleinsten positiven Divisionsreste von $a, 2a, 3a,..., \frac{p-1}{2}a$ bei Division durch $p$ in der Form 
	\begin{equation}\label{eq:7_6}
	r_1, r_2,..., r_s;\ p-r'_1, p-r'_2,...,p-r'_t
	\end{equation}
	darstellen mit $s+t=\frac{p-1}{2}$ und
	\begin{equation}\label{eq:7_7}
	r_1, r_2,..., r_s;\ r'_1, r'_2,...,r'_t\in \left\{1,2,...,\frac{p-1}{2}\right\}.
	\end{equation}
	Dabei sind die Reste in \eqref{eq:7_6} $\Mod p$ paarweise verschieden, denn $a\not\equiv 0~\ (p)$. Auch gibt es keine zwei Zahlen $j,k\in \left\{1,2,...,\frac{p-1}{2}\right\}$ mit $j\neq k$ und $j\cdot a\equiv-k\cdot a ~\ (p)$, da f\"ur diese $j+k\equiv 0 ~\ (p)$ mit $1<j+k<p$ gelten m\"usste, ein Widerspruch. Somit sind auch alle Reste in \eqref{eq:7_7} paarweise verschieden und m\"ussen wegen $s+t=\frac{p-1}{2}$ genau die Zahlen $1, 2, 3,...,\frac{p-1}{2}$ liefern.
	
	Es folgt f\"ur das Produkt aller Zahlen in \eqref{eq:7_6} modulo $p$:
	$$a^{\frac{p-1}{2}}\cdot \left(\frac{p-1}{2}\right)!\equiv (-1)^t \cdot \prod_{\alpha=1}^{s} r_\alpha\cdot \prod_{\beta=1}^{t} r'_\beta\equiv (-1)^t \cdot \left(\frac{p-1}{2}\right)!  \mod p,$$
	und da wir $\Mod p$ den Faktor $ \left(\frac{p-1}{2}\right)!$ k\"urzen d\"urfen:	
	\begin{equation}\label{eq:7_8}
	a^{\frac{p-1}{2}}\equiv (-1)^t \mod p.
	\end{equation}
	Aus dem Eulerschen Kriterium und \eqref{eq:7_8} folgt der erste Teil der Behauptung.\\
	
	Wir zeigen die Kongruenzformel f\"ur $t$ $\Mod 2$:	
	\begin{equation}\label{eq:7_9}
	\sum_{j=1}^{\frac{p-1}{2}}a\cdot j=a\frac{p^2-1}{8}=p\sum_{j=1}^{\frac{p-1}{2}}\left\lfloor\frac{ja}{p}\right\rfloor
	+\sum_{\alpha=1}^{s} r_\alpha+ \sum_{\beta=1}^{t} (p-r'_\beta)
	\end{equation}
	folgt mit der Bildung der Divisionsreste $r_\alpha, p-r'_\beta$.\\
	Unter Beachtung von $p-r'_\beta\equiv 1+r'_\beta ~\ (2)$, $p\equiv 1 ~\ (2)$, sowie mit $$\sum_{\alpha=1}^{s} r_\alpha+ \sum_{\beta=1}^{t} r'_\beta=\sum_{j=1}^{\frac{p-1}{2}}j=\frac{p^2-1}{8}$$ folgt aus \eqref{eq:7_9}:	
	\begin{equation}\label{eq:7_10}
	a\frac{p^2-1}{8}\equiv \sum_{j=1}^{\frac{p-1}{2}}\left\lfloor\frac{ja}{p}\right\rfloor
	+\frac{p^2-1}{8}+t\mod 2.
	\end{equation}
	Addieren wir in \eqref{eq:7_10} modulo 2 auf beiden Seiten die Summe $\sum\limits_{j=1}^{\frac{p-1}{2}}\left\lfloor\frac{ja}{p}\right\rfloor$ und subtrahieren $\frac{p^2-1}{8}$, so folgt die Behauptung.
	\dokendProof

\begin{Folg}\label{folg:7_8}
	Es sei $p\geq 3$ eine Primzahl. 
	\begin{enumerate}[(a)]
		\item $(2|p)=1\Leftrightarrow p\equiv \pm 1~\ (8)$. Allgemein gilt $(2|p)=(-1)^{\frac{p^2-1}{8}}$.
		\item $(-2|p)=1\Leftrightarrow p\equiv 1, 3~\ (8)$. Allgemein gilt $(-2|p)=(-1)^{\frac{1}{8}(p-1)(p-3)}$.
	\end{enumerate}
	\dokendDef
\end{Folg}

{\bf Beweis:}~
\begin{enumerate}[(a)]
	\item Wir setzen $a=2$ in Satz~\ref{satz:7_7} und beachten $\left\lfloor\frac{2j}{p}\right\rfloor=0$ f\"ur $j=1,...,\frac{p-1}{2}$. Damit ist 
	$$(2|p)=(-1)^{\frac{p^2-1}{8}}=(-1)^{\frac{1}{8}(p-1)(p+1)}$$
	und
	$$(2|p)=1\quad \text{genau f\"ur } p\equiv \pm 1~\ (8).$$
	\item folgt aus (a) und Folgerung~\ref{folg:7_6}: F\"ur $(-2|p)=1$ ist entweder $(2|p)=1$ und $(-1|p)=1$ mit $p\equiv 1~\ (8)$, oder $(2|p)=-1$ und $(-1|p)=-1$ mit $p\equiv 3~\ (8)$. Insgesamt ist 
	$$(-2|p)=(-1|p)\cdot (2|p)=(-1)^{-\frac{p-1}{2}}\cdot(-1)^{\frac{p^2-1}{8}}
	=(-1)^{\frac{1}{8}(p-1)(p-3)}.$$
\end{enumerate}
\dokendProof

\begin{Satz}[Das Reziprozit\"atsgesetz von Gau{\ss}\index{Reziprozit\"atsgesetz von Gau{\ss}}\label{Reziprozitaetsgesetz von Gauss}]\label{satz:7_9}
	Sind $p,q\geq 3$ zwei verschiedene Primzahlen, so gilt $(p|q)\cdot (q|p)=(-1)^{\frac{p-1}{2}\frac{q-1}{2}}$.
	\dokendSatz
\end{Satz}
{\bf Beweis:}~
Wegen $p,q\geq 3$ sind $p, q$ ungerade, und wegen $p\neq q$ gilt $(p|q)\cdot (q|p)=\pm 1$. Satz~\ref{satz:7_7} liefert 
\begin{equation}\label{eq:7_11}
(p|q)\cdot (q|p)=(-1)^{\sum\limits_{j=1}^{\frac{p-1}{2}}\left\lfloor\frac{qj}{p}\right\rfloor
	+\sum\limits_{k=1}^{\frac{q-1}{2}}\left\lfloor\frac{pk}{q}\right\rfloor}.
\end{equation}
F\"ur die Menge $G:=\left\{(j,k): j\in \left\{1,...,\frac{p-1}{2}\right\},\: k\in\left\{1,...,\frac{q-1}{2}\right\}\right\}$
gilt $qj\neq pk$ f\"ur alle $(j,k)\in G$ mit
\begin{equation}\label{eq:7_12}
|G|=\frac{p-1}{2}\cdot\frac{q-1}{2}.
\end{equation}
Somit ist $G=G_p\cup G_q$ die Vereinigung der beiden disjunkten Mengen
$$G_p:=\left\{(j,k)\in G: pk<qj\right\},\quad G_q:=\left\{(j,k)\in G: qj<pk\right\}.$$
Es besteht $G_p$ aus allen $(j,k)\in\N\times\N$ mit $j\in \left\{1,...,\frac{p-1}{2}\right\}$ und $k\leq \frac{qj}{p}$, wobei $k=\frac{qj}{p}$ nicht auftritt, und entsprechend $G_q$ aus allen $(j,k)\in\N\times\N$ mit $k\in\left\{1,...,\frac{q-1}{2}\right\}$ und $j\leq\frac{pk}{q}$, wobei $j=\frac{pk}{q}$ nicht auftritt. Wir erhalten
\begin{equation}\label{eq:7_13}
|G_p|=\sum\limits_{j=1}^{\frac{p-1}{2}}\left\lfloor\frac{qj}{p}\right\rfloor,
\quad |G_q|=\sum\limits_{k=1}^{\frac{q-1}{2}}\left\lfloor\frac{pk}{q}\right\rfloor,
\quad |G|=|G_p|+|G_q|.
\end{equation}
Aus \eqref{eq:7_11}-\eqref{eq:7_13} folgt nun die Behauptung des Satzes.
\dokendProof

\begin{Zus}\label{zus:7_10}~\hspace*{0cm}\\\vspace{-1cm}
	\begin{enumerate}[1)]
		\item Der Wert des Legendre-Symbols\index{Legendre-Symbol}\label{Legendre-Symbol2} $(a,p)\equiv a^{\frac{p-1}{2}}~\ (p)$ entscheidet f\"ur Primzahl-Moduln $p\geq 3$ \"uber die L\"osbarkeit von
		\begin{equation}\label{eq:7_14}
		x^2\equiv a~\ (p):
		\end{equation}
		Bei $(a,p)=1$ ist \eqref{eq:7_14} mit $x\not\equiv 0~\ (p)$ l\"osbar, bei $(a,p)=0$ mit $x\equiv 0~\ (p)$ und bei $(a,p)=-1$ ist \eqref{eq:7_14} unl\"osbar.
		Es ist $(a|p)=(a'|p)$ f\"ur $a\equiv a'~\ (p)$.
		\item Es gilt $\left(\prod\limits_{j=1}^n a_j|p\right)=\prod\limits_{j=1}^n \left(a_j|p\right)$ f\"ur alle $a_1,...,a_n\in\Z$.
		
		\item $(-1|p)=(-1)^\frac{p-1}{2}$ sowie $(2|p)=(-1)^\frac{p^2-1}{8}$, $(-2|p)=(-1)^{\frac{1}{8}(p-1)(p-3)}$.\\
		
		\item F\"ur jede zwei Primzahlen $p, q\geq 3$ gilt das quadratische Reziprozit\"atsgesetz\index{quadratisches Reziprozit\"atsgesetz}\label{quadratisches Reziprozitaetsgesetz}:
		$$(q|p)=(-1)^{\frac{p-1}{2}\frac{q-1}{2}}\cdot (p|q)$$
	\end{enumerate}
\end{Zus}\dokendSatz

\begin{Beis}\label{beis:7_11}~\hspace*{0cm}\\\vspace{-1cm}
	\begin{enumerate}[(a)]
		\item F\"ur welche Primzahlen $p>3$ ist 3 ein Q-Rest, f\"ur welche ein Q-Nichtrest?\\
		{\bf L\"osung:}~Aus Zusammenfassung~\ref{zus:7_10} 4) folgt mit $q:=3$:
		$$(3|p)=(-1)^{\frac{p-1}{2}}\cdot (p|3).$$
		Dabei gilt
		$$(p|3)=\begin{cases}
		\ \ \: 1,&\text{falls } p\equiv 1~\ (3),\\
		-1,&\text{falls } p\equiv -1~\ (3).
		\end{cases}$$
		Hieraus folgt $(3|p)=1$ f\"ur $p\equiv 1~\ (4)$, $p\equiv 1~\ (3)$, d.h. f\"ur $p\equiv 1~\ (12)$, oder aber f\"ur $p\equiv -1~\ (4)$, $p\equiv -1~\ (3)$, d.h. f\"ur $p\equiv -1~\ (12)$. F\"ur $p\equiv \pm 5~\ (12)$ ist dagegen $(3|p)=-1$, und $x^2\equiv 3~\ (p)$ besitzt keine L\"osung.\\
		
		\item F\"ur welche Primzahlen $p>3$ ist $-3$ ein Q-Rest, f\"ur welche ein Q-Nichtrest?\\
		{\bf L\"osung:}~Aus Zusammenfassung~\ref{zus:7_10} 2) und 3) folgt $$(-3|p)=(-1|p)\cdot(3|p)=(-1)^{\frac{p-1}{2}}\cdot (3|p),$$
		und weiter nach (a):
		$$(-3|p)=(p|3).$$
		Somit ist $x^2\equiv -3~\ (p)$ f\"ur $p\equiv 1~\ (3)$ l\"osbar ($p>3$ vorausgesetzt). F\"ur diese $p$ ist $-3$ ein Q-Rest, dagegen ist $-3$ ein Q-Nichtrest f\"ur $p\equiv -1~\ (3)$.\\
		
		\item Tabellen: Im Folgenden ist $p\geq 3$ eine Primzahl
		\begin{enumerate}
			\item [1.1)] $a=1$ ist Q-Rest f\"ur $p$, und $x\equiv \pm 1~\ (p)$ die L\"osungen von $x^2\equiv  1~\ (p)$.\\
			
			\item [1.2)] $a=-1$ ist Q-Rest f\"ur $p$ $\Leftrightarrow$ $p\equiv 1~\ (4)$.
			\begin{center}
				\begin{tabular}{|l||c|c|c|c|c|c|c|c|c|} \hline
					~$p\equiv 1~\ (4)$ ~  & ~~ $5$ ~~  & ~~ $13$ ~~ & ~~ $17$ ~~& ~~ $29$ ~~& ~~ $37$\\
					\hline
					~L\"osungen $x$ von~  & $x\equiv \pm 2$ & $x\equiv \pm 5$ & $x\equiv \pm 4$ &$x\equiv \pm 12$  & $x\equiv \pm 6$\\  %\hline
					~$x^2\equiv  -1~\ (p)$~  &  $\mod 5$~  & $\mod 13$~ &  $\mod 17$~&$\!\mod 29$~  & \!$\mod 37$ \\\cline{2-6}
					\rule{0pt}{2ex} 	& $41$ & $53$ & $61$ & $73$ & $89$\\\cline{2-6}
					\rule{0pt}{2ex}&$x\equiv \pm 9$  & $x\equiv \pm 23$ & $x\equiv \pm 11$ & $x\equiv \pm 27$ & $x\equiv \pm 34$ \\  %\hline
					%					~von~ $x^2\equiv  -1~\ (p)$~  
					&$\mod 41$~  & $\mod 53$~  & $\mod 61$~& $\mod 73$~ & $\mod 89$~\\ \hline
				\end{tabular}\\
			\end{center}
			
			\item [2.1)] $a=2$ ist Q-Rest f\"ur $p$ $\Leftrightarrow$ $p\equiv \pm 1~\ (8)$.
			\begin{center}
				\begin{tabular}{|l||c|c|c|c|c|c|c|c|c|} \hline
					~$p\equiv \pm 1~\ (8)$ ~  & ~~ $7$ ~~  & ~~ $17$ ~~ & ~~ $23$ ~~& ~~ $31$ ~~& ~~ $41$\\
					\hline
					~L\"osungen $x$ von~  & $x\equiv \pm 3$ & $x\equiv \pm 6$ & $x\equiv \pm 5$ &$x\equiv \pm 8$  & $x\equiv \pm 17$\\  
					~$x^2\equiv  2~\ (p)$~  &  $\!\mod 7$~  & $\!\mod 17$~ &  $\!\mod 23$~&$\!\mod 31$~  & \!$\mod 41$ \\\cline{2-6}
					\rule{0pt}{2ex} 	& $47$ & $71$ & $73$ & $79$ & $89$\\\cline{2-6}
					\rule{0pt}{2ex}&$x\equiv \pm 7$  & $x\equiv \pm 12$ & $x\equiv \pm 32$ & $x\equiv \pm 9$ & $x\equiv \pm 25$ \\  
					&$\mod 47$~  & $\mod 71$~  & $\mod 73$~& $\mod 79$~ & $\mod 89$~\\ \hline
				\end{tabular}\\
			\end{center}
			
			\item [2.2)] $a=-2$ ist Q-Rest f\"ur $p$ $\Leftrightarrow$ $p\equiv 1,3~\ (8)$.
			\begin{center}
				\begin{tabular}{|l||c|c|c|c|c|c|c|c|c|} \hline
					~$p\equiv 1,3~\ (8)$ ~  & ~~ $3$ ~~  & ~~ $11$ ~~ & ~~ $17$ ~~& ~~ $19$ ~~& ~~ $41$\\
					\hline
					~L\"osungen $x$ von~  & $x\equiv \pm 1$ & $x\equiv \pm 3$ & $x\equiv \pm 7$ &$x\equiv \pm 6$  & $x\equiv \pm 11$\\  
					~$x^2\equiv  -2~\ (p)$~  &  $\!\mod 3$~  & $\!\mod 11$~ &  $\!\mod 17$~&$\!\mod 19$~  & \!$\mod 41$ \\\cline{2-6}
					\rule{0pt}{2ex} 	& $43$ & $59$ & $67$ & $73$ & $83$\\\cline{2-6}
					\rule{0pt}{2ex}&$x\equiv \pm 16$  & $x\equiv \pm 23$ & $x\equiv \pm 20$ & $x\equiv \pm 12$ & $x\equiv \pm 9$ \\  
					&$\mod 43$~  & $\mod 59$~  & $\mod 67$~& $\mod 73$~ & $\mod 83$~\\ \hline
				\end{tabular}\\
			\end{center}
			
			\item [3.1)] $a=3$ ist Q-Rest f\"ur $p$ $\Leftrightarrow$ $p\equiv \pm 1~\ (12)$.
			\begin{center}
				\begin{tabular}{|l||c|c|c|c|c|c|c|c|c|} \hline
					~$p\equiv \pm 1~\ (12)$ ~  & ~~ $11$ ~~  & ~~ $13$ ~~ & ~~ $23$ ~~& ~~ $37$ ~~& ~~ $47$\\
					\hline
					~L\"osungen $x$ von~  & $x\equiv \pm 5$ & $x\equiv \pm 4$ & $x\equiv \pm 7$ &$x\equiv \pm 15$  & $x\equiv \pm 12$\\  
					~$x^2\equiv  3~\ (p)$~  &  $\!\mod 11$~  & $\!\mod 13$~ &  $\!\mod 23$~&$\!\mod 37$~  & \!$\mod 47$ \\\cline{2-6}
					\rule{0pt}{2ex} 	& $59$ & $61$ & $71$ & $73$ & $83$\\\cline{2-6}
					\rule{0pt}{2ex}&$x\equiv \pm 11$  & $x\equiv \pm 8$ & $x\equiv \pm 28$ & $x\equiv \pm 21$ & $x\equiv \pm 13$ \\ 
					&$\mod 59$~  & $\mod 61$~  & $\mod 71$~& $\mod 73$~ & $\mod 83$~\\ \hline
				\end{tabular}\\
			\end{center}			
			
			\item [3.2)] $a=-3$ ist Q-Rest f\"ur $p$ $\Leftrightarrow$ $p\equiv 1~\ (3)$.
		
			\begin{center}
					\begin{tabular}{|l||c|c|c|c|c|c|c|c|c|} \hline
						~$p\equiv 1~\ (3)$ ~  & ~~ $7$ ~~  & ~~ $13$ ~~ & ~~ $19$ ~~& ~~ $31$ ~~& ~~ $37$\\
						\hline
						~L\"osungen $x$ von~  & $x\equiv \pm 2$ & $x\equiv \pm 6$ & $x\equiv \pm 4$ &$x\equiv \pm 11$  & $x\equiv \pm 16$\\ 
						~$x^2\equiv  -3~\ (p)$~  &  $\!\mod 7$~  & $\!\mod 13$~ &  $\!\mod 19$~&$\!\mod 31$~  & \!$\mod 37$ \\\cline{2-6}
						\rule{0pt}{2ex} 	& $43$ & $61$ & $67$ & $73$ & $79$\\\cline{2-6}
						\rule{0pt}{2ex}&$x\equiv \pm 13$  & $x\equiv \pm 27$ & $x\equiv \pm 8$ & $x\equiv \pm 17$ & $x\equiv \pm 32$ \\  
						&$\mod 43$~  & $\mod 61$~  & $\mod 67$~& $\mod 73$~ & $\mod 79$~\\ \hline
					\end{tabular}\\
				\end{center}		
		\end{enumerate}
	\end{enumerate}
	\dokendSatz
\end{Beis}

Die Verwendung des Legendre-Symbols und des quadratischen Reziprozit\"atsgesetzes 
hat  in der vorliegenden Form den Nachteil, dass auf die Primfaktorzerlegung 
zur\"uckgegriffen werden muss. Dies l\"asst sich durch Verwendung
des sogenannten Jacobi-Symbols $(P|Q)$ vermeiden, einer sinnvollen Erweiterung des Legendre-Symbols.
Wir orientieren uns an \cite[Band 46,\S 3.3]{nz}:\\ Es seien $P,Q$ teilerfremd, $Q=q_1q_2 \ldots q_s>0$ sei ungerade
und das Produkt nicht notwendig voneinander verschiedener Primzahlen $q_j$. 
Dann ist durch das Jacobi-Symbol\index{Jacobi-Symbol}\label{Jacobi-Symbol}
$\begin{displaystyle}(P|Q)=\prod \limits_{j=1}^{s} (P|q_j)\end{displaystyle}$ eine Verallgemeinerung
des Legendre-Symbols gegeben, dass der Beziehung $(P|Q)=(P'|Q)$ f\"ur $P \equiv P'(Q)$
gen\"ugt. Nun lassen sich f\"ur je zwei teilerfremde und ungerade nat\"urliche Zahlen
$P, Q$ sowohl das quadratische Reziprozit\"atsgesetz
$$
(P|Q) \cdot (Q|P) = (-1)^{\frac{P-1}{2}\frac{Q-1}{2}}
$$
als auch die erg\"anzenden Beziehungen
$$
(-1|Q)=(-1)^{\frac{Q-1}{2}}\,,\quad (2|Q)=(-1)^{\frac{Q^2-1}{8}}
$$
ohne gr\"o{\ss}ere M\"uhe auf das Jacobi-Symbol \"ubertragen. Damit l\"asst sich insbesondere das
Legendre-Symbol, das Aufschluss \"uber das quadratische Restverhalten einer Zahl bzgl. eines Primzahlmoduls gibt,
sehr effizient mit Hilfe eines dem Euklidischen Algorithmus \"ahnendeln Verfahrens berechnen, 
siehe \cite[Band 46,\S 3.3]{nz}.

	\section{Aufgaben}\label{cha:7_A}

\begin{Auf}[L\"osungen quadratischer Kongruenzen\index{quadratischer Kongruenz}\label{quadratischer Kongruenz2}]\label{auf:7_1}
	\hspace*{0cm}\\\vspace{-1cm}
	\begin{enumerate}[(a)]
		\item Man bestimme die Anzahl der L\"osungen zur Kongruenz $x^2 \equiv 1 ~\,(360)$.\\
		
		\item F\"ur eine Primzahl $p \equiv -1 ~\,(4)$ sei $a \in \Z$ nicht durch $p$ teilbar 
		und quadratischer Rest mod $p$. Man zeige, dass dann 
		$x^2 \equiv a ~\,(p)$ genau die folgenden beiden L\"osungen besitzt:
		$$ x_{1/2} \equiv \pm \, a^\frac{p+1}{4} ~\,(p)\,.$$
		
		\item Es sei $p \equiv 5\,~(8)$ eine Primzahl. Man zeige:
		Es ist $2$ quadratischer Nichtrest von $p$, aber $-1$ Quadratrest, und die Kongruenz 
		$x^2 \equiv -1 ~\,(p)$ besitzt genau die beiden L\"osungen
		$$ x_{1/2} \equiv \pm \, 2^\frac{p-1}{4} ~\,(p)\,.$$
		
		\item Man bestimme alle L\"osungen der Kongruenz
		$x^2 \equiv -1 ~\,(65)\,.$\\
		
		\item Man bestimme jeweils alle Primzahlen $p \geq 3$,
		f\"ur die $-5$, $5$, $-6$ bzw. $6$ quadratische Reste sind.
	\end{enumerate}
\end{Auf}
{\bf L\"osung:}\\

\begin{enumerate}[(a)]
	\item $x^2\equiv 1 ~\, (360)$ hat wegen $360=2^3\cdot 3^2 \cdot 5$ genau $$\min(4,2^{\max(3,1)-1})\cdot 2^2=4\cdot 2^2=16$$ L\"osungen, siehe Satz~\ref{satz:7_2}.\\
	
	\item F\"ur $p\equiv-1 ~\, (4)$ sei $(a|p)=1$. Nach dem Eulerschen Kriterium ist $a^{\frac{p-1}{2}}\equiv1 ~\, (p)$, und somit gilt $\left(\pm a^{\frac{p+1}{4}}\right)^2\equiv a^{\frac{p+1}{2}}\equiv a\cdot a^{\frac{p-1}{2}}\equiv a ~\, (p)$. 
	Es sind $x_{1,2} \equiv \pm a^{\frac{p+1}{4}} ~\, (p)$ die einzigen L\"osungen von $x ^2\equiv a ~\, (p)$, da $p \geq 3$ Primzahl ist.\\
	
	\item Es sei $p\equiv5 ~\,(8)$ eine Primzahl. Dann ist $(-1)^{\frac{p-1}{2}}=1$, und somit $-1$ Q-Rest f\"ur $p$ nach dem Eulerschen Kriterium. Nach Folgerung~\ref{folg:7_8}~(a) ist $(2|p)=-1$, also auch $2^\frac{p-1}{2}\equiv -1 ~\, (p)$. Wir erhalten 
	\begin{equation*}
	\left(\pm 2^{\frac{p-1}{4}}\right)^{2}=2^\frac{p-1}{2}\equiv -1 ~\, (p),
	\end{equation*}
	und da die Kongruenz $x^2\equiv -1 ~\, (p)\mod p$ nicht mehr als zwei L\"osungen besitzten kann, sind $x_{1,2}\equiv\pm 2^\frac{p-1}{4} ~\, (p)$ alle L\"osungen von $x^2\equiv -1 ~\, (p)$.\\
	
	\item Es ist $65=5\cdot 13$ mit den beiden Primzahlen $5$ und $13\equiv 5 ~\,(8)$. Nach (c) sind $y_1=2$, $y_2=-2$ L\"osungen von $y^2\equiv -1 ~\, (5)$ sowie $y_3=5\equiv -8~\,(13)$ und $y_4=-5\equiv 8~\,(13)$ L\"osungen von $y^2\equiv-1~\,(13)$. Nach Satz~\ref{satz:7_2} hat die quadratische Kongruenz $x^2\equiv-1~\,(65)$ genau vier L\"osungen $x_1$, $x_2$, $x_3$, $x_4$, die wir aus den folgenden vier simultanen Kongruenzen erhalten:
	\begin{enumerate}[1)]
		\item $x_1\equiv2~\,(5)$, $x_1\equiv5~\,(13)$ liefert $x_1\equiv-8~\,(65)$,\\
		
		\item $x_2\equiv2~\,(5)$, $x_2\equiv-5~\,(13)$ liefert $x_2\equiv-18~\,(65)$,\\
		
		\item $x_3\equiv-2~\,(5)$, $x_3\equiv5~\,(13)$ liefert $x_3\equiv18~\,(65)$,\\
		
		\item $x_4\equiv-2~\,(5)$, $x_4\equiv-5~\,(13)$ liefert $x_4\equiv8~\,(65)$.
	\end{enumerate}
	
	\item Nach dem Eulerschen Kriterium und dem quadratischen Reziprozit\"atsgesetz gilt f\"ur Primzahl $p\geq 3$:
	\begin{equation*}
	(-5|p)=(-1)^\frac{p-1}{2}\cdot (p|5),
	\end{equation*}
	und somit gilt $(-5|p)=1$ genau dann, wenn
	\begin{equation*}
	p\equiv 1~\,(4)\wedge p\equiv\pm 1~\,(5)\quad \text{oder aber}\quad
	p\equiv -1~\,(4)\wedge p\equiv\pm 2~\,(5)\quad
	\end{equation*}
    gilt, siehe Folgerung~\ref{folg:7_5}. Wir erhalten $-5$ als Q-Rest f\"ur 
	\begin{equation*}
	\begin{array}{clclcc}
	&\quad p\equiv 1~\,(4) & \wedge & p\equiv 1~\,(5), & \quad\text{d.h.} & p\equiv 1~\,(20), \\
\text{oder} & \quad p\equiv 1~\,(4) & \wedge & p\equiv -1~\,(5), & \quad\text{d.h.} & p\equiv 9~\,(20), \\
\text{oder} &\quad p\equiv -1~\,(4) & \wedge & p\equiv 2~\,(5), &\quad \text{d.h.} & p\equiv 7~\,(20), \\
\text{oder} &\quad p\equiv -1~\,(4) & \wedge & p\equiv -2~\,(5), &\quad \text{d.h.} & p\equiv 3~\,(20). 
	\end{array}
	\end{equation*}
	Zusammengefasst gilt $(-5|p)=1$ genau dann, wenn $p\equiv 1,3,7,9~\,(20)$. Es ist
	\begin{equation*}
	(5|p)=(p|5)=1\Leftrightarrow p\equiv\pm1~\,(5)\wedge p\equiv 1~\,(2)\Leftrightarrow  p\equiv\pm1~\,(10).
	\end{equation*}
	F\"ur die Q-Reste $\pm 6$ beachten wir $(\pm6|p)=(2|p)\cdot(\pm3|p)$, siehe Folgerung~\ref{folg:7_5}, also gilt $(\pm6|p)=1$ genau f\"ur 
	\begin{equation*}
	(2|p)=(-1)^{\frac{1}{8}(p-1)(p+1)}=1\wedge (\pm3|p)=1
	\end{equation*}
	oder
	\begin{equation*}
	(2|p)=(-1)^{\frac{1}{8}(p-1)(p+1)}=-1\wedge (\pm3|p)=-1.
	\end{equation*}
	\underline{\textbf{Q-Rest $-6$}}: Wir verwenden Folgerung~\ref{folg:7_8} und Beispiel~\ref{beis:7_11}~(b).
		\begin{enumerate}[1)]
			\item $p\equiv1~\,(8)\wedge p\equiv1~\,(3)$ liefert $p\equiv1~\,(24)$,\\
			
			\item $p\equiv-1~\,(8)\wedge p\equiv1~\,(3)$ liefert $p\equiv7~\,(24)$,\\
			
			\item $p\equiv3~\,(8)\wedge p\equiv-1~\,(3)$ liefert $p\equiv11~\,(24)$,\\
			
			\item $p\equiv-3~\,(8)\wedge p\equiv-1~\,(3)$ liefert $p\equiv5~\,(24)$.
		\end{enumerate}
	\underline{\textbf{Q-Rest $6$}}: Wir verwenden Folgerung~\ref{folg:7_8} und Beispiel~\ref{beis:7_11}~(a).
		\begin{enumerate}[1)]
			\item $p\equiv1~\,(8)\wedge p\equiv1~\,(12)$ liefert $p\equiv1~\,(24)$,\\
			
			\item $p\equiv-1~\,(8)\wedge p\equiv -1~\,(12)$ liefert $p\equiv-1~\,(24)$,\\
			
			\item $p\equiv3~\,(8)\wedge p\equiv-5~\,(12)$ liefert $p\equiv-5~\,(24)$,\\
			
			\item $p\equiv-3~\,(8)\wedge p\equiv5~\,(12)$ liefert $p\equiv5~\,(24)$.
		\end{enumerate}
		
		{\bf Zusammenfassung:} 
		
		F\"ur jede Primzahl $p\geq 3$ gilt:
		\begin{enumerate}[(i)]
			\item $(-5|p)=1 \Leftrightarrow p\equiv1,3,7,9~\,(20)$.
				\item  $(5|p)=1 \Leftrightarrow p\equiv\pm1~\,(10)$.	
				\item $(-6|p)=1 \Leftrightarrow p\equiv1,5,7,11~\,(24)$.
				\item $(6|p)=1 \Leftrightarrow p\equiv\pm1,\pm5~\,(24)$.
		\end{enumerate}

		Zu jedem dieser vier F\"alle pr\"asentieren wir in den folgenden Tabellen
		jeweils die ersten zehn Primzahlen mit den entsprechenden L\"osungen  der
		quadratischen Kongruenzen:\\
		
		\begin{center}
			\begin{tabular}{|l||c|c|c|c|c|c|c|c|c|} \hline
				~$p\equiv 1,3,7, 9~\ (20)$ ~~  & ~~ $3$ ~~  & ~~ $7$ ~~ & ~~ $23$ ~~& ~~ $29$ ~~& ~~ $41$\\
				\hline
				~L\"osungen $x$ von~  & $x\equiv \pm 1$ & $x\equiv \pm 3$ & $x\equiv \pm 8$ &$x\equiv \pm 13$  & $x\equiv \pm 6$\\  %\hline
				~$x^2\equiv  -5~\ (p)$~  &  $\!\mod 3$~  & $\!\mod 7$~ &  $\!\mod 23$~&$\!\mod 29$~  & \!$\mod 41$ \\\cline{2-6}
				\rule{0pt}{2ex} 	& $43$ & $47$ & $61$ & $67$ & $83$\\\cline{2-6}
				\rule{0pt}{2ex}&$x\equiv \pm 9$  & $x\equiv \pm 18$ & $x\equiv \pm 19$ & $x\equiv \pm 14$ & $x\equiv \pm 24$ \\  
				&$\mod 43$~  & $\mod 47$~  & $\mod 61$~& $\mod 67$~ & $\mod 83$~\\ \hline
			\end{tabular}\\
		\end{center}

		\begin{center}
		\begin{tabular}{|l||c|c|c|c|c|c|c|c|c|} \hline
		~$p\equiv \pm 1~\ (10)$ ~~~~~~~~~~  & ~~ $11$ ~~  & ~~ $19$ ~~ & ~~ $29$ ~~& ~~ $31$ ~~& ~~ $41$\\
			\hline
			~L\"osungen $x$ von~  & $x\equiv \pm 4$ & $x\equiv \pm 9$ & $x\equiv \pm 11$ &$x\equiv \pm 6$  & $x\equiv \pm 13$\\  
			~$x^2\equiv  5~\ (p)$~  &  $\!\mod 11$~  & $\!\mod 19$~ &  $\!\mod 29$~&$\!\mod 31$~  & \!$\mod 41$ \\\cline{2-6}
			\rule{0pt}{2ex} 	& $59$ & $61$ & $71$ & $79$ & $89$\\\cline{2-6}
			\rule{0pt}{2ex}&$x\equiv \pm 8$  & $x\equiv \pm 26$ & $x\equiv \pm 17$ & $x\equiv \pm 20$ & $x\equiv \pm 19$ \\ 
			&$\mod 59$~  & $\mod 61$~  & $\mod 71$~& $\mod 79$~ & $\mod 89$~\\ \hline
		\end{tabular}\\
	\end{center}
			
	\begin{center}
		\begin{tabular}{|l||c|c|c|c|c|c|c|c|c|} \hline
			~$p\equiv 1,5,7,11~\,(24)$ ~  & ~~ $5$ ~~  & ~~ $7$ ~~ & ~~ $11$ ~~& ~~ $29$ ~~& ~~ $31$\\
			\hline
			~L\"osungen $x$ von~  & $x\equiv \pm 2$ & $x\equiv \pm 1$ & $x\equiv \pm 4$ &$x\equiv \pm 9$  & $x\equiv \pm 5$\\ 
			~$x^2\equiv  -6~\ (p)$~  &  $\!\mod 5$~  & $\!\mod 7$~ &  $\!\mod 11$~&$\!\mod 29$~  & \!$\mod 31$ \\\cline{2-6}
			\rule{0pt}{2ex} 	& $53$ & $59$ & $73$ & $79$ & $83$\\\cline{2-6}
			\rule{0pt}{2ex}&$x\equiv \pm 10$  & $x\equiv \pm 17$ & $x\equiv \pm 33$ & $x\equiv \pm 28$ & $x\equiv \pm 34$ \\
			&$\mod 53$~  & $\mod 59$~  & $\mod 73$~& $\mod 79$~ & $\mod 83$~\\ \hline
		\end{tabular}\\
	\end{center}

		\begin{center}
			\begin{tabular}{|l||c|c|c|c|c|c|c|c|c|} \hline
					~$p\equiv \pm1,\pm5~\,(24)$ ~~~~  & ~~ $5$ ~~  & ~~ $19$ ~~ & ~~ $23$ ~~& ~~ $29$ ~~& ~~ $43$\\
					\hline
					~L\"osungen $x$ von~  & $x\equiv \pm 1$ & $x\equiv \pm 5$ & $x\equiv \pm 11$ &$x\equiv \pm 8$  & $x\equiv \pm 7$\\ 
					~$x^2\equiv  6~\ (p)$~  &  $\!\mod 5$~  & $\!\mod 19$~ &  $\!\mod 23$~&$\!\mod 29$~  & \!$\mod 43$ \\\cline{2-6}
					\rule{0pt}{2ex} 	& $47$ & $53$ & $67$ & $71$ & $73$\\\cline{2-6}
					\rule{0pt}{2ex}&$x\equiv \pm 10$  & $x\equiv \pm 18$ & $x\equiv \pm 26$ & $x\equiv \pm 19$ & $x\equiv \pm 15$ \\
					&$\mod 47$~  & $\mod 53$~  & $\mod 67$~& $\mod 71$~ & $\mod 73$~\\ \hline
				\end{tabular}\\
			\end{center}
\end{enumerate}

\begin{Auf}[Quadratische Reste\index{quadratische Reste}\label{quadratische Reste3} Fermatscher Primzahlen\index{Fermatsche Primzahl}\label{Fermatsche Primzahl}]\label{auf:7_2}
\hspace*{0cm}\\\vspace{-1cm}	
	\begin{enumerate}[(a)]
		\item F\"ur $m \in \N$  sei $p=2^m+1$ eine Primzahl. Man zeige, dass dann
		$m$ eine Potenz von 2 sein muss, d.h. es muss $m=2^n$ f\"ur ein $n \in \N_0$ gelten.\\
		
		\item Es sei $p$ eine Primzahl der Gestalt $p = 2^m+1$ mit $m \in \N$ 
		und $a$ eine nicht durch $p$ teilbare ganze Zahl. Man zeige, dass 
		$a$ genau dann quadratischer Rest mod $p$ ist, wenn $a$ keine Primitivwurzel mod $p$
		ist.
		
	\end{enumerate}
	{\bf Bemerkung:} Die Primzahlen der Gestalt $p = 2^{2^n}+1$ mit $n \in \N_0$
	werden auch Fermatsche Primzahlen genannt. Nach Gau{\ss} ist
	f\"ur diese Primzahlen $p$ das regelm\"assige $p$-Eck 
	allein mit Zirkel und Lineal konstruierbar.
\end{Auf}
{\bf L\"osung:}\\

\begin{enumerate}[(a)]
	\item Es sei $p=2^m+1$ mit $m\in\N$ Primzahl. Dann wird $p$ auch Fermatsche Primzahl genannt. Ist $\lambda >1$ eine ungerade nat\"urliche Zahl und $k\in\N$ beliebig, so ist 
	\begin{equation*}
	2^{\lambda k}+1=(2^k+1)\cdot\sum\limits_{j=0}^{\lambda-1}(-1)^{j} 2^{jk}
	\end{equation*}
	wegen $1<2^k+1<2^{\lambda k}+1$ eine nichttriviale Zerlegung von $2^{\lambda k}+1$, so dass $m$ f\"ur die Fermatsche Primzahl $p$ eine Potenz von $2$ sein mu{\ss}, d.h. $m=2^n$ f\"ur ein $n\in\N_0$.\\
	
	\item Es sei $p=2^m+1$ eine Fermatsche Primzahl und $a$ eine Primitivwurzel $\Mod p$. Wir haben $\varphi(p)=p-1=2^m$, und die Kongruenz $a^{\frac{p-1}{2}}\equiv1 ~\,(p)$ kann nicht erf\"ullt sein. Nach Satz~\ref{satz:7_4} ist $a$ Q-Nichtrest $\Mod p$. Es gibt genau $\frac{p-1}{2}=2^{m-1}$ Q-Nichtreste $\Mod p$, siehe Folgerung~\ref{folg:7_5}, und genau $\varphi(\varphi(p))=\varphi(2^m)=2^{m-1}$ Primitivwurzeln $\Mod p$. Damit sind die Q-Nichtreste $\Mod p$ genau die Primitivwurzeln $\Mod p$.
\end{enumerate}

\begin{Auf}[Anwendung des quadratischen Reziprozit\"atsgesetzes zur Berechnung eines quadratischen Restsymboles]\label{auf:7_3}

Man berechne $(-6|101)$.
\end{Auf}
{\bf L\"osung:}
Da das Legendre-Symbol $(\cdot|p)$ f\"ur Primzahlen $p\geq 3$ vollst\"andig multiplikativ ist, errechnen wir zun\"achst mit $-6=(-1)\cdot 2\cdot 3$:
\begin{equation}\label{7_3:eq:1}
(-6|101)=(-1|101)\cdot(2|101)\cdot(3|101).
\end{equation}
Unter Beachtung von $(-1|p)=(-1)^\frac{p-1}{2}$ (Folgerung~\ref{folg:7_6}) und $(2|p)=(-1)^{\frac{1}{8}(p^2-1)}$ (Folgerung~\ref{folg:7_8}) erhalten wir hier f\"ur $p:=101$:
\begin{equation}\label{7_3:eq:2}
(-1|101)=1,\quad (2|101)=-1.
\end{equation}
Wir m\"ussen in \eqref{7_3:eq:1} nur noch $(3|101)$ berechnen:

Nach dem Reziprozit\"atsgesetz ist $(q|p)=(p|q)$, wenn eine der beiden Primzahlen $p\neq q$ modulo $4$ den Divisionsrest $1$ hat, und sonst $(q|p)=-(p|q)$. Hier ist $p=101\equiv1~\ (4)$, $q=3$, also mit $101\equiv-1~\ (3)$: 
\begin{equation}\label{7_3:eq:3}
(3|101)=(101|3)=(-1|3)=-1.
\end{equation}
Aus \eqref{7_3:eq:1} bis \eqref{7_3:eq:3} folgt
\begin{equation*}
(-6|101)=1\cdot(-1)\cdot(-1)=+1.
\end{equation*}
Es ist $-6$ quadratischer Rest $\Mod 101$.

	\chapter{Quadratische Formen}\label{cha:8}
	In diesem Kapitel betrachten wir bin\"are quadratische Formen, wobei wir uns auf die indefiniten Formen und deren Zusammenhang mit der Kettenbruchentwicklung reell quadratischer Irrationalzahlen konzentrieren. Dabei legen wir den Schwerpunkt auf die Entwicklung eines effizienten Reduktionsverfahrens f\"ur indefinite quadratische Formen,
	das Hand in Hand mit der Kettenbruchentwicklung quadratischer Irrationalzahlen geht, die wir diesen Formen
	zuordnen. Ein wichtiges Nebenprodukt dieser Untersuchungen ist, dass genau die reell quadratischen
	Irrationalzahlen eine Kettenbruchentwicklung besitzen, die in eine Periode m\"undet. Auch
	gehen wir mit Hilfe des hier entwickelten Verfahrens erste Schritte zur Beantwortung der Frage,
	wie sich die \"aquivalenten quadratischen Formen ineinander transformieren lassen.
	Weiterf\"uhrende Themen werden aufgrund ihres Umfangs, wenn \"uberhaupt, nur angerissen.
	Als Begleitlekt\"ure f\"ur die tiefergehenden Studien zu quadratischen Zahlk\"orpern
	empfehlen wir das Lehrb\"ucher von Niven und Zuckerman \cite[Band 47, \S 7.7-7.9]{nz} 
	bzw. \cite{nzm} sowie von Halter-Koch \cite{koch} und Steuding \cite{st}. 
	
        \section{Quadratische Formen und reell quadratische Irrationalzahlen}\label{cha:8A}
F\"ur gegebene Koeffizienten $a,b,c\in\Z$ nennen wir
\begin{equation}\label{eq:8_1}
D=b^2-4ac
\end{equation}
die Diskriminante\index{Diskriminante}\label{Diskriminante} einer nicht verschwindenden quadratischen Form\index{quadratische Form}	\label{quadratische Form2}
\begin{equation}\label{eq:8_2}
F(x,y)=ax^2+bxy+cy^2.
\end{equation}	
Es bestehen die Darstellungen
\begin{equation} \label{eq:8_3}
\left.\begin{tabular}{l}
	$4a F(x,y)=(2ax+by)^2-Dy^2$,\\
	$4c F(x,y)=(2cy+bx)^2-Dx^2$.	 
	\end{tabular}\right\}
\end{equation}

Die quadratische Form $F$ hei{\ss}t definit\index{definit quadratische Form}\label{definit quadratische Form}, wenn $D<0$ ist.
Ein Blick auf die Diskriminante in (\ref{eq:8_1}) lehrt, dass dies nur m\"oglich ist, wenn sowohl $a$ als auch $c$ von Null verschieden sind und zudem 
entweder beide positiv oder beide negativ sind.\\
	
Ist bei $D<0$ noch $a>0$, $c>0$, so nennen wir $F$ positiv definit\index{positiv definit quadratische Form}\label{positiv definit quadratische Form},
weil dann nach~(\ref{eq:8_3}) die Form $F$ nur positive Werte annimmt, wenn man in~(\ref{eq:8_2})
f\"ur $x,y\in \Z$, abgesehen von $x=y=0$, beliebige Werte einsetzt. 
Entsprechend hei{\ss}t $F$ f\"ur $D<0$ und $a<0$,~$c<0$ negativ definit.\\

Nun betrachten wir den Fall $D=f^2 \geq 0$ mit einem $f\in \N_0$. Dann folgt aus 
der ersten Gleichung von~(\ref{eq:8_3}):
\begin{equation} \label{eq:8_4}
a F(x,y)=\left(ax+\dfrac{b-f}{2} y\right) \left(ax+\dfrac{b+f}{2} y\right) 
\end{equation}
mit $b\equiv f$ $(2)$, so dass die Gr\"o{\ss}en $\dfrac{b\pm f}{2}$ ganzzahlig sind. 
F\"ur $a=0$ erhalten 	wir eine Zerlegung $F(x,y)=(bx+cy)\cdot y$. Nun setzen wir $a\neq 0$ voraus.
Dann ist $a=\lambda \cdot\lambda'$ mit den ganzzahligen Gr\"o{\ss}en 
	$$
	\lambda=\ggT\left(a,\frac{b-f}{2}\right), \quad \lambda'=\frac{a}{\lambda}.
	$$	

Die ganzen Zahlen $\dfrac{a}{\lambda}$, $\dfrac{b-f}{2\lambda}$ sind teilerfremd, und es gilt
$$
ac=\frac{b-f}{2}\cdot \frac{b+f}{2}, \quad 
\frac{a}{\lambda}c \left| \frac{b-f}{2\lambda}\cdot\frac{b+f}{2}\right.,
$$
folglich auch $\lambda' \left| \dfrac{b+f}{2}\right.$ neben $\lambda' \left|\  a\right.$.\\

Wir erhalten aus (\ref{eq:8_4}), dass auch im Fall $a\neq 0$ die Form
$$
F(x,y)=\left(\lambda'x+\frac{b-f}{2\lambda}y\right)\cdot\left(\lambda x+\frac{b+f}{2\lambda'}y\right)
$$
in das Produkt zweier Linearformen mit ganzzahligen Koeffizienten zerf\"allt.\\

Nun gelte umgekehrt $F(x,y)=\left(\lambda'x+\mu'y\right)\left(\lambda x+\mu y\right)$
mit irgendwelchen Zahlen $\lambda, \lambda', \mu, \mu'\in\Z$.
Dann errechnet man f\"ur diese Form die Diskriminante $D=f^2$ 
mit $f:=\left|\lambda'\mu-\lambda\mu'\right|\in\N_0$.
Diesen Fall schlie{\ss}en wir im Folgenden aus.\\

Jetzt mu{\ss} noch f\"ur die allgemeine Form $F$ in (\ref{eq:8_2}) 
der verbleibende Fall \linebreak
$D=b^2-4ac>0$ betrachtet werden, wobei $D$ keine Quadratzahl ist.
In diesem Falle nennen wir $F$ eine indefinite Form \index{indefinite quadratische Form}\label{indefinite quadratische Form}.
Dann ist $a\neq 0$ und $c\neq 0$. Bei indefiniten Formen werden wir im Folgenden stets 
stillschweigend voraussetzen, dass $D$ keine Quadratzahl ist.\\

Wir schicken eine Definition voraus, die f\"ur Formen mit beliebiger Diskriminante gilt:
\begin{Def}% Definition 8.1
\label{def:8_1}
	Die quadratische Form $F$ in (\ref{eq:8_2})  hei{\ss}t primitiv\index{primitive quadratische Form}\label{primitive quadratische Form}, 
	wenn $a,b,c$ teilerfremd sind, d.h. der gr\"o{\ss}te gemeinsame Teiler 
	von $a,b,c$ hat den Wert $1$.
	\dokendDef
\end{Def}

{\it Bemerkung}:~
Die weiterreichende paarweise Teilerfremdheit von $a,b,c$ wird 
in Definition~\ref{def:8_1} nicht verlangt.

\begin{DefSatz}% Definition und Satz 8.2
 \label{defsatz:8_2}
	Der indefiniten quadratischen Form $F(x,y)=ax^2+bxy+cy^2$ mit Diskriminante $D=b^2-4ac$
	ordnen wir ihre reell quadratische Irrationalzahl\index{quadratische Irrationalzahl}\label{quadratische Irrationalzahl3} $X(F):=\dfrac{\sqrt{D}-b}{2a}$ zu.
	Ist dann $F$ primitiv, so ist $F$ durch $X(F)$ eindeutig bestimmt. 
	\dokendDef
\end{DefSatz}
{\bf Beweis:}~
Da $D>0$ keine Quadratzahl ist, wird $X(F)$ eine reell quadratische Irrationalzahl.
Wir setzen nun $F$ als primitiv voraus, und betrachten eine weitere 
primitive und indefinite Form
$$
F'(x,y)=a'x^2+b'xy+c'y^2
$$
mit Diskriminante $D'>0$, $D'$ ebenfalls keine Quadratzahl, so dass $X(F)=X(F')$ gilt.
Wir erhalten
\begin{equation}\label{eq:8_5}
 \frac{\sqrt{D'}-b'}{2 a'}=\frac{\sqrt{D}-b}{2 a}
\end{equation}
bzw.   $\, a'b-ab'=a'\sqrt{D}-a\sqrt{D'}$, und hieraus durch Quadrieren:
\begin{equation*}
(a'b-ab')^2=a'^2 D+a^2 D'-2aa'\sqrt{DD'}\\
\end{equation*}
sowie
\begin{equation*}
\sqrt{DD'}=\frac{1}{2aa'}\left[a'^2 D+a^2D'-(a'b-ab')^2\right].
\end{equation*}	

Somit gibt es ein rationales $q>0$ mit $\sqrt{D'}=q\sqrt{D}$, und es folgt aus (\ref{eq:8_5}):
$$
\frac{b'}{a'}-\frac{b}{a}=\left(\frac{q}{a'}-\frac{1}{a}\right)\sqrt{D}.
$$
Aus
$$
\frac{2a}{b-\sqrt{D}}=\frac{\sqrt{D}+b}{2c}=\frac{2a'}{b'-\sqrt{D'}}=\frac{\sqrt{D'}+b'}{2c'}
$$
folgt aber auch
$$
\frac{b'}{c'}-\frac{b}{c}=\left(\frac{1}{c}-\frac{q}{c'}\right)\sqrt{D},
$$
und da $\sqrt{D}$ irrational ist:
\begin{equation}\label{eq:8_6}
\frac{b'}{a'}=\frac{b}{a},
\end{equation}
\begin{equation}\label{eq:8_7} 
\frac{b'}{c'}=\frac{b}{c},
\end{equation}
\begin{equation}\label{eq:8_8}
q=\frac{a'}{a}=\frac{c'}{c}, \quad \frac{a'}{c'}=\frac{a}{c}.
\end{equation}

Wir erinnern an $q>0$ und schreiben $q=\dfrac{\alpha}{\alpha'}$ 
mit teilerfremden Zahlen $\alpha, \alpha'\in\N$. Wir erhalten aus 
(\ref{eq:8_6}), (\ref{eq:8_7}), (\ref{eq:8_8}):
$$
a'=qa, \quad b'=qb, \quad c'=qc,
$$
und schlie{\ss}lich $\alpha=\alpha\cdot \ggT(a,b,c)=\alpha'\cdot \ggT(a',b',c')=\alpha'$,
da $F$ und $F'$ primitiv sind. Somit ist $\alpha=\alpha'=q=1$ und $F=F'$.
\dokendProof

\begin{DefSatz}[Transformation der Formen\index{Transformation der Formen}\label{Transformation der Formen}]% Definition und Satz 8.3
 \label{defsatz:8_3}
	Gegeben sind eine Matrix 
	$A=\left( \begin{matrix}
		\alpha & \beta  \\
		\gamma & \delta\\
		\end{matrix} \right)\in GL(2,\Z)$, also
	   $\alpha,\beta,\gamma,\delta\in\Z$ und
		\linebreak 
%   	\parfillskip=1pt
		 $\alpha\delta-\beta\gamma=\pm 1$, sowie f\"ur $a,b,c\in\Z$ eine quadratische Form 
		$$F(x,y)=ax^2+bxy+cy^2$$
	mit nichtquadratischer Diskriminante.
		Hiermit definieren wir die transformierte Form $F^A$ gem\"a{\ss}
		$$
		F^A(x,y)=F(\alpha x+\beta y,\gamma x+\delta y),
		$$
		also $F^A=F'$ mit $F'(x,y)=a'x^2+b'xy+c'y^2$ und 
		\begin{equation*}
		 \begin{split}
		  & a'=a\alpha^2+b\alpha\gamma+c\gamma^2=F(\alpha,\gamma),\\
		  & b'=2a\alpha\beta+b(\alpha\delta+\beta\gamma)+2c\gamma\delta,\\
		  & c'=a\beta^2+b\beta\delta+c\delta^2=F(\beta,\delta).
		 \end{split}
		\end{equation*}
		 Dann haben $F$ und $F^A$ dieselbe Diskriminante, und es gilt 
		 $$\ggT(a,b,c)=\ggT(a',b',c').$$
	\dokendDef
\end{DefSatz}
{\bf Beweis:}~
Es gilt die Darstellung
$$
2 F(x,y)=\left( \begin{matrix}
		x & \\
		y & 
		\end{matrix} \right)^T 
	\left( 	\begin{matrix}
		2a & b  \\
		b & 2c
		\end{matrix} \right)
	\left( 	\begin{matrix}
		x & \\
		y & 
		\end{matrix} \right),
$$		
wobei die der Form $F$ zugeordnete Matrix $\left( \begin{matrix}
		2a & b  \\
		b & 2c
		\end{matrix} \right)$
		symmetrisch ist, und
$$
D=-\Det\left(  \begin{matrix}
		2a & b  \\
		b & 2c
		\end{matrix} \right)
 =b^2-4ac
$$
die Diskriminante von $F$ ist.
Nun gilt entsprechend 
\begin{equation}\label{eq:8_9}
  \begin{split}
    2F^A(x,y)&=\left(A\left(\begin{matrix}
		    x &  \\
		    y & 
		    \end{matrix} \right) \right)^T
		    \left( \begin{matrix}
		     2a & b  \\
		     b & 2c
		    \end{matrix} \right)
		    A\left(\begin{matrix}
		    x &  \\
		    y & 
		    \end{matrix} \right)\\
		    &=\left(\begin{matrix}
		    x &  \\
		    y & 
		    \end{matrix} \right)^T A^T
		    \left( \begin{matrix}
		     2a & b  \\
		     b & 2c
		    \end{matrix} \right)
		    A\left(\begin{matrix}
		    x &  \\
		    y & 
		    \end{matrix} \right),
  \end{split}
\end{equation}
d. h. der Form $F^A$ ist die symmetrische Matrix $\displaystyle A^T
		    \left( \begin{matrix}
		     2a & b  \\
		     b & 2c
		    \end{matrix} \right)
		    A$	zugeordnet.

Wir bezeichnen die Diskriminante von $F'$ mit $D'$, und erhalten
\begin{equation*}
 \begin{split}
  D'=&-\Det\left(A^T \left( \begin{matrix}
		     2a & b  \\
		     b & 2c
		    \end{matrix} \right) A\right)\\
    =&	-(\pm1)^2 \Det 	  \left( \begin{matrix}
		     2a & b  \\
		     b & 2c
		    \end{matrix} \right)\\
    =& b^2-4ac=D.
 \end{split}
\end{equation*}

Die Berechnung der Matrix-Eintr\"age von $A^T \left( \begin{matrix}
		     2a & b  \\
		     b & 2c
		    \end{matrix} \right) A$
liefert ferner die angegebenen Formeln f\"ur $a',b',c'$, aus denen sofort 
$\ggT(a,b,c)\left|\ggT(a',b',c')\right.$ folgt. Aus $F=F'^{A^{-1}}$ mit $A^{-1}\in GL(2,\Z)$ folgt
umgekehrt $\ggT(a',b',c')\left|\ggT(a,b,c)\right.$, und wir erhalten $\ggT(a,b,c)=\ggT(a',b',c')$.
\dokendProof\\

Ist $F(x,y)=ax^2+bxy+cy^2$ mit $a,b,c\in \Z$ eine quadratische Form, dann gilt
\begin{equation}\label{eq:8_10}
\left(F^A\right)^B=F^{A\cdot B} \quad \forall A,B \in GL(2,\Z),
\end{equation}
denn gem\"a{\ss} (\ref{eq:8_9}) ist den beiden Formen $\left(F^A\right)^B$, $F^{A\cdot B}$
die Matrix 
$$
M=(AB)^T\left( \begin{matrix}
		     2a & b  \\
		     b & 2c
		    \end{matrix} \right) AB
$$ mit
$$
2\left(F^A\right)^B (x,y)=2F^{A\cdot B}(x,y)
			 =\left(\begin{matrix}
				 x &  \\
				 y & 
				 \end{matrix} \right)^T M \left(\begin{matrix}
				 x &  \\
				 y & 
				 \end{matrix} \right)
$$
zugeordnet. 
Wir erinnern an  
$$
SL(2,\Z)=\left\{\left( \begin{matrix}
		     \alpha & \beta  \\
		     \gamma & \delta
		    \end{matrix}\right):\alpha,\beta,\gamma,\delta\in\Z, ~\
		    \alpha\delta-\beta\gamma=1\right\}.		  
$$
Man schreibt $F\sim G$ f\"ur zwei quadratische Formen $F$ und $G$, wenn es sogar ein $A\in SL(2,\Z)$
gibt mit $G=F^A$. F\"ur eine sogenannte uneigentliche Transformation\index{uneigentliche Transformation}\label{uneigentliche Transformation} $T\in GL(2,\Z)$
mit $\Det\ T=-1$ ist die Bedingung $F\sim F^T$ i.A. verletzt.\\

Aus (\ref{eq:8_10}) folgt nun, dass durch $\sim$ eine \"Aquivalenzrelation auf der Menge aller nicht 
verschwindenden quadratischen Formen mit ganzzahligen Koeffizienten definiert ist. 
Bezeichnen wir die \"Aquivalenzklasse, der $F$ angeh\"ort, mit $[F]$, dann 
haben nach Satz~\ref{defsatz:8_3} alle Formen 
$G(x,y)=a'x^2+b'xy+c'y^2$ aus $[F]$ dieselbe Diskriminante $D=b^2-4ac$, und es gilt 
$\ggT(a',b',c')=\ggT(a,b,c)$.\\

Ist insbesondere $F$ primitiv, so auch jede weitere Form aus $[F]$, so dass wir auch die Formenklasse $[F]$
primitiv nennen d\"urfen.\\

\begin{DefSatz}[Automorphe Transformationen, Pellsche Gleichung\index{Automorphe Transformation}\label{Auto}]
% Definition und Satz 8.3b, neu hinzugekommen
 \label{defsatz:8_3b}
	Gegeben sind eine Matrix 
	$A=\left( \begin{matrix}
		\alpha & \beta  \\
		\gamma & \delta\\
		\end{matrix} \right)\in SL(2,\Z)$, also
	   $\alpha,\beta,\gamma,\delta\in\Z$ und \\
		 $\alpha\delta-\beta\gamma=1$, sowie f\"ur $a,b,c \in \Z$ eine quadratische Form 
		$$F(x,y)=ax^2+bxy+cy^2$$ mit nichtquadratischer Diskriminante
		$D=b^2-4ac$. Wir nennen $A$ eine automorphe Transformation von $F$, wenn $F^A=F$ gilt.
		Die Hauptform $H_D$ zur Diskriminante $D$ definieren wir folgendermassen:
		F\"ur $D=4m \equiv 0 (4)$ setzen wir $H_D(x,y)=x^2-my^2$, bzw. f\"ur $D=4m+1 \equiv 1 (4)$
		sei $H_D(x,y)=x^2+xy-my^2$. Ist dann $F$ primitiv, so gelten f\"ur jede Matrix 
		$A \in \Z^{2 \times 2}$ die folgenden Aussagen:
		\begin{itemize}
		\item[(a)] Im Falle $D=4m \equiv 0 (4)$ ist $A$ 
		genau dann automorphe Transformation f\"ur $F$, wenn sich $A$ mit einer ganzzahligen L\"osung
		$x_0,y_0 \in \Z$ der Gleichung $H_D(x_0,y_0)=1$ auf folgende Weise darstellen l\"asst:
		\begin{equation*}
		A=\left( \begin{matrix}
		x_0-\frac{b}{2}y_0 & -cy_0  \\
		ay_0 & x_0+\frac{b}{2}y_0\\
		\end{matrix} \right)\,.
		\end{equation*}
		\item[(b)] Im Falle $D=4m+1 \equiv 1 (4)$ ist $A$ 
		genau dann automorphe Transformation f\"ur $F$, wenn sich $A$ mit einer ganzzahligen L\"osung
		$x_0,y_0 \in \Z$ der Gleichung $H_D(x_0,y_0)=1$ auf folgende Weise darstellen l\"asst:
		\begin{equation*}
		A=\left( \begin{matrix}
		x_0-\frac{b-1}{2}y_0 & -cy_0  \\
		ay_0 & x_0+\frac{b+1}{2}y_0\\
		\end{matrix} \right)\,.
		\end{equation*}
		\item[(c)] Sowohl f\"ur $D=4m \equiv 0 (4)$ als auch  f\"ur $D=4m+1 \equiv 1 (4)$ ist $A$ 
		genau dann automorphe Transformation f\"ur $F$, wenn sich $A$ mit einer ganzzahligen L\"osung
		$t,u \in \Z$ der {\it Pellschen Gleichung} $t^2-Du^2=4$
    \index{Pellsche Gleichung}\label{Pell} auf folgende Weise darstellen 
		l\"asst:
		\begin{equation*}
		A=\left( \begin{matrix}
		\frac{1}{2}(t-bu )& -cu  \\
		au & \frac{1}{2}(t+bu)\\
		\end{matrix} \right)\,.
		\end{equation*}
		\end{itemize}
		\dokendDef
\end{DefSatz}
{\bf Beweis:} F\"ur eine automorphe Transformation $A$ von $F$ gilt nach 
Definition und Satz \ref{Transformation der Formen}:
\begin{equation*}
\left( \begin{matrix}
		\alpha & \gamma  \\
		\beta & \delta
		\end{matrix} \right)
\left( \begin{matrix}
		2a & b  \\
		b & 2c
		\end{matrix} \right)
\left( \begin{matrix}
		\alpha & \beta  \\
		\gamma & \delta
		\end{matrix} \right)
=	\left( \begin{matrix}
		2a & b  \\
		b & 2c
		\end{matrix} \right)\,.
		\end{equation*}
		Hieraus folgt durch Multiplikation mit der inversen Matrix auf der linken Seite
\begin{equation*}
\left( \begin{matrix}
		2a & b  \\
		b & 2c
		\end{matrix} \right)
\left( \begin{matrix}
		\alpha & \beta  \\
		\gamma & \delta
		\end{matrix} \right)
=	\left( \begin{matrix}
		\delta & -\gamma  \\
		-\beta & \alpha
		\end{matrix} \right)
  \left( \begin{matrix}
		2a & b  \\
		b & 2c
		\end{matrix} \right)\,,
\end{equation*}
und hieraus durch Vergleich der Koeffizienten der letzten beiden Produktmatrizen:
\begin{equation}\label{auto3}
a(\delta-\alpha)=b\gamma\,, \quad a\beta+c\gamma=0\,, \quad c(\delta-\alpha)=-b\beta\,.
\end{equation}		
Diese drei Gleichungen sind also notwendig und hinreichend daf\"ur, dass $A$
automorphe Transformation von $F$ ist, allerdings unter der Annahme $A \in SL(2,\Z)$.
Aus den ersten beiden Gleichungen folgt $a | b \gamma$ und $a | c \gamma$.
Da $F$ primitiv ist, sind $a$ und $\text{ggT}(b,c)$ teilerfremd, und $a$ muss
bereits ein Teiler von $\gamma$ sein. Wir erhalten daher mit \eqref{auto3}
eine ganze Zahl $y_0$ mit
\begin{equation}\label{autoy0}
\gamma=ay_0\,, \quad \beta=-cy_0\,, \quad \delta-\alpha=by_0\,.
\end{equation}		

Wir unterscheiden nun zwei F\"alle.\\

A) Es sei  $D=4m \equiv 0 (4)$\,. Dann ist $b$ gerade.
Zun\"achst nehmen wir an, dass $A$ automorphe Transformation ist. Da $F$ primitiv ist, k\"onnen
nicht auch noch $a$ und $c$ gerade sein, und aus \eqref{auto3} folgt $\delta \equiv \alpha (2)$.
Wir definieren damit die ganze Zahl $x_0 = \frac{1}{2}(\alpha+\delta)$\,,
und erhalten aus \eqref{autoy0} die gew\"unschte Darstellung von $A$ der
Teilaussage (a) des Satzes. Die Darstellung von $A$ in (c) folgt dann f\"ur $t=2x_0$ und $u=y_0$,
wobei $\alpha\delta-\beta\gamma=1$ garantiert, dass jeweils die Gleichungen
$H_D(x_0,y_0)=1$ bzw. $t^2-Du^2=4$ erf\"ullt sind.

Nun gelte umgekehrt $H_D(x_0,y_0)=1$ mit ganzen Zahlen $x_0,y_0$.  F\"ur $t=2x_0$ und $u=y_0$
gilt dann auch die Pellsche Gleichung, und wir erhalten aus \eqref{auto3}, dass
		\begin{equation*}
		A=\left( \begin{matrix}
		x_0-\frac{b}{2}y_0 & -cy_0  \\
		ay_0 & x_0+\frac{b}{2}y_0\\
		\end{matrix} \right)=
		\left( \begin{matrix}
		\frac{1}{2}(t-bu )& -cu  \\
		au & \frac{1}{2}(t+bu)\\
		\end{matrix} \right)
		\end{equation*}
automorphe Transformation von $F$ ist.

B) Nun sei  $D=4m+1 \equiv 1 (4)$\,. Dann ist $b$ ungerade. Wir nehmen an, dass $A$ automorphe Transformation ist. Dann definieren wir die ganze Zahl
$$
 x_0 = \alpha+\frac{b-1}{2}y_0\,,
$$
und erhalten aus \eqref{autoy0} die gew\"unschte Darstellung von $A$. Da $A$ Determinante $1$ besitzt,
ist zudem die Gleichung $H_D(x_0,y_0)=1$ erf\"ullt. Nun setzen wir
$t=2x_0+y_0$, $u=y_0$, und erhalten aus der ersten Gleichung von \eqref{eq:8_3},
dort mit $H_D$ anstelle von $F$,
dass $t,u \in \Z$ L\"osungen der Pellschen Gleichung
$$
4H_D(x_0,y_0)=t^2-Du^2=4
$$
sind, welche die Matrixdarstellung f\"ur $A$ in (c) erf\"ullen.

Abschliessend nehmen wir $H_D(x_0,y_0)=1$ mit ganzen Zahlen $x_0,y_0$ an.  F\"ur $t=2x_0+y_0$ 
und $u=y_0$ gilt dann auch die Pellsche Gleichung, und wir erhalten aus \eqref{auto3}, dass
		\begin{equation*}
		A=\left( \begin{matrix}
		x_0-\frac{b-1}{2}y_0 & -cy_0  \\
		ay_0 & x_0+\frac{b+1}{2}y_0\\
		\end{matrix} \right)=
		\left( \begin{matrix}
		\frac{1}{2}(t-bu )& -cu  \\
		au & \frac{1}{2}(t+bu)\\
		\end{matrix} \right)
		\end{equation*}
automorphe Transformation von $F$ ist.
Damit ist der Satz bewiesen.
\dokendProof\\

Die Gleichung $x^2-my^2=1$ im Teil (a) des vorigen Satzes wird ebenfalls Pellsche Gleichung genannt. Den L\"osungen dieser Pellschen Gleichungen entsprechen also umkehrbar eindeutig die automorphen Transformationen der primitiven Formen mit nichtquadratischer Determinante.
Aus der Teilaussage (c) dieses Satzes ergibt sich auch, dass es f\"ur
positiv (bzw. negativ) definite Formen nur jeweils endlich viele automorphe Transformationen 
gibt, genauer gilt der\\

\begin{Satz}[Automorphismen positiv definiter Formen]
% Satz 8.3c, neu hinzugekommen
 \label{satz:8_3c}
	Gegeben sei f\"ur $a,b,c \in \Z$ eine primitve quadratische Form 
		$$F(x,y)=ax^2+bxy+cy^2$$ mit negativer Diskriminante
		$D=b^2-4ac<0$. Dann gelten die folgenden Aussagen:
		\begin{itemize}
		\item[(a)] Im Falle $D=-3$ sind die einzigen automorphen Transformationen von $F$
		gegeben durch die Matrizen
		\begin{equation*}
		A=\pm \left( \begin{matrix}
		1 & 0 \\
		0 & 1\\
		\end{matrix} \right)\,, \quad
			A=\pm \left( \begin{matrix}
		\frac{1-b}{2} & -c \\
		a & \frac{1+b}{2}\\
		\end{matrix} \right)\,, \quad
			A=\pm \left( \begin{matrix}
		\frac{1+b}{2} & c \\
		-a & \frac{1-b}{2}\\
		\end{matrix} \right)\,.
		\end{equation*}
		\item[(b)] Im Falle $D=-4$ sind die einzigen automorphen Transformationen von $F$
		gegeben durch die Matrizen
			\begin{equation*}
		A=\pm \left( \begin{matrix}
		1 & 0 \\
		0 & 1\\
		\end{matrix} \right)\,, \quad
			A=\pm \left( \begin{matrix}
		-\frac{b}{2} & -c \\
		a & \frac{b}{2}\\
		\end{matrix} \right)\,.
		\end{equation*}
		\item[(c)] F\"ur $D<-4$ besitzt $F$ nur die beiden trivialen Automorphismen zu
		\begin{equation*}
				A=\pm \left( \begin{matrix}
		1 & 0 \\
		0 & 1\\
		\end{matrix} \right)\,.
		\end{equation*}
		\end{itemize}
		\dokendDef
\end{Satz}
{\bf Beweis: }
Dies folgt mit den L\"osungen der Pellschen Gleichung $t^2-Du^2=4$
aus Satz \ref{defsatz:8_3b}(c),
die f\"ur $D=-3$ durch $|t|=2$, $u=0$ bzw. $|t|=|u|=1$ gegeben sind,
f\"ur $D=-4$ durch $|t|=2$, $u=0$ sowie $t=0$, $|u|=1$ und endlich f\"ur
$D<-4$ durch $|t|=2$, $u=0$.
\dokendProof\\

Ist $F(x,y)=ax^2+bxy+cy^2$ eine quadratische Form, so schreiben wir auch k\"urzer $F=(a,b,c)$.
F\"ur den Rest dieses Abschnitts betrachten wir nur noch indefinite Formen $F=(a,b,c)$, die nach Satz~\ref{defsatz:8_2}
zu den quadratischen Irrationalzahlen $X(F)=\dfrac{\sqrt{D}-b}{2a}$ mit
\mbox{$D=b^2-4ac$} in enger Beziehung stehen. F\"ur diese Formen werden wir u.a. zeigen,
dass sie im Gegensatz zu den positiv definiten Formen \"uber 
unendlich viele Automorphismen verf\"ugen.\\

F\"ur $-F=(-a,-b,-c)$ erhalten wir die zu $X(F)$ quadratisch konjugierte Zahl 
$X(-F)=\dfrac{-\sqrt{D}-b}{2a}$.
Wir nennen daher sowohl die Formen $F$ und $-F$ als auch die Formenklassen $[F]$ und $[-F]$ 
zueinander konjugiert.\\

Wir definieren noch mit der Spiegelung $S:=\left(\begin{matrix*}[r]
		     -1 & 0  \\
		     0 & 1
		    \end{matrix*} \right)$ die zur Form $F=(a,b,c)$ 
uneigentlich konjugierte Form\index{uneigentlich konjugierte Formen}\label{uneigentlich konjugierte Formen}
\begin{equation}\label{eq:8_11}
 F_{-}:=-F^S, \quad \text{d.h.} \quad F_{-}=(a,b,c)_{-}=(-a,b,-c),
\end{equation}
sowie die zur Formenklasse $[F]$ uneigentlich konjugierte Klasse	\index{uneigentlich konjugierte Klasse}\label{uneigentlich konjugierte Klasse}
\begin{equation}\label{eq:8_12}
[F]_{-}:=[F_{-}].
\end{equation}

Genau wie die zu $F$ \"aquivalenten Formen besitzt jede Form $G=(a',b',c')\in[F]_{-}$ 
nach Satz~\ref{defsatz:8_3} dieselbe Diskriminante wie $F$, und es gilt
$$ \ggT(a',b',c')=\ggT(a,b,c).$$

\begin{Satz}\label{satz:8_4}% Satz 8.4
 Es sei $F=(a,b,c)$ eine indefinite Form. 
 \begin{enumerate}[(a)] 
	\item Genau dann ist $G\in [F]_{-}$, wenn es eine uneigentliche Transformation\index{uneigentliche Transformation}\label{uneigentliche Transformation2}\linebreak \mbox{$T\in GL(2,\Z)$} gibt,
	also $\Det\ T=-1$, mit $G=-F^T$.\vspace{0.1cm}
	\item
	$(F_{-})_{-}=F$ und $([F]_{-})_{-}=[F]$.
	\item
	Speziell f\"ur $G:=(-c,-b,-a)$ gilt
	$X(G)=\dfrac{1}{X(F)}$ sowie $G\in [F]_{-}$.
\end{enumerate}	
\dokendSatz
\end{Satz}
{\bf Beweis:}~ 
\begin{enumerate}[(a)]	
	\item
	$G\in [F]_{-}\Leftrightarrow G\sim F_{-}=-F^S\Leftrightarrow G=-F^{S A}$ 
	f\"ur ein $A\in SL(2,\Z)$.
	Die Transformation $T=SA$ ist uneigentlich, und umgekehrt l\"a{\ss}t sich jedes uneigentliche
	$T$ mit $A=ST\in SL(2,\Z)$ in der Form $T=SA$ schreiben.
	\item ist trivial, und f\"ur (c) beachten wir $X(F)=\dfrac{\sqrt{D}-b}{2a}$ f\"ur $D=b^2-4ac$ sowie 
	$\displaystyle \frac{1}{X(F)}=2a \frac{\sqrt{D}+b}{D-b^2}=\frac{\sqrt{D}+b}{-2c}=X(G)$
	f\"ur 
	$$G=(-c,-b,-a)=-(a,b,c)^C=-F^C$$
	mit $C:=\left(\begin{matrix}
		     0 & 1  \\
		     1 & 0
		    \end{matrix} \right)$ und $\Det\ C=-1$. Gem\"a{\ss} der Teilaussage (a) 
	folgt $G\in[F]_{-}$.
\end{enumerate}			    
\dokendProof

\begin{Def}% Definition 8.5
	Zwei Irrationalzahlen $x$, $x'$  hei{\ss}en strikt \"aquivalent, 
	wenn es eine Matrix \linebreak \mbox{$\left(\begin{matrix}
		     \alpha & \beta  \\
		     \gamma & \delta
		    \end{matrix} \right)\in SL(2,\Z)$}
		    mit $x=\dfrac{\alpha x'+\beta}{\gamma x'+\delta}$ gibt.
	\dokendDef
\end{Def}

{\it Bemerkung}:~
Man best\"atigt m\"uhelos, dass hierdurch eine \"Aquivalenzrelation auf der Menge aller
Irrationalzahlen gegeben ist.

\begin{Satz}\label{satz:8_6}% Satz 8.6
 Es sei $F=(a,b,c)$ eine indefinite Form und $x:=X(F)=\dfrac{\sqrt{D}-b}{2a}$. Genau dann sind $x$ und $x'$
 strikt \"aquivalent gem\"a{\ss} $x=\dfrac{\alpha x'+\beta}{\gamma x'+\delta}$ mit 
 $A=\left(\begin{matrix}
		     \alpha & \beta  \\
		     \gamma & \delta
		    \end{matrix} \right)\in SL(2,\Z)$, wenn $x'=X(F^A)$ gilt.
\hfill\dokendSatz		    
\end{Satz}

{\it Bemerkung}:~
Ist zudem $F$ primitiv, dann auch $F^A$ nach Satz~6.3. Dann entsprechen gem\"a{\ss} 
Satz~\ref{defsatz:8_2} und Satz~\ref{satz:8_6} strikt \"aquivalenten quadratischen
Irrationalzahlen $x,x'$ genau \"aquivalente, primitive und indefinite Formen $F,F'$.\\

{\bf Beweis:}~
Es ist $x=\dfrac{\sqrt{D}-b}{2a}=\dfrac{\alpha x'+\beta}{\gamma x'+\delta}$ 
\"aquivalent zu
$$x'=\dfrac{\delta x-\beta}{\alpha-\gamma x}=\dfrac{\delta \sqrt{D}-(\delta b+2a\beta)}
{-\gamma \sqrt{D}+(\gamma b+2a\alpha)}$$
 unter Beachtung von $\alpha\delta-\beta\gamma=1$.
 Um von dem letzten Bruch den Nenner ganzzahlig zu machen, erweitern wir ihn mit 
 $\gamma \sqrt{D}+(\gamma b+2a\alpha)$, wobei noch
 $$
 [\delta \sqrt{D}-(\delta b+2a\beta)][\gamma \sqrt{D}+(\gamma b+2a\alpha)]
 =2a[\sqrt{D}-(2a\alpha\beta+b(\alpha\delta+\beta\gamma)+2c\gamma\delta)]
 $$ 
 sowie
$$
(\gamma b+2a\alpha)^2-\gamma^2(b^2-4ac)=4a(a\alpha^2+b\alpha\gamma+c\gamma^2)
$$
zu beachten ist. Aus Satz~\ref{defsatz:8_3} folgt nun mit $x=X(F)$ die behauptete \"Aquivalenz wegen 
$$
x'=\frac{\sqrt{D}-(2a\alpha\beta+b(\alpha\delta+\beta\gamma)+2c\gamma\delta)}
{2(a\alpha^2+b\alpha\gamma+c\gamma^2)}=X(F')
$$
f\"ur die transformierte Form $F'=F^{\left(\begin{smallmatrix}
		     \alpha & \beta  \\
		     \gamma & \delta
		    \end{smallmatrix} \right)}$.
\dokendProof 
 
\begin{Beis}\label{beis:8_7}% Beispiel 8.7
Die indefinite Form $F=(-17,-29,-7)$ und ihre Transformierte 
$$F'=F^{\left(\begin{smallmatrix*}[r] -1 & -3 \\ 1 & 2 \end{smallmatrix*}\right)}=(5,15,-7)$$ 
haben die Diskriminante $D=365$ mit $x'=X(F')=\dfrac{\sqrt{365}-15}{10}$ und 
$$\displaystyle x=\frac{-x'-3}{x'+2}=-\frac{\sqrt{365}-15+30}{\sqrt{365}-15+20}
=\frac{\sqrt{365}-(-29)}{(-34)}=X(F).$$
\dokendSatz
\end{Beis}

\begin{Satz}\label{satz:8_8}% Satz 8.8
Ist $x=X(F)$ eine quadratische Irrationalzahl mit der indefiniten quadratischen Form $F=(a,b,c)$ 
und $q\in\Z$, so ist $X(F')=x-q$ f\"ur 
$$F'=F^{\left(\begin{smallmatrix*} 1 & q \\ 0 & 1 \end{smallmatrix*}\right)}
=(a,b+2aq,aq^2+bq+c).$$
\dokendSatz
\end{Satz}
 
{\bf Beweis:}~
Wir setzen $\alpha=\delta=1$, $\gamma=0$,~\ $\beta=q$ in Satz~\ref{satz:8_6} 
und beachten die Transformationsformeln f\"ur $F'$ aus Satz~\ref{defsatz:8_3}.
\dokendProof\\
 
Formen $F$ und $F'$ wie in Satz~\ref{satz:8_8} werden auch parallel genannt, wenn sie sich 
durch eine Transformation $T=\left(\begin{matrix} 1 & q \\ 0 & 1 \end{matrix}\right)$
ineinander \"uberf\"uhren lassen.\\

Im Hinblick auf Satz~\ref{satz:8_6} k\"onnte man geneigt sein, strikt \"aquivalente 
Irrationalzahlen 
\begin{equation}\label{eq:8_13}
 x=\frac{\alpha x'+\beta}{\gamma x'+\delta} \quad \text{und} \quad x'
\end{equation}
einfach nur als \"aquivalent zu bezeichnen.\\

Doch w\"urde dies die in der Theorie der Kettenbr\"uche \"ubliche Konvention verletzen, nach der
die Irrationalzahlen $x$, $x'$ bereits f\"ur eine Transformation 
$$T=\left(\begin{matrix} \alpha & \beta \\ \gamma & \delta \end{matrix}\right)\in GL(2,\Z)$$
als \"aquivalent bezeichnet werden.\\

Nach einem wohlbekannten Resultat aus der Lehre der Kettenbr\"uche, das man etwa in dem Lehrbuch
von G.H.~Hardy und E.M.~Wright ``An introduction to the theory of numbers'', \cite[Theorem~175]{hw}, findet,
sind zwei Irrationalzahlen $x,x'$ (nicht notwendigerweise quadratisch) genau dann \"aquivalent,
wenn sich ihre beiden Kettenbruchentwicklungen nur um jeweils endlich viele Anfangsglieder unterscheiden.\\

Wir betrachten daher noch als Erg\"anzung zum Satz~\ref{satz:8_6} in (\ref{eq:8_13}) zwei
quadratische Irrationalzahlen $x=X(F)$, $x'=X(F')$, die durch eine Transformation 
$T=\left(\begin{matrix} \alpha & \beta \\ \gamma & \delta \end{matrix}\right)$ mit $\Det\ T=-1$
verbunden sind. Die indefiniten Formen $F$, $F'$ d\"urfen wir hierbei als primitiv voraussetzen,
so dass sie sich aus $x$, $x'$ eindeutig ergeben.\\

Dann sind die beiden Irrationalzahlen
\begin{equation}\label{eq:8_14}
 x=\frac{\frac{1}{x'}\beta+\alpha}{\frac{1}{x'}\delta+\gamma} \quad \text{und} \quad \frac{1}{x'} 
\end{equation}
verm\"oge der Transformation 
$\left(\begin{matrix} \beta & \alpha \\ \delta & \gamma \end{matrix}\right)\in SL(2,\Z)$
wieder strikt \"aquivalent, und nach Satz~\ref{satz:8_6} sowie Satz~\ref{satz:8_4} (c)
geh\"oren die Formen $F$, $F'$ zu uneigentlich konjugierten Formenklassen, denn es gilt
\begin{equation}\label{eq:8_15}
 F'=-F^{\left(\begin{smallmatrix} \alpha & \beta \\ \gamma & \delta \end{smallmatrix}\right)},
 \quad [F']=[F]_{-}
\end{equation}
wegen $\Det \left(\begin{matrix} \alpha & \beta \\ \gamma & \delta \end{matrix}\right)=-1$.

       \section{Kettenbruchentwicklung reell quadratischer Irrationalzahlen} \label{cha:8B}
 
\begin{Satz}\label{satz:8_9}% Satz 8.9
Es sei $G=(a,b,c)$ eine indefinite Form mit Diskriminante $D=b^2-4ac$ und 
%\linebreak
\mbox{$f:=\lfloor\sqrt{D}\rfloor$}.
Hierf\"ur definieren wir die K-Nachfolgeform $G'$ zu $G$ gem\"a{\ss} 
$$G'=(a',b',c')$$ mit
$$
a'=-c, \quad b'=-2cq-b, \quad c'=q(-cq-b)-a
$$
und 
\begin{equation*}
 \begin{split}
  q:=\left\lfloor\frac{1}{X(G)}\right\rfloor=\left\{\begin{tabular}{l}
	$\left\lfloor\dfrac{f+b}{-2c}\right\rfloor,\quad c<0$,\\\\
	$\left\lfloor-\dfrac{f+b+1}{2c}\right\rfloor, \quad c>0$.
	\end{tabular}\right.
 \end{split}
\end{equation*}
Dann gilt 
$$X(G')=\frac{1}{X(G)}-\left\lfloor\frac{1}{X(G)}\right\rfloor\in (0,1)$$ 
mit $G'=-G^{\left(\begin{smallmatrix} 0 & 1 \\ 1 & q \end{smallmatrix}\right)}$ und $[G']=[G]_{-}$.
\hfill\dokendSatz
\end{Satz}
 
{\bf Beweis:}~
Die Transformationsformel 
$G'=-G^{\left(\begin{smallmatrix} 0 & 1 \\ 1 & q \end{smallmatrix}\right)}$ folgt sofort 
aus Satz~\ref{defsatz:8_3}, also ist $[G']=[G]_{-}$ nach Satz~\ref{satz:8_4} (a). 
Aus der Darstellung
$$
G'=(-c,-b,-a)^{\left(\begin{smallmatrix} 1 & q \\ 0 & 1 \end{smallmatrix}\right)}
$$
und den S\"atzen Satz~\ref{satz:8_4} (c) sowie Satz~\ref{satz:8_8} folgt nun auch die Beziehung
$$ X(G')=\frac{1}{X(G)}-\left\lfloor\frac{1}{X(G)}\right\rfloor\in (0,1).$$
Um die Darstellung f\"ur $\displaystyle q=\left\lfloor\frac{1}{X(G)}\right\rfloor$ zu beweisen, verwenden
wir die Beziehung $\displaystyle \left\lfloor\frac{\xi}{n}\right\rfloor=
\left\lfloor\frac{\left\lfloor\xi\right\rfloor}{n}\right\rfloor$, die f\"ur alle $\xi\in\R$, $n\in\N$
gilt. Unter Beachtung der Fallunterscheidung f\"ur $c<0$ bzw. $c>0$ folgt die Darstellung 
aus $\displaystyle\frac{1}{X(G)}=\frac{\sqrt{D}+b}{-2c}$.
\dokendProof\\ 
 
Nun setzen wir $0<X(G)<1$ in Satz~\ref{satz:8_9} voraus. Die Bildung der K-Nachfolgeform von $G$ kann dann
als Anwendung eines Kettenbruchschrittes auf die quadratische Irrationalzahl $\dfrac{1}{X(G)}>1$
interpretiert werden:
$$
\frac{1}{X(G')}=\frac{1}{\frac{1}{X(G)}-\left\lfloor\frac{1}{X(G)}\right\rfloor}>1.
$$
Das Pr\"afix ``K'' steht hierbei f\"ur ``Kettenbruch''.

\begin{Def}\label{def:8_10}% Definition 8.10
% \hspace*{0cm}\\\vspace{-1cm}
	Die indefinite Form $G=(a,b,c)$ hei{\ss}t K-reduziert,\index{K-reduzierte indefinite Form}
	\label{K-reduzierte indefinite Form}
	wenn  
	f\"ur $D=b^2-4ac$,\linebreak \mbox{$f=\lfloor\sqrt{D}\rfloor$} folgendes gilt: $a>0$, $b>0$, 
	$f-\text{min}(2a,2|c|)<b\leq f$.
	Sie hei{\ss}t reduziert,\index{reduzierte indefinite Form}\label{reduzierte indefinite Form} wenn $G$ oder $G_{-}$ K-reduziert ist.
	\dokendDef
\end{Def} 
 
\begin{Bem}\label{bem:8_11}% Bemerkung 8.11
 Die gegen\"uber der K-Reduziertheit schw\"achere Reduziertheit der indefiniten Form $G=(a,b,c)$
 l\"a{\ss}t sich wie folgt charakterisieren:
 $$
 b>0, \quad f-\text{min}(2|a|,2|c|)<b\leq f.
 $$
\dokendBem
\end{Bem}

\begin{Satz}\label{satz:8_12}% Satz 8.12
Es sei $G=(a,b,c)$ eine indefinite Form, und damit insbesondere ihre Diskriminante $D=b^2-4ac$ 
keine Quadratzahl. Wir setzen $f:=\lfloor\sqrt{D}\rfloor$.
\begin{enumerate}[(a)]
	\item  Die folgenden drei Aussagen sind \"aquivalent:
	\begin{enumerate}[(i)]
		\item G ist K-reduziert,
		\item $a>0$, $c<0$, $|a+c|<b$,
		\item $a>0$, $b>0$, $c<0$ und $a-c \leq f$.
	\end{enumerate}	
	\item Die folgenden drei Aussagen sind \"aquivalent:
	\begin{enumerate}[(i)]
		\item G ist reduziert,
		\item $a c<0$, $|a+c|<b$,
		\item $b>0$, $ac<0$ und $|a|+|c| \leq f$.
	\end{enumerate}	
\end{enumerate}	
\hfill\dokendSatz
\end{Satz}

{\bf Beweis:}~
Wir k\"onnen generell $b\in\N$ voraussetzen. Wir erinnern auch daran, 
dass $D=b^2-4ac>0$ keine Quadratzahl ist, so dass f\"ur $f:=\lfloor\sqrt{D}\rfloor$ gilt:
\begin{equation}\label{eq:8_16}
 f<\sqrt{D}<f+1, \quad a \cdot c \neq 0.
\end{equation}
Die Bedingung $a c<0$ ist wegen $D=b^2-4ac$ zu $D>b^2$ und wegen (\ref{eq:8_16}) zu $b\leq f$ 
\"aquivalent. Wir d\"urfen daher zur Charakterisierung der K-Reduziertheit der indefiniten Form
$G=(a,b,c)$ schon vorab
\begin{equation}\label{eq:8_17}
 a>0, \quad b>0, \quad c<0
\end{equation}
annehmen, und wir haben insbesondere
\begin{equation}\label{eq:8_18}
 b \leq f.
\end{equation}
Mit (\ref{eq:8_16}), (\ref{eq:8_17}) erh\"alt man die beiden \"Aquivalenzumformungen
\begin{equation*}
  \begin{split}
    f-2a<b
   \Leftrightarrow  \sqrt{D}<b+2a
   \Leftrightarrow  b^2-4ac<b^2+4ab+4a^2
   \Leftrightarrow -(a+c)<b
  \end{split}
\end{equation*}
sowie
\begin{equation*}
  \begin{split}
    f+2c<b &\
   \Leftrightarrow  \sqrt{D}<b-2c
   \Leftrightarrow  b^2-4ac<b^2-4bc+4c^2\\
   &\ \Leftrightarrow  bc<(a+c)c
   \Leftrightarrow  a+c<b.\\
  \end{split}
\end{equation*}
Wir erhalten aus (\ref{eq:8_16}), (\ref{eq:8_17}) die \"Aquivalenz
\begin{equation}\label{eq:8_19}
  f-\text{min}(2a,2|c|)<b \Leftrightarrow |a+c|<b.
\end{equation}
Wiederum mit (\ref{eq:8_16}), (\ref{eq:8_17}) k\"onnen wir die letzte Ungleichung in (\ref{eq:8_19})
wie folgt umformulieren:
\begin{equation*}
  \begin{split}
    |a+c|<b &\
   \Leftrightarrow  a^2+2ac+c^2<b^2
   \Leftrightarrow  a^2-2ac+c^2<D\\
  &\ \Leftrightarrow  (a-c)^2<D
   \Leftrightarrow  a-c\leq f. \\  
  \end{split}
\end{equation*}
Das entsprechende Kriterium f\'ur Reduziertheit ergibt sich sofort aus dem f\"ur K-Reduziertheit.
\dokendProof

\begin{Bem}\label{bem:8_13}% Bemerkung 8.13
 Die Charakterisierung der Reduziertheit von $G$ in Satz~\ref{satz:8_12} (b) erfordert in (ii) und
 (iii) die Bedingung $a\cdot c<0$, wie das Beispiel der nicht reduzierten Form $G=(1,3,1)$
 mit $D=5$, $f=2$ lehrt.
 \dokendBem
\end{Bem}

\begin{Satz}\label{satz:8_14}% Satz 8.14
F\"ur die indefinite Form $G=(a,b,c)$ sei $0<X(G)<1$. Mit $D=b^2-4ac$ und \mbox{$f=\lfloor\sqrt{D}\rfloor$}
gelte $|b|\leq f$. Es sei $G'$ die K-Nachfolgeform zu $G$. Dann ist $G'$ eine K-reduzierte Form.
\dokendSatz
\end{Satz}

{\bf Beweis:}~
Nach Voraussetzung ist
\begin{equation}\label{eq:8_20}
 0<\frac{\sqrt{D}-b}{2a}<1,
\end{equation}
\begin{equation}\label{eq:8_21}
 |b|<\sqrt{D}.
\end{equation}
Aus (\ref{eq:8_20}) folgt 
$$
\frac{2a}{\sqrt{D}-b}=\frac{2a(\sqrt{D}+b)}{-4ac}=\frac{\sqrt{D}+b}{-2c}>1
$$
und nach (\ref{eq:8_21}) ist $\sqrt{D}\pm b>0$. Somit ist 
\begin{equation}\label{eq:8_22}
a'=-c>0, \quad a>0.
\end{equation}
Nach Satz~\ref{satz:8_9} ist
$$
0<X(G')=\frac{\sqrt{D}-b'}{2a'}<1,
$$
und unter Beachtung von $a'>0$ folgt hieraus
\begin{equation}\label{eq:8_23}
b'<\sqrt{D}<2a'+b'.
\end{equation}
Nun ist $b'=2a'q-b\geq 2a'-b>2a'-\sqrt{D}$ mit der nat\"urlichen Zahl $\displaystyle q=\left\lfloor\frac{1}{X(G)}\right\rfloor$.
Zusammen mit der rechten Ungleichung in (\ref{eq:8_23}) folgt
\begin{equation}\label{eq:8_24}
2a'-b'<\sqrt{D}<2a'+b'.
\end{equation}
Dies ist nur f\"ur $b'>0$ m\"oglich, und wir erhalten mit der linken Ungleichung in~(\ref{eq:8_23}):
\begin{equation}\label{eq:8_25}
0<b'<\sqrt{D}.
\end{equation}
Die rechte Ungleichung von (\ref{eq:8_24}) schreiben wir in der Form
\begin{equation}\label{eq:8_26}
\sqrt{D}-2a'<b',
\end{equation}
und aus der linken folgern wir noch $\displaystyle \frac{\sqrt{D}+b'}{2a'}=\frac{-2c'}{\sqrt{D}-b'}>1$.
Zusammen mit (\ref{eq:8_25}) haben wir nun
\begin{equation}\label{eq:8_27}
c'<0, \quad \sqrt{D}+2c'<b'.
\end{equation}
Schlie{\ss}lich beachten wir $f<\sqrt{D}<f+1$ f\"ur die nat\"urliche Zahl $f$, und erhalten:
\begin{equation*}
a'>0 \quad \text{aus (\ref{eq:8_22})}, \quad 0<b'\leq f \quad \text{aus (\ref{eq:8_25})}
\end{equation*}
sowie
\begin{equation*}
f-\text{min}(2a',2|c'|)<b'\quad \text{aus (\ref{eq:8_26}) und (\ref{eq:8_27})}.
\end{equation*}
Somit ist $G'=(a',b',c')$ eine K-reduzierte Form.
\dokendProof

\begin{Satz}\label{satz:8_15}% Satz 8.15
Es seien $G, \tilde{G}, G'$ indefinite und K-reduzierte Formen und $G'$ die K-Nachfolgeform sowohl von $G$ 
als auch von $\tilde{G}$. Dann ist $G=\tilde{G}$.
\dokendSatz
\end{Satz}

{\bf Beweis:}~
Wir setzen 
$$
G=(a,b,c), \quad \tilde{G}=(\tilde{a},\tilde{b},\tilde{c}), \quad G'=(a',b',c').
$$
Dann gilt $a'=-c=-\tilde{c}$ und insbesondere $c=\tilde{c}$. Folglich gelten die Kongruenzen
$$
b'\equiv -b(2c), \quad b'\equiv -\tilde{b}(2c),
$$
und somit ist $b\equiv \tilde{b}(2c)$. Wir m\"ussen nur noch $b\equiv \tilde{b}$ zeigen, da mit $c=\tilde{c}$ 
und der Gleichheit der Diskriminanten von $G$ und $\tilde{G}$ mit der von $G'$ auch $a=\tilde{a}$ folgt.\\

Nun gilt wegen der K-Reduziertheit von $G$ und $\tilde{G}$:
$$
f+2c<b\leq f \quad \text{und} \quad f+2c<\tilde{b}\leq f.
$$
Da die Zahlen $k\in \Z$ mit $f+2c<k\leq f$ ein vollst\"andiges Restsystem $\text{mod } 2c$ bilden, folgt aus 
$b\equiv \tilde{b}(2c)$ in der Tat $b=\tilde{b}$.
\dokendProof

\begin{Satz}\label{satz:8_16}% Satz 8.16
Es sei $G=(a,b,c)$ eine indefinite Form mit Diskriminante $D=b^2-4ac$ und \linebreak
$f=\lfloor\sqrt{D}\rfloor$.
Es sei $0<X(G)<1$ und $|b|>f$. F\"ur die K-Nachfolgeform \linebreak $G'=(a',b',c')$ von $G$ sei $|b'|>f$.
Dann gilt $a\cdot c>0$,~ $a'\cdot c'>0$ sowie
$$ |a'+c'|<\frac{1}{2}|a+c|.$$
\dokendSatz
\end{Satz}
{\bf Beweis:}~
F\"ur $q:=\left\lfloor\dfrac{1}{X(G)}\right\rfloor\in\N$ folgt $G=-G'^{\left(\begin{smallmatrix*}[r]
		     -q & 1  \\
		     1 & 0
		    \end{smallmatrix*} \right)}$ aus Satz~\ref{satz:8_9}, und hieraus 
$c=-a'$, $a=-(a'q^2-b'q+c')$ bzw. 
\begin{equation}\label{eq:8_28}
a+c=-(a'+c')-(a'q^2-b'q).
\end{equation}
Aus $|b|>f$ folgt $|b|\geq f+1>\sqrt{D}$, und somit $b^2>b^2-4ac$, d.h. $a\cdot c>0$: 
\begin{equation}\label{eq:8_29}
|b|>\sqrt{D}, \quad a\cdot c>0.
\end{equation}
Entsprechend folgt aus $|b'|>f$: 
\begin{equation}\label{eq:8_30}
|b'|>\sqrt{D}, \quad a'\cdot c'>0.
\end{equation}

{\bf Fall 1:}~ $a'<0$, und somit auch $c'<0$ nach (\ref{eq:8_30}). Aus Satz~\ref{satz:8_9} folgt 
$0<X(G')<1$, also
$$\frac{1}{X(G')}=\frac{2a'}{\sqrt{D}-b'}=\frac{2|a'|}{b'-\sqrt{D}}=\frac{\sqrt{D}+b'}{-2c'}>1.$$
Wir erhalten
$$b'>\sqrt{D}, \quad \sqrt{D}+b'>-2c',$$
und hieraus 
$$2b'>b'+\sqrt{D}>-2c',$$
also 
\begin{equation}\label{eq:8_31}
b'>-c'>0.
\end{equation}
Aus $a'=-c<0$ erhalten wir $c>0$, $a>0$ aufgrund der zweiten Ungleichung in~(\ref{eq:8_29}).
Somit folgt aus (\ref{eq:8_28}) und (\ref{eq:8_31}) im Fall 1:
\begin{equation*}
  \begin{split}
    |a+c|=a+c&\geq-(a'+c')-a'+b'\\
	     &>-(a'+c')-a'-c'\\
	     &=2|a'+c'|.
   \end{split}
\end{equation*}

{\bf Fall 2:}~ $a'>0$, und somit auch $c'>0$ nach (\ref{eq:8_30}). Nach Satz~\ref{satz:8_9} ist
$0<X(G')<1$ mit
$$\frac{1}{X(G')}=\frac{2a'}{\sqrt{D}-b'}=\frac{-\left(\sqrt{D}+b'\right)}{2c'}>1.$$
Wir erhalten hier $-\left(\sqrt{D}+b'\right)>2c'$ sowie
\begin{equation}\label{eq:8_32}
-b'>\sqrt{D}+2c'>c'>0.
\end{equation}
Aus $a'=-c>0$ erhalten wir $c<0$, $a<0$ aufgrund der zweiten Ungleichung in~(\ref{eq:8_29}).
Somit folgt aus (\ref{eq:8_28}) und (\ref{eq:8_32}) auch im Fall 2:
\begin{equation*}
  \begin{split}
    |a+c|=-(a+c)&=a'+c'+a'q^2-b'q\\
		&\geq a'+c'+a'-b'\\
		&>a'+c'+a'+c'\\
		&=2|a'+c'|.
   \end{split}
\end{equation*}
\dokendSatz

Die Form $F=(a,b,c)$ sei indefinit, und $D=b^2-4ac>0$ keine Quadratzahl. Dann ist 
\begin{equation}\label{eq:8_33}
x_0:=X(F)=\frac{\sqrt{D}-b}{2a}
\end{equation}
Irrationalzahl.
%\begin{equation}\label{eq:8_34}
%  \begin{split}
%    q_0:=\left\lfloor X(F)\right\rfloor=\left\{\begin{tabular}{l} 
%	$\left\lfloor\dfrac{f-b}{2a}\right\rfloor,\quad a>0$,\\\\
%	$\left\lfloor\dfrac{b-(f+1)}{-2a}\right\rfloor, \quad a<0$,
%	\end{tabular}\right.
%  \end{split}
%\end{equation}
Zu $F$ definieren wir mit $q_0:=\left\lfloor X(F)\right\rfloor$ die Parallelform
\begin{equation}\label{eq:8_34}
G_1=F^{\left(\begin{smallmatrix}
		1 & q_0  \\
		0 & 1
		\end{smallmatrix}\right)}
   =\left(a,b+2aq_0,c+q_0\left(b+aq_0\right)\right).
\end{equation}

Beginnend mit $j=1$ berechnen wir nun schrittweise zu $G_j=(a_j,b_j,c_j)$ die \mbox{K-Nachfolgeform}
$G_{j+1}=(a_{j+1},b_{j+1},c_{j+1})$ und setzen
\begin{equation}\label{eq:8_35}
x_j:=\frac{1}{X(G_j)}, \quad q_j:=\lfloor x_j \rfloor \quad \text{f\"ur } j\in \N.
\end{equation}
Die S\"atze~\ref{satz:8_8} und \ref{satz:8_9} liefern dann
\begin{equation}\label{eq:8_36}
x_{j+1}=\frac{1}{x_j-q_j} \quad \forall j\in \N_0,
\end{equation}
und nach dem erweiterten Euklidischen Algorithmus ist
\begin{equation}\label{eq:8_37}
x_j=\langle q_j,q_{j+1},q_{j+2},...\rangle \quad \forall j\in \N_0.
\end{equation}
Mit diesen Notationen und Bezeichnungsweisen gilt nun der

\begin{Satz}\label{satz:8_17}% Satz 8.17
Die Folge $(G_j)_{j\in\N}$ der indefiniten Formen $G_j$ m\"undet in eine Periode, die aus lauter K-reduzierten
Formen besteht. Die Periode beginnt stets mit der ersten \mbox{K-reduzierten} Form, die in der Folge $(G_j)_{j\in\N}$
auftritt.
\dokendSatz
\end{Satz}
{\bf Beweis:}~ F\"ur die Form $G_1$ in (\ref{eq:8_34}) gilt
$$ 0<X(G_1)=X(F)-\lfloor X(F)\rfloor<1$$
nach Satz~\ref{satz:8_8}.\\

Da $G_{j+1}$ f\"ur alle $j\in\N$ die K-Nachfolgeform von $G_j$ ist, folgt mit Satz~\ref{satz:8_9}:
\begin{equation}\label{eq:8_38}
0<X(G_j)<1 \quad \forall j\in\N.
\end{equation}

Alle Formen $G_j=(a_j,b_j,c_j)$ besitzen dieselbe Diskriminante $D=b_j^2-4a_jc_j$, und wir setzen wieder
$f:=\lfloor\sqrt{D}\rfloor$. Nun kann wegen (\ref{eq:8_38}) und Satz~\ref{satz:8_16} nicht $|b_j|>f$
f\"ur alle $j\in\N$ gelten.
Folglich gibt es einen Index $j_0\in\N$ mit $|b_{j_0}|\leq f$, und nach~(\ref{eq:8_38}) sowie nach 
Satz~\ref{satz:8_14} sind alle auf $G_{j_0}$ folgende Formen K-reduziert.\\

Zur festen Diskriminante $D$ gibt es aber nach Satz~\ref{satz:8_12} (a) nur endlich viele K-reduzierte
Formen, so dass die Folge $(G_j)_j\in\N$ in eine Periode m\"undet, die aus lauter K-reduzierten Formen 
besteht. Es sei $j_*\in\N$ der erste Index, ab dem $G_{j_*}, G_{j_*+1}, G_{j_*+2},...$ usw.
K-reduziert ist. Dann gibt es Zahlen $s,t\in\N_0$ mit $s<t$ und $G_{j_*+s}=G_{j_*+t}$.
Wir w\"ahlen $s$ minimal und f\"uhren die Annahme $s\geq 1$ zum Widerspruch:\\

Es ist $G_{j_*+s}$ die K-Nachfolgeform sowohl von $G_{j_*+s-1}$ als auch von $G_{j_*+t-1}$, und
alle drei Formen $G_{j_*+s}$, $G_{j_*+s-1}$, $G_{j_*+t-1}$ sind wegen $s\geq 1$ auch K-reduziert.
Satz~\ref{satz:8_15} liefert $G_{j_*+s-1}=G_{j_*+t-1}$, was der Minimalit\"at von $s$ widerspricht.
Somit ist $s=0$, und die Periode K-reduzierter Formen beginnt wie behauptet mit $G_{j_*}$.
~\dokendProof

    \section{Reduktion indefiniter quadratischer Formen} \label{cha:8C}
Hier fassen wir die Formeln aus dem vorigen Abschnitt noch einmal zu einem leicht zu implementierenden Rechenschema zusammen. Die Form \mbox{$F=(a,b,c)$}
sei indefinit, und $D=b^2-4ac>0$ keine Quadratzahl. Dann ist
\begin{equation}\label{eq:8_40}
X(F)=\frac{\sqrt{D}-b}{2a}
\end{equation}	
Irrationalzahl mit
\begin{equation}\label{eq:8_41}
  \begin{split}
    q_0:=\left\lfloor X(F)\right\rfloor=\left\{\begin{tabular}{l}
	$\left\lfloor\dfrac{f-b}{2a}\right\rfloor,\quad a>0$,\\\\
	$\left\lfloor\dfrac{b-(f+1)}{-2a}\right\rfloor, \quad a<0$,
	\end{tabular}\right.
  \end{split}
\end{equation}
f\"ur $f:=\lfloor \sqrt{D} \rfloor$. Zu $F$ definieren wir die Parallelform
\begin{equation}\label{eq:8_42}
G_1:=F^{\left(\begin{smallmatrix}
		1 & q_0  \\
		0 & 1
		\end{smallmatrix}\right)}
   =\left(a,b+2aq_0,c+q_0\left(b+aq_0\right)\right).
\end{equation}
Dabei gilt
\begin{equation} \label{eq:8_43}
X(G_1)=X(F)-q_0.
\end{equation}

Beginnend mit $j=1$ berechnen wir schrittweise zu $G_j=(a_j,b_j,c_j)$ die 
K-Nach\-fol\-ge\-form
$G_{j+1}=(a_{j+1},b_{j+1},c_{j+1})$, d.h.
\begin{equation} \label{eq:8_44}
\left\{\begin{tabular}{l}
	$a_{j+1}=-c_j$,\\
	$b_{j+1}=-2c_j q_j-b_j$,\\
	$c_{j+1}=q_j\cdot(-c_jq_j-b_j)-a_j$	 
	\end{tabular}\right.
\end{equation}
mit der Gr\"o{\ss}e
\begin{equation}\label{eq:8_45}
  \begin{split}
    q_j:=\left\lfloor\frac{1}{X(G_j)}\right\rfloor=\left\{\begin{tabular}{l}
	$\left\lfloor\dfrac{f+b_j}{-2c_j}\right\rfloor,\quad c_j<0$,\\\\
	$\left\lfloor-\dfrac{f+b_j+1}{2c_j}\right\rfloor, \quad c_j>0$.
	\end{tabular}\right.
  \end{split}
\end{equation}

Dann gilt f\"ur alle $j\in\N$:
\begin{equation} \label{eq:8_46}
X(G_{j+1})=\frac{1}{X(G_j)}-\left\lfloor\frac{1}{X(G_j)}\right\rfloor=\frac{1}{X(G_j)}-q_j.
\end{equation}
Mit den $G_j$ definieren wir f\"ur alle $j\in\N$ die Formen
\begin{equation}\label{eq:8_47}
  \begin{split}
    F_j:=\left\{\begin{tabular}{l}
	$G_j,\quad $ falls $j$ ungerade ist,\\
	$(-a_j,b_j,-c_j), \quad $ falls $j$ gerade ist. 
	\end{tabular}\right.
  \end{split}
\end{equation}
Mit der Matrix $S=\begin{pmatrix*}[r]
		-1 & 0  \\
		0 & 1\\
		\end{pmatrix*}$ kann man die $F_j$ einheitlich f\"ur $j\in\N$ in der Form 
\begin{equation} \label{eq:8_48}
F_j=\left((-1)^{j+1}a_j,b_j,(-1)^{j+1}c_j\right)=(-1)^{j+1}G_j^{\;S^{\; j+1}}
\end{equation}
schreiben.\\

Schlie{\ss}lich definieren wir die Transformationsmatrizen $T_j\in SL(2,\Z)$ f\"ur $j\in\N$ rekursiv gem\"a{ss}
\begin{equation} \label{eq:8_49}
T_1=\begin{pmatrix*}
		1 & q_0  \\
		0 & 1\\
		\end{pmatrix*}, \quad		
T_{j+1}=T_j\ A_j \quad \text{mit} \quad
A_j:=\begin{pmatrix*}
		0 & (-1)^{j+1}  \\
		(-1)^j & q_j\\
		\end{pmatrix*}.
\end{equation}
Dann gilt f\"ur alle $j\in\N$:
\begin{equation} \label{eq:8_50}
F_j=F^{T_j}, \quad F_{j+1}=F_j^{A_j},
\end{equation}
wobei $F_j$ und $F$ eigentlich \"aquivalent sind.

Wir setzen $x_0:=X(F)=\dfrac{\sqrt{D}-b}{2a}$ sowie $x_k:=\dfrac{1}{X(G_k)} \ \forall k\in\N$,
wenden den erweiterten Euklidischen Algorithmus auf die beiden Eingabewerte $x_0,1$ an und erhalten:
\begin{equation} \label{eq:8_51}
x_k=\langle q_k,q_{k+1},q_{k+2},...\rangle, \quad x_{k+1}=\frac{1}{x_k-q_k}, \quad q_k=\lfloor x_k \rfloor 
\quad \forall k\in\N_0,
\end{equation}
\begin{equation} \label{eq:8_52}
\frac{\sqrt{D}-b}{2a}=x_0=\langle q_0,q_1,q_2,...\rangle.
\end{equation}

Hiermit konstruieren wir ein Zahlenschema mit $6$ Spalten:			 

\begin{center}
	\begin{tabular*}{\textwidth}{ |p{0.8cm}||p{2.8cm}|p{1.8cm}|p{1.8cm}|p{1.8cm}|p{1.85cm}| }		\hline		
		\multicolumn{1}{|c||}{$j$} & \multicolumn{1}{c|}{ $G_j$} & \multicolumn{1}{c|}{ $q_j$}	
		& \multicolumn{1}{c|}{$F_j$} & \multicolumn{1}{c|}{ $T_j$} & \multicolumn{1}{c|}{ $A_j$}\\
		\hline		
		Index  $j\in\N$
		& F\"ur $j=1$ berechnet aus $F$ mit (\ref{eq:8_41}), (\ref{eq:8_42})
		und f\"ur \mbox{$j\geq 2$} jeweils mit (\ref{eq:8_44}) aus der Vorg\"angerzeile berechnet
		&  Berechnung mit (\ref{eq:8_45}) aus der zweiten Spalte
		&  Berechnung mit (\ref{eq:8_48}) aus den Spalten f\"ur $j$ und $G_j$
		&  Berechnung mit (\ref{eq:8_49});
			f\"ur $j\geq 2$ unter Verwendung der Vorg\"angerzeile
		& Aus der ersten und dritten Spalte gem\"a{\ss} (\ref{eq:8_49}).\\ \hline
	\end{tabular*}
\end{center}

Gesamt\"ubersicht der Berechnungsvorschriften zur Reduktion indefiniter quadratischer Formen.

\begin{center}
	
	\begin{tabular}{|l|l|} \hline
		
		~ Berechnung von $q_0$ und $G_1$ ~ &  ~ Berechnung von $q_j$ und $G_{j+1}$  ~ \\ 
		~ aus $F=(a, b, c)$:  ~ & ~ aus $G_j=(a_j, b_j, c_j)$ f\"ur $j\in\N$: ~ \\ \hline
		\rule{0pt}{4ex}	~ $q_0=\left\lfloor\dfrac{f-b}{2a}\right\rfloor$ f\"ur $a>0$ bzw. ~ &  ~ $q_j=\left\lfloor\dfrac{f+b_j}{-2c_j}\right\rfloor$ f\"ur $c_j<0$ bzw. ~  \\ 
		\rule{0pt}{4ex}	~ $q_0=\left\lfloor\dfrac{b-(f+1)}{-2a}\right\rfloor$ f\"ur $a<0$, ~ &  ~ $q_j=\left\lfloor-\dfrac{f+b_j+1}{2c_j}\right\rfloor$ f\"ur $c_j>0$, ~  \\
		\rule{0pt}{4ex}	~ $G_1=(a, b+2aq_0, c+q_0(b+a q_0))$, ~ & ~ $a_{j+1}=-c_j$, $b_{j+1}=-2c_jq_j-b_j$, ~ \\
		\rule{0pt}{3ex}	~ mit $f=\left\lfloor\sqrt{D}\right\rfloor$, $D=b^2-4ac>0$. ~ & ~ $c_{j+1}=q_j(-c_jq_j-b_j)-a_j$. ~  \\[0.2cm] \hline 
	\end{tabular}
\end{center}
\vspace{0.25cm}
\begin{center}
	
	\begin{tabular}{|l|l|} \hline
		~ Berechnung von $F_j$ aus $j$ und $G_j$: ~ & ~  Berechnung von $T_j$: ~  \\ \hline
		\rule{0pt}{4ex}	~ $F_j=((-1)^{j+1}a_j, b_j, (-1)^{j+1} c_j)$ ~ & ~ $T_1=\begin{pmatrix}
		1 & q_0 \\ 0 &  1
		\end{pmatrix}$, $T_{j+1}=T_j \cdot A_j$ und ~ \\ 
		\rule{0pt}{4ex}	~ f\"ur $j\in\N$ und $G_j=(a_j, b_j, c_j)$. ~ & ~  $A_j=\begin{pmatrix}
		0 & (-1)^{j+1} \\ (-1)^j & q_j
		\end{pmatrix}$ f\"ur $j\in\N$. ~ \\ \hline

	\end{tabular}
\end{center}
\vspace{0.25cm}
	
{\bf Beispiel:}~ $F(x,y)=-17x^2-29xy-7y^2$ liefert $a=-17$, $b=-29$, $c=-7$, \linebreak
\mbox{$D=b^2-4ac=365$} und $f=19$.
Wir haben $x_0=X(F)=-\dfrac{\sqrt{365}+29}{34}$, und\linebreak aus~(\ref{eq:8_41}) folgt $q_0=-2$.
Wir erhalten f\"ur $j=1$ aus (\ref{eq:8_42}), (\ref{eq:8_47}):
$$ G_1=F_1=(-17,39,-17).$$

\begin{center}
	\begin{tabular}{|l||c|c|c|c|c|c|c|c|c|c||c|} \hline
		~ $j$ ~  & ~~ $G_j$ ~~  & ~~ $q_j$ ~~ & ~~ $F_j$ ~~ & ~~ $T_j$ ~~ & $A_j$ \\
		\hline
		~ 1 ~  & -17, 39, -17 & 1 & -17, 39, -17 &\rule{0pt}{4ex} $\begin{array}{rr} 1 & -2\\0 & 1\end{array}$  & $\begin{array}{cc} 0 & 1\\-1 & 1\end{array}$  \\  \hline
		~ 2 ~  &  17, -5, -5  & 1 &  -17, -5, 5   & \rule{0pt}{4ex}$\begin{array}{rr} 2 & -1\\-1 & 1\end{array}$  & $\begin{array}{cc} 0 & -1\\1 & 1\end{array}$  \\ \hline
		~ 3 ~  &  5, 15, -7   & 2 &  5, 15, -7   & \rule{0pt}{4ex}$\begin{array}{rr} -1 & -3\\1 & 2\end{array}$  & $\begin{array}{cc} 0 & 1\\-1 & 2\end{array}$   \\ \hline
		~ 4 ~  &  7, 13, -7   & 2 &  -7, 13, 7   & \rule{0pt}{4ex}$\begin{array}{rr} 3 & -7\\-2 & 5\end{array}$  & $\begin{array}{cc} 0 & -1\\1 & 2\end{array}$   \\ \hline
		~ 5 ~  &  7, 15, -5   & 3 &  7, 15, -5   &\rule{0pt}{4ex} $\begin{array}{rr} -7 & -17\\5 & 12\end{array}$  & $\begin{array}{cc} 0 & 1\\-1 & 3\end{array}$      \\ \hline
		~ 6 ~  &  5, 15, -7   & 2 &  -5, 15, 7   & \rule{0pt}{4ex}$\begin{array}{rr} 17 & -58\\-12 & 41\end{array}$  & $\begin{array}{cc} 0 & -1\\1 & 2\end{array}$      \\ \hline
		~ 7 ~  &  7, 13, -7   & 2 &  7, 13, -7   & \rule{0pt}{4ex}$\begin{array}{rr} -58 & -133\\41 & 94\end{array}$  & $\begin{array}{cc} 0 & 1\\-1 & 2\end{array}$      \\ \hline
		~ 8 ~  &  7, 15, -5   & 3 &  -7, 15, 5   &\rule{0pt}{4ex} $\begin{array}{rr} 133 & -324\\-94 & 229\end{array}$  & $\begin{array}{cc} 0 & -1\\1 & 3\end{array}$      \\ \hline
		~ 9 ~  &  5, 15, -7   & 2 &  5, 15, -7   & \rule{0pt}{4ex}$\begin{array}{rr} -324 & -1105\\229 & 781\end{array}$  & $\begin{array}{cc} 0 & 1\\-1 & 2\end{array}$      \\ \hline
		
	\end{tabular}
\end{center}
\vspace{0.25cm}
Wir haben die Form $F=(-17, -29, -7)$,
$$X(F)=x_0=-\frac{\sqrt{365}+29}{34}=\langle-2, 1, 1,\overline{ 2, 2, 3}\rangle,$$
$F_3=F^{T_3}=F^{T_9}$, so dass $T_9 T_3^{-1}=\left(\begin{array}{rr} 457 & 133\\-323 & -94\end{array}\right)$ automorphe Substitution f\"ur $F$ ist.\\

Zum Vergleich mit der Reduktion der Form $F$ f\"uhren wir nun noch die vollkommen analoge Kettenbruchentwicklung von $X(F)$ durch:\\

\textbf{Erweiterter Euklidischer Algorithmus zur Berechnung der \boldmath{$x_j$}:}
\begin{center}
	\begin{tabular}{|c||c|c|c|c|c|} \hline
		
		~ $j$ ~ & ~ $x_j$ ~ & ~~ $q_j$ ~~ & ~~ $s_j$ ~~ & ~~ $t_j$ ~~ & $x_j=\langle q_j,q_{j+1},q_{j+2},... \rangle$ \\
		\hline
		~ 0 ~  & \rule{0pt}{5ex}\centering $-\dfrac{\sqrt{365}+29}{34}$ &-2 & 1    & 0  & $\langle-2, 1, 1,\overline{ 2, 2, 3}\rangle=-1.41485215...$  \\ [0.3cm] \hline
		~~ 1 ~~   & \rule{0pt}{5ex}$\dfrac{\sqrt{365}+39}{34}$  & 1 & -2   & 1  & $\langle1, 1,\overline{ 2, 2, 3}\rangle=1.708969799...$  \\[0.3cm] \hline
		~ 2 ~   & \rule{0pt}{5ex}$\dfrac{\sqrt{365}-5}{10}$   & 1 & -1   & 1  & $\langle1,\overline{ 2, 2, 3}\rangle=1.4104973...$  \\ [0.3cm]\hline
		~ 3 ~   & \rule{0pt}{5ex}$\dfrac{\sqrt{365}+15}{14}$  & 2 & -3   & 2  & $\langle\overline{ 2, 2, 3}\rangle=2.4360695...$  \\[0.3cm] \hline
		~ 4 ~   & \rule{0pt}{5ex}$\dfrac{\sqrt{365}+13}{14}$  & 2 & -7   & 5  & $\langle\overline{ 2, 3, 2}\rangle=2.293212...$  \\ [0.3cm]\hline
		~ 5 ~   &\rule{0pt}{5ex} $\dfrac{\sqrt{365}+15}{10}$  & 3 & -17  & 12 & $\langle\overline{ 3, 2, 2}\rangle=3.4104973...$  \\[0.3cm] \hline
		~ 6 ~   &\rule{0pt}{5ex} $\dfrac{\sqrt{365}+15}{14}$  & 2 & -58  & 41 & $\langle\overline{ 2, 2, 3}\rangle=2.4360695...$  \\ [0.3cm]\hline
		~ 7 ~   & \rule{0pt}{5ex}$\dfrac{\sqrt{365}+13}{14}$  & 2 & -133 & 94 & $\langle\overline{ 2, 3, 2}\rangle=2.293212...$  \\[0.3cm] \hline
		~ 8 ~   &\rule{0pt}{5ex} $\dfrac{\sqrt{365}+15}{10}$  & 3 & -324 & 229 & $\langle\overline{ 3, 2, 2}\rangle=3.4104973...$  \\[0.3cm] \hline
		
	\end{tabular}
\end{center}
\vspace{0.25cm}
$F=(-17, -29, -7)$, $x_0=X(F)$, und f\"ur $j\in\N_0$:
\begin{equation*}
\begin{split}
x_{j+1}&=\dfrac{1}{X(G_{j+1})}=\dfrac{1}{x_j-q_j}\quad \text{mit}\quad q_j=\lfloor x_j \rfloor.\\
s_0&=1,\quad s_1=q_0,\quad s_{j+1}=s_{j-1}+s_j\cdot q_j,\\
t_0&=0,\ \quad t_1=1,\enspace\quad t_{j+1}=t_{j-1}+t_j\cdot q_j
\end{split}
\end{equation*}
f\"ur alle $j\in\N$.\\

Nun heben wir zwei wichtige Resultate hervor, die
eine direkte Folge unseres Reduktionsverfahrens f\"ur
indefinite Formen sind. So stellt die n\"achste einfache Folgerung
aus der Darstellung \eqref{eq:3_15} der allgemeinen Kettenbr\"uche
$x_j=\langle q_j,q_{j+1},q_{j+2} \ldots \rangle$ aus Lektion \ref{cha:3A}
und Satz \ref{satz:8_17} eines unserer Hauptergebnisse dar:

\begin{Satz}\label{satz:kettenperiode}
Genau die quadratischen Irrationalzahlen besitzen eine Kettenbruchentwicklung, die (ggf. nach einer
endlichen Vorperiode) in eine Periode einm\"undet.
\hfill\dokendSatz
\end{Satz}

Mit dem folgenden Ergebnis schliesst sich auch der Kreis, der in Satz \ref{defsatz:8_3b}
und Satz \ref{satz:8_3c} seinen Ursprung hat:

\begin{Satz}\label{satz:autoinfty}
Jede indefinite, primitive Form $F=(a,b,c)$ mit Diskriminante\footnote{wir erinnern wieder daran, dass wir nur Formen mit nichtquadratischer Diskriminante betrachten.} $D=b^2-4ac>0$ besitzt
unendlich viele automorphe Transformationen\index{Automorphe Transformation} \label{Auto2}. Insbesondere besitzt die Pellsche Gleichung\index{Pellsche Gleichung}\label{Pell2}
\begin{equation*}
t^2-Du^2=4 \quad \text{bzw.}\quad H_D(x,y)=1
\end{equation*}
aus Satz \ref{defsatz:8_3b} jeweils unendlich viele ganzzahlige L\"osungen $t,u$ bzw. $x,y$.
\hfill\dokendSatz
\end{Satz}
{\bf Beweis:~} Da sich jede Form $F$ duch eine Kette \"aquivalenter Formen $F_j$ 
(in der vierten Spalte unseres Schemas) in eine reduzierte Form
\"uberf\"uhren l\"asst, k\"onnen wir annehmen, $F$ sei reduziert. Sobald die Kette 
der $F_j=(a_j,b_j,c_j)$ periodisch wird, alternieren die Vorzeichen der $a_j$.
Wir k\"onnen daher von vorneherein annehmen, dass $F=(a,b,c)$ mit $a>0$ sogar K-reduziert ist.

Wenden wir nun den Reduktionsmechanismus
dieses Abschnitts auf $F=F_1$ an, so entsteht eine reine Periode gerader L\"ange 
von Formen $F_1 \sim F_2 \sim \ldots \sim F_{2m}$ mit $m \in \N$ und $F_{2m+1}=F_{1}$,
und f\"ur alle $j \in \N$ gilt $F_{j+1}=F_j^{A_j}$. 
Nun ist das Produkt zweier aufeinanderfolgender Matrizen $A_j, A_{j+1} \in SL(2,\Z)$ 
aus der letzten Spalte des Schemas f\"ur ungerades $j$ gegeben durch
$$
A_j A_{j+1}= \begin{pmatrix}
0 & 1\\
-1 & q_j
\end{pmatrix}
\begin{pmatrix}
0 & -1\\
1 & q_{j+1}
\end{pmatrix}=
\begin{pmatrix}
1 & q_{j+1}\\
q_j & 1+q_j q_{j+1}
\end{pmatrix}\,.
$$
Die aufeinanderfolgenden Produkte der Matrizen $A_jA_{j+1}$ liefern beliebig gro{\ss}e Eintr\"age, wenn man f\"ur $j$ die Folge der ungeraden Zahlen durchl\"auft und somit auch beliebig oft
die volle Periode reduzierter Formen, die von $F=F_1$ ausgeht. Dies liefert unendlich viele automorphe Transformationen von $F$, und der Rest der Behauptung folgt sofort 
aus Satz \ref{defsatz:8_3b}. 
\dokendProof\\

Wir erw\"ahnen an dieser Stelle, dass f\"ur die Reduktion der indefiniten Formen auch
andere Verfahren in der Literatur beschrieben werden. So findet man etwa im Lehrbuch
von Scholz und Schoeneberg \cite[\S 31]{Scholzundco} das Verfahren
der sogenannten halbreduzierten rechten Nachbarformen, dass dieselben Perioden
reduzierter indefiniter Formen wie unser Schema liefert, sich aber 
bei der Reduktion der Formen in der Vorperiode 
unterscheiden kann. Die Form $F=(a,b,c)$ sei indefinit mit nichtquadratischer Diskriminante
$D=b^2-4ac>0$, und es sei wieder $f=\lfloor\sqrt{D} \rfloor$ (Vorsicht: bei 
\cite[\S 31]{Scholzundco} ist $f=\lceil \sqrt{D} \rceil$).
Die halbreduzierte rechte Nachbarform von $F$ ist dann 
$$
\mathcal{R}(F)=(c,2ct-b,a-bt+ct^2) \quad \text{mit} \quad 
t=\sign(c)\left\lfloor \frac{f+b}{2|c|} \right \rfloor\,.
$$
Beim Reduktionsverfahren von \cite[\S 31]{Scholzundco} betrachten wir die folgende
Kette \"aquivalenter Formen, die aus $F$ durch schrittweise Bildung der
halbreduzierten rechten Nachbarformen hervorgeht:
$$
 F,~\mathcal{R}(F),~ \mathcal{R}(\mathcal{R}(F)),~
\mathcal{R}(\mathcal{R}(\mathcal{R}(F))), \ldots \quad \text{usw.}
$$
Nun l\"asst sich \cite[Satz 81]{Scholzundco} folgendermassen auf unser Schema \"ubertragen:

\begin{Satz}\label{satz:autochar}
Sind $F, F'$ indefinit und reduziert und gilt $F' = F^A$
mit einem $A \in SL(2,\Z)$, so liegt $F'$ in der von $F$ ausgehenden Periode reduzierter Formen aus der vierten Spalte unseres Schemas.
\hfill\dokendSatz
\end{Satz}

Durch die abschlie{\ss}enden Internet-Recherchen bei der Fertigstellung dieses Buches 
ist uns noch ein drittes Reduktionsverfahren f\"ur
indefinite Formen bekannt geworden, das auf einer alternativen Art von
Kettenbruchentwicklung bzw. Formen-Reduziertheit basiert, 
siehe hierzu Zagier \cite[\S 13]{zagier}. \\

F\"ur positiv definite Formen gestaltet sich die Formen-Reduktion mit Hilfe eines effizienten Verfahrens
sowie die Bestimmung der \"Aquivalenzklassen wesentlich einfacher
als im indefiniten Fall, siehe hierzu die kompakte Darstellung \cite[\S 30]{Scholzundco}.\\
Zum Abschluss stellen wir nur eine interessante \"Ubungsaufgabe, deren Ausf\"uhrung wir dem geneigten Leser \"uberlassen m\"ochten.\\

{\bf Aufgabe:~} Zur Reduktion der indefiniten Formen implementiere man die
drei oben genannten Verfahren, n\"amlich das Kettenbruchverfahren dieses Abschnittes
sowie das Verfahren der halbreduzierten rechten Nachbarformen und
das in Zagier  \cite[\S 13]{zagier} beschriebene Verfahren. Hierauf 
vergleiche man diese Verfahren, indem man sie f\"ur gr\"ossere Werte von 
$n \in \N$, $n \geq 2$ auf folgende indefiniten Formen anwendet:
$$
F_n=((n+1)^2-2,-2(n^2+n-2),n^2-2)
$$
mit Diskriminante $D=8$ und $X(F_n)=\langle 0,1,n+1,\overline{2}\rangle$ bzw.
$$
\tilde{F}_n=(n,-n,-1)
$$
mit Diskriminante $D=n(n+4)$ und 
$\begin{displaystyle}X(\tilde{F}_n)
=\frac{1}{2}+\sqrt{\frac{1}{4}+\frac{1}{n}} = \langle \overline{1,n}\rangle\,.
\end{displaystyle}$ 

	\chapter{Anhang}\label{cha:anhang}
	\section{Logische Symbole, Mengen und Abbildungen}\label{cha:anhang1}

\begin{center}
		{\bf Logische Symbole\index{Logische Symbole}\label{Logische Symbole} der mathematischen Umgangssprache}
	\end{center}
	\begin{eqnarray*}
		(1) \quad & \neg \, A   \qquad  & {\mbox{nicht~$A$}}\,,\qquad  \\
		(2) \quad & A \, \wedge \, B \qquad  & {\mbox{$A$~und~$B$}}\,, \qquad \\
		(3) \quad & A \, \vee \, B  \qquad  & {\mbox{$A$~oder~$B$}}\,,\qquad \\ 
		(4) \quad & A \, \Rightarrow \, B  \qquad  & {\mbox{$A$~impliziert~$B$}}\,,\qquad \\
		(5) \quad & A \, \Leftrightarrow \, B  \qquad  & 
		{\mbox{$A$~und~$B$~sind~\"aquivalent}}\,,\qquad \\
		(6) \quad & \forall \, x ~ A(x)   \qquad  & {\mbox{f\"ur~alle~$x$~gilt~$A(x)$}}\,, \qquad \\
		(7) \quad & \exists \, x ~ B(x)   \qquad  & 
		{\mbox{es~gibt~ein~$x$~f\"ur~das~$B(x)$~gilt}}\,. \qquad \\  
	\end{eqnarray*}
	
	In (1)-(5) sind $A$, $B$ Aussagen\index{Aussage}\label{Aussage}, in (6) und (7) dagegen Aussageformen\index{Aussageform}\label{Aussageform},
	die von einer freien Variablen $x$ abh\"angen d\"urfen. Die Variable $x$
	entstammt dabei einer festen, vorgegebenen Grundmenge $M$, die oft
	nicht explizit in den Formeln mitgef\"uhrt wird. Ein Beipiel f\"ur (6)
	w\"are demnach $\forall x \in \Z \,:\, x^2 \geq 0$, oder einfach
	$\forall x \, ( x^2 \geq 0)$, nachdem zuvor die Grundmenge $M := \Z$
	festgelegt worden ist.\\

	\begin{center}
		{\bf Wahrheitstabellen\index{Wahrheitstabelle}\label{Wahrheitstabelle} f\"ur aussagenlogische Verkn\"upfungen\index{aussagenlogische Verkn\"upfung}\label{aussagenlogische Verknuepfung}}
	\end{center}
	Hier sind $\alpha$ und $\beta$ Aussagen mit dem Wahrheitsgehalt
	w=wahr oder f=falsch.
	\begin{center}
		\begin{tabular}{||l|l||c||c||c||c||c||} \hline
			$\alpha$ \quad & $\beta$ \quad & $\neg \alpha$ \quad & 
			$\alpha \, \wedge \,\beta$ \quad &
			$\alpha \vee \beta$ \quad & $\alpha \Rightarrow \beta$ \quad &
			$\alpha \Leftrightarrow \beta$ \\ \hline
			w \quad  & w \quad & f \quad & w \quad  & w \quad & w \quad  & w \\ \hline
			w \quad  & f \quad & f \quad & f \quad  & w \quad & f \quad  & f \\ \hline
			f \quad  & w \quad & w \quad & f \quad  & w \quad & w \quad  & f \\ \hline
			f \quad  & f \quad & w \quad & f \quad  & f \quad & w \quad  & w \\ \hline
		\end{tabular}	\\
	\end{center}

	\begin{center}
		{\bf Symbole der (nicht formalisierten) Mengenlehre}
	\end{center}
	
	Wir betrachten hier Teilmengen $K$, $L$ einer vorgegebenen Grundmenge $M$.
	\begin{align*}
	(1) \quad & x \in M \setminus \, K   \quad  & 
	x \notin K \quad & {\mbox{Komplement~von~$K$}}\,,\\
	(2) \quad & x \in K \cap \, L   \quad  & 
	x \in K \,\wedge\, x \in L \quad & {\mbox{Durchschnitt}}\,,\\
	(3) \quad & x \in K \cup \, L   \quad  & 
	x \in K \,\vee\, x \in L \quad & {\mbox{Vereinigung}}\,,\\
	(4) \quad & \forall x \,(x \in K \Rightarrow x \in L)   \quad  & 
	K \,\subseteq L \quad & {\mbox{Inklusion}}\,,\\
	(5) \quad & \forall x \,(x \in K \Leftrightarrow x \in L)   \quad  & 
	K = L \quad & {\mbox{Mengengleichheit}}\,.\\
	\end{align*}
	
	\begin{center}
		{\bf Wichtige Beispiele f\"ur Mengen\index{Menge}\label{Menge}}
	\end{center}
	\begin{itemize}
	\item[(1)] \quad $\N=\{1,2,3,\ldots\}$ ist die Menge der nat\"urlichen Zahlen.
	\item[(2)] \quad $\N_0=\{0,1,2,3,\ldots\}$ ist die Menge der nat\"urlichen Zahlen inklusive der Null.
        \item[(3)] \quad $\Z=\{0,\pm 1, \pm 2, \pm 3,\ldots\}$ ist die Menge der ganzen Zahlen.
       \item[(4)] \quad $\Q=\{a/b\,:\,a \in \Z\,,~b \in \N \}$ ist die Menge der rationalen Zahlen.
    \item[(5)] \quad $\R$ ist die Menge der reellen Zahlen.
  \item[(6)] \quad $\C=\{x+iy\,:\,x,y \in \R\}$ ist die Menge der komplexen Zahlen.
	\end{itemize}
	
	Besonders oft werden Intervalle als spezielle Teilmengen der reellen Zahlen in der Mathematik ben\"otigt:
	
	\begin{center}		
		{\bf Notationen f\"ur Intervalle\index{Intervall}\label{Intervall}}
	\end{center}
	{\bf Abgeschlossenes Intervall:} 
	$[a,b] := \{\, x \in \R \,:\, a \leq x  \leq b\,\}$\,.\\
	
	{\bf Offenes Intervall:} 
	$(a,b) := \{\, x \in \R \,:\, a < x  < b\,\}$\,.\\
	Die Menge $\R^+ := \{\, x \in \R \,:\, x >0 \,\}$
	ist ein ``unendliches'' offenes Intervall. \\
	
	{\bf Halboffene Intervalle:}\\ 
	$(a,b] := \{\, x \in \R \,:\, a < x  \leq b\,\}$\,, \quad
	$[a,b) := \{\, x \in \R \,:\, a \leq x  < b\,\}$\,.\\
	Die Menge $\R^+_0 := \{\, x \in \R \,:\, x \geq 0 \,\}$
	ist ein ``unendliches'' halboffenes Intervall. \\
	
Die Bildung kartesischer Produktmengen und deren Teilmengen ist ein besonders wichtiges Konstruktionsprinzip in der Mathematik, um aus gegebenen Mengen neue Mengen zu bilden und um Eigenschaften von komplexerer Struktur zu beschreiben:\\

	{\bf Kartesisches Produkt\index{Kartesisches Produkt}\label{Kartesisches Produkt} von $n$ Mengen und $n$-stellige Relationen\index{Relation}\label{Relation}}
	
	Sind $M_1$, $M_2$, ... , $M_n$ nichtleere Mengen, so ist ihr 
	kartesisches Produkt erkl\"art als Menge aller geordneter ``$n$-Tupel''
	$(x_1, ... , x_n)$ mit $x_1 \in M_1$, $x_2 \in M_2$, ... , $x_n \in M_n$, d.h.
	$$
	M_1 \times M_2 \times ... \times M_n := 
	\{\,(x_1, ... , x_n)\,:\, x_k \in M_k \mbox{~f\"ur~} k=1,...,n \,\}\,.
	$$
	Der $\R^n := \R  \times ... \times \R$ mit $n$ Faktoren ist ein wichtiges
	Beispiel. Eine Teilmenge des kartesischen Produktes
	$M_1 \times M_2 \times ... \times M_n$ hei{\ss}t 
	{\it n-stellige Relation}\,.\\
	
	{\bf Funktionen\index{Funktion}\label{Funktion} (auch Abbildungen\index{Abbildung}\label{Abbildung} genannt)}
	
	Es seien $A$, $B$ nichtleere Mengen. Eine {\it Funktion} bzw. {\it Abbildung}
	$f$ mit Definitionsbereich $A$ und Wertebereich $B$ ist eine Zuordnung,
	die jedem $x \in A$ genau einen Wert $y \in B$ zuordnet. Wir schreiben dann
	$y = f(x)$ und nennen $f(x)$ das Bild bzw. den Funktionswert von $x$.
	
	Formal gesehen sind Funktionen spezielle Teilmengen $G \subseteq A \times B$
	der kartesi\-schen Produktmenge $A \times B$ mit der Eigenschaft, 
	da{\ss} es zu jedem $x \in A$ genau ein Paar $(x,y) \in G$ gibt.
	Im Sprachgebrauch nennt man $G$ aber meistens den ``Graphen'' der Funktion $f$.
	\begin{itemize}
		\item Die Funktion $f$ hei{\ss}t {\it injektiv}, wenn 
		f\"ur alle $x, y \in A$ aus $f(x)=f(y)$ stets $x=y$ folgt.
		\item Die Funktion $f$ hei{\ss}t {\it surjektiv}, wenn es zu jedem $z \in B$
		mindestens ein $x \in A$ gibt mit $f(x)=z$.
		\item Eine injektive und surjektive Funktion $f$ wird auch {\it bijektiv}
		bzw. {\it Bijektion} genannt. Zu jeder bijektiven Funktion $f : A \to B$
		gibt es die sogenannte Umkehrabbildung $f^{-1} : B \to A$, wobei
		f\"ur jedes $y \in B$ der Wert $x = f^{-1}(y)$ der {\it Umkehrabbildung} 
		durch die Beziehung $f(x) =  f(f^{-1}(y)) = y$ eindeutig bestimmt ist.
		Es gilt $(f^{-1})^{-1}=f$.
	\end{itemize}
	
	{\bf Verkettung von Funktionen\index{Verkettung von Funktionen}\label{Verkettung von Funktionen}}
	
	Sind $A,B,B',C$ nichtleere Mengen mit $B \subseteq B'$
	und $h : A \to B$ bzw. $g : B' \to C$
	Abbildungen, so definiert ihre {\it Verkettung} oder 
	{\it Komposition} 
	eine neue Funktion $g \circ h : A \to C$ gem\"a{\ss}
	$$
	(g \circ h)(x) = g(h(x)) \quad \mbox{f\"ur ~ alle ~} x \in A\,.
	$$
	Sind $h : A \to B$, $g : B \to C$ und $f : C \to D$ Abbildungen,
	so sind die Verkettungen $f \circ (g \circ h), (f \circ g) \circ h : A \to D$
	definiert, und es gilt das {\it Assoziativgesetz}
	$f \circ (g \circ h) = (f \circ g) \circ h$\,.\\

	{\bf Beispiele f\"ur Funktionen und deren Verkettungen}
	
\begin{enumerate}[(a)]
	\item $f_1 : \R \to [-1,1]$\, mit \, $f_1(x) := \sin x$\, ist eine
	surjektive Funktion,\\ aber nicht injektiv.
	\item $f_2 : [-\frac{\pi}{2}, \frac{\pi}{2}] \to \R$ \, mit \, $f_2(x) := \sin x$ \, ist 
	injektiv, aber nicht surjektiv. 
	\item $f_3 : \R^+_0 \to \R^+_0$\, mit $f_3(x) := x^2$\,
	ist bijektiv mit Umkehrabbildung $f_3^{-1} : \R^+_0 \to \R^+_0$, $f_3^{-1}(x) = \sqrt{x}$\,.
	\item $f_4 : \R \to \R$\, mit $f_4(x) := x^2$\,
	ist weder injektiv noch surjektiv.
	\item $f_5 : \R \to \R^+$\, mit $f_5(x) := e^x$\,
	ist bijektiv mit Umkehrabbildung\\ $f_5^{-1} : \R^+ \to \R$, 
	$f_5^{-1}(x) = \ln x$.
\end{enumerate}

	Verkettungen wie $f_3 \circ f_1$ bzw. $f_2 \circ f_5$ 
	sind hier nicht m\"oglich, da weder
	$[-1,1] \subseteq  \R^+_0$ noch
	$\R^+ \subseteq [-\frac{\pi}{2}, \frac{\pi}{2}]$ gelten.
	Beispiele f\"ur ``erlaubte'' Verkettungen sind dagegen:
	\begin{enumerate}
		\item[(f)] $f_1 \circ f_3 : \R^+_0 \to [-1,1]$\, mit \, 
		$(f_1\circ f_3)(x) = \sin(x^2)$\,,
		\item[(g)]  $f_4 \circ f_1 : \R \to \R$\, mit \, 
		$(f_4\circ f_1)(x) = \sin^2 x$\,,
		\item[(h)]  $f_5 \circ f_2 : [-\frac{\pi}{2}, \frac{\pi}{2}] \to \R^+$\, mit \, 
		$(f_5\circ f_2)(x) = e^{\sin x}$\,,
		\item[(i)]  $f_5 \circ f_3^{-1} : \R^+_0 \to \R^+$\, mit \, 
		$(f_5\circ f_3^{-1})(x) = e^{\sqrt{x}}$\,.		
	\end{enumerate}	

	Die Verkettung bijektiver Abbildungen auf einer endlichen Tr\"agermenge f\"uhrt nun zu den Permutationsgruppen,
	die nicht nur in der linearen Algebra (Determinanten) sondern auch in der Zahlentheorie und Kombinatorik von Bedeutung sind:

\section{Permutationsgruppen}\label{cha:anhang2}\index{Permutationsgruppe}\label{Permutationsgruppe3}
	
	Permutationen sind bijektive Abbildungen einer Menge auf sich selbst.
	Bei unendlicher Tr\"agermenge nennt man sie auch Transformationen.
	Liegt dagegen eine endliche Tr\"agermenge mit $n \geq 1$ Elementen zugrunde,
	dann spricht man von Permutationen vom Grad $n$. Wir w\"ahlen
	im folgenden die feste Tr\"agermenge \mbox{$\N_n := \{\,1,2,...,n\,\}$}.\\
	
	{\bf Matrixdarstellung der Permutationen\index{Matrixdarstellung der Permutationen}\label{Matrixdarstellung der Permutationen}}
	
	Eine Permutation $f : \N_n \to \N_n$ l\"a{\ss}t sich wie folgt als
	Matrix schreiben:
	$$
	f = \begin{pmatrix}
	1    & 2    & ... & n\\
	f(1) & f(2) & ... & f(n)\\
	\end{pmatrix}\,.
	$$\\
\hspace*{0cm}\\\vspace{-1cm}
	
	{\bf Die Permutationsgruppe $\Sigma_n$}
	
	Sind $f,g : \N_n \to \N_n$ zwei beliebige Permutationen auf $\N_n$,
	so lassen sie sich gem\"a{\ss} $f \circ g : \N_n \to \N_n$ mit
	$(f \circ g)(x) := f(g(x))$ f\"ur alle $x \in \N_n$
	zu einer neuen Permutation $f \circ g$ verkn\"upfen.
	Damit wird die Menge $\Sigma_n = (\Sigma_n,\circ)$ 
	aller Permutationen auf $\N_n$ zu einer Gruppe,
	der sogenannten {\it Permutationsgruppe}
	$n$-ten Grades\index{Permutationsgruppe $n$-ten Grades}  
	\label{Permutationsgruppe n-ten Grades}.
	Sie besteht aus $n ! = 1 \cdot 2 \cdot ... \cdot n$ Permutationen
	\index{Fakult\"at} \label{Fakultaet2}.
	Bei dieser Verkn\"upfung ist nicht nur deshalb Vorsicht geboten, weil die
	Reihenfolge der ``Faktoren'' i.a. nicht vertauschbar ist, sondern auch
	deshalb, weil einige Autoren $f \circ g$ in der umgekehrten Reihenfolge
	$g(f)$ definieren! Dies h\"angt damit zusammen, da{\ss} bei unserer
	gel\"aufigeren Schreibweise die Funktionsauswertung zwar von
	``rechts nach links'' erfolgt, aber die Komposition von
	``links nach rechts'' aufgeschrieben wird. Dies kann als
	Diskrepanz empfunden werden.
	
	Das Einselement dieser Gruppe wird auch als {\it Identit\"at} $\Id$ 
	bzw. $\Id_n$ bezei\-chnet und hat die Darstellung
	$$
	\Id = \begin{pmatrix}
	1  & 2  & ... & n\\
	1  & 2  & ... & n\\
	\end{pmatrix}\,.
	$$
	Die zu $f$ inverse Permutation $f^{-1}$ \index{inverse Permutation}\label{inverse_Permutation} entsteht aus der
	Matrix von $f$ durch Vertauschung ihrer beiden Zeilen, d.h.
	$$
	f^{-1} = \begin{pmatrix}
	f(1) & f(2) & ... & f(n)\\
	1    & 2    & ... & n\\
	\end{pmatrix}\,.
	$$
	So erhalten wir etwa f\"ur $n=4$, d.h. $\N_n=\{\,1,2,3,4\,\}$,
	das Beispiel
	$$
	f = \begin{pmatrix}
	1 & 2  & 3 & 4\\
	2 & 4  & 1 & 3\\
	\end{pmatrix}\,, \quad
	f^{-1} = \begin{pmatrix}
	2 & 4  & 1 & 3\\
	1 & 2  & 3 & 4\\
	\end{pmatrix}
	= \begin{pmatrix}
	1 & 2  & 3 & 4 \\
	3 & 1  & 4 & 2 \\
	\end{pmatrix}\,.
	$$\\
	
\hspace*{0cm}\\\vspace{-1cm}
	
	{\bf Die Zyklenschreibweise f\"ur Permutationen\index{Zyklenschreibweise f\"ur Permutationen}\label{Zyklenschreibweise fuer Permutationen}}
	Neben der Matrixdarstellung gibt es aber auch noch die Zerlegung
	einer Permutation in {\it elementfremde Zyklen}\index{Zyklen einer Permutation}\label{Zyklen_einer_Permutation}. 
	Diese f\"uhrt auf eine weitere sehr wichtige Darstellung f\"ur Permutationen. 
	Wir betrachten als Beispiel die Permutationen $f,g : \N_6 \to \N_6$ mit
	$$
	f = \begin{pmatrix}
	1 & 2  & 3 & 4 & 5 & 6\\
	2 & 1  & 3 & 6 & 4 & 5\\
	\end{pmatrix}\,, \quad
	g = \begin{pmatrix}
	1 & 2  & 3 & 4 & 5 & 6\\
	2 & 3  & 4 & 1 & 6 & 5\\
	\end{pmatrix}\,.
	$$
	Die Permutation $f$ vertauscht die Ziffern 1,2 miteinander,
	hat die Ziffer 3 als sogenannten {\it Fixpunkt} und \"uberf\"uhrt
	die Ziffern 4,6,5 zyklisch ineinander in der angegebenen Reihenfolge
	$4 \to 6 \to 5 \to 4$. Entsprechend finden wir f\"ur $g$ die beiden Zyklen
	$1 \to 2 \to 3 \to 4 \to 1$ bzw. $5 \to 6 \to 5$. 
	
	Allgemein schreibt man einen Zyklus $k_1 \to k_2 \to ... \to k_m \to k_1$ 
	mit verschiedenen $k_1$,...,$k_m$
	in der Form $Z=(k_1, k_2, ... ,k_m)$.
	Mit $|Z|=m$ bezeichnen wir die L\"ange dieses Zyklus.

	F\"ur $f$ und $g$ haben wir somit die folgenden Zerlegungen 
	in elementfremde Zyklen gefunden:
	$$
	f = [(1,2)(3)(4,6,5)]\,, \quad
	g = [(1,2,3,4)(5,6)]\,.
	$$
        Fixpunkte, d.h. Zyklen der L\"ange 1, 
	l\"a{\ss}t man meistens weg und schreibt dann etwa
	$f= [(1,2)(4,6,5)]$\,, $\Id_6= [\,]$.\\

	Die Injektivit\"at der Permutationen auf $\N_n$ garantiert im allgemeinen Fall,
	da{\ss} sich jeder Zyklus wieder mit dem Element schlie{\ss}t,
	mit dem man begonnen hat. Jedes Element $k \in N_n$ besitzt n\"amlich bzgl. einer Abbildung
        $f \in \Sigma_n$ einen eindeutigen Vorg\"anger $f^{-1}(k)$, und somit nicht nur 
        einen eindeutigen Nachfolger $f(k)$. Daher gilt auch der folgende 
	
	\begin{Satz}\label{satz:zyklenform}
	Jede Permutation auf $\N_n$
	l\"a{\ss}t sich eindeutig in elementfremde Zyklen zerlegen. 
	\hfill \dokendSatz
        \end{Satz}
        
        	Die Zyklenzerlegung der Permutationen l\"a{\ss}t sich 
	graphisch gut illustrieren:
	
	\begin{figure}[h] 
		\begin{center}
			\epsfxsize=80mm
			\epsffile{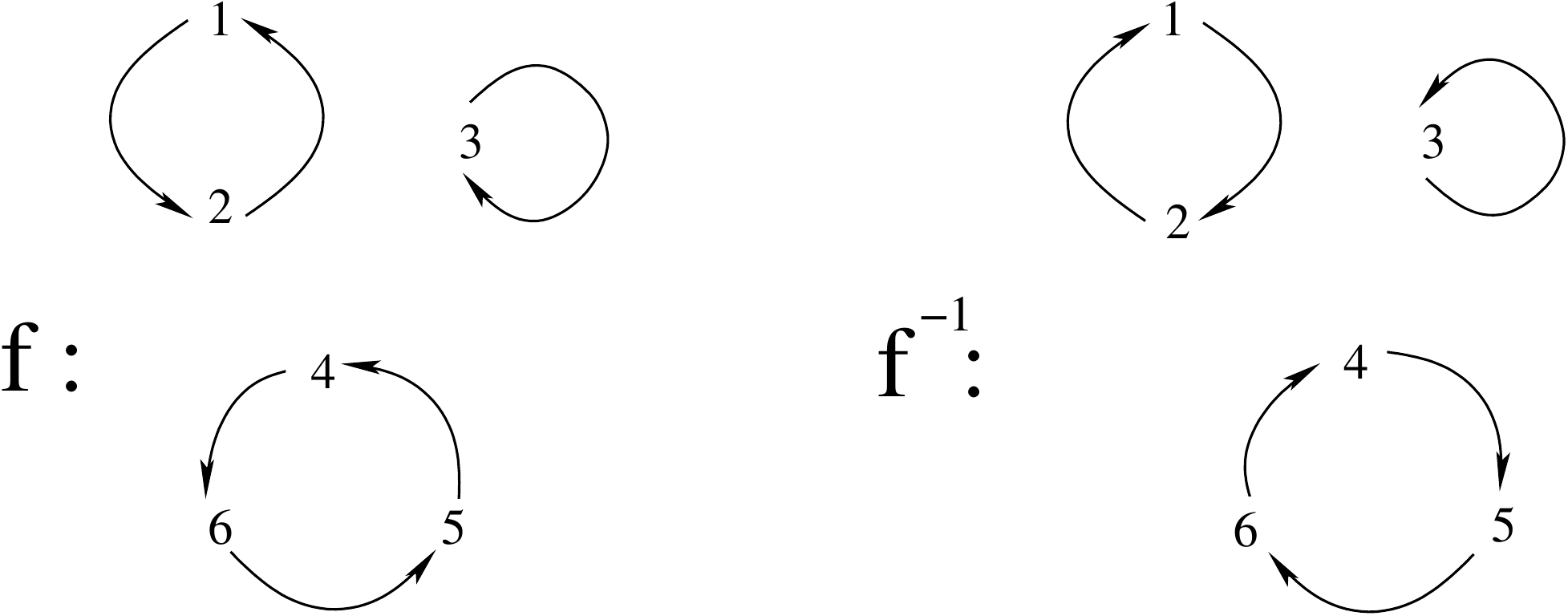}
		\end{center}
	\end{figure}

	Wir k\"onnen auch aus der Zyklen\-zerlegung
	sofort die Inversen bzw. die Kompositionen erhalten:
	\begin{align*}
	f^{-1} = [(2,1)(5,6,4)]\,, ~
	g^{-1} = [(4,3,2,1)(6,5)]\,, \\
	f \circ g = [(2,3,6,4)]\,, ~
	g \circ f = [(1,3,4,5)]\,.
	\end{align*}

	{\bf Zerlegung einer Permutationen in Transpositionen\index{Transposition}\label{Transposition2}}

	Eine Transposition ist eine Permutation der Form $[(a,b)]$,
	die nur zwei Ziffern $a \neq b$ miteinander vertauscht.
	F\"ur eine zyklische Permutation $[(n_1,n_2,...,n_r)]$
	mit der Zyklenl\"ange $r \geq 2$
	besteht die folgende Zerlegung in $(r-1)$ Transpositionen,
	die sich mittels vollst\"andiger Induktion zeigen l\"a{\ss}t:
	\begin{equation}\label{zykluszerlegung}
	[(n_1,n_2,...,n_r)] = [(n_1,n_r)] \circ ... \circ [(n_1,n_2)]
	\end{equation}
	Im folgenden sei $f : \N_n \to \N_n$ eine Permutation und $n \geq 2$.
	Da sich nach dem vorigen Satz $f$ in paarweise disjunkte (d.h. elementfremde) 
	Zyklen $Z_1$,...,$Z_s$ gem\"a{\ss} $f = [Z_1] \circ [Z_2] ... \circ [Z_s]$
	zerlegen l\"a{\ss}t und wir f\"ur $f \neq \Id$ 
	die Fixpunktzyklen aus dieser Zerlegung
	streichen k\"onnen, folgt in diesem Fall die Zerlegbarkeit von $f$ 
	in ein Produkt von Transpositionen. F\"ur $f = \Id$ k\"onnen wir dagegen
	wegen $n \geq 2$ die Zerlegung $\Id = [(1,2)] \circ [(1,2)]$ angeben.\\
		
	\begin{Def}[Gerade und ungerade Permutationen]\label{def:persignum} 
	\index{ungerade Permutation}\label{ungerade_Permutation}\index{gerade Permutation}\label{gerade_Permutation}
        Eine Permutation $f : \N_n \to \N_n$ hei{\ss}t {\it gerade},
	wenn sie sich in eine gerade Anzahl von Transpositionen \index{Transposition}\label{Transposition}
        faktorisieren l\"a{\ss}t.
	In diesem Falle schreiben wir $\mbox{sign}(f) = +1$. 
	Ist dagegen eine solche Zerlegung nicht
	m\"oglich, so hei{\ss}t die Permutation {\it ungerade}, und wir schreiben
	dann $\mbox{sign}(f) = -1$.
        \hfill \dokendDef
	\end{Def}

        Die Zerlegung einer Permutation in Transpositionen ist im allgemeinen
	nicht eindeutig. Umso wichtiger ist der folgende

	\begin{Satz}\label{satz:geradeungerade} Die Permutation $f : \N_n \to \N_n$ mit $n \geq 2$
	sei auf zwei verschiedene Arten in Transpositionen $T_k$, $T'_k$ zerlegt
	gem\"a{\ss}
	$$
	f = T_1 \circ ... \circ T_r = T'_1 \circ ... \circ T'_{r'}\,. 
	$$
	Dann sind $r$ und $r'$ entweder beide gerade oder beide ungerade.
        \hfill \dokendSatz
	\end{Satz}
	{\bf Beweis:}~
	Wir definieren das folgende Polynom\index{Polynom}\label{Polynom}:
	$$
	P(x_1,x_2,...,x_n) := \prod \limits_{1 \leq j < k \leq n} (x_k - x_j)\,.
	$$
	Nun geben wir zwei beliebige Zahlen $m > m'$ aus $\N_n$ vor und zerlegen
	dieses Polynom in f\"unf Faktoren gem\"a{\ss}
	\begin{align*}
	P(x_1,x_2,...,x_n)  &= (x_m - x_{m'})\,\cdot \,
	\prod \limits_{j < k \,\wedge \, j,k \notin \{m,m'\}} (x_k - x_j)\,\cdot\\
	& \prod \limits_{j > m} \Big\{ (x_j - x_m)(x_j - x_{m'}) \Big\} \cdot
	\prod \limits_{k < m'} \Big\{ (x_m - x_k)(x_{m'}-x_k) \Big\}\,\cdot\\
	&\prod \limits_{m' < k < m} \Big\{(x_m - x_k)(x_k-x_{m'})\Big\}\,.
	\end{align*}
	Produkte \"uber einen leeren Indexbereich sollen hierbei den Wert 1 haben.
	Vertauschen wir die Variablen $x_m$ und $x_{m'}$ in $P(x_1,x_2,...,x_n)$,
	so wechselt das Polynom nur sein Vorzeichen, da die vier mit $\prod$
	beginnenden Produkte hierbei unver\"andert bleiben, w\"ahrend der
	erste Faktor $(x_m - x_{m'})$ sein Vorzei\-chen wechselt. 
	
	Wir definieren f\"ur jedes $g \in \Sigma_n$ das Polynom
	$P_g(x_1,...,x_n) :=P(x_{g(1)},...,x_{g(n)})$ und beachten
	f\"ur alle $g,h \in \Sigma_n$ die Assoziativit\"at
	$$
	(P_g)_h = P_{g \circ h} \,.
	$$
	F\"ur die beliebige Transposition $T = [(m,m')]$ folgt nach dem oben gezeigten 
	$$
	P_T(x_1,...,x_n) = -P(x_1,...,x_n)\,.
	$$  
	Wenden wir die letzten beiden Beziehungen wiederholt auf die beiden Zerlegungen
	$f = T_1 \circ ... \circ T_r = T'_1 \circ ... \circ T'_{r'}$ an, so erhalten
	wir die folgende Gleichung, die unsere Ausgangsbehauptung beweist:
	$$
	P_f(x_1,...,x_n) = (-1)^{r}P(x_1,...,x_n) = (-1)^{r'}P(x_1,...,x_n)\,.
	$$  
	Speziell f\"ur $x_k:=k \in \N_n$ erhalten wir zudem
	$\mbox{sign}(f) = P_f(1,...,n)/P(1,...,n)$\,.
        \hfill \dokendProof\\

	Nun gilt der folgende wichtige
	\begin{Satz}\label{satz:permsigman}
	Wir betrachten die Permutationsgruppe $(\Sigma_n,\circ)$
	auf $\N_n$, $n \geq 2$.
%	\begin{itemize}
%		\item[(a)]~
%		F\"ur je zwei Permutationen $f,g \in \Sigma_n$ gilt
%		\begin{align*}
%		\mbox{sign}(f \circ g) = \mbox{sign}(f) \cdot \mbox{sign}(g)\,,\quad
%		\mbox{sign}(\Id) = 1\,, \quad
%		\mbox{sign}(f^{-1}) = \mbox{sign}(f)\,.
%		\end{align*}
%		\item[(b)] ~Die geraden Permutationen bilden eine Untergruppe von 
%		$(\Sigma_n,\circ)$, die sogenannte {\it alternierende Gruppe}
%                \index{alternierende Gruppe}
%		$(A_n,\circ)$, die aus $\frac12 n!$ Permutationen besteht.\\
%		\item[(c)] ~
%		Ist weiter $g$ die Anzahl der Zyklen von $f$ mit gerader L\"ange, so 
%                gilt f\"ur $\mbox{sign}(f)$ die Berechnungsvorschrift
%		$
%		\mbox{sign}(f) = (-1)^{g}.
%		$
%	\end{itemize}
\begin{enumerate}[(a)]
	\item F\"ur je zwei Permutationen $f,g \in \Sigma_n$ gilt
	\begin{align*}
	\mbox{sign}(f \circ g) = \mbox{sign}(f) \cdot \mbox{sign}(g)\,,\quad
	\mbox{sign}(\Id) = 1\,, \quad
	\mbox{sign}(f^{-1}) = \mbox{sign}(f)\,.
	\end{align*}
	\item Die geraden Permutationen bilden eine Untergruppe von 
	$(\Sigma_n,\circ)$, die sogenannte {\it alternierende Gruppe}
	\index{alternierende Gruppe}\label{alt_gruppe}
	$(A_n,\circ)$, die aus $\frac12 n!$ Permutationen besteht.\\
	\item Ist weiter $g$ die Anzahl der Zyklen\index{Zyklus}\label{Zyklus} von $f$ mit gerader L\"ange, so 
	gilt f\"ur $\mbox{sign}(f)$ die Berechnungsvorschrift
	$
	\mbox{sign}(f) = (-1)^{g}.
	$
\end{enumerate}
       \hfill \dokendSatz
        \end{Satz}

	{\bf Beweis:}~
	Die Teilaussage (a) ergibt sich aus Satz \ref{satz:geradeungerade},
        und (b) ist eine direkte Folge von (a). Wir zeigen die Berechnungsformel 
        f\"ur $\mbox{sign}(f)$: Ist $f$ vollst\"andig in seine paarweise
		disjunkten (d.h. elementfremden) Zyklen
		$Z_1$,...,$Z_s$ (mit oder ohne Einerzyklen) zerlegt und bezeichnet
		$|Z_k|$ die L\"ange des $k$-ten Zyklus, $k=1,...,s$, so haben wir in
                \eqref{zykluszerlegung}
                jeden Zyklus $Z_k$ als Produkt von $|Z_k|-1$ Transpositionen dargestellt.
                Folglich gilt die Beziehung
		$$
		\mbox{sign}(f) = (-1)^{m}~~
		\mbox{mit}~~ m := \sum \limits_{k=1}^{s}(|Z_k|-1)\,.
		$$
        Allein f\"ur die Zyklen $Z_k$ mit gerader L\"ange $|Z_k|$ ist $|Z_k|-1 \equiv 1\, (2)$,
       f\"ur die $Z_k$ mit ungerader L\"ange ist dagegen $|Z_k|-1 \equiv 0 \,(2)$\,. Somit ist
       $\mbox{sign}(f) = (-1)^{g}$\,.
        \hfill \dokendProof\\

	{\it Beispiel:}~
	Ist $f : \N_8 \to \N_8$ in der Zyklenform
	$f :=[(1,7,8)(2,5,4,3)(6)]$ gegeben, so ist $(2,5,4,3)$
        der einzige Zyklus gerader L\"ange von $f$ und 
        $
	\mbox{sign}(f) = (-1)^1 = -1\,.
	$

	\newpage
	\section{Primzahltabelle}\label{cha:anhang3}

\begin{figure}[hb]
	\centering
	\includegraphics[width=1\linewidth]{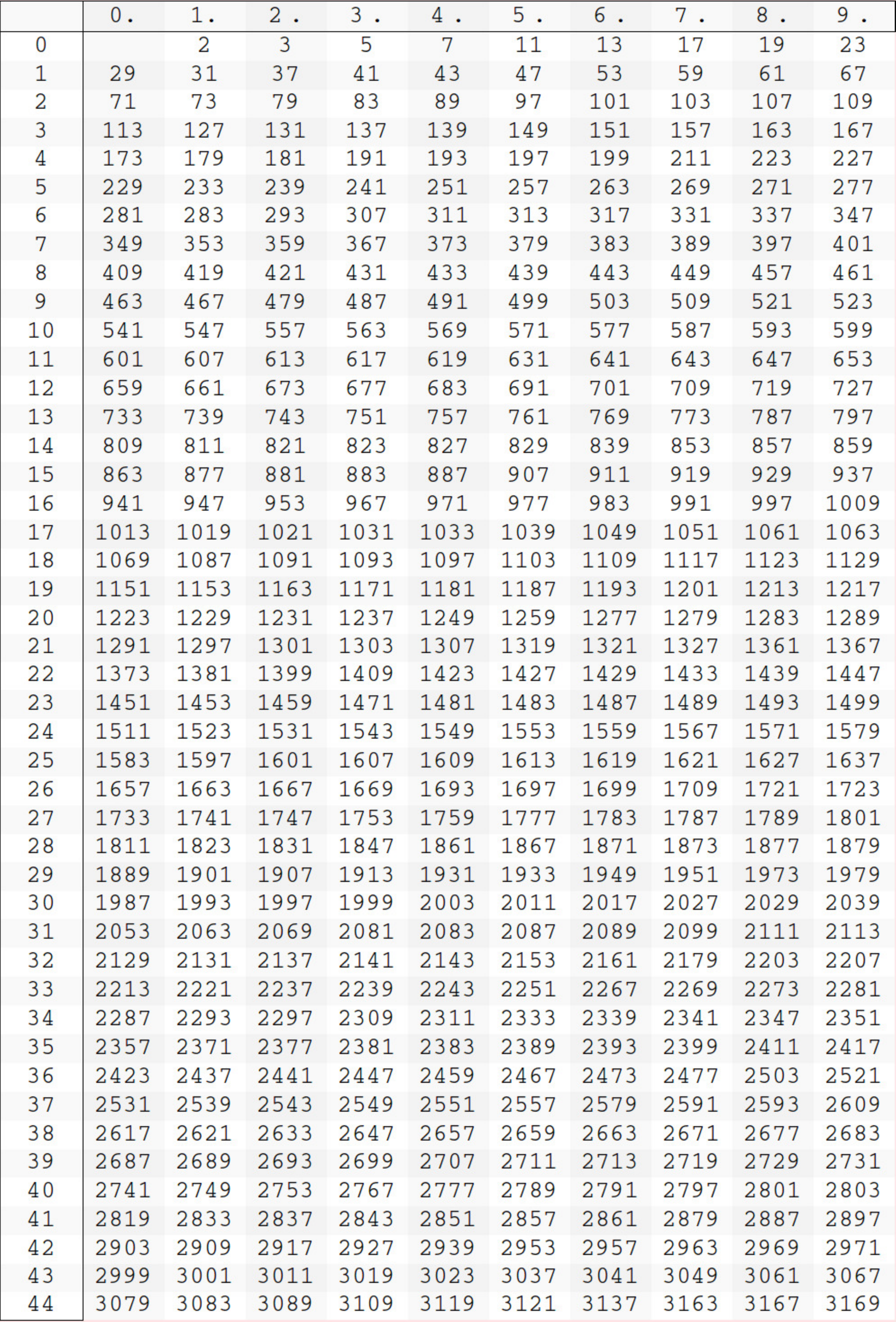}
\end{figure}
%	\newpage
\begin{figure}[hb]
	\centering
	\includegraphics[width=1\linewidth]{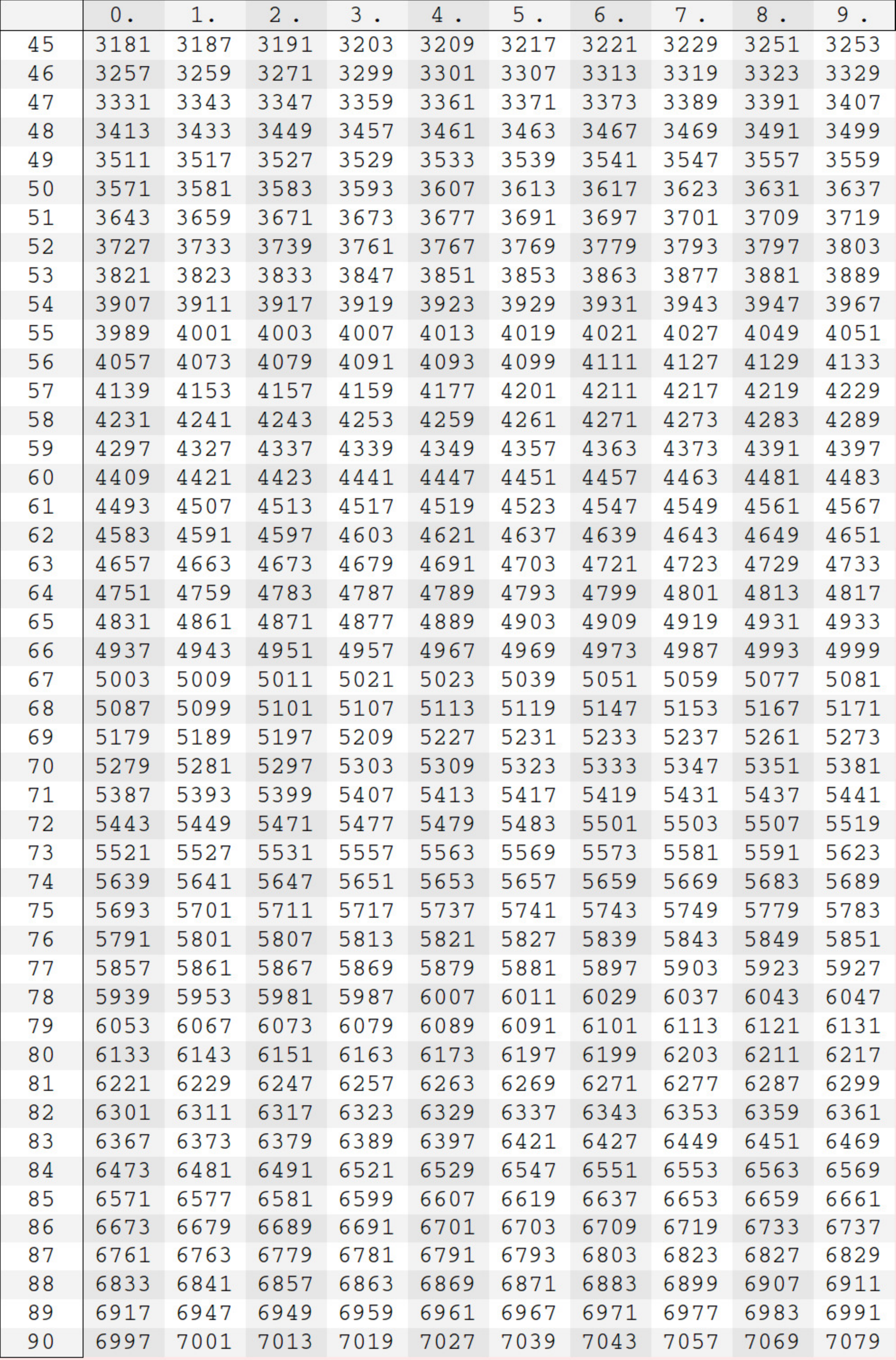}
\end{figure}
%	\newpage
\begin{figure}[hb]
	\centering
	\includegraphics[width=1\linewidth]{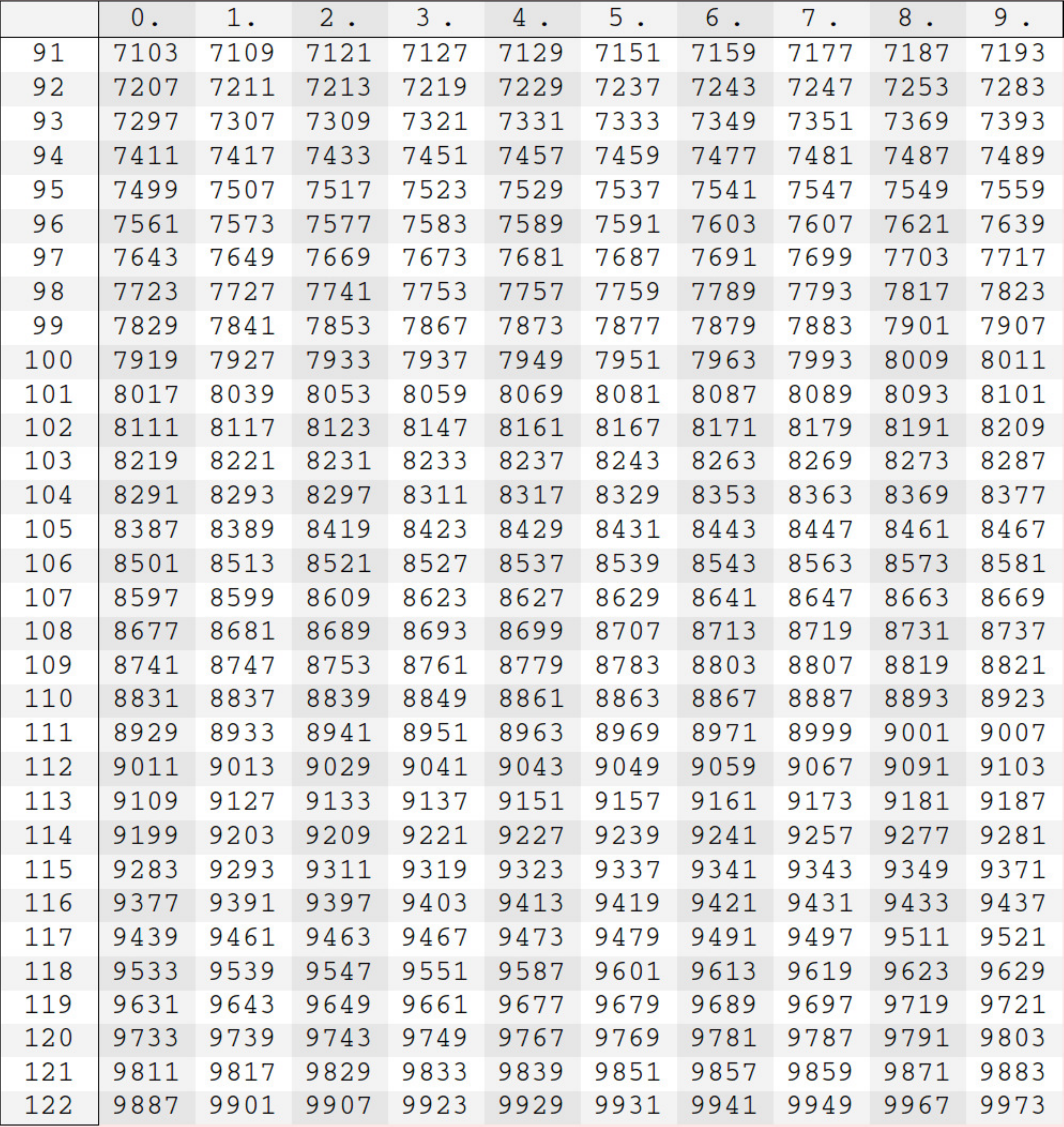}
\end{figure}

% die Tabelle mit Lighter[Gray,0.9] - ganz hell
%\begin{figure}[hb]
%	\centering
%	\includegraphics[width=1\linewidth]{prime_1}
%\end{figure}
% die Tabelle mit Lighter[Gray,0.8] - hell
%\begin{figure}[hb]
%	\centering
%	\includegraphics[width=1\linewidth]{prime_2}
%\end{figure}
	
	\backmatter
	%%%%%%%%%%%%%%%%%%%%%%%%%%%%%%%%%%%%%%%%%%%%%%%%%%%%%%%%%%%%%%%%%%%%%%
	%\include{glossary}
	%\include{solutions}
	%\printindex
	
	%%%%%%%%%%%%%%%%%%%%%%%%%%%%%%%%%%%%%%%%%%%%%%%%%%%%%%%%%%%%%%%%%%%%%%
	\backmatter
	\pagestyle{empty}
	\clearpage
	\renewcommand{\refname}{Literaturverzeichnis}
	%\addcontentsline{toc}{chapter}{Literaturverzeichnis}
	
	\renewcommand{\indexname}{Indexverzeichnis}
	\begin{theindex}
{\bf A}\nopagebreak%
 \indexspace\nopagebreak%
 \item Abbildung\idxquad \pageref{Abbildung}
  \item abelsche Gruppe\idxquad \pageref{abelsche_Gruppe}, \pageref{abelsche_Gruppe2}
  \item alternierende Gruppe\idxquad \pageref{alt_gruppe}
  \item Approximation einer Irrationalzahl\idxquad \pageref{Approximation einer Irrationalzahl}
  \item Approximationssatz f\"ur Farey-Br\"uche\idxquad \pageref{Approximationssatz fuer Farey-Brueche}
  \item Approximationssatz von Hurwitz\idxquad \pageref{Approximationssatz von Hurwitz}
  \item Assoziativgesetz\idxquad \pageref{Assoziativgesetz}
  \item Ausgabewerte des erweiterten Euklidischen Algorithmus\idxquad \pageref{Ausgabewerte des erweiterten Euklidischen Algorithmus}
  \item Aussage\idxquad \pageref{Aussage}
  \item Aussageform\idxquad \pageref{Aussageform}
  \item aussagenlogische Verkn\"upfung\idxquad \pageref{aussagenlogische Verknuepfung}
 \item Automorphe Transformation\idxquad \pageref{Auto},\pageref{Auto2}

  \indexspace
{\bf B}\nopagebreak%
\indexspace\nopagebreak%
\item $b$-adische Darstellung\idxquad \pageref{b-adische Darstellung}
\item beste rationale Approximation\idxquad \pageref{beste rationale Approximation}
\item Binetsche Formel\idxquad \pageref{Binetsche_Formel}, \pageref{Binetsche_Formel2}
 
 \indexspace
{\bf C}\nopagebreak%
\indexspace\nopagebreak%
\item Chinesischer Restsatz\idxquad \pageref{Chinesischer Restsatz}

  \indexspace
{\bf D}\nopagebreak%
 \indexspace\nopagebreak%
 \item definit quadratische Form\idxquad \pageref{definit quadratische Form}
 \item Dirichlet-Faltung\idxquad \pageref{Dirichlet-Faltung}
 \item Dirichletscher Approximationssatz\idxquad \pageref {Dirichletscher Approximationssatz}
 \item Diskriminante\idxquad \pageref{Diskriminante}
  \item Divisionskoeffizient\idxquad \pageref{Divisionskoeffizient}, \pageref{Divisionskoeffizient2}
  \item Divisionsrest\idxquad \pageref{Divisionsrest}, \pageref{Divisionskoeffizient2}

  \indexspace
{\bf E}\nopagebreak%
 \indexspace\nopagebreak%
 \item eindeutige Primfaktorzerlegung\idxquad \pageref{eindeutige Primfaktorzerlegung}
 \item Eingabewerte des erweiterten Euklidischen Algorithmus\idxquad \pageref{Eingabewerte des erweiterten Euklidischen Algorithmus}
  \item Einheit\idxquad \pageref{Einheit}
  \item Einselement einer Gruppe\idxquad \pageref{Einselement_einer_Gruppe}
   \item Einselement eines Ringes\idxquad \pageref{Einselement_eines_Ringes}
   \item endlicher Kettenbruch\idxquad \pageref{endlicher Kettenbruch}
   \item erweiterte Farey-Sequenz\idxquad \pageref{erweiterte Farey-Sequenz}
  \item erweiterter Euklidischer Algorithmus\idxquad \pageref{erweiterter_Euklidischer_Algorithmus}, \pageref{erweiterter_Euklidischer_Algorithmus2}, \pageref{erweiterter_Euklidischer_Algorithmus3}
  \item Euklidischer Algorithmus\idxquad  \pageref{Euklidischer_Algorithmus}, \pageref{Euklidischer_Algorithmus2}, \pageref{Euklidischer_Algorithmus3}
  \item Euklidischer Ring\idxquad \pageref{Euklidischer_Ring}
  \item{Euler, Leonard}\idxquad \pageref{Leonard Euler}
	\item Eulersche Funktion\idxquad \pageref{Eulersche Funktion}, \pageref{Eulersche Funktion2}
	\item Eulersches Kriterium\idxquad \pageref{Eulersches Kriterium}
  \item Exponent\idxquad \pageref{Exponent}

  \indexspace
{\bf F}\nopagebreak%
 \indexspace\nopagebreak%	
 \item Fakult\"at \idxquad \pageref{Fakultaet},\pageref{Fakultaet2}
 \item Faltungsgruppe der multiplikativen Funktionen\idxquad \pageref{Faltungsgruppe der multiplikativen Funktionen}
\item{Farey, John}\idxquad \pageref{John Farey}
 \item Farey-Sequenz\idxquad \pageref{Farey-Sequenz}
 \item Fermatsche Primzahl\idxquad \pageref{Fermatsche Primzahl}
 \item Fibonacci-Folge\idxquad \pageref{Fibonacci-Folge}, \pageref{Fibonacci-Folge2}
 \item Fibonacci-Zahlen\idxquad \pageref{Fibonacci-Zahlen}, \pageref{Fibonacci-Zahlen2}, \pageref{Fibonacci-Zahlen3}, \pageref{Fibonacci-Zahlen4}
  \item Fundamentalsatz der Arithmetik\idxquad \pageref{Fundamentalsatz_der_Arithmetik}, \pageref{Fundamentalsatz_der_Arithmetik2}, \pageref{Fundamentalsatz_der_Arithmetik3},
  \pageref{Fundamentalsatz_der_Arithmetik4}
  \item Funktion\idxquad \pageref{Funktion}

  \indexspace
{\bf G}\nopagebreak%
 \indexspace\nopagebreak%
\item Gau{\ss}, Carl Friedrich \idxquad \pageref{Carl Friedrich Gauss}
  \item Gau{\ss}-Klammer\idxquad \pageref{Gaussklammer}, \pageref{Gaussklammer2}
  \item Gau{\ss}sches Lemma\idxquad \pageref{Gausssches Lemma}
  \item gek\"urzter Bruch\idxquad \pageref{gekuerzter Bruch}
  \item gerade Permutation\idxquad \pageref{gerade_Permutation}
  \item gro{\ss}e Faltungsgruppe\idxquad \pageref{grosse Faltungsgruppe}
  \item gr\"o{\ss}ter gemeinsamer Teiler\idxquad \pageref{groesster_gemeinsamer _Teiler}, \pageref{groesster_gemeinsamer _Teiler2}
  \item Gruppe\idxquad \pageref{Gruppe}
\indexspace
{\bf H}\nopagebreak
\item Huygens, Christiaan\idxquad \pageref{Christiaan Huygens}
  \indexspace
{\bf I}\nopagebreak%
 \indexspace\nopagebreak%
  \item Identit\"at\idxquad \pageref{Identitaet}
  \item indefinite qadratische Form\idxquad \pageref{indefinite quadratische Form} 
	%\pageref{indefinite quadratische Form2}
	 \item Index\idxquad \pageref{Index}
  \item Induktionsprinzip\idxquad \pageref{Induktionsprinzip}
  \item Integrit\"atsbereich\idxquad \pageref{Integritaetsbereich}, \pageref{Integritaetsbereich2}, \pageref{Integritaetsbereich3}
  \item Intervall\idxquad \pageref{Intervall}
  \item inverse Permutation\idxquad \pageref{inverse_Permutation}
  \item inverses Element\idxquad \pageref{inverses_Element}
  \item Irrationalzahl\idxquad \pageref{Irrationalzahl}, \pageref{Irrationalzahl2}, \pageref{Irrationalzahl3}
  \item Isomorphismus\idxquad \pageref{Isomorphismus}

 \indexspace
 {\bf J}\nopagebreak%
 \indexspace\nopagebreak%
\item Jacobi-Symbol\idxquad \pageref{Jacobi-Symbol}

  \indexspace
{\bf K}\nopagebreak%
 \indexspace\nopagebreak%
 \item Kartesisches Produkt\idxquad \pageref{Kartesisches Produkt}
  \item Kettenbruch\idxquad \pageref{Kettenbruch}, \pageref{Kettenbruch2}, \pageref{Kettenbruch3}
  \item Kettenbruchentwicklung\idxquad \pageref{Kettenbruchentwicklung}, \pageref{Kettenbruchentwicklung2}
  \item kommutative Gruppe\idxquad \pageref{kommutative_Gruppe}
  \item kommutativer Ring\idxquad \pageref{kommutativer_Ring}, \pageref{kommutativer_Ring2}, \pageref{kommutativer_Ring3}
  \item Kongruenz\idxquad \pageref{Kongruenz}, \pageref{Kongruenz2}
  \item K\"orper\idxquad \pageref{Koerper}, \pageref{Koerper2}
  \item K-reduzierte indefinite Form \idxquad \pageref{K-reduzierte indefinite Form}

  \indexspace
{\bf L}\nopagebreak%
 \indexspace\nopagebreak%
 \item Lagrange, Joseph-Louis \idxquad \pageref{Joseph-Louis Lagrange}
 \item Legendre-Symbol\idxquad \pageref{Legendre-Symbol}, \pageref{Legendre-Symbol2}
  \item Linearform\idxquad \pageref{Linearform}
  \item Linksnebenklasse\idxquad \pageref {Linksnebenklasse}
  \item Logische Symbole\idxquad \pageref{Logische Symbole}

  \indexspace
  {\bf M}\nopagebreak%
  \indexspace\nopagebreak%
  \item Matrixdarstellung der Permutationen\idxquad \pageref{Matrixdarstellung der Permutationen}
  \item Matrizenmultiplikation\idxquad \pageref{Matrizenmultiplikation}, \pageref{Matrizenmultiplikation2}
  \item Mediant\idxquad \pageref{Mediant}
  \item Mediantensatz\idxquad \pageref{Mediantensatz}, \pageref{Mediantensatz2}, \pageref{Mediantensatz3}, \pageref{Mediantensatz4}
  \item Menge\idxquad \pageref{Menge}
  \item M\"obius-Funktion\idxquad \pageref{Moebius-Funktion}, \pageref{Moebius-Funktion2}
  \item M\"obiussche Umkehrformel\idxquad \pageref{Moebiussche Umkehrformel}, \pageref{Moebiussche Umkehrformel2}
  \item Modul\idxquad \pageref{Modul}
  \item multiplikativ\idxquad \pageref{multiplikativ}
  \item multiplikative Inverse\idxquad \pageref{multiplikative Inverse}

  \indexspace
{\bf N}\nopagebreak%
 \indexspace\nopagebreak%
  \item Nullelement\idxquad \pageref{Nullelement}
  \item Nullteiler\idxquad \pageref{Nullteiler}

  \indexspace
  {\bf O}\nopagebreak%
  \indexspace\nopagebreak%
  \item Ordnung\idxquad \pageref{Ordnung}

  \indexspace
{\bf P}\nopagebreak%
 \indexspace\nopagebreak%
\item Pellsche Gleichung \idxquad \pageref{Pell}, \pageref{Pell2}
 \item periodische Funktion\idxquad \pageref{periodische Funktion}
 \item periodische Kettenbruchentwicklung\idxquad \pageref{periodische Kettenbruchentwicklung}
  \item Permutationsgruppe\idxquad\pageref{Permutationsgruppe}, \pageref{Permutationsgruppe2}, \pageref{Permutationsgruppe3}
  \item Permutationsgruppe $n$-ten Grades\idxquad \pageref{Permutationsgruppe n-ten Grades}
	  \item Perron, Oskar \idxquad \pageref{Oskar Perron}
  \item Polynom\idxquad \pageref{Polynom}
  \item Polynom mit ganzzahligen Koeffizienten\idxquad \pageref{Polynom mit ganzzahligen Koeffizienten}
  \item positiv definit quadratische Form\idxquad \pageref{positiv definit quadratische Form}
  \item prime Restklassengruppe\idxquad \pageref{prime Restklassengruppe}
  \item Primelemente\idxquad \pageref{Primelement}, \pageref{Primelement2}
  \item Primfaktorzerlegung\idxquad \pageref{Primfaktorzerlegung}
  \item primitive quadratische Form\idxquad \pageref{primitive quadratische Form}
  \item Primitivwurzel\idxquad \pageref{Primitivwurzel}, \pageref{Primitivwurzel2}
  \item Primzahl\idxquad \pageref{Primzahl}, \pageref{Primzahl2}, \pageref{Primzahl3}
  \item Primzahl-Moduln\idxquad \pageref{Primzahl-Modul}
	 \item Primzahlpotenz\idxquad \pageref{Primzahlpotenz}
  \item Pythagoreisches Zahlentripel\idxquad \pageref{Pythagoreisches_Zahlentripel}

  \indexspace
  {\bf Q}\nopagebreak%
  \indexspace\nopagebreak%
\item quadratische Irrationalzahl\idxquad \pageref{quadratische Irrationalzahl}, \pageref{quadratische Irrationalzahl2}, \pageref{quadratische Irrationalzahl3}
\item quadratische Nichtreste\idxquad \pageref{quadratische Nichtreste}
\item quadratische Reste\idxquad \pageref{quadratische Reste}, \pageref{quadratische Reste2}, \pageref{quadratische Reste3}
\item quadratische Form\idxquad \pageref{quadratische Form}, \pageref{quadratische Form2}
\item quadratischer Kongruenz\idxquad \pageref{quadratischer Kongruenz}, \pageref{quadratischer Kongruenz2}
\item quadratisches Reziprozit\"atsgesetz\idxquad \pageref{quadratisches Reziprozitaetsgesetz}

  \indexspace
{\bf R}\nopagebreak%
 \indexspace\nopagebreak%
 \item rationale Bestapproximation\idxquad \pageref{rationale Bestapproximation}
 \item reduziertes Restsystem\idxquad \pageref{reduziertes Restsystem}
\item reduzierte indefinite Form \idxquad \pageref{reduzierte indefinite Form}
 \item reell quadratische Irrationalzahl\idxquad \pageref{reell quadratische Irrationalzahl}, \pageref{reell quadratische Irrationalzahl2}
 \item Relation\idxquad \pageref{Relation}
 \item Reziprozit\"atsgesetz von Gau{\ss}\idxquad \pageref{Reziprozitaetsgesetz von Gauss}
 \item Ring\idxquad \pageref{Ring}
 \item r\"uckl\"aufige Rekursion f\"ur Farey-Br\"uche\idxquad \pageref{ruecklaeufige Rekursion fuer Farey-Brueche}

  \indexspace
{\bf S}\nopagebreak%
 \indexspace\nopagebreak%
  \item Satz von Lagrange\idxquad \pageref{Satz_von_Lagrange}

  \indexspace
{\bf T}\nopagebreak%
 \indexspace\nopagebreak%
   \item Teiler\idxquad \pageref{Teiler}
  \item teilerfremde Zahlen\idxquad \pageref{teilerfremde_Zahlen}
  \item Transformation der Formen\idxquad \pageref{Transformation der Formen}
  \item Transposition\idxquad \pageref{Transposition}, \pageref{Transposition2}

  \indexspace
{\bf U}\nopagebreak%
 \indexspace\nopagebreak%
 \item uneigentlich konjugierte Formen\idxquad \pageref{uneigentlich konjugierte Formen}
 \item uneigentlich konjugierte Klasse\idxquad \pageref{uneigentlich konjugierte Klasse}
 \item uneigentliche Transformation\idxquad \pageref{uneigentliche Transformation}, \pageref{uneigentliche Transformation2}
 \item unendlicher Kettenbruch\idxquad \pageref{unendlicher Kettenbruch}, \pageref{unendlicher Kettenbruch2}
  \item ungerade Permutation\idxquad \pageref{ungerade_Permutation}
  \item Untergruppe\idxquad \pageref{Untergruppe}
  
  \indexspace
  {\bf V}\nopagebreak%
  \indexspace\nopagebreak%
  \item Verkettung von Funktionen\idxquad \pageref{Verkettung von Funktionen}
  \item vollst\"andig multiplikativ\idxquad \pageref{vollstaendig multiplikativ}, \pageref{vollstaendig multiplikativ2}, \pageref{vollstaendig multiplikativ3}
  \item vollst\"andig multiplikative Funktion\idxquad \pageref{vollstaendig multiplikative Funktion}
  \item vollst\"andige Induktion\idxquad \pageref{vollstaendige_Induktion}
    \item vollst\"andiges Restsystem\idxquad \pageref{vollstaendiges Restsystem}
    
  \indexspace
{\bf W}\nopagebreak%
  \indexspace\nopagebreak%    
  \item Wahrheitstabelle\idxquad \pageref{Wahrheitstabelle}
 \item Wilsonscher Satz\idxquad \pageref{Wilsonscher Satz}

  \indexspace
{\bf Z}\nopagebreak%
 \indexspace\nopagebreak%
 \item zahlentheoretische Funktion\idxquad \pageref{zahlentheoretische Funktion}
  \item Zyklen einer Permutation\idxquad \pageref{Zyklen_einer_Permutation}
  \item Zyklenschreibweise f\"ur Permutationen\idxquad \pageref{Zyklenschreibweise fuer Permutationen}
  \item zyklische Untergruppe\idxquad \pageref{zyklische_Untergruppe}
  \item Zyklus \pageref{Zyklus}

\end{theindex}

\end{document}